\documentclass[reqno,twoside]{amsart}

\usepackage{amsfonts}
\usepackage{amsmath}
\usepackage{amssymb}
\usepackage{enumerate}
\usepackage[english]{babel}
\usepackage{color}
\usepackage{fullpage}
\usepackage{microtype}
\usepackage{hyperref}
\usepackage[leftcaption]{sidecap}
\sidecaptionvpos{figure}{m}

\usepackage{xcolor}
\usepackage{tikz}
\usepackage{pgfplots}
\usepackage{pgfplotstable}
\usepackage{subfig}

\DisableLigatures{encoding = *, family = * }

\newcommand{\NN}{\mathbb{N}}
\newcommand{\ZZ}{\mathbb{Z}}
\newcommand{\RR}{\mathbb{R}}

\newcommand{\Ggh}{\mathcal{G}^h_g}
\newcommand{\Bf}[1]{\textbf{#1}}
\newcommand{\Bs}[1]{\boldsymbol{#1}}

\ifx\figforTeXisloaded\relax \else\global\let\figforTeXisloaded=\relax\fi
\catcode`\@=11
\ifx\ctr@ln@m\undefined\else%
    \immediate\write16{*** Fig4TeX WARNING : \string\ctr@ln@m\space already defined.}\fi
\def\ctr@ln@m#1{\ifx#1\undefined\else%
    \immediate\write16{*** Fig4TeX WARNING : \string#1 already defined.}\fi}
\ctr@ln@m\ctr@ld@f
\def\ctr@ld@f#1#2{\ctr@ln@m#2#1#2}
\ctr@ld@f\def\ctr@ln@w#1#2{\ctr@ln@m#2\csname#1\endcsname#2}
{\catcode`\/=0 \catcode`/\=12 /ctr@ld@f/gdef/BS@{\}}
\ctr@ld@f\def\ctr@lcsn@m#1{\expandafter\ifx\csname#1\endcsname\relax\else%
    \immediate\write16{*** Fig4TeX WARNING : \BS@\expandafter\string#1\space already defined.}\fi}
\ctr@ld@f\edef\colonc@tcode{\the\catcode`\:}
\ctr@ld@f\edef\semicolonc@tcode{\the\catcode`\;}
\ctr@ld@f\def\t@stc@tcodech@nge{{\let\c@tcodech@nged=\z@%
    \ifnum\colonc@tcode=\the\catcode`\:\else\let\c@tcodech@nged=\@ne\fi%
    \ifnum\semicolonc@tcode=\the\catcode`\;\else\let\c@tcodech@nged=\@ne\fi%
    \ifx\c@tcodech@nged\@ne%
    \immediate\write16{}
    \immediate\write16{!!!=============================================================!!!}
    \immediate\write16{ Fig4TeX WARNING:}
    \immediate\write16{ The category code of some characters has been changed, which will}
    \immediate\write16{ result in an error (message "Runaway argument?").}
    \immediate\write16{ This probably comes from another package that changed the category}
    \immediate\write16{ code after Fig4TeX was loaded. If that proves to be exact, the}
    \immediate\write16{ solution is to exchange the loading commands on top of your file}
    \immediate\write16{ so that Fig4TeX is loaded last. For example, in LaTeX, we should}
    \immediate\write16{ say :}
    \immediate\write16{\BS@ usepackage[french]{babel}}
    \immediate\write16{\BS@ usepackage{fig4tex}}
    \immediate\write16{!!!=============================================================!!!}
    \immediate\write16{}
    \fi}}
\ctr@ld@f\def\FigforTeX{F\kern-.05em i\kern-.05em g\kern-.1em\raise-.14em\hbox{4}\kern-.19em\TeX}
\ctr@ld@f\def\W@rnmesoldA#1{\W@rnmesold}
\ctr@ld@f\def\W@rnmesoldAB#1(#2){\W@rnmesold}
\ctr@ld@f\def\W@rnmesold{%
    \immediate\write16{}
    \immediate\write16{!!!=============================================================!!!}
    \immediate\write16{ Fig4TeX WARNING:}
    \immediate\write16{ The file to be compiled is not compatible with the current version}
    \immediate\write16{ of Fig4TeX. To fix that, upgrade the source file (mainly change \BS@ ps*}
    \immediate\write16{ macros by \BS@ fig* macros), or use fig4tex184.tex instead (\BS@ input fig4tex184}
    \immediate\write16{ or \BS@ usepackage{fig4tex184}).}
    \immediate\write16{!!!=============================================================!!!}
    \immediate\write16{}}
\ctr@ln@m\psbeginfig\let\psbeginfig\W@rnmesoldA
\ctr@ln@m\psset\let\psset\W@rnmesoldAB
\ctr@ln@m\pssetdefault\let\pssetdefault\W@rnmesoldAB
\ctr@ln@m\pssetupdate\let\pssetupdate\W@rnmesoldA
\ctr@ln@w{newdimen}\epsil@n\epsil@n=0.00005pt
\ctr@ln@w{newdimen}\Cepsil@n\Cepsil@n=0.005pt
\ctr@ln@w{newdimen}\dcq@\dcq@=254pt
\ctr@ln@w{newdimen}\PI@\PI@=3.141592pt
\ctr@ln@w{newdimen}\DemiPI@deg\DemiPI@deg=90pt
\ctr@ln@w{newdimen}\PI@deg\PI@deg=180pt
\ctr@ln@w{newdimen}\DePI@deg\DePI@deg=360pt
\ctr@ld@f\chardef\t@n=10
\ctr@ld@f\chardef\c@nt=100
\ctr@ld@f\chardef\@lxxiv=74
\ctr@ld@f\chardef\@xci=91
\ctr@ld@f\mathchardef\@nMnCQn=9949
\ctr@ld@f\chardef\@vi=6
\ctr@ld@f\chardef\@xxx=30
\ctr@ld@f\chardef\@lvi=56
\ctr@ld@f\chardef\@@lxxi=71
\ctr@ld@f\chardef\@lxxxv=85
\ctr@ld@f\mathchardef\@@mmmmlxviii=4068
\ctr@ld@f\mathchardef\@ccclx=360
\ctr@ld@f\mathchardef\@dccxx=720
\ctr@ln@w{newcount}\p@rtent \ctr@ln@w{newcount}\f@ctech \ctr@ln@w{newcount}\result@tent
\ctr@ln@w{newdimen}\v@lmin \ctr@ln@w{newdimen}\v@lmax \ctr@ln@w{newdimen}\v@leur
\ctr@ln@w{newdimen}\result@t\ctr@ln@w{newdimen}\result@@t
\ctr@ln@w{newdimen}\mili@u \ctr@ln@w{newdimen}\c@rre \ctr@ln@w{newdimen}\delt@
\ctr@ld@f\def\degT@rd{0.017453 }  
\ctr@ld@f\def\rdT@deg{57.295779 } 
\ctr@ln@m\v@leurseule
{\catcode`p=12 \catcode`t=12 \gdef\v@leurseule#1pt{#1}}
\ctr@ld@f\def\repdecn@mb#1{\expandafter\v@leurseule\the#1\space}
\ctr@ld@f\def\arct@n#1(#2,#3){{\v@lmin=#2\v@lmax=#3%
    \maxim@m{\mili@u}{-\v@lmin}{\v@lmin}\maxim@m{\c@rre}{-\v@lmax}{\v@lmax}%
    \delt@=\mili@u\m@ech\mili@u%
    \ifdim\c@rre>\@nMnCQn\mili@u\divide\v@lmax\tw@\c@lATAN\v@leur(\z@,\v@lmax)
    \else%
    \maxim@m{\mili@u}{-\v@lmin}{\v@lmin}\maxim@m{\c@rre}{-\v@lmax}{\v@lmax}%
    \m@ech\c@rre%
    \ifdim\mili@u>\@nMnCQn\c@rre\divide\v@lmin\tw@
    \maxim@m{\mili@u}{-\v@lmin}{\v@lmin}\c@lATAN\v@leur(\mili@u,\z@)%
    \else\c@lATAN\v@leur(\delt@,\v@lmax)\fi\fi%
    \ifdim\v@lmin<\z@\v@leur=-\v@leur\ifdim\v@lmax<\z@\advance\v@leur-\PI@%
    \else\advance\v@leur\PI@\fi\fi%
    \global\result@t=\v@leur}#1=\result@t}
\ctr@ld@f\def\m@ech#1{\ifdim#1>1.646pt\divide\mili@u\t@n\divide\c@rre\t@n\m@ech#1\fi}
\ctr@ld@f\def\c@lATAN#1(#2,#3){{\v@lmin=#2\v@lmax=#3\v@leur=\z@\delt@=\tw@ pt%
    \un@iter{0.785398}{\v@lmax<}%
    \un@iter{0.463648}{\v@lmax<}%
    \un@iter{0.244979}{\v@lmax<}%
    \un@iter{0.124355}{\v@lmax<}%
    \un@iter{0.062419}{\v@lmax<}%
    \un@iter{0.031240}{\v@lmax<}%
    \un@iter{0.015624}{\v@lmax<}%
    \un@iter{0.007812}{\v@lmax<}%
    \un@iter{0.003906}{\v@lmax<}%
    \un@iter{0.001953}{\v@lmax<}%
    \un@iter{0.000976}{\v@lmax<}%
    \un@iter{0.000488}{\v@lmax<}%
    \un@iter{0.000244}{\v@lmax<}%
    \un@iter{0.000122}{\v@lmax<}%
    \un@iter{0.000061}{\v@lmax<}%
    \un@iter{0.000030}{\v@lmax<}%
    \un@iter{0.000015}{\v@lmax<}%
    \global\result@t=\v@leur}#1=\result@t}
\ctr@ld@f\def\un@iter#1#2{%
    \divide\delt@\tw@\edef\dpmn@{\repdecn@mb{\delt@}}%
    \mili@u=\v@lmin%
    \ifdim#2\z@%
      \advance\v@lmin-\dpmn@\v@lmax\advance\v@lmax\dpmn@\mili@u%
      \advance\v@leur-#1pt%
    \else%
      \advance\v@lmin\dpmn@\v@lmax\advance\v@lmax-\dpmn@\mili@u%
      \advance\v@leur#1pt%
    \fi}
\ctr@ld@f\def\c@ssin#1#2#3{\expandafter\ifx\csname COS@\number#3\endcsname\relax\c@lCS{#3pt}%
    \expandafter\xdef\csname COS@\number#3\endcsname{\repdecn@mb\result@t}%
    \expandafter\xdef\csname SIN@\number#3\endcsname{\repdecn@mb\result@@t}\fi%
    \edef#1{\csname COS@\number#3\endcsname}\edef#2{\csname SIN@\number#3\endcsname}}
\ctr@ld@f\def\c@lCS#1{{\mili@u=#1\p@rtent=\@ne%
    \relax\ifdim\mili@u<\z@\red@ng<-\else\red@ng>+\fi\f@ctech=\p@rtent%
    \relax\ifdim\mili@u<\z@\mili@u=-\mili@u\f@ctech=-\f@ctech\fi\c@@lCS}}
\ctr@ld@f\def\c@@lCS{\v@lmin=\mili@u\c@rre=-\mili@u\advance\c@rre\DemiPI@deg\v@lmax=\c@rre%
    \mili@u\@@lxxi\mili@u\divide\mili@u\@@mmmmlxviii%
    \edef\v@larg{\repdecn@mb{\mili@u}}\mili@u=-\v@larg\mili@u%
    \edef\v@lmxde{\repdecn@mb{\mili@u}}%
    \c@rre\@@lxxi\c@rre\divide\c@rre\@@mmmmlxviii%
    \edef\v@largC{\repdecn@mb{\c@rre}}\c@rre=-\v@largC\c@rre%
    \edef\v@lmxdeC{\repdecn@mb{\c@rre}}%
    \fctc@s\mili@u\v@lmin\global\result@t\p@rtent\v@leur%
    \let\t@mp=\v@larg\let\v@larg=\v@largC\let\v@largC=\t@mp%
    \let\t@mp=\v@lmxde\let\v@lmxde=\v@lmxdeC\let\v@lmxdeC=\t@mp%
    \fctc@s\c@rre\v@lmax\global\result@@t\f@ctech\v@leur}
\ctr@ld@f\def\fctc@s#1#2{\v@leur=#1\relax\ifdim#2<\@lxxxv\p@\cosser@h\else\sinser@t\fi}
\ctr@ld@f\def\cosser@h{\advance\v@leur\@lvi\p@\divide\v@leur\@lvi%
    \v@leur=\v@lmxde\v@leur\advance\v@leur\@xxx\p@%
    \v@leur=\v@lmxde\v@leur\advance\v@leur\@ccclx\p@%
    \v@leur=\v@lmxde\v@leur\advance\v@leur\@dccxx\p@\divide\v@leur\@dccxx}
\ctr@ld@f\def\sinser@t{\v@leur=\v@lmxdeC\p@\advance\v@leur\@vi\p@%
    \v@leur=\v@largC\v@leur\divide\v@leur\@vi}
\ctr@ld@f\def\red@ng#1#2{\relax\ifdim\mili@u#1#2\DemiPI@deg\advance\mili@u#2-\PI@deg%
    \p@rtent=-\p@rtent\red@ng#1#2\fi}
\ctr@ld@f\def\pr@c@lCS#1#2#3{\ctr@lcsn@m{COS@\number#3 }%
    \expandafter\xdef\csname COS@\number#3\endcsname{#1}%
    \expandafter\xdef\csname SIN@\number#3\endcsname{#2}}
\pr@c@lCS{1}{0}{0}
\pr@c@lCS{0.7071}{0.7071}{45}\pr@c@lCS{0.7071}{-0.7071}{-45}
\pr@c@lCS{0}{1}{90}          \pr@c@lCS{0}{-1}{-90}
\pr@c@lCS{-1}{0}{180}        \pr@c@lCS{-1}{0}{-180}
\pr@c@lCS{0}{-1}{270}        \pr@c@lCS{0}{1}{-270}
\ctr@ld@f\def\invers@#1#2{{\v@leur=#2\maxim@m{\v@lmax}{-\v@leur}{\v@leur}%
    \f@ctech=\@ne\m@inv@rs%
    \multiply\v@leur\f@ctech\edef\v@lv@leur{\repdecn@mb{\v@leur}}%
    \p@rtentiere{\p@rtent}{\v@leur}\v@lmin=\p@\divide\v@lmin\p@rtent%
    \inv@rs@\multiply\v@lmax\f@ctech\global\result@t=\v@lmax}#1=\result@t}
\ctr@ld@f\def\m@inv@rs{\ifdim\v@lmax<\p@\multiply\v@lmax\t@n\multiply\f@ctech\t@n\m@inv@rs\fi}
\ctr@ld@f\def\inv@rs@{\v@lmax=-\v@lmin\v@lmax=\v@lv@leur\v@lmax%
    \advance\v@lmax\tw@ pt\v@lmax=\repdecn@mb{\v@lmin}\v@lmax%
    \delt@=\v@lmax\advance\delt@-\v@lmin\ifdim\delt@<\z@\delt@=-\delt@\fi%
    \ifdim\delt@>\epsil@n\v@lmin=\v@lmax\inv@rs@\fi}
\ctr@ld@f\def\minim@m#1#2#3{\relax\ifdim#2<#3#1=#2\else#1=#3\fi}
\ctr@ld@f\def\maxim@m#1#2#3{\relax\ifdim#2>#3#1=#2\else#1=#3\fi}
\ctr@ld@f\def\p@rtentiere#1#2{#1=#2\divide#1by65536 }
\ctr@ld@f\def\r@undint#1#2{{\v@leur=#2\divide\v@leur\t@n\p@rtentiere{\p@rtent}{\v@leur}%
    \v@leur=\p@rtent pt\global\result@t=\t@n\v@leur}#1=\result@t}
\ctr@ld@f\def\sqrt@#1#2{{\v@leur=#2%
    \minim@m{\v@lmin}{\p@}{\v@leur}\maxim@m{\v@lmax}{\p@}{\v@leur}%
    \f@ctech=\@ne\m@sqrt@\sqrt@@%
    \mili@u=\v@lmin\advance\mili@u\v@lmax\divide\mili@u\tw@\multiply\mili@u\f@ctech%
    \global\result@t=\mili@u}#1=\result@t}
\ctr@ld@f\def\m@sqrt@{\ifdim\v@leur>\dcq@\divide\v@leur\c@nt\v@lmax=\v@leur%
    \multiply\f@ctech\t@n\m@sqrt@\fi}
\ctr@ld@f\def\sqrt@@{\mili@u=\v@lmin\advance\mili@u\v@lmax\divide\mili@u\tw@%
    \c@rre=\repdecn@mb{\mili@u}\mili@u%
    \ifdim\c@rre<\v@leur\v@lmin=\mili@u\else\v@lmax=\mili@u\fi%
    \delt@=\v@lmax\advance\delt@-\v@lmin\ifdim\delt@>\epsil@n\sqrt@@\fi}
\ctr@ld@f\def\extrairelepremi@r#1\de#2{\expandafter\lepremi@r#2@#1#2}
\ctr@ld@f\def\lepremi@r#1,#2@#3#4{\def#3{#1}\def#4{#2}\ignorespaces}
\ctr@ld@f\def\@cfor#1:=#2\do#3{%
  \edef\@fortemp{#2}%
  \ifx\@fortemp\empty\else\@cforloop#2,\@nil,\@nil\@@#1{#3}\fi}
\ctr@ln@m\@nextwhile
\ctr@ld@f\def\@cforloop#1,#2\@@#3#4{%
  \def#3{#1}%
  \ifx#3\Fig@nnil\let\@nextwhile=\Fig@fornoop\else#4\relax\let\@nextwhile=\@cforloop\fi%
  \@nextwhile#2\@@#3{#4}}

\ctr@ld@f\def\@ecfor#1:=#2\do#3{%
  \def\@@cfor{\@cfor#1:=}%
  \edef\@@@cfor{#2}%
  \expandafter\@@cfor\@@@cfor\do{#3}}
\ctr@ld@f\def\Fig@nnil{\@nil}
\ctr@ld@f\def\Fig@fornoop#1\@@#2#3{}
\ctr@ln@m\list@@rg
\ctr@ld@f\def\trtlis@rg#1#2{\def\list@@rg{#1}%
    \@ecfor\p@rv@l:=\list@@rg\do{\expandafter#2\p@rv@l|}}
\ctr@ld@f\def\trtlis@rgtok#1{\let@xte={}\let\n@xt\addt@t@xt\addt@t@xt #1}
\ctr@ln@m\M@cro
\ctr@ln@m\n@xt
\ctr@ld@f\def\addt@t@xt#1{\if#1|\let\n@xt\relax\else%
    \if#1,\expandafter\M@cro\the\let@xte|\let@xte={}%
    \else\let@xte=\expandafter{\the\let@xte #1}\fi\fi\n@xt}
\ctr@ln@w{newbox}\b@xvisu
\ctr@ln@w{newtoks}\let@xte
\ctr@ln@w{newif}\ifitis@K
\ctr@ln@w{newcount}\s@mme
\ctr@ln@w{newcount}\l@mbd@un \ctr@ln@w{newcount}\l@mbd@de
\ctr@ln@w{newcount}\superc@ntr@l\superc@ntr@l=\@ne        
\ctr@ln@w{newcount}\typec@ntr@l\typec@ntr@l=\superc@ntr@l 
\ctr@ln@w{newdimen}\v@lX  \ctr@ln@w{newdimen}\v@lY  \ctr@ln@w{newdimen}\v@lZ
\ctr@ln@w{newdimen}\v@lXa \ctr@ln@w{newdimen}\v@lYa \ctr@ln@w{newdimen}\v@lZa
\ctr@ln@w{newdimen}\unit@\unit@=\p@ 
\ctr@ld@f\def\unit@util{pt}
\ctr@ld@f\def\ptT@ptps{0.996264 }
\ctr@ld@f\def\ptpsT@pt{1.00375 }
\ctr@ld@f\def\ptT@unit@{1} 
\ctr@ld@f\def\setunit@#1{\def\unit@util{#1}\setunit@@#1:\invers@{\result@t}{\unit@}%
    \edef\ptT@unit@{\repdecn@mb\result@t}}
\ctr@ld@f\def\setunit@@#1#2:{\ifcat#1a\unit@=\@ne#1#2\else\unit@=#1#2\fi}
\ctr@ld@f\def\d@fm@cdim#1#2{{\v@leur=#2\v@leur=\ptT@unit@\v@leur\xdef#1{\repdecn@mb\v@leur}}}
\ctr@ln@w{newif}\ifBdingB@x\BdingB@xtrue
\ctr@ln@w{newdimen}\c@@rdXmin \ctr@ln@w{newdimen}\c@@rdYmin  
\ctr@ln@w{newdimen}\c@@rdXmax \ctr@ln@w{newdimen}\c@@rdYmax
\ctr@ld@f\def\b@undb@x#1#2{\ifBdingB@x%
    \relax\ifdim#1<\c@@rdXmin\global\c@@rdXmin=#1\fi%
    \relax\ifdim#2<\c@@rdYmin\global\c@@rdYmin=#2\fi%
    \relax\ifdim#1>\c@@rdXmax\global\c@@rdXmax=#1\fi%
    \relax\ifdim#2>\c@@rdYmax\global\c@@rdYmax=#2\fi\fi}
\ctr@ld@f\def\b@undb@xP#1{{\Figg@tXY{#1}\b@undb@x{\v@lX}{\v@lY}}}
\ctr@ld@f\def\ellBB@x#1;#2,#3(#4,#5,#6){{\s@uvc@ntr@l\et@tellBB@x%
    \setc@ntr@l{2}\figptell-2::#1;#2,#3(#4,#6)\b@undb@xP{-2}%
    \figptell-2::#1;#2,#3(#5,#6)\b@undb@xP{-2}%
    \c@ssin{\C@}{\S@}{#6}\v@lmin=\C@ pt\v@lmax=\S@ pt%
    \mili@u=#3\v@lmin\delt@=#2\v@lmax\arct@n\v@leur(\delt@,\mili@u)%
    \mili@u=-#3\v@lmax\delt@=#2\v@lmin\arct@n\c@rre(\delt@,\mili@u)%
    \v@leur=\rdT@deg\v@leur\advance\v@leur-\DePI@deg%
    \c@rre=\rdT@deg\c@rre\advance\c@rre-\DePI@deg%
    \v@lmin=#4pt\v@lmax=#5pt%
    \loop\ifdim\v@leur<\v@lmax\ifdim\v@leur>\v@lmin%
    \edef\@ngle{\repdecn@mb\v@leur}\figptell-2::#1;#2,#3(\@ngle,#6)%
    \b@undb@xP{-2}\fi\advance\v@leur\PI@deg\repeat%
    \loop\ifdim\c@rre<\v@lmax\ifdim\c@rre>\v@lmin%
    \edef\@ngle{\repdecn@mb\c@rre}\figptell-2::#1;#2,#3(\@ngle,#6)%
    \b@undb@xP{-2}\fi\advance\c@rre\PI@deg\repeat%
    \resetc@ntr@l\et@tellBB@x}\ignorespaces}
\ctr@ld@f\def\initb@undb@x{\c@@rdXmin=\maxdimen\c@@rdYmin=\maxdimen%
    \c@@rdXmax=-\maxdimen\c@@rdYmax=-\maxdimen}
\ctr@ld@f\def\c@ntr@lnum#1{%
    \relax\ifnum\typec@ntr@l=\@ne%
    \ifnum#1<\z@%
    \immediate\write16{*** Forbidden point number (#1). Abort.}\end\fi\fi%
    \set@bjc@de{#1}}
\ctr@ln@m\objc@de
\ctr@ld@f\def\set@bjc@de#1{\edef\objc@de{@BJ\ifnum#1<\z@ M\romannumeral-#1\else\romannumeral#1\fi}}
\s@mme=\m@ne\loop\ifnum\s@mme>-19
  \set@bjc@de{\s@mme}\ctr@lcsn@m\objc@de\ctr@lcsn@m{\objc@de T}
\advance\s@mme\m@ne\repeat
\s@mme=\@ne\loop\ifnum\s@mme<6
  \set@bjc@de{\s@mme}\ctr@lcsn@m\objc@de\ctr@lcsn@m{\objc@de T}
\advance\s@mme\@ne\repeat
\ctr@ld@f\def\setc@ntr@l#1{\ifnum\superc@ntr@l>#1\typec@ntr@l=\superc@ntr@l%
    \else\typec@ntr@l=#1\fi}
\ctr@ld@f\def\resetc@ntr@l#1{\global\superc@ntr@l=#1\setc@ntr@l{#1}}
\ctr@ld@f\def\s@uvc@ntr@l#1{\edef#1{\the\superc@ntr@l}}
\ctr@ln@m\c@lproscal
\ctr@ld@f\def\c@lproscalDD#1[#2,#3]{{\Figg@tXY{#2}%
    \edef\Xu@{\repdecn@mb{\v@lX}}\edef\Yu@{\repdecn@mb{\v@lY}}\Figg@tXY{#3}%
    \global\result@t=\Xu@\v@lX\global\advance\result@t\Yu@\v@lY}#1=\result@t}
\ctr@ld@f\def\c@lproscalTD#1[#2,#3]{{\Figg@tXY{#2}\edef\Xu@{\repdecn@mb{\v@lX}}%
    \edef\Yu@{\repdecn@mb{\v@lY}}\edef\Zu@{\repdecn@mb{\v@lZ}}%
    \Figg@tXY{#3}\global\result@t=\Xu@\v@lX\global\advance\result@t\Yu@\v@lY%
    \global\advance\result@t\Zu@\v@lZ}#1=\result@t}
\ctr@ld@f\def\c@lprovec#1{%
    \det@rmC\v@lZa(\v@lX,\v@lY,\v@lmin,\v@lmax)%
    \det@rmC\v@lXa(\v@lY,\v@lZ,\v@lmax,\v@leur)%
    \det@rmC\v@lYa(\v@lZ,\v@lX,\v@leur,\v@lmin)%
    \Figv@ctCreg#1(\v@lXa,\v@lYa,\v@lZa)}
\ctr@ld@f\def\det@rm#1[#2,#3]{{\Figg@tXY{#2}\Figg@tXYa{#3}%
    \delt@=\repdecn@mb{\v@lX}\v@lYa\advance\delt@-\repdecn@mb{\v@lY}\v@lXa%
    \global\result@t=\delt@}#1=\result@t}
\ctr@ld@f\def\det@rmC#1(#2,#3,#4,#5){{\global\result@t=\repdecn@mb{#2}#5%
    \global\advance\result@t-\repdecn@mb{#3}#4}#1=\result@t}
\ctr@ld@f\def\getredf@ctDD#1(#2,#3){{\maxim@m{\v@lXa}{-#2}{#2}\maxim@m{\v@lYa}{-#3}{#3}%
    \maxim@m{\v@lXa}{\v@lXa}{\v@lYa}
    \ifdim\v@lXa>\@xci pt\divide\v@lXa\@xci%
    \p@rtentiere{\p@rtent}{\v@lXa}\advance\p@rtent\@ne\else\p@rtent=\@ne\fi%
    \global\result@tent=\p@rtent}#1=\result@tent\ignorespaces}
\ctr@ld@f\def\getredf@ctTD#1(#2,#3,#4){{\maxim@m{\v@lXa}{-#2}{#2}\maxim@m{\v@lYa}{-#3}{#3}%
    \maxim@m{\v@lZa}{-#4}{#4}\maxim@m{\v@lXa}{\v@lXa}{\v@lYa}%
    \maxim@m{\v@lXa}{\v@lXa}{\v@lZa}
    \ifdim\v@lXa>\@lxxiv pt\divide\v@lXa\@lxxiv%
    \p@rtentiere{\p@rtent}{\v@lXa}\advance\p@rtent\@ne\else\p@rtent=\@ne\fi%
    \global\result@tent=\p@rtent}#1=\result@tent\ignorespaces}
\ctr@ln@m\getredf@ctB
\ctr@ld@f\def\getredf@ctBDD#1{\getredf@ctDD#1(\v@lX,\v@lY)}
\ctr@ld@f\def\getredf@ctBTD#1{\getredf@ctTD#1(\v@lX,\v@lY,\v@lZ)}
\ctr@ld@f\def\FigptintercircB@zDD#1:#2:#3,#4[#5,#6,#7,#8]{{\s@uvc@ntr@l\et@tfigptintercircB@zDD%
    \setc@ntr@l{2}\figvectPDD-1[#5,#8]\Figg@tXY{-1}\getredf@ctDD\f@ctech(\v@lX,\v@lY)%
    \mili@u=#4\unit@\divide\mili@u\f@ctech\c@rre=\repdecn@mb{\mili@u}\mili@u%
    \figptBezierDD-5::#3[#5,#6,#7,#8]%
    \v@lmin=#3\p@\v@lmax=\v@lmin\advance\v@lmax0.1\p@%
    \loop\edef\T@{\repdecn@mb{\v@lmax}}\figptBezierDD-2::\T@[#5,#6,#7,#8]%
    \figvectPDD-1[-5,-2]\n@rmeucCDD{\delt@}{-1}\ifdim\delt@<\c@rre\v@lmin=\v@lmax%
    \advance\v@lmax0.1\p@\repeat%
    \loop\mili@u=\v@lmin\advance\mili@u\v@lmax%
    \divide\mili@u\tw@\edef\T@{\repdecn@mb{\mili@u}}\figptBezierDD-2::\T@[#5,#6,#7,#8]%
    \figvectPDD-1[-5,-2]\n@rmeucCDD{\delt@}{-1}\ifdim\delt@>\c@rre\v@lmax=\mili@u%
    \else\v@lmin=\mili@u\fi\v@leur=\v@lmax\advance\v@leur-\v@lmin%
    \ifdim\v@leur>\epsil@n\repeat\figptcopyDD#1:#2/-2/%
    \resetc@ntr@l\et@tfigptintercircB@zDD}\ignorespaces}
\ctr@ln@m\figptinterlines
\ctr@ld@f\def\inters@cDD#1:#2[#3,#4;#5,#6]{{\s@uvc@ntr@l\et@tinters@cDD%
    \setc@ntr@l{2}\vecunit@{-1}{#4}\vecunit@{-2}{#6}%
    \Figg@tXY{-1}\setc@ntr@l{1}\Figg@tXYa{#3}%
    \edef\A@{\repdecn@mb{\v@lX}}\edef\B@{\repdecn@mb{\v@lY}}%
    \v@lmin=\B@\v@lXa\advance\v@lmin-\A@\v@lYa%
    \Figg@tXYa{#5}\setc@ntr@l{2}\Figg@tXY{-2}%
    \edef\C@{\repdecn@mb{\v@lX}}\edef\D@{\repdecn@mb{\v@lY}}%
    \v@lmax=\D@\v@lXa\advance\v@lmax-\C@\v@lYa%
    \delt@=\A@\v@lY\advance\delt@-\B@\v@lX%
    \invers@{\v@leur}{\delt@}\edef\v@ldelta{\repdecn@mb{\v@leur}}%
    \v@lXa=\A@\v@lmax\advance\v@lXa-\C@\v@lmin%
    \v@lYa=\B@\v@lmax\advance\v@lYa-\D@\v@lmin%
    \v@lXa=\v@ldelta\v@lXa\v@lYa=\v@ldelta\v@lYa%
    \setc@ntr@l{1}\Figp@intregDD#1:{#2}(\v@lXa,\v@lYa)%
    \resetc@ntr@l\et@tinters@cDD}\ignorespaces}
\ctr@ld@f\def\inters@cTD#1:#2[#3,#4;#5,#6]{{\s@uvc@ntr@l\et@tinters@cTD%
    \setc@ntr@l{2}\figvectNVTD-1[#4,#6]\figvectNVTD-2[#6,-1]\figvectPTD-1[#3,#5]%
    \r@pPSTD\v@leur[-2,-1,#4]\edef\v@lcoef{\repdecn@mb{\v@leur}}%
    \figpttraTD#1:{#2}=#3/\v@lcoef,#4/\resetc@ntr@l\et@tinters@cTD}\ignorespaces}
\ctr@ld@f\def\r@pPSTD#1[#2,#3,#4]{{\Figg@tXY{#2}\edef\Xu@{\repdecn@mb{\v@lX}}%
    \edef\Yu@{\repdecn@mb{\v@lY}}\edef\Zu@{\repdecn@mb{\v@lZ}}%
    \Figg@tXY{#3}\v@lmin=\Xu@\v@lX\advance\v@lmin\Yu@\v@lY\advance\v@lmin\Zu@\v@lZ%
    \Figg@tXY{#4}\v@lmax=\Xu@\v@lX\advance\v@lmax\Yu@\v@lY\advance\v@lmax\Zu@\v@lZ%
    \invers@{\v@leur}{\v@lmax}\global\result@t=\repdecn@mb{\v@leur}\v@lmin}%
    #1=\result@t}
\ctr@ln@m\n@rminf
\ctr@ld@f\def\n@rminfDD#1#2{{\Figg@tXY{#2}\maxim@m{\v@lX}{\v@lX}{-\v@lX}%
    \maxim@m{\v@lY}{\v@lY}{-\v@lY}\maxim@m{\global\result@t}{\v@lX}{\v@lY}}%
    #1=\result@t}
\ctr@ld@f\def\n@rminfTD#1#2{{\Figg@tXY{#2}\maxim@m{\v@lX}{\v@lX}{-\v@lX}%
    \maxim@m{\v@lY}{\v@lY}{-\v@lY}\maxim@m{\v@lZ}{\v@lZ}{-\v@lZ}%
    \maxim@m{\v@lX}{\v@lX}{\v@lY}\maxim@m{\global\result@t}{\v@lX}{\v@lZ}}%
    #1=\result@t}
\ctr@ln@m\n@rmeucC
\ctr@ld@f\def\n@rmeucCDD#1#2{\Figg@tXY{#2}\divide\v@lX\f@ctech\divide\v@lY\f@ctech%
    #1=\repdecn@mb{\v@lX}\v@lX\v@lX=\repdecn@mb{\v@lY}\v@lY\advance#1\v@lX}
\ctr@ld@f\def\n@rmeucCTD#1#2{\Figg@tXY{#2}%
    \divide\v@lX\f@ctech\divide\v@lY\f@ctech\divide\v@lZ\f@ctech%
    #1=\repdecn@mb{\v@lX}\v@lX\v@lX=\repdecn@mb{\v@lY}\v@lY\advance#1\v@lX%
    \v@lX=\repdecn@mb{\v@lZ}\v@lZ\advance#1\v@lX}
\ctr@ln@m\n@rmeucSV
\ctr@ld@f\def\n@rmeucSVDD#1#2{{\Figg@tXY{#2}%
    \v@lXa=\repdecn@mb{\v@lX}\v@lX\v@lYa=\repdecn@mb{\v@lY}\v@lY%
    \advance\v@lXa\v@lYa\sqrt@{\global\result@t}{\v@lXa}}#1=\result@t}
\ctr@ld@f\def\n@rmeucSVTD#1#2{{\Figg@tXY{#2}\v@lXa=\repdecn@mb{\v@lX}\v@lX%
    \v@lYa=\repdecn@mb{\v@lY}\v@lY\v@lZa=\repdecn@mb{\v@lZ}\v@lZ%
    \advance\v@lXa\v@lYa\advance\v@lXa\v@lZa\sqrt@{\global\result@t}{\v@lXa}}#1=\result@t}
\ctr@ln@m\n@rmeuc
\ctr@ld@f\def\n@rmeucDD#1#2{{\Figg@tXY{#2}\getredf@ctDD\f@ctech(\v@lX,\v@lY)%
    \divide\v@lX\f@ctech\divide\v@lY\f@ctech%
    \v@lXa=\repdecn@mb{\v@lX}\v@lX\v@lYa=\repdecn@mb{\v@lY}\v@lY%
    \advance\v@lXa\v@lYa\sqrt@{\global\result@t}{\v@lXa}%
    \global\multiply\result@t\f@ctech}#1=\result@t}
\ctr@ld@f\def\n@rmeucTD#1#2{{\Figg@tXY{#2}\getredf@ctTD\f@ctech(\v@lX,\v@lY,\v@lZ)%
    \divide\v@lX\f@ctech\divide\v@lY\f@ctech\divide\v@lZ\f@ctech%
    \v@lXa=\repdecn@mb{\v@lX}\v@lX%
    \v@lYa=\repdecn@mb{\v@lY}\v@lY\v@lZa=\repdecn@mb{\v@lZ}\v@lZ%
    \advance\v@lXa\v@lYa\advance\v@lXa\v@lZa\sqrt@{\global\result@t}{\v@lXa}%
    \global\multiply\result@t\f@ctech}#1=\result@t}
\ctr@ln@m\vecunit@
\ctr@ld@f\def\vecunit@DD#1#2{{\Figg@tXY{#2}\getredf@ctDD\f@ctech(\v@lX,\v@lY)%
    \divide\v@lX\f@ctech\divide\v@lY\f@ctech%
    \Figv@ctCreg#1(\v@lX,\v@lY)\n@rmeucSV{\v@lYa}{#1}%
    \invers@{\v@lXa}{\v@lYa}\edef\v@lv@lXa{\repdecn@mb{\v@lXa}}%
    \v@lX=\v@lv@lXa\v@lX\v@lY=\v@lv@lXa\v@lY%
    \Figv@ctCreg#1(\v@lX,\v@lY)\multiply\v@lYa\f@ctech\global\result@t=\v@lYa}}
\ctr@ld@f\def\vecunit@TD#1#2{{\Figg@tXY{#2}\getredf@ctTD\f@ctech(\v@lX,\v@lY,\v@lZ)%
    \divide\v@lX\f@ctech\divide\v@lY\f@ctech\divide\v@lZ\f@ctech%
    \Figv@ctCreg#1(\v@lX,\v@lY,\v@lZ)\n@rmeucSV{\v@lYa}{#1}%
    \invers@{\v@lXa}{\v@lYa}\edef\v@lv@lXa{\repdecn@mb{\v@lXa}}%
    \v@lX=\v@lv@lXa\v@lX\v@lY=\v@lv@lXa\v@lY\v@lZ=\v@lv@lXa\v@lZ%
    \Figv@ctCreg#1(\v@lX,\v@lY,\v@lZ)\multiply\v@lYa\f@ctech\global\result@t=\v@lYa}}
\ctr@ld@f\def\vecunitC@TD[#1,#2]{\Figg@tXYa{#1}\Figg@tXY{#2}%
    \advance\v@lX-\v@lXa\advance\v@lY-\v@lYa\advance\v@lZ-\v@lZa\c@lvecunitTD}
\ctr@ld@f\def\vecunitCV@TD#1{\Figg@tXY{#1}\c@lvecunitTD}
\ctr@ld@f\def\c@lvecunitTD{\getredf@ctTD\f@ctech(\v@lX,\v@lY,\v@lZ)%
    \divide\v@lX\f@ctech\divide\v@lY\f@ctech\divide\v@lZ\f@ctech%
    \v@lXa=\repdecn@mb{\v@lX}\v@lX%
    \v@lYa=\repdecn@mb{\v@lY}\v@lY\v@lZa=\repdecn@mb{\v@lZ}\v@lZ%
    \advance\v@lXa\v@lYa\advance\v@lXa\v@lZa\sqrt@{\v@lYa}{\v@lXa}%
    \invers@{\v@lXa}{\v@lYa}\edef\v@lv@lXa{\repdecn@mb{\v@lXa}}%
    \v@lX=\v@lv@lXa\v@lX\v@lY=\v@lv@lXa\v@lY\v@lZ=\v@lv@lXa\v@lZ}
\ctr@ln@m\figgetangle
\ctr@ld@f\def\figgetangleDD#1[#2,#3,#4]{\ifGR@cri{\s@uvc@ntr@l\et@tfiggetangleDD\setc@ntr@l{2}%
    \figvectPDD-1[#2,#3]\figvectPDD-2[#2,#4]\vecunit@{-1}{-1}%
    \c@lproscalDD\delt@[-2,-1]\figvectNVDD-1[-1]\c@lproscalDD\v@leur[-2,-1]%
    \arct@n\v@lmax(\delt@,\v@leur)\v@lmax=\rdT@deg\v@lmax%
    \ifdim\v@lmax<\z@\advance\v@lmax\DePI@deg\fi\xdef#1{\repdecn@mb{\v@lmax}}%
    \resetc@ntr@l\et@tfiggetangleDD}\ignorespaces\fi}
\ctr@ld@f\def\figgetangleTD#1[#2,#3,#4,#5]{\ifGR@cri{\s@uvc@ntr@l\et@tfiggetangleTD\setc@ntr@l{2}%
    \figvectPTD-1[#2,#3]\figvectPTD-2[#2,#5]\figvectNVTD-3[-1,-2]%
    \figvectPTD-2[#2,#4]\figvectNVTD-4[-3,-1]%
    \vecunit@{-1}{-1}\c@lproscalTD\delt@[-2,-1]\c@lproscalTD\v@leur[-2,-4]%
    \arct@n\v@lmax(\delt@,\v@leur)\v@lmax=\rdT@deg\v@lmax%
    \ifdim\v@lmax<\z@\advance\v@lmax\DePI@deg\fi\xdef#1{\repdecn@mb{\v@lmax}}%
    \resetc@ntr@l\et@tfiggetangleTD}\ignorespaces\fi}    
\ctr@ld@f\def\figgetdist#1[#2,#3]{\ifGR@cri{\s@uvc@ntr@l\et@tfiggetdist\setc@ntr@l{2}%
    \figvectP-1[#2,#3]\n@rmeuc{\v@lX}{-1}\v@lX=\ptT@unit@\v@lX\xdef#1{\repdecn@mb{\v@lX}}%
    \resetc@ntr@l\et@tfiggetdist}\ignorespaces\fi}
\ctr@ld@f\def\figget#1=#2[#3]{\keln@mun#1|%
    \def\n@mref{a}\ifx\l@debut\n@mref\figgetangle#2[#3]\else
    \def\n@mref{d}\ifx\l@debut\n@mref\figgetdist#2[#3]\else
    \W@rnmeskwd{figget}{#1}\fi\fi\ignorespaces}
\ctr@ld@f\def\Figg@tT#1{\c@ntr@lnum{#1}%
    {\expandafter\expandafter\expandafter\extr@ctT\csname\objc@de\endcsname:%
     \ifnum\B@@ltxt=\z@\ptn@me{#1}\else\csname\objc@de T\endcsname\fi}}
\ctr@ld@f\def\extr@ctT#1,#2,#3/#4:{\def\B@@ltxt{#3}}
\ctr@ld@f\def\Figg@tXY#1{\c@ntr@lnum{#1}%
    \expandafter\expandafter\expandafter\extr@ctC\csname\objc@de\endcsname:}
\ctr@ln@m\extr@ctC
\ctr@ld@f\def\extr@ctCDD#1/#2,#3,#4:{\v@lX=#2\v@lY=#3}
\ctr@ld@f\def\extr@ctCTD#1/#2,#3,#4:{\v@lX=#2\v@lY=#3\v@lZ=#4}
\ctr@ld@f\def\Figg@tXYa#1{\c@ntr@lnum{#1}%
    \expandafter\expandafter\expandafter\extr@ctCa\csname\objc@de\endcsname:}
\ctr@ln@m\extr@ctCa
\ctr@ld@f\def\extr@ctCaDD#1/#2,#3,#4:{\v@lXa=#2\v@lYa=#3}
\ctr@ld@f\def\extr@ctCaTD#1/#2,#3,#4:{\v@lXa=#2\v@lYa=#3\v@lZa=#4}
\ctr@ln@m\t@xt@
\ctr@ld@f\def\figinit#1{\t@stc@tcodech@nge\initpr@lim\Figinit@#1,:\initpss@ttings\ignorespaces}
\ctr@ld@f\def\Figinit@#1,#2:{\setunit@{#1}\def\t@xt@{#2}\ifx\t@xt@\empty\else\Figinit@@#2:\fi}
\ctr@ld@f\def\Figinit@@#1#2:{\if#12 \else\Figs@tproj{#1}\initTD@\fi}
\ctr@ln@w{newif}\ifTr@isDim
\ctr@ld@f\def\UnD@fined{UNDEFINED}
\ctr@ln@m\@utoFN
\ctr@ln@m\@utoFInDone
\ctr@ln@m\disob@unit
\ctr@ld@f\def\initpr@lim{\initb@undb@x\figsetmark{}\figsetptname{$A_{##1}$}\def\Sc@leFact{1}%
    \initDD@\figsetroundcoord{yes}\GR@critrue\expandafter\setupd@te\D@FTupdate:%
    \edef\disob@unit{\UnD@fined}\edef\t@rgetpt{\UnD@fined}\gdef\@utoFInDone{1}\gdef\@utoFN{0}}
\ctr@ld@f\def\initDD@{\Tr@isDimfalse%
    \ifPDFm@ke%
     \let\Ps@rcerc=\Ps@rcercBz%
     \let\Ps@rell=\Ps@rellBz%
    \fi
    \let\c@lDCUn=\c@lDCUnDD%
    \let\c@lDCDeux=\c@lDCDeuxDD%
    \let\c@ldefproj=\relax%
    \let\c@lproscal=\c@lproscalDD%
    \let\c@lprojSP=\relax%
    \let\extr@ctC=\extr@ctCDD%
    \let\extr@ctCa=\extr@ctCaDD%
    \let\extr@ctCF=\extr@ctCFDD%
    \let\Figp@intreg=\Figp@intregDD%
    \let\Figpts@xes=\Figpts@xesDD%
    \let\getredf@ctB=\getredf@ctBDD%
    \let\n@rmeucSV=\n@rmeucSVDD\let\n@rmeuc=\n@rmeucDD\let\n@rmeucC\n@rmeucCDD\let\n@rminf=\n@rminfDD%
    \let\pr@dMatV=\pr@dMatVDD%
    \let\Q@@xes=\Q@@xesDD%
    \let\vecunit@=\vecunit@DD%
    \let\figcoord=\figcoordDD%
    \let\figgetangle=\figgetangleDD%
    \let\figpt=\figptDD%
    \let\figptBezier=\figptBezierDD%
    \let\figptbary=\figptbaryDD%
    \let\figptcirc=\figptcircDD%
    \let\figptcircumcenter=\figptcircumcenterDD%
    \let\figptcopy=\figptcopyDD%
    \let\figptcurvcenter=\figptcurvcenterDD%
    \let\figptell=\figptellDD%
    \let\figptendnormal=\figptendnormalDD%
    \let\figptinterlineplane=\figptinterlineplaneDD%
    \let\figptinterlines=\inters@cDD%
    \let\figptorthocenter=\figptorthocenterDD%
    \let\figptorthoprojline=\figptorthoprojlineDD%
    \let\figptorthoprojplane=\figptorthoprojplaneDD%
    \let\figptrot=\figptrotDD%
    \let\figptscontrol=\figptscontrolDD%
    \let\figptsintercirc=\figptsintercircDD%
    \let\figptsinterlinell=\figptsinterlinellDD%
    \let\figptsorthoprojline=\figptsorthoprojlineDD%
    \let\figptorthoprojplane=\figptorthoprojplaneDD%
    \let\figptsrot=\figptsrotDD%
    \let\figptssym=\figptssymDD%
    \let\figptstra=\figptstraDD%
    \let\figptsym=\figptsymDD%
    \let\figpttraC=\figpttraCDD%
    \let\figpttra=\figpttraDD%
    \let\figptvisilimSL=\figptvisilimSLDD%
    \let\figsetobdist=\figsetobdistDD%
    \let\figsettarget=\figsettargetDD%
    \let\figsetview=\figsetviewDD%
    \let\figvectDBezier=\figvectDBezierDD%
    \let\figvectN=\figvectNDD%
    \let\figvectNV=\figvectNVDD%
    \let\figvectP=\figvectPDD%
    \let\figvectU=\figvectUDD%
    \let\figdrawarccircP=\Q@arccircPDD%
    \let\figdrawarccirc=\Q@arccircDD%
    \let\figdrawarcell=\Q@arcellDD%
    \let\figdrawarcellPA=\Q@arcellPADD%
    \let\figdrawarrowBezier=\Q@arrowBezierDD%
    \let\figdrawarrowcircP=\Q@arrowcircPDD%
    \let\figdrawarrowcirc=\Q@arrowcircDD%
    \let\figdrawarrowhead=\Q@arrowheadDD%
    \let\figdrawarrow=\Q@arrowDD%
    \let\figdrawBezier=\Q@BezierDD%
    \let\figdrawcirc=\Q@circDD%
    \let\figdrawcurve=\Q@curveDD%
    \let\figdrawnormal=\Q@normalDD%
    }
\ctr@ld@f\def\initTD@{\Tr@isDimtrue\initb@undb@xTD\newt@rgetptfalse\newdis@bfalse%
    \let\c@lDCUn=\c@lDCUnTD%
    \let\c@lDCDeux=\c@lDCDeuxTD%
    \let\c@ldefproj=\c@ldefprojTD%
    \let\c@lproscal=\c@lproscalTD%
    \let\extr@ctC=\extr@ctCTD%
    \let\extr@ctCa=\extr@ctCaTD%
    \let\extr@ctCF=\extr@ctCFTD%
    \let\Figp@intreg=\Figp@intregTD%
    \let\Figpts@xes=\Figpts@xesTD%
    \let\getredf@ctB=\getredf@ctBTD%
    \let\n@rmeucSV=\n@rmeucSVTD\let\n@rmeuc=\n@rmeucTD\let\n@rmeucC\n@rmeucCTD\let\n@rminf=\n@rminfTD%
    \let\pr@dMatV=\pr@dMatVTD%
    \let\Q@@xes=\Q@@xesTD%
    \let\vecunit@=\vecunit@TD%
    \let\figcoord=\figcoordTD%
    \let\figgetangle=\figgetangleTD%
    \let\figpt=\figptTD%
    \let\figptBezier=\figptBezierTD%
    \let\figptbary=\figptbaryTD%
    \let\figptcirc=\figptcircTD%
    \let\figptcircumcenter=\figptcircumcenterTD%
    \let\figptcopy=\figptcopyTD%
    \let\figptcurvcenter=\figptcurvcenterTD%
    \let\figptinterlineplane=\figptinterlineplaneTD%
    \let\figptinterlines=\inters@cTD%
    \let\figptorthocenter=\figptorthocenterTD%
    \let\figptorthoprojline=\figptorthoprojlineTD%
    \let\figptorthoprojplane=\figptorthoprojplaneTD%
    \let\figptrot=\figptrotTD%
    \let\figptscontrol=\figptscontrolTD%
    \let\figptsintercirc=\figptsintercircTD%
    \let\figptsorthoprojline=\figptsorthoprojlineTD%
    \let\figptsorthoprojplane=\figptsorthoprojplaneTD%
    \let\figptsrot=\figptsrotTD%
    \let\figptssym=\figptssymTD%
    \let\figptstra=\figptstraTD%
    \let\figptsym=\figptsymTD%
    \let\figpttraC=\figpttraCTD%
    \let\figpttra=\figpttraTD%
    \let\figptvisilimSL=\figptvisilimSLTD%
    \let\figsetobdist=\figsetobdistTD%
    \let\figsettarget=\figsettargetTD%
    \let\figsetview=\figsetviewTD%
    \let\figvectDBezier=\figvectDBezierTD%
    \let\figvectN=\figvectNTD%
    \let\figvectNV=\figvectNVTD%
    \let\figvectP=\figvectPTD%
    \let\figvectU=\figvectUTD%
    \let\figdrawarccircP=\Q@arccircPTD%
    \let\figdrawarccirc=\Q@arccircTD%
    \let\figdrawarcell=\Q@arcellTD%
    \let\figdrawarcellPA=\Q@arcellPATD%
    \let\figdrawarrowBezier=\Q@arrowBezierTD%
    \let\figdrawarrowcircP=\Q@arrowcircPTD%
    \let\figdrawarrowcirc=\Q@arrowcircTD%
    \let\figdrawarrowhead=\Q@arrowheadTD%
    \let\figdrawarrow=\Q@arrowTD%
    \let\figdrawBezier=\Q@BezierTD%
    \let\figdrawcirc=\Q@circTD%
    \let\figdrawcurve=\Q@curveTD%
    }
\ctr@ld@f\def\un@v@ilable#1{\immediate\write16{*** The macro #1 is not available in the current context.}}
\ctr@ld@f\def\figinsert#1{{\def\t@xt@{#1}\relax%
    \ifx\t@xt@\empty\ifnum\@utoFInDone>\z@\Figinsert@\DefGIfilen@me,:\fi%
    \else\expandafter\FiginsertNu@#1 :\fi}\ignorespaces}
\ctr@ld@f\def\FiginsertNu@#1 #2:{\def\t@xt@{#1}\relax\ifx\t@xt@\empty\def\t@xt@{#2}%
    \ifx\t@xt@\empty\ifnum\@utoFInDone>\z@\Figinsert@\DefGIfilen@me,:\fi%
    \else\FiginsertNu@#2:\fi\else\expandafter\FiginsertNd@#1 #2:\fi}
\ctr@ld@f\def\FiginsertNd@#1#2:{\ifcat#1a\Figinsert@#1#2,:\else%
    \ifnum\@utoFInDone>\z@\Figinsert@\DefGIfilen@me,#1#2,:\fi\fi}
\ctr@ln@m\Sc@leFact
\ctr@ld@f\def\Figinsert@#1,#2:{\def\t@xt@{#2}\ifx\t@xt@\empty\xdef\Sc@leFact{1}\else%
    \X@rgdeux@#2\xdef\Sc@leFact{\@rgdeux}\fi%
    \Figdisc@rdLTS{#1}{\t@xt@}\@psfgetbb{\t@xt@}%
    \v@lX=\@psfllx\p@\v@lX=\ptpsT@pt\v@lX\v@lX=\Sc@leFact\v@lX%
    \v@lY=\@psflly\p@\v@lY=\ptpsT@pt\v@lY\v@lY=\Sc@leFact\v@lY%
    \b@undb@x{\v@lX}{\v@lY}%
    \v@lX=\@psfurx\p@\v@lX=\ptpsT@pt\v@lX\v@lX=\Sc@leFact\v@lX%
    \v@lY=\@psfury\p@\v@lY=\ptpsT@pt\v@lY\v@lY=\Sc@leFact\v@lY%
    \b@undb@x{\v@lX}{\v@lY}%
    \ifPDFm@ke\Figinclud@PDF{\t@xt@}{\Sc@leFact}\else%
    \v@lX=\c@nt pt\v@lX=\Sc@leFact\v@lX\edef\F@ct{\repdecn@mb{\v@lX}}%
    \ifx\TeXturesonMacOSltX\special{postscriptfile #1 vscale=\F@ct\space hscale=\F@ct}%
    \else\includegraphics{#1}\fi\fi%
    \message{[\t@xt@]}\ignorespaces}
\ctr@ld@f\def\Figdisc@rdLTS#1#2{\expandafter\Figdisc@rdLTS@#1 :#2}
\ctr@ld@f\def\Figdisc@rdLTS@#1 #2:#3{\def#3{#1}\relax\ifx#3\empty\expandafter\Figdisc@rdLTS@#2:#3\fi}
\ctr@ld@f\def\figinsertE#1{\FiginsertE@#1,:\ignorespaces}
\ctr@ld@f\def\FiginsertE@#1,#2:{{\def\t@xt@{#2}\ifx\t@xt@\empty\xdef\Sc@leFact{1}\else%
    \X@rgdeux@#2\xdef\Sc@leFact{\@rgdeux}\fi%
    \Figdisc@rdLTS{#1}{\t@xt@}\pdfximage{\t@xt@}%
    \setbox\Gb@x=\hbox{\pdfrefximage\pdflastximage}%
    \v@lX=\z@\v@lY=-\Sc@leFact\dp\Gb@x\b@undb@x{\v@lX}{\v@lY}%
    \advance\v@lX\Sc@leFact\wd\Gb@x\advance\v@lY\Sc@leFact\dp\Gb@x%
    \advance\v@lY\Sc@leFact\ht\Gb@x\b@undb@x{\v@lX}{\v@lY}%
    \v@lX=\Sc@leFact\wd\Gb@x\pdfximage width \v@lX {\t@xt@}%
    \rlap{\pdfrefximage\pdflastximage}\message{[\t@xt@]}}\ignorespaces}
\ctr@ld@f\def\X@rgdeux@#1,{\edef\@rgdeux{#1}}
\ctr@ln@m\figpt
\ctr@ld@f\def\figptDD#1:#2(#3,#4){\ifGR@cri\c@ntr@lnum{#1}%
    {\v@lX=#3\unit@\v@lY=#4\unit@\Fig@dmpt{#2}{\z@}}\ignorespaces\fi}
\ctr@ld@f\def\Fig@dmpt#1#2{\def\t@xt@{#1}\ifx\t@xt@\empty\def\B@@ltxt{\z@}%
    \else\expandafter\gdef\csname\objc@de T\endcsname{#1}\def\B@@ltxt{\@ne}\fi%
    \expandafter\xdef\csname\objc@de\endcsname{\ifitis@vect@r\C@dCl@svect%
    \else\C@dCl@spt\fi,\z@,\B@@ltxt/\the\v@lX,\the\v@lY,#2}}
\ctr@ld@f\def\C@dCl@spt{P}
\ctr@ld@f\def\C@dCl@svect{V}
\ctr@ln@m\c@@rdYZ
\ctr@ln@m\c@@rdY
\ctr@ld@f\def\figptTD#1:#2(#3,#4){\ifGR@cri\c@ntr@lnum{#1}%
    \def\c@@rdYZ{#4,0,0}\extrairelepremi@r\c@@rdY\de\c@@rdYZ%
    \extrairelepremi@r\c@@rdZ\de\c@@rdYZ%
    {\v@lX=#3\unit@\v@lY=\c@@rdY\unit@\v@lZ=\c@@rdZ\unit@\Fig@dmpt{#2}{\the\v@lZ}%
    \b@undb@xTD{\v@lX}{\v@lY}{\v@lZ}}\ignorespaces\fi}
\ctr@ln@m\Figp@intreg
\ctr@ld@f\def\Figp@intregDD#1:#2(#3,#4){\c@ntr@lnum{#1}%
    {\result@t=#4\v@lX=#3\v@lY=\result@t\Fig@dmpt{#2}{\z@}}\ignorespaces}
\ctr@ld@f\def\Figp@intregTD#1:#2(#3,#4){\c@ntr@lnum{#1}%
    \def\c@@rdYZ{#4,\z@,\z@}\extrairelepremi@r\c@@rdY\de\c@@rdYZ%
    \extrairelepremi@r\c@@rdZ\de\c@@rdYZ%
    {\v@lX=#3\v@lY=\c@@rdY\v@lZ=\c@@rdZ\Fig@dmpt{#2}{\the\v@lZ}%
    \b@undb@xTD{\v@lX}{\v@lY}{\v@lZ}}\ignorespaces}
\ctr@ln@m\figptBezier
\ctr@ld@f\def\figptBezierDD#1:#2:#3[#4,#5,#6,#7]{\ifGR@cri{\s@uvc@ntr@l\et@tfigptBezierDD%
    \FigptBezier@#3[#4,#5,#6,#7]\Figp@intregDD#1:{#2}(\v@lX,\v@lY)%
    \resetc@ntr@l\et@tfigptBezierDD}\ignorespaces\fi}
\ctr@ld@f\def\figptBezierTD#1:#2:#3[#4,#5,#6,#7]{\ifGR@cri{\s@uvc@ntr@l\et@tfigptBezierTD%
    \FigptBezier@#3[#4,#5,#6,#7]\Figp@intregTD#1:{#2}(\v@lX,\v@lY,\v@lZ)%
    \resetc@ntr@l\et@tfigptBezierTD}\ignorespaces\fi}
\ctr@ld@f\def\FigptBezier@#1[#2,#3,#4,#5]{\setc@ntr@l{2}%
    \edef\T@{#1}\v@leur=\p@\advance\v@leur-#1pt\edef\UNmT@{\repdecn@mb{\v@leur}}%
    \figptcopy-4:/#2/\figptcopy-3:/#3/\figptcopy-2:/#4/\figptcopy-1:/#5/%
    \l@mbd@un=-4 \l@mbd@de=-\thr@@\p@rtent=\m@ne\c@lDecast%
    \l@mbd@un=-4 \l@mbd@de=-\thr@@\p@rtent=-\tw@\c@lDecast%
    \l@mbd@un=-4 \l@mbd@de=-\thr@@\p@rtent=-\thr@@\c@lDecast\Figg@tXY{-4}}
\ctr@ln@m\c@lDCUn
\ctr@ld@f\def\c@lDCUnDD#1#2{\Figg@tXY{#1}\v@lX=\UNmT@\v@lX\v@lY=\UNmT@\v@lY%
    \Figg@tXYa{#2}\advance\v@lX\T@\v@lXa\advance\v@lY\T@\v@lYa%
    \Figp@intregDD#1:(\v@lX,\v@lY)}
\ctr@ld@f\def\c@lDCUnTD#1#2{\Figg@tXY{#1}\v@lX=\UNmT@\v@lX\v@lY=\UNmT@\v@lY\v@lZ=\UNmT@\v@lZ%
    \Figg@tXYa{#2}\advance\v@lX\T@\v@lXa\advance\v@lY\T@\v@lYa\advance\v@lZ\T@\v@lZa%
    \Figp@intregTD#1:(\v@lX,\v@lY,\v@lZ)}
\ctr@ld@f\def\c@lDecast{\relax\ifnum\l@mbd@un<\p@rtent\c@lDCUn{\l@mbd@un}{\l@mbd@de}%
    \advance\l@mbd@un\@ne\advance\l@mbd@de\@ne\c@lDecast\fi}
\ctr@ld@f\def\figptmap#1:#2=#3/#4/#5/{\ifGR@cri{\s@uvc@ntr@l\et@tfigptmap%
    \setc@ntr@l{2}\figvectP-1[#4,#3]\Figg@tXY{-1}%
    \pr@dMatV/#5/\figpttra#1:{#2}=#4/1,-1/%
    \resetc@ntr@l\et@tfigptmap}\ignorespaces\fi}
\ctr@ln@m\pr@dMatV
\ctr@ld@f\def\pr@dMatVDD/#1,#2;#3,#4/{\v@lXa=#1\v@lX\advance\v@lXa#2\v@lY%
    \v@lYa=#3\v@lX\advance\v@lYa#4\v@lY\Figv@ctCreg-1(\v@lXa,\v@lYa)}
\ctr@ld@f\def\pr@dMatVTD/#1,#2,#3;#4,#5,#6;#7,#8,#9/{%
    \v@lXa=#1\v@lX\advance\v@lXa#2\v@lY\advance\v@lXa#3\v@lZ%
    \v@lYa=#4\v@lX\advance\v@lYa#5\v@lY\advance\v@lYa#6\v@lZ%
    \v@lZa=#7\v@lX\advance\v@lZa#8\v@lY\advance\v@lZa#9\v@lZ%
    \Figv@ctCreg-1(\v@lXa,\v@lYa,\v@lZa)}
\ctr@ln@m\figptbary
\ctr@ld@f\def\figptbaryDD#1:#2[#3;#4]{\ifGR@cri{\edef\list@num{#3}\extrairelepremi@r\p@int\de\list@num%
    \s@mme=\z@\@ecfor\c@ef:=#4\do{\advance\s@mme\c@ef}%
    \edef\listec@ef{#4,0}\extrairelepremi@r\c@ef\de\listec@ef%
    \Figg@tXY{\p@int}\divide\v@lX\s@mme\divide\v@lY\s@mme%
    \multiply\v@lX\c@ef\multiply\v@lY\c@ef%
    \@ecfor\p@int:=\list@num\do{\extrairelepremi@r\c@ef\de\listec@ef%
           \Figg@tXYa{\p@int}\divide\v@lXa\s@mme\divide\v@lYa\s@mme%
           \multiply\v@lXa\c@ef\multiply\v@lYa\c@ef%
           \advance\v@lX\v@lXa\advance\v@lY\v@lYa}%
    \Figp@intregDD#1:{#2}(\v@lX,\v@lY)}\ignorespaces\fi}
\ctr@ld@f\def\figptbaryTD#1:#2[#3;#4]{\ifGR@cri{\edef\list@num{#3}\extrairelepremi@r\p@int\de\list@num%
    \s@mme=\z@\@ecfor\c@ef:=#4\do{\advance\s@mme\c@ef}%
    \edef\listec@ef{#4,0}\extrairelepremi@r\c@ef\de\listec@ef%
    \Figg@tXY{\p@int}\divide\v@lX\s@mme\divide\v@lY\s@mme\divide\v@lZ\s@mme%
    \multiply\v@lX\c@ef\multiply\v@lY\c@ef\multiply\v@lZ\c@ef%
    \@ecfor\p@int:=\list@num\do{\extrairelepremi@r\c@ef\de\listec@ef%
           \Figg@tXYa{\p@int}\divide\v@lXa\s@mme\divide\v@lYa\s@mme\divide\v@lZa\s@mme%
           \multiply\v@lXa\c@ef\multiply\v@lYa\c@ef\multiply\v@lZa\c@ef%
           \advance\v@lX\v@lXa\advance\v@lY\v@lYa\advance\v@lZ\v@lZa}%
    \Figp@intregTD#1:{#2}(\v@lX,\v@lY,\v@lZ)}\ignorespaces\fi}
\ctr@ld@f\def\figptbaryR#1:#2[#3;#4]{\ifGR@cri{%
    \v@leur=\z@\@ecfor\c@ef:=#4\do{\maxim@m{\v@lmax}{\c@ef pt}{-\c@ef pt}%
    \ifdim\v@lmax>\v@leur\v@leur=\v@lmax\fi}%
    \ifdim\v@leur<\p@\f@ctech=\@M\else\ifdim\v@leur<\t@n\p@\f@ctech=\@m\else%
    \ifdim\v@leur<\c@nt\p@\f@ctech=\c@nt\else\ifdim\v@leur<\@m\p@\f@ctech=\t@n\else%
    \f@ctech=\@ne\fi\fi\fi\fi%
    \def\listec@ef{0}%
    \@ecfor\c@ef:=#4\do{\sc@lec@nvRI{\c@ef pt}\edef\listec@ef{\listec@ef,\the\s@mme}}%
    \extrairelepremi@r\c@ef\de\listec@ef\figptbary#1:#2[#3;\listec@ef]}\ignorespaces\fi}
\ctr@ld@f\def\sc@lec@nvRI#1{\v@leur=#1\p@rtentiere{\s@mme}{\v@leur}\advance\v@leur-\s@mme\p@%
    \multiply\v@leur\f@ctech\p@rtentiere{\p@rtent}{\v@leur}%
    \multiply\s@mme\f@ctech\advance\s@mme\p@rtent}
\ctr@ln@m\figptcirc
\ctr@ld@f\def\figptcircDD#1:#2:#3;#4(#5){\ifGR@cri{\s@uvc@ntr@l\et@tfigptcircDD%
    \c@lptellDD#1:{#2}:#3;#4,#4(#5)\resetc@ntr@l\et@tfigptcircDD}\ignorespaces\fi}
\ctr@ld@f\def\figptcircTD#1:#2:#3,#4,#5;#6(#7){\ifGR@cri{\s@uvc@ntr@l\et@tfigptcircTD%
    \setc@ntr@l{2}\c@lExtAxes#3,#4,#5(#6)\figptellP#1:{#2}:#3,-4,-5(#7)%
    \resetc@ntr@l\et@tfigptcircTD}\ignorespaces\fi}
\ctr@ln@m\figptcircumcenter
\ctr@ld@f\def\figptcircumcenterDD#1:#2[#3,#4,#5]{\ifGR@cri{\s@uvc@ntr@l\et@tfigptcircumcenterDD%
    \setc@ntr@l{2}\figvectNDD-5[#3,#4]\figptbaryDD-3:[#3,#4;1,1]%
                  \figvectNDD-6[#4,#5]\figptbaryDD-4:[#4,#5;1,1]%
    \resetc@ntr@l{2}\inters@cDD#1:{#2}[-3,-5;-4,-6]%
    \resetc@ntr@l\et@tfigptcircumcenterDD}\ignorespaces\fi}
\ctr@ld@f\def\figptcircumcenterTD#1:#2[#3,#4,#5]{\ifGR@cri{\s@uvc@ntr@l\et@tfigptcircumcenterTD%
    \setc@ntr@l{2}\figvectNTD-1[#3,#4,#5]%
    \figvectPTD-3[#3,#4]\figvectNVTD-5[-1,-3]\figptbaryTD-3:[#3,#4;1,1]%
    \figvectPTD-4[#4,#5]\figvectNVTD-6[-1,-4]\figptbaryTD-4:[#4,#5;1,1]%
    \resetc@ntr@l{2}\inters@cTD#1:{#2}[-3,-5;-4,-6]%
    \resetc@ntr@l\et@tfigptcircumcenterTD}\ignorespaces\fi}
\ctr@ln@m\figptcopy
\ctr@ld@f\def\figptcopyDD#1:#2/#3/{\ifGR@cri{\Figg@tXY{#3}%
    \Figp@intregDD#1:{#2}(\v@lX,\v@lY)}\ignorespaces\fi}
\ctr@ld@f\def\figptcopyTD#1:#2/#3/{\ifGR@cri{\Figg@tXY{#3}%
    \Figp@intregTD#1:{#2}(\v@lX,\v@lY,\v@lZ)}\ignorespaces\fi}
\ctr@ln@m\figptcurvcenter
\ctr@ld@f\def\figptcurvcenterDD#1:#2:#3[#4,#5,#6,#7]{\ifGR@cri{\s@uvc@ntr@l\et@tfigptcurvcenterDD%
    \setc@ntr@l{2}\c@lcurvradDD#3[#4,#5,#6,#7]\edef\Sprim@{\repdecn@mb{\result@t}}%
    \figptBezierDD-1::#3[#4,#5,#6,#7]\figpttraDD#1:{#2}=-1/\Sprim@,-5/%
    \resetc@ntr@l\et@tfigptcurvcenterDD}\ignorespaces\fi}
\ctr@ld@f\def\figptcurvcenterTD#1:#2:#3[#4,#5,#6,#7]{\ifGR@cri{\s@uvc@ntr@l\et@tfigptcurvcenterTD%
    \setc@ntr@l{2}\figvectDBezierTD -5:1,#3[#4,#5,#6,#7]%
    \figvectDBezierTD -6:2,#3[#4,#5,#6,#7]\vecunit@TD{-5}{-5}%
    \edef\Sprim@{\repdecn@mb{\result@t}}\figvectNVTD-1[-6,-5]%
    \figvectNVTD-5[-5,-1]\c@lproscalTD\v@leur[-6,-5]%
    \invers@{\v@leur}{\v@leur}\v@leur=\Sprim@\v@leur\v@leur=\Sprim@\v@leur%
    \figptBezierTD-1::#3[#4,#5,#6,#7]\edef\Sprim@{\repdecn@mb{\v@leur}}%
    \figpttraTD#1:{#2}=-1/\Sprim@,-5/\resetc@ntr@l\et@tfigptcurvcenterTD}\ignorespaces\fi}
\ctr@ld@f\def\c@lcurvradDD#1[#2,#3,#4,#5]{{\figvectDBezierDD -5:1,#1[#2,#3,#4,#5]%
    \figvectDBezierDD -6:2,#1[#2,#3,#4,#5]\vecunit@DD{-5}{-5}%
    \edef\Sprim@{\repdecn@mb{\result@t}}\figvectNVDD-5[-5]\c@lproscalDD\v@leur[-6,-5]%
    \invers@{\v@leur}{\v@leur}\v@leur=\Sprim@\v@leur\v@leur=\Sprim@\v@leur%
    \global\result@t=\v@leur}}
\ctr@ln@m\figptell
\ctr@ld@f\def\figptellDD#1:#2:#3;#4,#5(#6,#7){\ifGR@cri{\s@uvc@ntr@l\et@tfigptell%
    \c@lptellDD#1::#3;#4,#5(#6)\figptrotDD#1:{#2}=#1/#3,#7/%
    \resetc@ntr@l\et@tfigptell}\ignorespaces\fi}
\ctr@ld@f\def\c@lptellDD#1:#2:#3;#4,#5(#6){\c@ssin{\C@}{\S@}{#6}\v@lmin=\C@ pt\v@lmax=\S@ pt%
    \v@lmin=#4\v@lmin\v@lmax=#5\v@lmax%
    \edef\Xc@mp{\repdecn@mb{\v@lmin}}\edef\Yc@mp{\repdecn@mb{\v@lmax}}%
    \setc@ntr@l{2}\figvectC-1(\Xc@mp,\Yc@mp)\figpttraDD#1:{#2}=#3/1,-1/}
\ctr@ld@f\def\figptellP#1:#2:#3,#4,#5(#6){\ifGR@cri{\s@uvc@ntr@l\et@tfigptellP%
    \setc@ntr@l{2}\figvectP-1[#3,#4]\figvectP-2[#3,#5]%
    \v@leur=#6pt\c@lptellP{#3}{-1}{-2}\figptcopy#1:{#2}/-3/%
    \resetc@ntr@l\et@tfigptellP}\ignorespaces\fi}
\ctr@ln@m\@ngle
\ctr@ld@f\def\c@lptellP#1#2#3{\edef\@ngle{\repdecn@mb\v@leur}\c@ssin{\C@}{\S@}{\@ngle}%
    \figpttra-3:=#1/\C@,#2/\figpttra-3:=-3/\S@,#3/}
\ctr@ln@m\figptendnormal
\ctr@ld@f\def\figptendnormalDD#1:#2:#3,#4[#5,#6]{\ifGR@cri{\s@uvc@ntr@l\et@tfigptendnormal%
    \Figg@tXYa{#5}\Figg@tXY{#6}%
    \advance\v@lX-\v@lXa\advance\v@lY-\v@lYa%
    \setc@ntr@l{2}\Figv@ctCreg-1(\v@lX,\v@lY)\vecunit@{-1}{-1}\Figg@tXY{-1}%
    \delt@=#3\unit@\maxim@m{\delt@}{\delt@}{-\delt@}\edef\l@ngueur{\repdecn@mb{\delt@}}%
    \v@lX=\l@ngueur\v@lX\v@lY=\l@ngueur\v@lY%
    \delt@=\p@\advance\delt@-#4pt\edef\l@ngueur{\repdecn@mb{\delt@}}%
    \figptbaryR-1:[#5,#6;#4,\l@ngueur]\Figg@tXYa{-1}%
    \advance\v@lXa\v@lY\advance\v@lYa-\v@lX%
    \setc@ntr@l{1}\Figp@intregDD#1:{#2}(\v@lXa,\v@lYa)\resetc@ntr@l\et@tfigptendnormal}%
    \ignorespaces\fi}
\ctr@ld@f\def\figptexcenter#1:#2[#3,#4,#5]{\ifGR@cri{\let@xte={-}
    \Figptexinsc@nter#1:#2[#3,#4,#5]}\ignorespaces\fi}
\ctr@ld@f\def\figptincenter#1:#2[#3,#4,#5]{\ifGR@cri{\let@xte={}
    \Figptexinsc@nter#1:#2[#3,#4,#5]}\ignorespaces\fi}
\ctr@ld@f
\ctr@ld@f\def\Figptexinsc@nter#1:#2[#3,#4,#5]{%
    \figgetdist\LA@[#4,#5]\figgetdist\LB@[#3,#5]\figgetdist\LC@[#3,#4]%
    \figptbaryR#1:{#2}[#3,#4,#5;\the\let@xte\LA@,\LB@,\LC@]}
\ctr@ln@m\figptinterlineplane
\ctr@ld@f\def\figptinterlineplaneDD{\un@v@ilable{figptinterlineplane}}
\ctr@ld@f\def\figptinterlineplaneTD#1:#2[#3,#4;#5,#6]{\ifGR@cri{\s@uvc@ntr@l\et@tfigptinterlineplane%
    \setc@ntr@l{2}\figvectPTD-1[#3,#5]\vecunit@TD{-2}{#6}%
    \r@pPSTD\v@leur[-2,-1,#4]\edef\v@lcoef{\repdecn@mb{\v@leur}}%
    \figpttraTD#1:{#2}=#3/\v@lcoef,#4/\resetc@ntr@l\et@tfigptinterlineplane}\ignorespaces\fi}
\ctr@ln@m\figptorthocenter
\ctr@ld@f\def\figptorthocenterDD#1:#2[#3,#4,#5]{\ifGR@cri{\s@uvc@ntr@l\et@tfigptorthocenterDD%
    \setc@ntr@l{2}\figvectNDD-3[#3,#4]\figvectNDD-4[#4,#5]%
    \resetc@ntr@l{2}\inters@cDD#1:{#2}[#5,-3;#3,-4]%
    \resetc@ntr@l\et@tfigptorthocenterDD}\ignorespaces\fi}
\ctr@ld@f\def\figptorthocenterTD#1:#2[#3,#4,#5]{\ifGR@cri{\s@uvc@ntr@l\et@tfigptorthocenterTD%
    \setc@ntr@l{2}\figvectNTD-1[#3,#4,#5]%
    \figvectPTD-2[#3,#4]\figvectNVTD-3[-1,-2]%
    \figvectPTD-2[#4,#5]\figvectNVTD-4[-1,-2]%
    \resetc@ntr@l{2}\inters@cTD#1:{#2}[#5,-3;#3,-4]%
    \resetc@ntr@l\et@tfigptorthocenterTD}\ignorespaces\fi}
\ctr@ln@m\figptorthoprojline
\ctr@ld@f\def\figptorthoprojlineDD#1:#2=#3/#4,#5/{\ifGR@cri{\s@uvc@ntr@l\et@tfigptorthoprojlineDD%
    \setc@ntr@l{2}\figvectPDD-3[#4,#5]\figvectNVDD-4[-3]\resetc@ntr@l{2}%
    \inters@cDD#1:{#2}[#3,-4;#4,-3]\resetc@ntr@l\et@tfigptorthoprojlineDD}\ignorespaces\fi}
\ctr@ld@f\def\figptorthoprojlineTD#1:#2=#3/#4,#5/{\ifGR@cri{\s@uvc@ntr@l\et@tfigptorthoprojlineTD%
    \setc@ntr@l{2}\figvectPTD-1[#4,#3]\figvectPTD-2[#4,#5]\vecunit@TD{-2}{-2}%
    \c@lproscalTD\v@leur[-1,-2]\edef\v@lcoef{\repdecn@mb{\v@leur}}%
    \figpttraTD#1:{#2}=#4/\v@lcoef,-2/\resetc@ntr@l\et@tfigptorthoprojlineTD}\ignorespaces\fi}
\ctr@ln@m\figptorthoprojplane
\ctr@ld@f\def\figptorthoprojplaneDD{\un@v@ilable{figptorthoprojplane}}
\ctr@ld@f\def\figptorthoprojplaneTD#1:#2=#3/#4,#5/{\ifGR@cri{\s@uvc@ntr@l\et@tfigptorthoprojplane%
    \setc@ntr@l{2}\figvectPTD-1[#3,#4]\vecunit@TD{-2}{#5}%
    \c@lproscalTD\v@leur[-1,-2]\edef\v@lcoef{\repdecn@mb{\v@leur}}%
    \figpttraTD#1:{#2}=#3/\v@lcoef,-2/\resetc@ntr@l\et@tfigptorthoprojplane}\ignorespaces\fi}
\ctr@ld@f\def\figpthom#1:#2=#3/#4,#5/{\ifGR@cri{\s@uvc@ntr@l\et@tfigpthom%
    \setc@ntr@l{2}\figvectP-1[#4,#3]\figpttra#1:{#2}=#4/#5,-1/%
    \resetc@ntr@l\et@tfigpthom}\ignorespaces\fi}
\ctr@ld@f\def\figptinv#1:#2=#3/#4,#5/{\ifGR@cri{\s@uvc@ntr@l\et@tfigptinv%
    \setc@ntr@l{2}\figvectP-1[#4,#3]\Figg@tXY{-1}%
    \getredf@ctB\f@ctech\n@rmeucC{\delt@}{-1}%
    \delt@=\ptT@unit@\delt@\delt@=\ptT@unit@\delt@%
    \invers@{\delt@}{\delt@}\multiply\f@ctech\f@ctech\divide\delt@\f@ctech%
    \delt@=#5\delt@\edef\v@lcoef{\repdecn@mb{\delt@}}\figpttra#1:{#2}=#4/\v@lcoef,-1/%
    \resetc@ntr@l\et@tfigptinv}\ignorespaces\fi}
\ctr@ln@m\figptrot
\ctr@ld@f\def\figptrotDD#1:#2=#3/#4,#5/{\ifGR@cri{\s@uvc@ntr@l\et@tfigptrotDD%
    \c@ssin{\C@}{\S@}{#5}\setc@ntr@l{2}\figvectPDD-1[#4,#3]\Figg@tXY{-1}%
    \v@lXa=\C@\v@lX\advance\v@lXa-\S@\v@lY%
    \v@lYa=\S@\v@lX\advance\v@lYa\C@\v@lY%
    \Figv@ctCreg-1(\v@lXa,\v@lYa)\figpttraDD#1:{#2}=#4/1,-1/%
    \resetc@ntr@l\et@tfigptrotDD}\ignorespaces\fi}
\ctr@ld@f\def\figptrotTD#1:#2=#3/#4,#5,#6/{\ifGR@cri{\s@uvc@ntr@l\et@tfigptrotTD%
    \c@ssin{\C@}{\S@}{#5}%
    \setc@ntr@l{2}\figptorthoprojplaneTD-3:=#4/#3,#6/\figvectPTD-2[-3,#3]%
    \n@rmeucTD\v@leur{-2}\ifdim\v@leur<\Cepsil@n\Figg@tXYa{#3}\else%
    \edef\v@lcoef{\repdecn@mb{\v@leur}}\figvectNVTD-1[#6,-2]%
    \Figg@tXYa{-1}\v@lXa=\v@lcoef\v@lXa\v@lYa=\v@lcoef\v@lYa\v@lZa=\v@lcoef\v@lZa%
    \v@lXa=\S@\v@lXa\v@lYa=\S@\v@lYa\v@lZa=\S@\v@lZa\Figg@tXY{-2}%
    \advance\v@lXa\C@\v@lX\advance\v@lYa\C@\v@lY\advance\v@lZa\C@\v@lZ%
    \Figg@tXY{-3}\advance\v@lXa\v@lX\advance\v@lYa\v@lY\advance\v@lZa\v@lZ\fi%
    \Figp@intregTD#1:{#2}(\v@lXa,\v@lYa,\v@lZa)\resetc@ntr@l\et@tfigptrotTD}\ignorespaces\fi}
\ctr@ln@m\figptsym
\ctr@ld@f\def\figptsymDD#1:#2=#3/#4,#5/{\ifGR@cri{\s@uvc@ntr@l\et@tfigptsymDD%
    \resetc@ntr@l{2}\figptorthoprojlineDD-5:=#3/#4,#5/\figvectPDD-2[#3,-5]%
    \figpttraDD#1:{#2}=#3/2,-2/\resetc@ntr@l\et@tfigptsymDD}\ignorespaces\fi}
\ctr@ld@f\def\figptsymTD#1:#2=#3/#4,#5/{\ifGR@cri{\s@uvc@ntr@l\et@tfigptsymTD%
    \resetc@ntr@l{2}\figptorthoprojplaneTD-3:=#3/#4,#5/\figvectPTD-2[#3,-3]%
    \figpttraTD#1:{#2}=#3/2,-2/\resetc@ntr@l\et@tfigptsymTD}\ignorespaces\fi}
\ctr@ln@m\figpttra
\ctr@ld@f\def\figpttraDD#1:#2=#3/#4,#5/{\ifGR@cri{\Figg@tXYa{#5}\v@lXa=#4\v@lXa\v@lYa=#4\v@lYa%
    \Figg@tXY{#3}\advance\v@lX\v@lXa\advance\v@lY\v@lYa%
    \Figp@intregDD#1:{#2}(\v@lX,\v@lY)}\ignorespaces\fi}
\ctr@ld@f\def\figpttraTD#1:#2=#3/#4,#5/{\ifGR@cri{\Figg@tXYa{#5}\v@lXa=#4\v@lXa\v@lYa=#4\v@lYa%
    \v@lZa=#4\v@lZa\Figg@tXY{#3}\advance\v@lX\v@lXa\advance\v@lY\v@lYa%
    \advance\v@lZ\v@lZa\Figp@intregTD#1:{#2}(\v@lX,\v@lY,\v@lZ)}\ignorespaces\fi}
\ctr@ln@m\figpttraC
\ctr@ld@f\def\figpttraCDD#1:#2=#3/#4,#5/{\ifGR@cri{\v@lXa=#4\unit@\v@lYa=#5\unit@%
    \Figg@tXY{#3}\advance\v@lX\v@lXa\advance\v@lY\v@lYa%
    \Figp@intregDD#1:{#2}(\v@lX,\v@lY)}\ignorespaces\fi}
\ctr@ld@f\def\figpttraCTD#1:#2=#3/#4,#5,#6/{\ifGR@cri{\v@lXa=#4\unit@\v@lYa=#5\unit@\v@lZa=#6\unit@%
    \Figg@tXY{#3}\advance\v@lX\v@lXa\advance\v@lY\v@lYa\advance\v@lZ\v@lZa%
    \Figp@intregTD#1:{#2}(\v@lX,\v@lY,\v@lZ)}\ignorespaces\fi}
\ctr@ld@f\def\figptsaxes#1:#2(#3){\ifGR@cri{\an@lys@xes#3,:\ifx\t@xt@\empty%
    \ifTr@isDim\Figpts@xes#1:#2(0,#3,0,#3,0,#3)\else\Figpts@xes#1:#2(0,#3,0,#3)\fi%
    \else\Figpts@xes#1:#2(#3)\fi}\ignorespaces\fi}
\ctr@ln@m\Figpts@xes
\ctr@ld@f\def\Figpts@xesDD#1:#2(#3,#4,#5,#6){%
    \s@mme=#1\figpttraC\the\s@mme:$x$=#2/#4,0/%
    \advance\s@mme\@ne\figpttraC\the\s@mme:$y$=#2/0,#6/}
\ctr@ld@f\def\Figpts@xesTD#1:#2(#3,#4,#5,#6,#7,#8){%
    \s@mme=#1\figpttraC\the\s@mme:$x$=#2/#4,0,0/%
    \advance\s@mme\@ne\figpttraC\the\s@mme:$y$=#2/0,#6,0/%
    \advance\s@mme\@ne\figpttraC\the\s@mme:$z$=#2/0,0,#8/}
\ctr@ld@f\def\figptsmap#1=#2/#3/#4/{\ifGR@cri{\s@uvc@ntr@l\et@tfigptsmap%
    \setc@ntr@l{2}\def\list@num{#2}\s@mme=#1%
    \@ecfor\p@int:=\list@num\do{\figvectP-1[#3,\p@int]\Figg@tXY{-1}%
    \pr@dMatV/#4/\figpttra\the\s@mme:=#3/1,-1/\advance\s@mme\@ne}%
    \resetc@ntr@l\et@tfigptsmap}\ignorespaces\fi}
\ctr@ln@m\figptscontrol
\ctr@ld@f\def\figptscontrolDD#1[#2,#3,#4,#5]{\ifGR@cri{\s@uvc@ntr@l\et@tfigptscontrolDD\setc@ntr@l{2}%
    \v@lX=\z@\v@lY=\z@\Figtr@nptDD{-5}{#2}\Figtr@nptDD{2}{#5}%
    \divide\v@lX\@vi\divide\v@lY\@vi%
    \Figtr@nptDD{3}{#3}\Figtr@nptDD{-1.5}{#4}\Figp@intregDD-1:(\v@lX,\v@lY)%
    \v@lX=\z@\v@lY=\z@\Figtr@nptDD{2}{#2}\Figtr@nptDD{-5}{#5}%
    \divide\v@lX\@vi\divide\v@lY\@vi\Figtr@nptDD{-1.5}{#3}\Figtr@nptDD{3}{#4}%
    \s@mme=#1\advance\s@mme\@ne\Figp@intregDD\the\s@mme:(\v@lX,\v@lY)%
    \figptcopyDD#1:/-1/\resetc@ntr@l\et@tfigptscontrolDD}\ignorespaces\fi}
\ctr@ld@f\def\figptscontrolTD#1[#2,#3,#4,#5]{\ifGR@cri{\s@uvc@ntr@l\et@tfigptscontrolTD\setc@ntr@l{2}%
    \v@lX=\z@\v@lY=\z@\v@lZ=\z@\Figtr@nptTD{-5}{#2}\Figtr@nptTD{2}{#5}%
    \divide\v@lX\@vi\divide\v@lY\@vi\divide\v@lZ\@vi%
    \Figtr@nptTD{3}{#3}\Figtr@nptTD{-1.5}{#4}\Figp@intregTD-1:(\v@lX,\v@lY,\v@lZ)%
    \v@lX=\z@\v@lY=\z@\v@lZ=\z@\Figtr@nptTD{2}{#2}\Figtr@nptTD{-5}{#5}%
    \divide\v@lX\@vi\divide\v@lY\@vi\divide\v@lZ\@vi\Figtr@nptTD{-1.5}{#3}\Figtr@nptTD{3}{#4}%
    \s@mme=#1\advance\s@mme\@ne\Figp@intregTD\the\s@mme:(\v@lX,\v@lY,\v@lZ)%
    \figptcopyTD#1:/-1/\resetc@ntr@l\et@tfigptscontrolTD}\ignorespaces\fi}
\ctr@ld@f\def\Figtr@nptDD#1#2{\Figg@tXYa{#2}\v@lXa=#1\v@lXa\v@lYa=#1\v@lYa%
    \advance\v@lX\v@lXa\advance\v@lY\v@lYa}
\ctr@ld@f\def\Figtr@nptTD#1#2{\Figg@tXYa{#2}\v@lXa=#1\v@lXa\v@lYa=#1\v@lYa\v@lZa=#1\v@lZa%
    \advance\v@lX\v@lXa\advance\v@lY\v@lYa\advance\v@lZ\v@lZa}
\ctr@ld@f\def\figptscontrolcurve#1,#2[#3]{\ifGR@cri{\s@uvc@ntr@l\et@tfigptscontrolcurve%
    \def\list@num{#3}\extrairelepremi@r\Ak@\de\list@num%
    \extrairelepremi@r\Ai@\de\list@num\extrairelepremi@r\Aj@\de\list@num%
    \s@mme=#1\figptcopy\the\s@mme:/\Ai@/%
    \setc@ntr@l{2}\figvectP -1[\Ak@,\Aj@]%
    \@ecfor\Ak@:=\list@num\do{\advance\s@mme\@ne\figpttra\the\s@mme:=\Ai@/\curv@roundness,-1/%
       \figvectP -1[\Ai@,\Ak@]\advance\s@mme\@ne\figpttra\the\s@mme:=\Aj@/-\curv@roundness,-1/%
       \advance\s@mme\@ne\figptcopy\the\s@mme:/\Aj@/%
       \edef\Ai@{\Aj@}\edef\Aj@{\Ak@}}\advance\s@mme-#1\divide\s@mme\thr@@%
       \xdef#2{\the\s@mme}%
    \resetc@ntr@l\et@tfigptscontrolcurve}\ignorespaces\fi}
\ctr@ln@m\figptsintercirc
\ctr@ld@f\def\figptsintercircDD#1[#2,#3;#4,#5]{\ifGR@cri{\s@uvc@ntr@l\et@tfigptsintercircDD%
    \setc@ntr@l{2}\let\c@lNVintc=\c@lNVintcDD\Figptsintercirc@#1[#2,#3;#4,#5]%
    \resetc@ntr@l\et@tfigptsintercircDD}\ignorespaces\fi}
\ctr@ld@f\def\figptsintercircTD#1[#2,#3;#4,#5;#6]{\ifGR@cri{\s@uvc@ntr@l\et@tfigptsintercircTD%
    \setc@ntr@l{2}\let\c@lNVintc=\c@lNVintcTD\vecunitC@TD[#2,#6]%
    \Figv@ctCreg-3(\v@lX,\v@lY,\v@lZ)\Figptsintercirc@#1[#2,#3;#4,#5]%
    \resetc@ntr@l\et@tfigptsintercircTD}\ignorespaces\fi}
\ctr@ld@f\def\Figptsintercirc@#1[#2,#3;#4,#5]{\figvectP-1[#2,#4]%
    \vecunit@{-1}{-1}\delt@=\result@t\f@ctech=\result@tent%
    \s@mme=#1\advance\s@mme\@ne\figptcopy#1:/#2/\figptcopy\the\s@mme:/#4/%
    \ifdim\delt@=\z@\else%
    \v@lmin=#3\unit@\v@lmax=#5\unit@\v@leur=\v@lmin\advance\v@leur\v@lmax%
    \ifdim\v@leur>\delt@%
    \v@leur=\v@lmin\advance\v@leur-\v@lmax\maxim@m{\v@leur}{\v@leur}{-\v@leur}%
    \ifdim\v@leur<\delt@%
    \divide\v@lmin\f@ctech\divide\v@lmax\f@ctech\divide\delt@\f@ctech%
    \v@lmin=\repdecn@mb{\v@lmin}\v@lmin\v@lmax=\repdecn@mb{\v@lmax}\v@lmax%
    \invers@{\v@leur}{\delt@}\advance\v@lmax-\v@lmin%
    \v@lmax=-\repdecn@mb{\v@leur}\v@lmax\advance\delt@\v@lmax\delt@=.5\delt@%
    \v@lmax=\delt@\multiply\v@lmax\f@ctech%
    \edef\t@ille{\repdecn@mb{\v@lmax}}\figpttra-2:=#2/\t@ille,-1/%
    \delt@=\repdecn@mb{\delt@}\delt@\advance\v@lmin-\delt@%
    \sqrt@{\v@leur}{\v@lmin}\multiply\v@leur\f@ctech\edef\t@ille{\repdecn@mb{\v@leur}}%
    \c@lNVintc\figpttra#1:=-2/-\t@ille,-1/\figpttra\the\s@mme:=-2/\t@ille,-1/\fi\fi\fi}
\ctr@ld@f\def\c@lNVintcDD{\Figg@tXY{-1}\Figv@ctCreg-1(-\v@lY,\v@lX)} 
\ctr@ld@f\def\c@lNVintcTD{{\Figg@tXY{-3}\v@lmin=\v@lX\v@lmax=\v@lY\v@leur=\v@lZ%
    \Figg@tXY{-1}\c@lprovec{-3}\vecunit@{-3}{-3}
    \Figg@tXY{-1}\v@lmin=\v@lX\v@lmax=\v@lY%
    \v@leur=\v@lZ\Figg@tXY{-3}\c@lprovec{-1}}} 
\ctr@ln@m\figptsinterlinell
\ctr@ld@f\def\figptsinterlinellDD#1[#2,#3,#4,#5;#6,#7]{\ifGR@cri{\s@uvc@ntr@l\et@tfigptsinterlinellDD%
    \figptcopy#1:/#6/\s@mme=#1\advance\s@mme\@ne\figptcopy\the\s@mme:/#7/%
    \v@lmin=#3\unit@\v@lmax=#4\unit@
    \setc@ntr@l{2}\figptbaryDD-4:[#6,#7;1,1]\figptsrotDD-3=-4,#7/#2,-#5/
    \Figg@tXY{-3}\Figg@tXYa{#2}\advance\v@lX-\v@lXa\advance\v@lY-\v@lYa
    \figvectP-1[-3,-2]\Figg@tXYa{-1}\figvectP-3[-4,#7]\Figptsint@rLE{#1}
    \resetc@ntr@l\et@tfigptsinterlinellDD}\ignorespaces\fi}
\ctr@ld@f\def\figptsinterlinellP#1[#2,#3,#4;#5,#6]{\ifGR@cri{\s@uvc@ntr@l\et@tfigptsinterlinellP%
    \figptcopy#1:/#5/\s@mme=#1\advance\s@mme\@ne\figptcopy\the\s@mme:/#6/\setc@ntr@l{2}%
    \figvectP-1[#2,#3]\vecunit@{-1}{-1}\v@lmin=\result@t
    \figvectP-2[#2,#4]\vecunit@{-2}{-2}\v@lmax=\result@t
    \figptbary-4:[#5,#6;1,1]
    \figvectP-3[#2,-4]\c@lproscal\v@lX[-3,-1]\c@lproscal\v@lY[-3,-2]
    \figvectP-3[-4,#6]\c@lproscal\v@lXa[-3,-1]\c@lproscal\v@lYa[-3,-2]
    \Figptsint@rLE{#1}\resetc@ntr@l\et@tfigptsinterlinellP}\ignorespaces\fi}
\ctr@ld@f\def\Figptsint@rLE#1{%
    \getredf@ctDD\f@ctech(\v@lmin,\v@lmax)%
    \getredf@ctDD\p@rtent(\v@lX,\v@lY)\ifnum\p@rtent>\f@ctech\f@ctech=\p@rtent\fi%
    \getredf@ctDD\p@rtent(\v@lXa,\v@lYa)\ifnum\p@rtent>\f@ctech\f@ctech=\p@rtent\fi%
    \divide\v@lmin\f@ctech\divide\v@lmax\f@ctech\divide\v@lX\f@ctech\divide\v@lY\f@ctech%
    \divide\v@lXa\f@ctech\divide\v@lYa\f@ctech%
    \c@rre=\repdecn@mb\v@lXa\v@lmax\mili@u=\repdecn@mb\v@lYa\v@lmin%
    \getredf@ctDD\f@ctech(\c@rre,\mili@u)%
    \c@rre=\repdecn@mb\v@lX\v@lmax\mili@u=\repdecn@mb\v@lY\v@lmin%
    \getredf@ctDD\p@rtent(\c@rre,\mili@u)\ifnum\p@rtent>\f@ctech\f@ctech=\p@rtent\fi%
    \divide\v@lmin\f@ctech\divide\v@lmax\f@ctech\divide\v@lX\f@ctech\divide\v@lY\f@ctech%
    \divide\v@lXa\f@ctech\divide\v@lYa\f@ctech%
    \v@lmin=\repdecn@mb{\v@lmin}\v@lmin\v@lmax=\repdecn@mb{\v@lmax}\v@lmax%
    \edef\G@xde{\repdecn@mb\v@lmin}\edef\P@xde{\repdecn@mb\v@lmax}%
    \c@rre=-\v@lmax\v@leur=\repdecn@mb\v@lY\v@lY\advance\c@rre\v@leur\c@rre=\G@xde\c@rre%
    \v@leur=\repdecn@mb\v@lX\v@lX\v@leur=\P@xde\v@leur\advance\c@rre\v@leur
    \v@lmin=\repdecn@mb\v@lYa\v@lmin\v@lmax=\repdecn@mb\v@lXa\v@lmax%
    \mili@u=\repdecn@mb\v@lX\v@lmax\advance\mili@u\repdecn@mb\v@lY\v@lmin
    \v@lmax=\repdecn@mb\v@lXa\v@lmax\advance\v@lmax\repdecn@mb\v@lYa\v@lmin
    \ifdim\v@lmax>\epsil@n%
    \maxim@m{\v@leur}{\c@rre}{-\c@rre}\maxim@m{\v@lmin}{\mili@u}{-\mili@u}%
    \maxim@m{\v@leur}{\v@leur}{\v@lmin}\maxim@m{\v@lmin}{\v@lmax}{-\v@lmax}%
    \maxim@m{\v@leur}{\v@leur}{\v@lmin}\p@rtentiere{\p@rtent}{\v@leur}\advance\p@rtent\@ne%
    \divide\c@rre\p@rtent\divide\mili@u\p@rtent\divide\v@lmax\p@rtent%
    \delt@=\repdecn@mb{\mili@u}\mili@u\v@leur=\repdecn@mb{\v@lmax}\c@rre%
    \advance\delt@-\v@leur\ifdim\delt@<\z@\else\sqrt@\delt@\delt@%
    \invers@\v@lmax\v@lmax\edef\Uns@rAp{\repdecn@mb\v@lmax}%
    \v@leur=-\mili@u\advance\v@leur-\delt@\v@leur=\Uns@rAp\v@leur%
    \edef\t@ille{\repdecn@mb\v@leur}\figpttra#1:=-4/\t@ille,-3/\s@mme=#1\advance\s@mme\@ne%
    \v@leur=-\mili@u\advance\v@leur\delt@\v@leur=\Uns@rAp\v@leur%
    \edef\t@ille{\repdecn@mb\v@leur}\figpttra\the\s@mme:=-4/\t@ille,-3/\fi\fi}
\ctr@ln@m\figptsorthoprojline
\ctr@ld@f\def\figptsorthoprojlineDD#1=#2/#3,#4/{\ifGR@cri{\s@uvc@ntr@l\et@tfigptsorthoprojlineDD%
    \setc@ntr@l{2}\figvectPDD-3[#3,#4]\figvectNVDD-4[-3]\resetc@ntr@l{2}%
    \def\list@num{#2}\s@mme=#1\@ecfor\p@int:=\list@num\do{%
    \inters@cDD\the\s@mme:[\p@int,-4;#3,-3]\advance\s@mme\@ne}%
    \resetc@ntr@l\et@tfigptsorthoprojlineDD}\ignorespaces\fi}
\ctr@ld@f\def\figptsorthoprojlineTD#1=#2/#3,#4/{\ifGR@cri{\s@uvc@ntr@l\et@tfigptsorthoprojlineTD%
    \setc@ntr@l{2}\figvectPTD-2[#3,#4]\vecunit@TD{-2}{-2}%
    \def\list@num{#2}\s@mme=#1\@ecfor\p@int:=\list@num\do{%
    \figvectPTD-1[#3,\p@int]\c@lproscalTD\v@leur[-1,-2]%
    \edef\v@lcoef{\repdecn@mb{\v@leur}}\figpttraTD\the\s@mme:=#3/\v@lcoef,-2/%
    \advance\s@mme\@ne}\resetc@ntr@l\et@tfigptsorthoprojlineTD}\ignorespaces\fi}
\ctr@ln@m\figptsorthoprojplane
\ctr@ld@f\def\figptsorthoprojplaneDD{\un@v@ilable{figptsorthoprojplane}}
\ctr@ld@f\def\figptsorthoprojplaneTD#1=#2/#3,#4/{\ifGR@cri{\s@uvc@ntr@l\et@tfigptsorthoprojplane%
    \setc@ntr@l{2}\vecunit@TD{-2}{#4}%
    \def\list@num{#2}\s@mme=#1\@ecfor\p@int:=\list@num\do{\figvectPTD-1[\p@int,#3]%
    \c@lproscalTD\v@leur[-1,-2]\edef\v@lcoef{\repdecn@mb{\v@leur}}%
    \figpttraTD\the\s@mme:=\p@int/\v@lcoef,-2/\advance\s@mme\@ne}%
    \resetc@ntr@l\et@tfigptsorthoprojplane}\ignorespaces\fi}
\ctr@ld@f\def\figptshom#1=#2/#3,#4/{\ifGR@cri{\s@uvc@ntr@l\et@tfigptshom%
    \setc@ntr@l{2}\def\list@num{#2}\s@mme=#1%
    \@ecfor\p@int:=\list@num\do{\figvectP-1[#3,\p@int]%
    \figpttra\the\s@mme:=#3/#4,-1/\advance\s@mme\@ne}%
    \resetc@ntr@l\et@tfigptshom}\ignorespaces\fi}
\ctr@ld@f\def\figptsinv#1=#2/#3,#4/{\ifGR@cri{\s@uvc@ntr@l\et@tfigptsinv%
    \setc@ntr@l{2}\def\list@num{#2}\s@mme=#1%
    \@ecfor\p@int:=\list@num\do{\figvectP-1[#3,\p@int]\Figg@tXY{-1}%
    \getredf@ctB\f@ctech\n@rmeucC{\delt@}{-1}%
    \delt@=\ptT@unit@\delt@\delt@=\ptT@unit@\delt@%
    \invers@{\delt@}{\delt@}\multiply\f@ctech\f@ctech\divide\delt@\f@ctech%
    \delt@=#4\delt@\edef\v@lcoef{\repdecn@mb{\delt@}}\figpttra\the\s@mme:=#3/\v@lcoef,-1/%
    \advance\s@mme\@ne}\resetc@ntr@l\et@tfigptsinv}\ignorespaces\fi}
\ctr@ln@m\figptsrot
\ctr@ld@f\def\figptsrotDD#1=#2/#3,#4/{\ifGR@cri{\s@uvc@ntr@l\et@tfigptsrotDD%
    \c@ssin{\C@}{\S@}{#4}\setc@ntr@l{2}\def\list@num{#2}\s@mme=#1%
    \@ecfor\p@int:=\list@num\do{\figvectPDD-1[#3,\p@int]\Figg@tXY{-1}%
    \v@lXa=\C@\v@lX\advance\v@lXa-\S@\v@lY%
    \v@lYa=\S@\v@lX\advance\v@lYa\C@\v@lY%
    \Figv@ctCreg-1(\v@lXa,\v@lYa)\figpttraDD\the\s@mme:=#3/1,-1/\advance\s@mme\@ne}%
    \resetc@ntr@l\et@tfigptsrotDD}\ignorespaces\fi}
\ctr@ld@f\def\figptsrotTD#1=#2/#3,#4,#5/{\ifGR@cri{\s@uvc@ntr@l\et@tfigptsrotTD%
    \c@ssin{\C@}{\S@}{#4}%
    \setc@ntr@l{2}\def\list@num{#2}\s@mme=#1%
    \@ecfor\p@int:=\list@num\do{\figptorthoprojplaneTD-3:=#3/\p@int,#5/%
    \figvectPTD-2[-3,\p@int]%
    \figvectNVTD-1[#5,-2]\n@rmeucTD\v@leur{-2}\edef\v@lcoef{\repdecn@mb{\v@leur}}%
    \Figg@tXYa{-1}\v@lXa=\v@lcoef\v@lXa\v@lYa=\v@lcoef\v@lYa\v@lZa=\v@lcoef\v@lZa%
    \v@lXa=\S@\v@lXa\v@lYa=\S@\v@lYa\v@lZa=\S@\v@lZa\Figg@tXY{-2}%
    \advance\v@lXa\C@\v@lX\advance\v@lYa\C@\v@lY\advance\v@lZa\C@\v@lZ%
    \Figg@tXY{-3}\advance\v@lXa\v@lX\advance\v@lYa\v@lY\advance\v@lZa\v@lZ%
    \Figp@intregTD\the\s@mme:(\v@lXa,\v@lYa,\v@lZa)\advance\s@mme\@ne}%
    \resetc@ntr@l\et@tfigptsrotTD}\ignorespaces\fi}
\ctr@ln@m\figptssym
\ctr@ld@f\def\figptssymDD#1=#2/#3,#4/{\ifGR@cri{\s@uvc@ntr@l\et@tfigptssymDD%
    \setc@ntr@l{2}\figvectPDD-3[#3,#4]\Figg@tXY{-3}\Figv@ctCreg-4(-\v@lY,\v@lX)%
    \resetc@ntr@l{2}\def\list@num{#2}\s@mme=#1%
    \@ecfor\p@int:=\list@num\do{\inters@cDD-5:[#3,-3;\p@int,-4]\figvectPDD-2[\p@int,-5]%
    \figpttraDD\the\s@mme:=\p@int/2,-2/\advance\s@mme\@ne}%
    \resetc@ntr@l\et@tfigptssymDD}\ignorespaces\fi}
\ctr@ld@f\def\figptssymTD#1=#2/#3,#4/{\ifGR@cri{\s@uvc@ntr@l\et@tfigptssymTD%
    \setc@ntr@l{2}\vecunit@TD{-2}{#4}\def\list@num{#2}\s@mme=#1%
    \@ecfor\p@int:=\list@num\do{\figvectPTD-1[\p@int,#3]%
    \c@lproscalTD\v@leur[-1,-2]\v@leur=2\v@leur\edef\v@lcoef{\repdecn@mb{\v@leur}}%
    \figpttraTD\the\s@mme:=\p@int/\v@lcoef,-2/\advance\s@mme\@ne}%
    \resetc@ntr@l\et@tfigptssymTD}\ignorespaces\fi}
\ctr@ln@m\figptstra
\ctr@ld@f\def\figptstraDD#1=#2/#3,#4/{\ifGR@cri{\Figg@tXYa{#4}\v@lXa=#3\v@lXa\v@lYa=#3\v@lYa%
    \def\list@num{#2}\s@mme=#1\@ecfor\p@int:=\list@num\do{\Figg@tXY{\p@int}%
    \advance\v@lX\v@lXa\advance\v@lY\v@lYa%
    \Figp@intregDD\the\s@mme:(\v@lX,\v@lY)\advance\s@mme\@ne}}\ignorespaces\fi}
\ctr@ld@f\def\figptstraTD#1=#2/#3,#4/{\ifGR@cri{\Figg@tXYa{#4}\v@lXa=#3\v@lXa\v@lYa=#3\v@lYa%
    \v@lZa=#3\v@lZa\def\list@num{#2}\s@mme=#1\@ecfor\p@int:=\list@num\do{\Figg@tXY{\p@int}%
    \advance\v@lX\v@lXa\advance\v@lY\v@lYa\advance\v@lZ\v@lZa%
    \Figp@intregTD\the\s@mme:(\v@lX,\v@lY,\v@lZ)\advance\s@mme\@ne}}\ignorespaces\fi}
\ctr@ln@m\figptvisilimSL
\ctr@ld@f\def\figptvisilimSLDD{\un@v@ilable{figptvisilimSL}}
\ctr@ld@f\def\figptvisilimSLTD#1:#2[#3,#4;#5,#6]{\ifGR@cri{\s@uvc@ntr@l\et@tfigptvisilimSLTD%
    \setc@ntr@l{2}\figvectP-1[#3,#4]\n@rminf{\delt@}{-1}%
    \ifcase\CUR@proj\v@lX=\cxa@\p@\v@lY=-\p@\v@lZ=\cxb@\p@
    \Figv@ctCreg-2(\v@lX,\v@lY,\v@lZ)\figvectP-3[#5,#6]\figvectNV-1[-2,-3]%
    \or\figvectP-1[#5,#6]\vecunitCV@TD{-1}\v@lmin=\v@lX\v@lmax=\v@lY
    \v@leur=\v@lZ\v@lX=\cza@\p@\v@lY=\czb@\p@\v@lZ=\czc@\p@\c@lprovec{-1}%
    \or\c@ley@pt{-2}\figvectN-1[#5,#6,-2]\fi
    \edef\Ai@{#3}\edef\Aj@{#4}\figvectP-2[#5,\Ai@]\c@lproscal\v@leur[-1,-2]%
    \ifdim\v@leur>\z@\p@rtent=\@ne\else\p@rtent=\m@ne\fi%
    \figvectP-2[#5,\Aj@]\c@lproscal\v@leur[-1,-2]%
    \ifdim\p@rtent\v@leur>\z@\figptcopy#1:#2/#3/%
    \message{*** \BS@ figptvisilimSL: points are on the same side.}\else%
    \figptcopy-3:/#3/\figptcopy-4:/#4/%
    \loop\figptbary-5:[-3,-4;1,1]\figvectP-2[#5,-5]\c@lproscal\v@leur[-1,-2]%
    \ifdim\p@rtent\v@leur>\z@\figptcopy-3:/-5/\else\figptcopy-4:/-5/\fi%
    \divide\delt@\tw@\ifdim\delt@>\epsil@n\repeat%
    \figptbary#1:#2[-3,-4;1,1]\fi\resetc@ntr@l\et@tfigptvisilimSLTD}\ignorespaces\fi}
\ctr@ld@f\def\c@ley@pt#1{\t@stp@r\ifitis@K\v@lX=\cza@\p@\v@lY=\czb@\p@\v@lZ=\czc@\p@%
    \Figv@ctCreg-1(\v@lX,\v@lY,\v@lZ)\Figp@intreg-2:(\wd\Bt@rget,\ht\Bt@rget,\dp\Bt@rget)%
    \figpttra#1:=-2/-\disob@intern,-1/\else\end\fi}
\ctr@ld@f\def\t@stp@r{\itis@Ktrue\ifnewt@rgetpt\else\itis@Kfalse%
    \message{*** \BS@ figptvisilimXX: target point undefined.}\fi\ifnewdis@b\else%
    \itis@Kfalse\message{*** \BS@ figptvisilimXX: observation distance undefined.}\fi%
    \ifitis@K\else\message{*** This macro must be called after \BS@ figdrawbegin or after
    having set the missing parameter(s) with \BS@ figset proj()}\fi}
\ctr@ld@f\def\figscan#1(#2,#3){{\s@uvc@ntr@l\et@tfigscan\@psfgetbb{#1}\if@psfbbfound\else%
    \def\@psfllx{0}\def\@psflly{20}\def\@psfurx{540}\def\@psfury{640}\fi\figscan@{#2}{#3}%
    \resetc@ntr@l\et@tfigscan}\ignorespaces}
\ctr@ld@f\def\figscan@#1#2{%
    \unit@=\@ne bp\setc@ntr@l{2}\figsetmark{}%
    \def\minst@p{20pt}%
    \v@lX=\@psfllx\p@\v@lX=\Sc@leFact\v@lX\r@undint\v@lX\v@lX%
    \v@lY=\@psflly\p@\v@lY=\Sc@leFact\v@lY\ifdim\v@lY>\z@\r@undint\v@lY\v@lY\fi%
    \delt@=\@psfury\p@\delt@=\Sc@leFact\delt@%
    \advance\delt@-\v@lY\v@lXa=\@psfurx\p@\v@lXa=\Sc@leFact\v@lXa\v@leur=\minst@p%
    \edef\valv@lY{\repdecn@mb{\v@lY}}\edef\LgTr@it{\the\delt@}%
    \loop\ifdim\v@lX<\v@lXa\edef\valv@lX{\repdecn@mb{\v@lX}}%
    \figptDD -1:(\valv@lX,\valv@lY)\figwriten -1:\hbox{\vrule height\LgTr@it}(0)%
    \ifdim\v@leur<\minst@p\else\figsetmark{\raise-8bp\hbox{$\scriptscriptstyle\triangle$}}%
    \figwrites -1:\@ffichnb{0}{\valv@lX}(6)\v@leur=\z@\figsetmark{}\fi%
    \advance\v@leur#1pt\advance\v@lX#1pt\repeat%
    \def\minst@p{10pt}%
    \v@lX=\@psfllx\p@\v@lX=\Sc@leFact\v@lX\ifdim\v@lX>\z@\r@undint\v@lX\v@lX\fi%
    \v@lY=\@psflly\p@\v@lY=\Sc@leFact\v@lY\r@undint\v@lY\v@lY%
    \delt@=\@psfurx\p@\delt@=\Sc@leFact\delt@%
    \advance\delt@-\v@lX\v@lYa=\@psfury\p@\v@lYa=\Sc@leFact\v@lYa\v@leur=\minst@p%
    \edef\valv@lX{\repdecn@mb{\v@lX}}\edef\LgTr@it{\the\delt@}%
    \loop\ifdim\v@lY<\v@lYa\edef\valv@lY{\repdecn@mb{\v@lY}}%
    \figptDD -1:(\valv@lX,\valv@lY)\figwritee -1:\vbox{\hrule width\LgTr@it}(0)%
    \ifdim\v@leur<\minst@p\else\figsetmark{$\triangleright$\kern4bp}%
    \figwritew -1:\@ffichnb{0}{\valv@lY}(6)\v@leur=\z@\figsetmark{}\fi%
    \advance\v@leur#2pt\advance\v@lY#2pt\repeat}
\ctr@ld@f
\ctr@ld@f\def\figscan@E#1(#2,#3){{\s@uvc@ntr@l\et@tfigscan@E%
    \Figdisc@rdLTS{#1}{\t@xt@}\pdfximage{\t@xt@}%
    \setbox\Gb@x=\hbox{\pdfrefximage\pdflastximage}%
    \edef\@psfllx{0}\v@lY=-\dp\Gb@x\edef\@psflly{\repdecn@mb{\v@lY}}%
    \edef\@psfurx{\repdecn@mb{\wd\Gb@x}}%
    \v@lY=\dp\Gb@x\advance\v@lY\ht\Gb@x\edef\@psfury{\repdecn@mb{\v@lY}}%
    \figscan@{#2}{#3}\resetc@ntr@l\et@tfigscan@E}\ignorespaces}
\ctr@ld@f\def\figshowpts[#1,#2]{{\figsetmark{$\bullet$}\figsetptname{\bf ##1}%
    \p@rtent=#2\relax\ifnum\p@rtent<\z@\p@rtent=\z@\fi%
    \s@mme=#1\relax\ifnum\s@mme<\z@\s@mme=\z@\fi%
    \loop\ifnum\s@mme<\p@rtent\pt@rvect{\s@mme}%
    \ifitis@K\figwriten{\the\s@mme}:(4pt)\fi\advance\s@mme\@ne\repeat%
    \pt@rvect{\s@mme}\ifitis@K\figwriten{\the\s@mme}:(4pt)\fi}\ignorespaces}
\ctr@ld@f\def\pt@rvect#1{\set@bjc@de{#1}%
    \expandafter\expandafter\expandafter\inqpt@rvec\csname\objc@de\endcsname:}
\ctr@ld@f\def\inqpt@rvec#1#2:{\if#1\C@dCl@spt\itis@Ktrue\else\itis@Kfalse\fi}
\ctr@ld@f\def\figshowsettings{{%
    \immediate\write16{====================================================================}%
    \immediate\write16{ Current settings are (DDV means "with dynamic default value"):}%
    \immediate\write16{ --- GENERAL ---}%
    \immediate\write16{Scale factor and Unit = \unit@util\space (\the\unit@)
     \space -> \BS@ figinit{ScaleFactorUnit}}%
    \immediate\write16{Update mode = \ifGRupdatem@de yes\else no\fi
     \space-> \BS@ figset(update=yes/no) or \BS@ figsetdefault(update=yes/no)}%
    \immediate\write16{ --- WRITING ---}%
    \immediate\write16{Implicit point name = \ptn@me{i} \space-> \BS@ figset write(ptname={Name})}%
    \immediate\write16{Point marker = \the\c@nsymb \space -> \BS@ figset write(mark=Mark)}%
    \immediate\write16{Print rounded coordinates = \ifr@undcoord yes\else no\fi
     \space-> \BS@ figset write(roundcoord=yes/no)}%
    \immediate\write16{ --- GRAPHICAL (general) ---}%
    \immediate\write16{Color = \CUR@color \space-> \BS@ figset(color=ColorDefinition)}%
    \immediate\write16{Filling mode = \iffillm@de yes\else no\fi
     \space-> \BS@ figset(fillmode=yes/no)}%
    \immediate\write16{Line join = \CUR@join \space-> \BS@ figset(join=miter/round/bevel)}%
    \immediate\write16{Line style = \CUR@dash \space-> \BS@ figset(dash=Index/Pattern)}%
    \immediate\write16{Line width = \CUR@width
     \space-> \BS@ figset(width=real in PostScript units)}%
    \immediate\write16{ --- GRAPHICAL (specific) ---}%
    \immediate\write16{Altitude (all the following attributes are DDV):}%
    \immediate\write16{ Base line color =
     \ifx\DDV@blcolor\D@FTref general color\else\DDV@blcolor\fi
     \space-> \BS@ figset altitude(blcolor=ColorDefinition)}%
    \immediate\write16{ Base line style =
     \ifx\DDV@bldash\D@FTref general style\else\DDV@bldash\fi
     \space-> \BS@ figset altitude(bldash=Index/Pattern)}%
    \immediate\write16{ Base line width =
     \ifx\DDV@blwidth\D@FTref general width\else\DDV@blwidth\fi
     \space-> \BS@ figset altitude(blwidth=real in PostScript units)}%
    \immediate\write16{ Square line color =
     \ifx\DDV@sqcolor\D@FTref general color\else\DDV@sqcolor\fi
     \space-> \BS@ figset altitude(sqcolor=ColorDefinition)}%
    \immediate\write16{ Square line style =
     \ifx\DDV@sqdash\D@FTref general style\else\DDV@sqdash\fi
     \space-> \BS@ figset altitude(sqdash=Index/Pattern)}%
    \immediate\write16{ Square line width =
     \ifx\DDV@sqwidth\D@FTref general width\else\DDV@sqwidth\fi
     \space-> \BS@ figset altitude(sqwidth=real in PostScript units)}%
    \immediate\write16{Arrowhead:}%
    \immediate\write16{ (half-)Angle = \@rrowheadangle
     \space-> \BS@ figset arrowhead(angle=real in degrees)}%
    \immediate\write16{ Filling mode = \if@rrowhfill yes\else no\fi
     \space-> \BS@ figset arrowhead(fillmode=yes/no)}%
    \immediate\write16{ "Outside" = \if@rrowhout yes\else no\fi
     \space-> \BS@ figset arrowhead(out=yes/no)}%
    \immediate\write16{ Length = \@rrowheadlength
     \if@rrowratio\space(not active)\else\space(active)\fi
     \space-> \BS@ figset arrowhead(length=real in user coord.)}%
    \immediate\write16{ Ratio = \@rrowheadratio
     \if@rrowratio\space(active)\else\space(not active)\fi
     \space-> \BS@ figset arrowhead(ratio=real in [0,1])}%
    \immediate\write16{Curve:}%
    \immediate\write16{ Roundness = \curv@roundness
     \space-> \BS@ figset curve(roundness=real in [0,0.5])}%
    \immediate\write16{Flow chart:}%
    \immediate\write16{ Arrow position = \@rrowp@s
     \space-> \BS@ figset flowchart(arrowposition=real in [0,1])}%
    \immediate\write16{ Arrow reference point = \ifcase\@rrowr@fpt start\else end\fi
     \space-> \BS@ figset flowchart(arrowrefpt = start/end)}%
    \immediate\write16{ Background color = \fcbgc@lor
     \space-> \BS@ figset flowchart(bgcolor=ColorDefinition)}%
    \immediate\write16{ Line type = \ifcase\fclin@typ@ curve\else polygon\fi
     \space-> \BS@ figset flowchart(line=polygon/curve)}%
    \immediate\write16{ Padding = (\Xp@dd, \Yp@dd)
     \space-> \BS@ figset flowchart(padding = real in user coord.)}%
    \immediate\write16{\space\space\space\space(or
     \BS@ figset flowchart(xpadding=real, ypadding=real) )}%
    \immediate\write16{ Radius = \fclin@r@d
     \space-> \BS@ figset flowchart(radius=positive real in user coord.)}%
    \immediate\write16{ Shape = \fcsh@pe
     \space-> \BS@ figset flowchart(shape = rectangle, ellipse or lozenge)}%
    \immediate\write16{ Thickness color (DDV) = 
     \ifx\DDV@thickcolor\D@FTref general color\else\DDV@thickcolor\fi
     \space-> \BS@ figset flowchart(thickcolor=ColorDefinition)}%
    \immediate\write16{ Thickness = \thickn@ss
     \space-> \BS@ figset flowchart(thickness = real in user coord.)}%
    \immediate\write16{Mesh:}%
    \immediate\write16{ Diagonal = \c@ntrolmesh
     \space-> \BS@ figset mesh(diag=integer in {-1,0,1})}%
    \immediate\write16{ Lines color (DDV) =
     \ifx\DDV@meshcolor\D@FTref general color\else\DDV@meshcolor\fi
     \space-> \BS@ figset mesh(color=ColorDefinition)}%
    \immediate\write16{ Lines style (DDV) =
     \ifx\DDV@meshdash\D@FTref general style\else\DDV@meshdash\fi
     \space-> \BS@ figset mesh(dash=Index/Pattern)}%
    \immediate\write16{ Lines width (DDV) =
     \ifx\DDV@meshwidth\D@FTref general width\else\DDV@meshwidth\fi
     \space-> \BS@ figset mesh(width=real in PostScript units)}%
    \immediate\write16{Trimesh:}%
    \immediate\write16{ Lines color (DDV) =
     \ifx\DDV@tmeshcolor\D@FTref general color\else\DDV@tmeshcolor\fi
     \space-> \BS@ figset trimesh(color=ColorDefinition)}%
    \immediate\write16{ Lines style (DDV) =
     \ifx\DDV@tmeshdash\D@FTref general style\else\DDV@tmeshdash\fi
     \space-> \BS@ figset trimesh(dash=Index/Pattern)}%
    \immediate\write16{ Lines width (DDV) =
     \ifx\DDV@tmeshwidth\D@FTref general width\else\DDV@tmeshwidth\fi
     \space-> \BS@ figset trimesh(width=real in PostScript units)}%
    \ifTr@isDim%
    \immediate\write16{ --- 3D to 2D PROJECTION ---}%
    \immediate\write16{Projection : \typ@proj \space-> \BS@ figinit{ScaleFactorUnit, ProjType}}%
    \immediate\write16{Longitude (psi) = \v@lPsi \space-> \BS@ figset proj(psi=real in degrees)}%
    \ifcase\CUR@proj\immediate\write16{Depth coeff. (Lambda)
     \space = \v@lTheta \space-> \BS@ figset proj(lambda=real in [0,1])}%
    \else\immediate\write16{Latitude (theta)
     \space = \v@lTheta \space-> \BS@ figset proj(theta=real in degrees)}%
    \fi%
    \ifnum\CUR@proj=\tw@%
    \immediate\write16{Observation distance = \disob@unit
     \space-> \BS@ figset proj(dist=real in user coord.)}%
    \immediate\write16{Target point = \t@rgetpt \space-> \BS@ figset proj(targetpt=pt number)}%
     \v@lX=\ptT@unit@\wd\Bt@rget\v@lY=\ptT@unit@\ht\Bt@rget\v@lZ=\ptT@unit@\dp\Bt@rget%
    \immediate\write16{ Its coordinates are
     (\repdecn@mb{\v@lX}, \repdecn@mb{\v@lY}, \repdecn@mb{\v@lZ})}%
    \fi%
    \fi%
    \immediate\write16{====================================================================}%
    \ignorespaces}}
\ctr@ln@w{newif}\ifitis@vect@r
\ctr@ld@f\def\figvectC#1(#2,#3){{\itis@vect@rtrue\figpt#1:(#2,#3)}\ignorespaces}
\ctr@ld@f\def\Figv@ctCreg#1(#2,#3){{\itis@vect@rtrue\Figp@intreg#1:(#2,#3)}\ignorespaces}
\ctr@ln@m\figvectDBezier
\ctr@ld@f\def\figvectDBezierDD#1:#2,#3[#4,#5,#6,#7]{\ifGR@cri{\s@uvc@ntr@l\et@tfigvectDBezierDD%
    \FigvectDBezier@#2,#3[#4,#5,#6,#7]\v@lX=\c@ef\v@lX\v@lY=\c@ef\v@lY%
    \Figv@ctCreg#1(\v@lX,\v@lY)\resetc@ntr@l\et@tfigvectDBezierDD}\ignorespaces\fi}
\ctr@ld@f\def\figvectDBezierTD#1:#2,#3[#4,#5,#6,#7]{\ifGR@cri{\s@uvc@ntr@l\et@tfigvectDBezierTD%
    \FigvectDBezier@#2,#3[#4,#5,#6,#7]\v@lX=\c@ef\v@lX\v@lY=\c@ef\v@lY\v@lZ=\c@ef\v@lZ%
    \Figv@ctCreg#1(\v@lX,\v@lY,\v@lZ)\resetc@ntr@l\et@tfigvectDBezierTD}\ignorespaces\fi}
\ctr@ld@f\def\FigvectDBezier@#1,#2[#3,#4,#5,#6]{\setc@ntr@l{2}%
    \edef\T@{#2}\v@leur=\p@\advance\v@leur-#2pt\edef\UNmT@{\repdecn@mb{\v@leur}}%
    \ifnum#1=\tw@\def\c@ef{6}\else\def\c@ef{3}\fi%
    \figptcopy-4:/#3/\figptcopy-3:/#4/\figptcopy-2:/#5/\figptcopy-1:/#6/%
    \l@mbd@un=-4 \l@mbd@de=-\thr@@\p@rtent=\m@ne\c@lDecast%
    \ifnum#1=\tw@\c@lDCDeux{-4}{-3}\c@lDCDeux{-3}{-2}\c@lDCDeux{-4}{-3}\else%
    \l@mbd@un=-4 \l@mbd@de=-\thr@@\p@rtent=-\tw@\c@lDecast%
    \c@lDCDeux{-4}{-3}\fi\Figg@tXY{-4}}
\ctr@ln@m\c@lDCDeux
\ctr@ld@f\def\c@lDCDeuxDD#1#2{\Figg@tXY{#2}\Figg@tXYa{#1}%
    \advance\v@lX-\v@lXa\advance\v@lY-\v@lYa\Figp@intregDD#1:(\v@lX,\v@lY)}
\ctr@ld@f\def\c@lDCDeuxTD#1#2{\Figg@tXY{#2}\Figg@tXYa{#1}\advance\v@lX-\v@lXa%
    \advance\v@lY-\v@lYa\advance\v@lZ-\v@lZa\Figp@intregTD#1:(\v@lX,\v@lY,\v@lZ)}
\ctr@ln@m\figvectN
\ctr@ld@f\def\figvectNDD#1[#2,#3]{\ifGR@cri{\Figg@tXYa{#2}\Figg@tXY{#3}%
    \advance\v@lX-\v@lXa\advance\v@lY-\v@lYa%
    \Figv@ctCreg#1(-\v@lY,\v@lX)}\ignorespaces\fi}
\ctr@ld@f\def\figvectNTD#1[#2,#3,#4]{\ifGR@cri{\vecunitC@TD[#2,#4]\v@lmin=\v@lX\v@lmax=\v@lY%
    \v@leur=\v@lZ\vecunitC@TD[#2,#3]\c@lprovec{#1}}\ignorespaces\fi}
\ctr@ln@m\figvectNV
\ctr@ld@f\def\figvectNVDD#1[#2]{\ifGR@cri{\Figg@tXY{#2}\Figv@ctCreg#1(-\v@lY,\v@lX)}\ignorespaces\fi}
\ctr@ld@f\def\figvectNVTD#1[#2,#3]{\ifGR@cri{\vecunitCV@TD{#3}\v@lmin=\v@lX\v@lmax=\v@lY%
    \v@leur=\v@lZ\vecunitCV@TD{#2}\c@lprovec{#1}}\ignorespaces\fi}
\ctr@ln@m\figvectP
\ctr@ld@f\def\figvectPDD#1[#2,#3]{\ifGR@cri{\Figg@tXYa{#2}\Figg@tXY{#3}%
    \advance\v@lX-\v@lXa\advance\v@lY-\v@lYa%
    \Figv@ctCreg#1(\v@lX,\v@lY)}\ignorespaces\fi}
\ctr@ld@f\def\figvectPTD#1[#2,#3]{\ifGR@cri{\Figg@tXYa{#2}\Figg@tXY{#3}%
    \advance\v@lX-\v@lXa\advance\v@lY-\v@lYa\advance\v@lZ-\v@lZa%
    \Figv@ctCreg#1(\v@lX,\v@lY,\v@lZ)}\ignorespaces\fi}
\ctr@ln@m\figvectU
\ctr@ld@f\def\figvectUDD#1[#2]{\ifGR@cri{\n@rmeuc\v@leur{#2}\invers@\v@leur\v@leur%
    \delt@=\repdecn@mb{\v@leur}\unit@\edef\v@ldelt@{\repdecn@mb{\delt@}}%
    \Figg@tXY{#2}\v@lX=\v@ldelt@\v@lX\v@lY=\v@ldelt@\v@lY%
    \Figv@ctCreg#1(\v@lX,\v@lY)}\ignorespaces\fi}
\ctr@ld@f\def\figvectUTD#1[#2]{\ifGR@cri{\n@rmeuc\v@leur{#2}\invers@\v@leur\v@leur%
    \delt@=\repdecn@mb{\v@leur}\unit@\edef\v@ldelt@{\repdecn@mb{\delt@}}%
    \Figg@tXY{#2}\v@lX=\v@ldelt@\v@lX\v@lY=\v@ldelt@\v@lY\v@lZ=\v@ldelt@\v@lZ%
    \Figv@ctCreg#1(\v@lX,\v@lY,\v@lZ)}\ignorespaces\fi}
\ctr@ld@f\def\figvisu#1#2#3{\c@ldefproj\initb@undb@x\xdef\figforTeXFigno{\figforTeXnextFigno}%
    \s@mme=\figforTeXnextFigno\advance\s@mme\@ne\xdef\figforTeXnextFigno{\number\s@mme}%
    \setbox\b@xvisu=\hbox{\ifnum\@utoFN>\z@\figinsert{}\gdef\@utoFInDone{0}\fi\ignorespaces#3}%
    \gdef\@utoFInDone{1}\gdef\@utoFN{0}%
    \v@lXa=-\c@@rdYmin\v@lYa=\c@@rdYmax\advance\v@lYa-\c@@rdYmin%
    \v@lX=\c@@rdXmax\advance\v@lX-\c@@rdXmin%
    \setbox#1=\hbox{#2}\v@lY=-\v@lX\maxim@m{\v@lX}{\v@lX}{\wd#1}%
    \advance\v@lY\v@lX\divide\v@lY\tw@\advance\v@lY-\c@@rdXmin%
    \setbox#1=\vbox{\parindent\z@\hsize=\v@lX\vskip\v@lYa%
    \rlap{\hskip\v@lY\smash{\raise\v@lXa\box\b@xvisu}}%
    \def\t@xt@{#2}\ifx\t@xt@\empty\else\medskip\centerline{#2}\fi}\wd#1=\v@lX}
\ctr@ld@f\def\figDecrementFigno{{\xdef\figforTeXnextFigno{\figforTeXFigno}%
    \s@mme=\figforTeXFigno\advance\s@mme\m@ne\xdef\figforTeXFigno{\number\s@mme}}}
\ctr@ln@w{newbox}\Bt@rget\setbox\Bt@rget=\null
\ctr@ln@w{newbox}\BminTD@\setbox\BminTD@=\null
\ctr@ln@w{newbox}\BmaxTD@\setbox\BmaxTD@=\null
\ctr@ln@w{newif}\ifnewt@rgetpt\ctr@ln@w{newif}\ifnewdis@b
\ctr@ld@f\def\b@undb@xTD#1#2#3{%
    \relax\ifdim#1<\wd\BminTD@\global\wd\BminTD@=#1\fi%
    \relax\ifdim#2<\ht\BminTD@\global\ht\BminTD@=#2\fi%
    \relax\ifdim#3<\dp\BminTD@\global\dp\BminTD@=#3\fi%
    \relax\ifdim#1>\wd\BmaxTD@\global\wd\BmaxTD@=#1\fi%
    \relax\ifdim#2>\ht\BmaxTD@\global\ht\BmaxTD@=#2\fi%
    \relax\ifdim#3>\dp\BmaxTD@\global\dp\BmaxTD@=#3\fi}
\ctr@ld@f\def\c@ldefdisob{{\ifdim\wd\BminTD@<\maxdimen\v@leur=\wd\BmaxTD@\advance\v@leur-\wd\BminTD@%
    \delt@=\ht\BmaxTD@\advance\delt@-\ht\BminTD@\maxim@m{\v@leur}{\v@leur}{\delt@}%
    \delt@=\dp\BmaxTD@\advance\delt@-\dp\BminTD@\maxim@m{\v@leur}{\v@leur}{\delt@}%
    \v@leur=5\v@leur\else\v@leur=800pt\fi\c@ldefdisob@{\v@leur}}}
\ctr@ln@m\disob@intern
\ctr@ln@m\disob@
\ctr@ln@m\divf@ctproj
\ctr@ld@f\def\c@ldefdisob@#1{{\v@leur=#1\ifdim\v@leur<\p@\v@leur=800pt\fi%
    \xdef\disob@intern{\repdecn@mb{\v@leur}}%
    \delt@=\ptT@unit@\v@leur\xdef\disob@unit{\repdecn@mb{\delt@}}%
    \f@ctech=\@ne\loop\ifdim\v@leur>\t@n pt\divide\v@leur\t@n\multiply\f@ctech\t@n\repeat%
    \xdef\disob@{\repdecn@mb{\v@leur}}\xdef\divf@ctproj{\the\f@ctech}}%
    \global\newdis@btrue}
\ctr@ln@m\t@rgetpt
\ctr@ld@f\def\c@ldeft@rgetpt{\newt@rgetpttrue\def\t@rgetpt{CenterBoundBox}{%
    \delt@=\wd\BmaxTD@\advance\delt@-\wd\BminTD@\divide\delt@\tw@%
    \v@leur=\wd\BminTD@\advance\v@leur\delt@\global\wd\Bt@rget=\v@leur%
    \delt@=\ht\BmaxTD@\advance\delt@-\ht\BminTD@\divide\delt@\tw@%
    \v@leur=\ht\BminTD@\advance\v@leur\delt@\global\ht\Bt@rget=\v@leur%
    \delt@=\dp\BmaxTD@\advance\delt@-\dp\BminTD@\divide\delt@\tw@%
    \v@leur=\dp\BminTD@\advance\v@leur\delt@\global\dp\Bt@rget=\v@leur}}
\ctr@ln@m\c@ldefproj
\ctr@ld@f\def\c@ldefprojTD{\ifnewt@rgetpt\else\c@ldeft@rgetpt\fi\ifnewdis@b\else\c@ldefdisob\fi}
\ctr@ld@f\def\c@lprojcav{
    \v@lZa=\cxa@\v@lY\advance\v@lX\v@lZa%
    \v@lZa=\cxb@\v@lY\v@lY=\v@lZ\advance\v@lY\v@lZa\ignorespaces}
\ctr@ln@m\v@lcoef
\ctr@ld@f\def\c@lprojrea{
    \advance\v@lX-\wd\Bt@rget\advance\v@lY-\ht\Bt@rget\advance\v@lZ-\dp\Bt@rget%
    \v@lZa=\cza@\v@lX\advance\v@lZa\czb@\v@lY\advance\v@lZa\czc@\v@lZ%
    \divide\v@lZa\divf@ctproj\advance\v@lZa\disob@ pt\invers@{\v@lZa}{\v@lZa}%
    \v@lZa=\disob@\v@lZa\edef\v@lcoef{\repdecn@mb{\v@lZa}}%
    \v@lXa=\cxa@\v@lX\advance\v@lXa\cxb@\v@lY\v@lXa=\v@lcoef\v@lXa%
    \v@lY=\cyb@\v@lY\advance\v@lY\cya@\v@lX\advance\v@lY\cyc@\v@lZ%
    \v@lY=\v@lcoef\v@lY\v@lX=\v@lXa\ignorespaces}
\ctr@ld@f\def\c@lprojort{
    \v@lXa=\cxa@\v@lX\advance\v@lXa\cxb@\v@lY%
    \v@lY=\cyb@\v@lY\advance\v@lY\cya@\v@lX\advance\v@lY\cyc@\v@lZ%
    \v@lX=\v@lXa\ignorespaces}
\ctr@ld@f\def\Figptpr@j#1:#2/#3/{{\Figg@tXY{#3}\superc@lprojSP%
    \Figp@intregDD#1:{#2}(\v@lX,\v@lY)}\ignorespaces}
\ctr@ln@m\figsetobdist
\ctr@ld@f\def\figsetobdistDD{\un@v@ilable{figsetobdist}}
\ctr@ld@f\def\figsetobdistTD(#1){{\ifCUR@PS\W@rnmesIgn{figset proj(dist=...)}%
    \else\v@leur=#1\unit@\c@ldefdisob@{\v@leur}\fi}\ignorespaces}
\ctr@ln@m\c@lprojSP
\ctr@ln@m\CUR@proj
\ctr@ln@m\typ@proj
\ctr@ln@m\superc@lprojSP
\ctr@ld@f\def\Figs@tproj#1{%
    \if#13 \def@ultproj\else\if#1c\def@ultproj%
    \else\if#1o\xdef\CUR@proj{1}\xdef\typ@proj{orthogonal}%
         \figsetviewTD(\def@ultpsi,\def@ulttheta)%
         \global\let\c@lprojSP=\c@lprojort\global\let\superc@lprojSP=\c@lprojort%
    \else\if#1r\xdef\CUR@proj{2}\xdef\typ@proj{realistic}%
         \figsetviewTD(\def@ultpsi,\def@ulttheta)%
         \global\let\c@lprojSP=\c@lprojrea\global\let\superc@lprojSP=\c@lprojrea%
    \else\def@ultproj\message{*** Unknown projection. Cavalier projection assumed.}%
    \fi\fi\fi\fi}
\ctr@ld@f\def\def@ultproj{\xdef\CUR@proj{0}\xdef\typ@proj{cavalier}\figsetviewTD(\def@ultpsi,0.5)%
         \global\let\c@lprojSP=\c@lprojcav\global\let\superc@lprojSP=\c@lprojcav}
\ctr@ln@m\figsettarget
\ctr@ld@f\def\figsettargetDD{\un@v@ilable{figsettarget}}
\ctr@ld@f\def\figsettargetTD[#1]{{\ifCUR@PS\W@rnmesIgn{figset proj(targetpt=...)}%
    \else\global\newt@rgetpttrue\xdef\t@rgetpt{#1}\Figg@tXY{#1}\global\wd\Bt@rget=\v@lX%
    \global\ht\Bt@rget=\v@lY\global\dp\Bt@rget=\v@lZ\fi}\ignorespaces}
\ctr@ln@m\figsetview
\ctr@ld@f\def\figsetviewDD{\un@v@ilable{figsetview}}
\ctr@ld@f\def\figsetviewTD(#1){\ifCUR@PS\W@rnmesIgn{figset proj(Psi|Theta|Lambda=...)}%
     \else\Figsetview@#1,:\fi\ignorespaces}
\ctr@ld@f\def\Figsetview@#1,#2:{{\xdef\v@lPsi{#1}\def\t@xt@{#2}%
    \ifx\t@xt@\empty\def\@rgdeux{\v@lTheta}\else\X@rgdeux@#2\fi%
    \c@ssin{\costhet@}{\sinthet@}{#1}\v@lmin=\costhet@ pt\v@lmax=\sinthet@ pt%
    \ifcase\CUR@proj%
    \v@leur=\@rgdeux\v@lmin\xdef\cxa@{\repdecn@mb{\v@leur}}%
    \v@leur=\@rgdeux\v@lmax\xdef\cxb@{\repdecn@mb{\v@leur}}\v@leur=\@rgdeux pt%
    \relax\ifdim\v@leur>\p@\message{*** Lambda too large ! See \BS@ figset proj() !}\fi%
    \else%
    \v@lmax=-\v@lmax\xdef\cxa@{\repdecn@mb{\v@lmax}}\xdef\cxb@{\costhet@}%
    \ifx\t@xt@\empty\edef\@rgdeux{\def@ulttheta}\fi\c@ssin{\C@}{\S@}{\@rgdeux}%
    \v@lmax=-\S@ pt%
    \v@leur=\v@lmax\v@leur=\costhet@\v@leur\xdef\cya@{\repdecn@mb{\v@leur}}%
    \v@leur=\v@lmax\v@leur=\sinthet@\v@leur\xdef\cyb@{\repdecn@mb{\v@leur}}%
    \xdef\cyc@{\C@}\v@lmin=-\C@ pt%
    \v@leur=\v@lmin\v@leur=\costhet@\v@leur\xdef\cza@{\repdecn@mb{\v@leur}}%
    \v@leur=\v@lmin\v@leur=\sinthet@\v@leur\xdef\czb@{\repdecn@mb{\v@leur}}%
    \xdef\czc@{\repdecn@mb{\v@lmax}}\fi%
    \xdef\v@lTheta{\@rgdeux}}}
\ctr@ld@f\def\def@ultpsi{40}
\ctr@ld@f\def\def@ulttheta{25}
\ctr@ln@m\l@debut
\ctr@ln@m\n@mref
\ctr@ld@f\def\Figsetpr@j#1=#2|{\keln@mtr#1|%
    \def\n@mref{dep}\ifx\l@debut\n@mref\Figsetd@p{#2}\else
    \def\n@mref{dis}\ifx\l@debut\n@mref%
     \ifnum\CUR@proj=\tw@\figsetobdist(#2)\else\Figset@rr\fi\else
    \def\n@mref{lam}\ifx\l@debut\n@mref\Figsetd@p{#2}\else
    \def\n@mref{lat}\ifx\l@debut\n@mref\Figsetth@{#2}\else
    \def\n@mref{lon}\ifx\l@debut\n@mref\figsetview(#2)\else
    \def\n@mref{psi}\ifx\l@debut\n@mref\figsetview(#2)\else
    \def\n@mref{tar}\ifx\l@debut\n@mref%
     \ifnum\CUR@proj=\tw@\figsettarget[#2]\else\Figset@rr\fi\else
    \def\n@mref{the}\ifx\l@debut\n@mref\Figsetth@{#2}\else
    \W@rnmesAttr{figset proj}{#1}\fi\fi\fi\fi\fi\fi\fi\fi}
\ctr@ld@f\def\Figsetd@p#1{\ifnum\CUR@proj=\z@\figsetview(\v@lPsi,#1)\else\Figset@rr\fi}
\ctr@ld@f\def\Figsetth@#1{\ifnum\CUR@proj=\z@\Figset@rr\else\figsetview(\v@lPsi,#1)\fi}
\ctr@ld@f\def\Figset@rr{\message{*** \BS@ figset proj(): Attribute "\n@mref" ignored, incompatible
    with current projection}}
\ctr@ld@f\def\initb@undb@xTD{\wd\BminTD@=\maxdimen\ht\BminTD@=\maxdimen\dp\BminTD@=\maxdimen%
    \wd\BmaxTD@=-\maxdimen\ht\BmaxTD@=-\maxdimen\dp\BmaxTD@=-\maxdimen}
\ctr@ln@w{newbox}\Gb@x      
\ctr@ln@w{newbox}\Gb@xSC    
\ctr@ln@w{newtoks}\c@nsymb  
\ctr@ln@w{newif}\ifr@undcoord\ctr@ln@w{newif}\ifunitpr@sent
\ctr@ld@f\def\unssqrttw@{0.707106 }
\ctr@ld@f\def\figAst{\raise-1.15ex\hbox{$\ast$}}
\ctr@ld@f\def\figBullet{\raise-1.15ex\hbox{$\bullet$}}
\ctr@ld@f\def\figCirc{\raise-1.15ex\hbox{$\circ$}}
\ctr@ld@f\def\figDiamond{\raise-1.15ex\hbox{$\diamond$}}%
\ctr@ld@f\def\boxit#1#2{\leavevmode\hbox{\vrule\vbox{\hrule\vglue#1%
    \vtop{\hbox{\kern#1{#2}\kern#1}\vglue#1\hrule}}\vrule}}
\ctr@ld@f
\ctr@ld@f
\ctr@ld@f\def\c@nterpt{\ignorespaces%
    \kern-.5\wd\Gb@xSC%
    \raise-.5\ht\Gb@xSC\rlap{\hbox{\raise.5\dp\Gb@xSC\hbox{\copy\Gb@xSC}}}%
    \kern .5\wd\Gb@xSC\ignorespaces}
\ctr@ld@f\def\b@undb@xSC#1#2{{\v@lXa=#1\v@lYa=#2%
    \v@leur=\ht\Gb@xSC\advance\v@leur\dp\Gb@xSC%
    \advance\v@lXa-.5\wd\Gb@xSC\advance\v@lYa-.5\v@leur\b@undb@x{\v@lXa}{\v@lYa}%
    \advance\v@lXa\wd\Gb@xSC\advance\v@lYa\v@leur\b@undb@x{\v@lXa}{\v@lYa}}}
\ctr@ln@m\Dist@n
\ctr@ln@m\l@suite
\ctr@ld@f\def\@keldist#1#2{\edef\Dist@n{#2}\y@tiunit{\Dist@n}%
    \ifunitpr@sent#1=\Dist@n\else#1=\Dist@n\unit@\fi}
\ctr@ld@f\def\y@tiunit#1{\unitpr@sentfalse\expandafter\y@tiunit@#1:}
\ctr@ld@f\def\y@tiunit@#1#2:{\ifcat#1a\unitpr@senttrue\else\def\l@suite{#2}%
    \ifx\l@suite\empty\else\y@tiunit@#2:\fi\fi}
\ctr@ln@m\figcoord
\ctr@ld@f\def\figcoordDD#1{{\v@lX=\ptT@unit@\v@lX\v@lY=\ptT@unit@\v@lY%
    \ifr@undcoord\ifcase#1\v@leur=0.5pt\or\v@leur=0.05pt\or\v@leur=0.005pt%
    \or\v@leur=0.0005pt\else\v@leur=\z@\fi%
    \ifdim\v@lX<\z@\advance\v@lX-\v@leur\else\advance\v@lX\v@leur\fi%
    \ifdim\v@lY<\z@\advance\v@lY-\v@leur\else\advance\v@lY\v@leur\fi\fi%
    (\@ffichnb{#1}{\repdecn@mb{\v@lX}},\ifmmode\else\thinspace\fi%
    \@ffichnb{#1}{\repdecn@mb{\v@lY}})}}
\ctr@ld@f\def\@ffichnb#1#2{{\def\@@ffich{\@ffich#1(}\edef\n@mbre{#2}%
    \expandafter\@@ffich\n@mbre)}}
\ctr@ld@f\def\@ffich#1(#2.#3){{#2\ifnum#1>\z@.\fi\def\dig@ts{#3}\s@mme=\z@%
    \loop\ifnum\s@mme<#1\expandafter\@ffichdec\dig@ts:\advance\s@mme\@ne\repeat}}
\ctr@ld@f\def\@ffichdec#1#2:{\relax#1\def\dig@ts{#20}}
\ctr@ld@f\def\figcoordTD#1{{\v@lX=\ptT@unit@\v@lX\v@lY=\ptT@unit@\v@lY\v@lZ=\ptT@unit@\v@lZ%
    \ifr@undcoord\ifcase#1\v@leur=0.5pt\or\v@leur=0.05pt\or\v@leur=0.005pt%
    \or\v@leur=0.0005pt\else\v@leur=\z@\fi%
    \ifdim\v@lX<\z@\advance\v@lX-\v@leur\else\advance\v@lX\v@leur\fi%
    \ifdim\v@lY<\z@\advance\v@lY-\v@leur\else\advance\v@lY\v@leur\fi%
    \ifdim\v@lZ<\z@\advance\v@lZ-\v@leur\else\advance\v@lZ\v@leur\fi\fi%
    (\@ffichnb{#1}{\repdecn@mb{\v@lX}},\ifmmode\else\thinspace\fi%
     \@ffichnb{#1}{\repdecn@mb{\v@lY}},\ifmmode\else\thinspace\fi%
     \@ffichnb{#1}{\repdecn@mb{\v@lZ}})}}
\ctr@ld@f\def\figsetroundcoord#1{\expandafter\Figsetr@undcoord#1:\ignorespaces}
\ctr@ld@f\def\Figsetr@undcoord#1#2:{\if#1n\r@undcoordfalse\else\r@undcoordtrue\fi}
\ctr@ld@f\def\Figsetwr@te#1=#2|{\keln@mun#1|%
    \def\n@mref{m}\ifx\l@debut\n@mref\figsetmark{#2}\else
    \def\n@mref{p}\ifx\l@debut\n@mref\figsetptname{#2}\else
    \def\n@mref{r}\ifx\l@debut\n@mref\figsetroundcoord{#2}\else
    \W@rnmesAttr{figset write}{#1}\fi\fi\fi}
\ctr@ld@f\def\figsetmark#1{\c@nsymb={#1}\setbox\Gb@xSC=\hbox{\the\c@nsymb}\ignorespaces}
\ctr@ln@m\ptn@me
\ctr@ld@f\def\figsetptname#1{\def\ptn@me##1{#1}\ignorespaces}
\ctr@ld@f\def\FigWrit@L#1:#2(#3,#4){\ignorespaces\@keldist\v@leur{#3}\@keldist\delt@{#4}%
    \C@rp@r@m\def\list@num{#1}\@ecfor\p@int:=\list@num\do{\FigWrit@pt{\p@int}{#2}}}
\ctr@ld@f\def\FigWrit@pt#1#2{\FigWp@r@m{#1}{#2}\Vc@rrect\figWp@si%
    \ifdim\wd\Gb@xSC>\z@\b@undb@xSC{\v@lX}{\v@lY}\fi\figWBB@x}
\ctr@ld@f\def\FigWp@r@m#1#2{\Figg@tXY{#1}%
    \setbox\Gb@x=\hbox{\def\t@xt@{#2}\ifx\t@xt@\empty\Figg@tT{#1}\else#2\fi}\c@lprojSP}
\ctr@ld@f\let\Vc@rrect=\relax
\ctr@ld@f\let\C@rp@r@m=\relax
\ctr@ld@f\def\figwrite[#1]#2{{\ignorespaces\def\list@num{#1}\@ecfor\p@int:=\list@num\do{%
    \setbox\Gb@x=\hbox{\def\t@xt@{#2}\ifx\t@xt@\empty\Figg@tT{\p@int}\else#2\fi}%
    \Figwrit@{\p@int}}}\ignorespaces}
\ctr@ld@f\def\Figwrit@#1{\Figg@tXY{#1}\c@lprojSP%
    \rlap{\kern\v@lX\raise\v@lY\hbox{\unhcopy\Gb@x}}\v@leur=\v@lY%
    \advance\v@lY\ht\Gb@x\b@undb@x{\v@lX}{\v@lY}\advance\v@lX\wd\Gb@x%
    \v@lY=\v@leur\advance\v@lY-\dp\Gb@x\b@undb@x{\v@lX}{\v@lY}}
\ctr@ld@f\def\figwritec[#1]#2{{\ignorespaces\def\list@num{#1}%
    \@ecfor\p@int:=\list@num\do{\Figwrit@c{\p@int}{#2}}}\ignorespaces}
\ctr@ld@f\def\Figwrit@c#1#2{\FigWp@r@m{#1}{#2}%
    \rlap{\kern\v@lX\raise\v@lY\hbox{\rlap{\kern-.5\wd\Gb@x%
    \raise-.5\ht\Gb@x\hbox{\raise.5\dp\Gb@x\hbox{\unhcopy\Gb@x}}}}}%
    \v@leur=\ht\Gb@x\advance\v@leur\dp\Gb@x%
    \advance\v@lX-.5\wd\Gb@x\advance\v@lY-.5\v@leur\b@undb@x{\v@lX}{\v@lY}%
    \advance\v@lX\wd\Gb@x\advance\v@lY\v@leur\b@undb@x{\v@lX}{\v@lY}}
\ctr@ld@f\def\figwritep[#1]{{\ignorespaces\def\list@num{#1}\setbox\Gb@x=\hbox{\c@nterpt}%
    \@ecfor\p@int:=\list@num\do{\Figwrit@{\p@int}}}\ignorespaces}
\ctr@ld@f\def\figwritew#1:#2(#3){\figwritegcw#1:{#2}(#3,0pt)}
\ctr@ld@f\def\figwritee#1:#2(#3){\figwritegce#1:{#2}(#3,0pt)}
\ctr@ld@f\def\figwriten#1:#2(#3){{\def\Vc@rrect{\v@lZ=\v@leur\advance\v@lZ\dp\Gb@x}%
    \Figwrit@NS#1:{#2}(#3)}\ignorespaces}
\ctr@ld@f\def\figwrites#1:#2(#3){{\def\Vc@rrect{\v@lZ=-\v@leur\advance\v@lZ-\ht\Gb@x}%
    \Figwrit@NS#1:{#2}(#3)}\ignorespaces}
\ctr@ld@f\def\Figwrit@NS#1:#2(#3){\let\figWp@si=\FigWp@siNS\let\figWBB@x=\FigWBB@xNS%
    \FigWrit@L#1:{#2}(#3,0pt)}
\ctr@ld@f\def\FigWp@siNS{\rlap{\kern\v@lX\raise\v@lY\hbox{\rlap{\kern-.5\wd\Gb@x%
    \raise\v@lZ\hbox{\unhcopy\Gb@x}}\c@nterpt}}}
\ctr@ld@f\def\FigWBB@xNS{\advance\v@lY\v@lZ%
    \advance\v@lY-\dp\Gb@x\advance\v@lX-.5\wd\Gb@x\b@undb@x{\v@lX}{\v@lY}%
    \advance\v@lY\ht\Gb@x\advance\v@lY\dp\Gb@x%
    \advance\v@lX\wd\Gb@x\b@undb@x{\v@lX}{\v@lY}}
\ctr@ld@f\def\figwritenw#1:#2(#3){{\let\figWp@si=\FigWp@sigW\let\figWBB@x=\FigWBB@xgWE%
    \def\C@rp@r@m{\v@leur=\unssqrttw@\v@leur\delt@=\v@leur%
    \ifdim\delt@=\z@\delt@=\epsil@n\fi}\let@xte={-}\FigWrit@L#1:{#2}(#3,0pt)}\ignorespaces}
\ctr@ld@f\def\figwritesw#1:#2(#3){{\let\figWp@si=\FigWp@sigW\let\figWBB@x=\FigWBB@xgWE%
    \def\C@rp@r@m{\v@leur=\unssqrttw@\v@leur\delt@=-\v@leur%
    \ifdim\delt@=\z@\delt@=-\epsil@n\fi}\let@xte={-}\FigWrit@L#1:{#2}(#3,0pt)}\ignorespaces}
\ctr@ld@f\def\figwritene#1:#2(#3){{\let\figWp@si=\FigWp@sigE\let\figWBB@x=\FigWBB@xgWE%
    \def\C@rp@r@m{\v@leur=\unssqrttw@\v@leur\delt@=\v@leur%
    \ifdim\delt@=\z@\delt@=\epsil@n\fi}\let@xte={}\FigWrit@L#1:{#2}(#3,0pt)}\ignorespaces}
\ctr@ld@f\def\figwritese#1:#2(#3){{\let\figWp@si=\FigWp@sigE\let\figWBB@x=\FigWBB@xgWE%
    \def\C@rp@r@m{\v@leur=\unssqrttw@\v@leur\delt@=-\v@leur%
    \ifdim\delt@=\z@\delt@=-\epsil@n\fi}\let@xte={}\FigWrit@L#1:{#2}(#3,0pt)}\ignorespaces}
\ctr@ld@f\def\figwritegw#1:#2(#3,#4){{\let\figWp@si=\FigWp@sigW\let\figWBB@x=\FigWBB@xgWE%
    \let@xte={-}\FigWrit@L#1:{#2}(#3,#4)}\ignorespaces}
\ctr@ld@f\def\figwritege#1:#2(#3,#4){{\let\figWp@si=\FigWp@sigE\let\figWBB@x=\FigWBB@xgWE%
    \let@xte={}\FigWrit@L#1:{#2}(#3,#4)}\ignorespaces}
\ctr@ld@f\def\FigWp@sigW{\v@lXa=\z@\v@lYa=\ht\Gb@x\advance\v@lYa\dp\Gb@x%
    \ifdim\delt@>\z@\relax%
    \rlap{\kern\v@lX\raise\v@lY\hbox{\rlap{\kern-\wd\Gb@x\kern-\v@leur%
          \raise\delt@\hbox{\raise\dp\Gb@x\hbox{\unhcopy\Gb@x}}}\c@nterpt}}%
    \else\ifdim\delt@<\z@\relax\v@lYa=-\v@lYa%
    \rlap{\kern\v@lX\raise\v@lY\hbox{\rlap{\kern-\wd\Gb@x\kern-\v@leur%
          \raise\delt@\hbox{\raise-\ht\Gb@x\hbox{\unhcopy\Gb@x}}}\c@nterpt}}%
    \else\v@lXa=-.5\v@lYa%
    \rlap{\kern\v@lX\raise\v@lY\hbox{\rlap{\kern-\wd\Gb@x\kern-\v@leur%
          \raise-.5\ht\Gb@x\hbox{\raise.5\dp\Gb@x\hbox{\unhcopy\Gb@x}}}\c@nterpt}}%
    \fi\fi}
\ctr@ld@f\def\FigWp@sigE{\v@lXa=\z@\v@lYa=\ht\Gb@x\advance\v@lYa\dp\Gb@x%
    \ifdim\delt@>\z@\relax%
    \rlap{\kern\v@lX\raise\v@lY\hbox{\c@nterpt\kern\v@leur%
          \raise\delt@\hbox{\raise\dp\Gb@x\hbox{\unhcopy\Gb@x}}}}%
    \else\ifdim\delt@<\z@\relax\v@lYa=-\v@lYa%
    \rlap{\kern\v@lX\raise\v@lY\hbox{\c@nterpt\kern\v@leur%
          \raise\delt@\hbox{\raise-\ht\Gb@x\hbox{\unhcopy\Gb@x}}}}%
    \else\v@lXa=-.5\v@lYa%
    \rlap{\kern\v@lX\raise\v@lY\hbox{\c@nterpt\kern\v@leur%
          \raise-.5\ht\Gb@x\hbox{\raise.5\dp\Gb@x\hbox{\unhcopy\Gb@x}}}}%
    \fi\fi}
\ctr@ld@f\def\FigWBB@xgWE{\advance\v@lY\delt@%
    \advance\v@lX\the\let@xte\v@leur\advance\v@lY\v@lXa\b@undb@x{\v@lX}{\v@lY}%
    \advance\v@lX\the\let@xte\wd\Gb@x\advance\v@lY\v@lYa\b@undb@x{\v@lX}{\v@lY}}
\ctr@ld@f\def\figwritegcw#1:#2(#3,#4){{\let\figWp@si=\FigWp@sigcW\let\figWBB@x=\FigWBB@xgcWE%
    \let@xte={-}\FigWrit@L#1:{#2}(#3,#4)}\ignorespaces}
\ctr@ld@f\def\figwritegce#1:#2(#3,#4){{\let\figWp@si=\FigWp@sigcE\let\figWBB@x=\FigWBB@xgcWE%
    \let@xte={}\FigWrit@L#1:{#2}(#3,#4)}\ignorespaces}
\ctr@ld@f\def\FigWp@sigcW{\rlap{\kern\v@lX\raise\v@lY\hbox{\rlap{\kern-\wd\Gb@x\kern-\v@leur%
     \raise-.5\ht\Gb@x\hbox{\raise\delt@\hbox{\raise.5\dp\Gb@x\hbox{\unhcopy\Gb@x}}}}%
     \c@nterpt}}}
\ctr@ld@f\def\FigWp@sigcE{\rlap{\kern\v@lX\raise\v@lY\hbox{\c@nterpt\kern\v@leur%
    \raise-.5\ht\Gb@x\hbox{\raise\delt@\hbox{\raise.5\dp\Gb@x\hbox{\unhcopy\Gb@x}}}}}}
\ctr@ld@f\def\FigWBB@xgcWE{\v@lZ=\ht\Gb@x\advance\v@lZ\dp\Gb@x%
    \advance\v@lX\the\let@xte\v@leur\advance\v@lY\delt@\advance\v@lY.5\v@lZ%
    \b@undb@x{\v@lX}{\v@lY}%
    \advance\v@lX\the\let@xte\wd\Gb@x\advance\v@lY-\v@lZ\b@undb@x{\v@lX}{\v@lY}}
\ctr@ld@f\def\figwritebn#1:#2(#3){{\def\Vc@rrect{\v@lZ=\v@leur}\Figwrit@NS#1:{#2}(#3)}\ignorespaces}
\ctr@ld@f\def\figwritebs#1:#2(#3){{\def\Vc@rrect{\v@lZ=-\v@leur}\Figwrit@NS#1:{#2}(#3)}\ignorespaces}
\ctr@ld@f\def\figwritebw#1:#2(#3){{\let\figWp@si=\FigWp@sibW\let\figWBB@x=\FigWBB@xbWE%
    \let@xte={-}\FigWrit@L#1:{#2}(#3,0pt)}\ignorespaces}
\ctr@ld@f\def\figwritebe#1:#2(#3){{\let\figWp@si=\FigWp@sibE\let\figWBB@x=\FigWBB@xbWE%
    \let@xte={}\FigWrit@L#1:{#2}(#3,0pt)}\ignorespaces}
\ctr@ld@f\def\FigWp@sibW{\rlap{\kern\v@lX\raise\v@lY\hbox{\rlap{\kern-\wd\Gb@x\kern-\v@leur%
          \hbox{\unhcopy\Gb@x}}\c@nterpt}}}
\ctr@ld@f\def\FigWp@sibE{\rlap{\kern\v@lX\raise\v@lY\hbox{\c@nterpt\kern\v@leur%
          \hbox{\unhcopy\Gb@x}}}}
\ctr@ld@f\def\FigWBB@xbWE{\v@lZ=\ht\Gb@x\advance\v@lZ\dp\Gb@x%
    \advance\v@lX\the\let@xte\v@leur\advance\v@lY\ht\Gb@x\b@undb@x{\v@lX}{\v@lY}%
    \advance\v@lX\the\let@xte\wd\Gb@x\advance\v@lY-\v@lZ\b@undb@x{\v@lX}{\v@lY}}
\ctr@ln@w{newread}\frf@g  \ctr@ln@w{newwrite}\fwf@g
\ctr@ln@w{newif}\ifCUR@PS
\ctr@ln@w{newif}\ifGR@cri
\ctr@ln@w{newif}\ifUse@llipse
\ctr@ln@w{newif}\ifGRdebugm@de \GRdebugm@defalse 
\ctr@ln@w{newif}\ifPDFm@ke
\ifx\pdfliteral\undefined\else\ifnum\pdfoutput>\z@\PDFm@ketrue\fi\fi
\ctr@ld@f\def\initPDF@rDVI{%
\ifPDFm@ke
 \let\figscan=\figscan@E
 \let\newGr@FN=\newGr@FNPDF
 \ctr@ld@f\def\c@mcurveto{c}
 \ctr@ld@f\def\c@mfill{f}
 \ctr@ld@f\def\c@mgsave{q}
 \ctr@ld@f\def\c@mgrestore{Q}
 \ctr@ld@f\def\c@mlineto{l}
 \ctr@ld@f\def\c@mmoveto{m}
 \ctr@ld@f\def\c@msetgray{g}     \ctr@ld@f\def\c@msetgrayStroke{G}
 \ctr@ld@f\def\c@msetcmykcolor{k}\ctr@ld@f\def\c@msetcmykcolorStroke{K}
 \ctr@ld@f\def\c@msetrgbcolor{rg}\ctr@ld@f\def\c@msetrgbcolorStroke{RG}
 \ctr@ld@f\def\d@fprimarC@lor{\CUR@color\space\CUR@colorc@md%
               \space\CUR@color\space\CUR@colorc@mdStroke}
 \ctr@ld@f\def\c@msetdash{d}
 \ctr@ld@f\def\c@msetlinejoin{j}
 \ctr@ld@f\def\c@msetlinewidth{w}
 \ctr@ld@f\def\f@gclosestroke{\immediate\write\fwf@g{s}}
 \ctr@ld@f\def\f@gfill{\immediate\write\fwf@g{\fillc@md}}
 \ctr@ld@f\def\f@gnewpath{}
 \ctr@ld@f\def\f@gstroke{\immediate\write\fwf@g{S}}
\else
 \let\figinsertE=\figinsert
 \let\newGr@FN=\newGr@FNDVI
 \ctr@ld@f\def\c@mcurveto{curveto}
 \ctr@ld@f\def\c@mfill{fill}
 \ctr@ld@f\def\c@mgsave{gsave}
 \ctr@ld@f\def\c@mgrestore{grestore}
 \ctr@ld@f\def\c@mlineto{lineto}
 \ctr@ld@f\def\c@mmoveto{moveto}
 \ctr@ld@f\def\c@msetgray{setgray}          \ctr@ld@f\def\c@msetgrayStroke{}
 \ctr@ld@f\def\c@msetcmykcolor{setcmykcolor}\ctr@ld@f\def\c@msetcmykcolorStroke{}
 \ctr@ld@f\def\c@msetrgbcolor{setrgbcolor}  \ctr@ld@f\def\c@msetrgbcolorStroke{}
 \ctr@ld@f\def\d@fprimarC@lor{\CUR@color\space\CUR@colorc@md}
 \ctr@ld@f\def\c@msetdash{setdash}
 \ctr@ld@f\def\c@msetlinejoin{setlinejoin}
 \ctr@ld@f\def\c@msetlinewidth{setlinewidth}
 \ctr@ld@f\def\f@gclosestroke{\immediate\write\fwf@g{closepath\space stroke}}
 \ctr@ld@f\def\f@gfill{\immediate\write\fwf@g{\fillc@md}}
 \ctr@ld@f\def\f@gnewpath{\immediate\write\fwf@g{newpath}}
 \ctr@ld@f\def\f@gstroke{\immediate\write\fwf@g{stroke}}
\fi}
\ctr@ld@f\def\c@pypsfile#1#2{\c@pyfil@{\immediate\write#1}{#2}}
\ctr@ld@f\def\Figinclud@PDF#1#2{\openin\frf@g=#1\pdfliteral{q #2 0 0 #2 0 0 cm}%
    \c@pyfil@{\pdfliteral}{\frf@g}\pdfliteral{Q}\closein\frf@g}
\ctr@ln@w{newif}\ifmored@ta
\ctr@ln@m\bl@nkline
\ctr@ld@f\def\c@pyfil@#1#2{\def\bl@nkline{\par}{\catcode`\%=12
    \loop\ifeof#2\mored@tafalse\else\mored@tatrue\immediate\read#2 to\tr@c
    \ifx\tr@c\bl@nkline\else#1{\tr@c}\fi\fi\ifmored@ta\repeat}}
\ctr@ld@f\def\keln@mun#1#2|{\def\l@debut{#1}\def\l@suite{#2}}
\ctr@ld@f\def\keln@mde#1#2#3|{\def\l@debut{#1#2}\def\l@suite{#3}}
\ctr@ld@f\def\keln@mtr#1#2#3#4|{\def\l@debut{#1#2#3}\def\l@suite{#4}}
\ctr@ld@f\def\keln@mqu#1#2#3#4#5|{\def\l@debut{#1#2#3#4}\def\l@suite{#5}}
\ctr@ld@f\let\@psffilein=\frf@g 
\ctr@ln@w{newif}\if@psffileok    
\ctr@ln@w{newif}\if@psfbbfound   
\ctr@ln@w{newif}\if@psfverbose   
\@psfverbosetrue
\ctr@ln@m\@psfllx \ctr@ln@m\@psflly
\ctr@ln@m\@psfurx \ctr@ln@m\@psfury
\ctr@ln@m\resetcolonc@tcode
\ctr@ld@f\def\@psfgetbb#1{\global\@psfbbfoundfalse%
\global\def\@psfllx{0}\global\def\@psflly{0}%
\global\def\@psfurx{30}\global\def\@psfury{30}%
\openin\@psffilein=#1\relax
\ifeof\@psffilein\errmessage{I couldn't open #1, will ignore it}\else
   \edef\resetcolonc@tcode{\catcode`\noexpand\:\the\catcode`\:\relax}%
   {\@psffileoktrue \chardef\other=12
    \def\do##1{\catcode`##1=\other}\dospecials \catcode`\ =10 \resetcolonc@tcode
    \loop
       \read\@psffilein to \@psffileline
       \ifeof\@psffilein\@psffileokfalse\else
          \expandafter\@psfaux\@psffileline:. \\%
       \fi
   \if@psffileok\repeat
   \if@psfbbfound\else
    \if@psfverbose\message{No bounding box comment in #1; using defaults}\fi\fi
   }\closein\@psffilein\fi}%
\ctr@ln@m\@psfbblit
\ctr@ln@m\@psfpercent
{\catcode`\%=12 \global\let\@psfpercent=
\ctr@ln@m\@psfaux
\long\def\@psfaux#1#2:#3\\{\ifx#1\@psfpercent
   \def\testit{#2}\ifx\testit\@psfbblit
      \@psfgrab #3 . . . \\%
      \@psffileokfalse
      \global\@psfbbfoundtrue
   \fi\else\ifx#1\par\else\@psffileokfalse\fi\fi}%
\ctr@ld@f\def\@psfempty{}%
\ctr@ld@f\def\@psfgrab #1 #2 #3 #4 #5\\{%
\global\def\@psfllx{#1}\ifx\@psfllx\@psfempty
      \@psfgrab #2 #3 #4 #5 .\\\else
   \global\def\@psflly{#2}%
   \global\def\@psfurx{#3}\global\def\@psfury{#4}\fi}%
\ctr@ld@f\def\PSwrit@cmd#1#2#3{{\Figg@tXY{#1}\c@lprojSP\b@undb@x{\v@lX}{\v@lY}%
    \v@lX=\ptT@ptps\v@lX\v@lY=\ptT@ptps\v@lY%
    \immediate\write#3{\repdecn@mb{\v@lX}\space\repdecn@mb{\v@lY}\space#2}}}
\ctr@ld@f\def\PSwrit@cmdS#1#2#3#4#5{{\Figg@tXY{#1}\c@lprojSP\b@undb@x{\v@lX}{\v@lY}%
    \global\result@t=\v@lX\global\result@@t=\v@lY%
    \v@lX=\ptT@ptps\v@lX\v@lY=\ptT@ptps\v@lY%
    \immediate\write#3{\repdecn@mb{\v@lX}\space\repdecn@mb{\v@lY}\space#2}}%
    \edef#4{\the\result@t}\edef#5{\the\result@@t}}
\ctr@ld@f\def\update@ttr#1#2#3{\Figdisc@rdLTS{#3}{\n@mref}%
    \ifx\n@mref\D@FTref#2{#1}\else#2{#3}\fi}
\ctr@ld@f\def\D@FTref{default}
\ctr@ld@f\def\W@rnmesAttr#1#2{%
    \immediate\write16{*** Unknown attribute: \BS@ #1(..., #2=...)}}
\ctr@ld@f\def\W@rnmeskwd#1#2{%
    \immediate\write16{*** Unknown keyword #2 in \BS@ #1}}
\ctr@ld@f\def\W@rnmesIgn#1{\immediate\write16{*** \BS@ #1 is ignored inside a
     \BS@ figdrawbegin-\BS@ figdrawend block.}}
\ctr@ld@f\def\Psset@lti#1=#2|{\keln@mtr#1|%
    \def\n@mref{blc}\ifx\l@debut\n@mref\update@ttr\D@FTref\P@setblcolor{#2}\else
    \def\n@mref{bld}\ifx\l@debut\n@mref\update@ttr\D@FTref\P@setbldash{#2}\else
    \def\n@mref{blw}\ifx\l@debut\n@mref\update@ttr\D@FTref\P@setblwidth{#2}\else
    \def\n@mref{sqc}\ifx\l@debut\n@mref\update@ttr\D@FTref\P@setsqcolor{#2}\else
    \def\n@mref{sqd}\ifx\l@debut\n@mref\update@ttr\D@FTref\P@setsqdash{#2}\else
    \def\n@mref{sqw}\ifx\l@debut\n@mref\update@ttr\D@FTref\P@setsqwidth{#2}\else
    \W@rnmesAttr{figset altitude}{#1}\fi\fi\fi\fi\fi\fi}
\ctr@ln@m\DDV@blcolor
\ctr@ld@f\def\P@setblcolor#1{\edef\DDV@blcolor{#1}}
\ctr@ln@m\DDV@bldash
\ctr@ld@f\def\P@setbldash#1{\edef\DDV@bldash{#1}}
\ctr@ln@m\DDV@blwidth
\ctr@ld@f\def\P@setblwidth#1{\edef\DDV@blwidth{#1}}
\ctr@ln@m\DDV@sqcolor
\ctr@ld@f\def\P@setsqcolor#1{\edef\DDV@sqcolor{#1}}
\ctr@ln@m\DDV@sqdash
\ctr@ld@f\def\P@setsqdash#1{\edef\DDV@sqdash{#1}}
\ctr@ln@m\DDV@sqwidth
\ctr@ld@f\def\P@setsqwidth#1{\edef\DDV@sqwidth{#1}}
\ctr@ld@f\def\figdrawaltitude#1[#2,#3,#4]{{\ifCUR@PS\ifGR@cri%
    \PSc@mment{altitude Square Dim=#1, Triangle=[#2 / #3,#4]}%
    \s@uvc@ntr@l\et@tpsaltitude\resetc@ntr@l{2}\figptorthoprojline-5:=#2/#3,#4/%
    \figvectP -1[#3,#4]\n@rminf{\v@leur}{-1}\vecunit@{-3}{-1}%
    \figvectP -1[-5,#3]\n@rminf{\v@lmin}{-1}\figvectP -2[-5,#4]\n@rminf{\v@lmax}{-2}%
    \ifdim\v@lmin<\v@lmax\s@mme=#3\else\v@lmax=\v@lmin\s@mme=#4\fi%
    \figvectP -4[-5,#2]\vecunit@{-4}{-4}\delt@=#1\unit@%
    \edef\t@ille{\repdecn@mb{\delt@}}\figpttra-1:=-5/\t@ille,-3/%
    \figptstra-3=-5,-1/\t@ille,-4/\figdrawline[#2,-5]%
    \Pss@tspecifSt{color=\DDV@sqcolor,dash=\DDV@sqdash,width=\DDV@sqwidth}%
    \figdrawline[-1,-2,-3]%
    \Psrest@reSt{color=\DDV@sqcolor,dash=\DDV@sqdash,width=\DDV@sqwidth}%
    \ifdim\v@leur<\v@lmax%
    \Pss@tspecifSt{color=\DDV@blcolor,dash=\DDV@bldash,width=\DDV@blwidth}%
    \figdrawline[-5,\the\s@mme]%
    \Psrest@reSt{color=\DDV@blcolor,dash=\DDV@bldash,width=\DDV@blwidth}%
    \fi\PSc@mment{End altitude}\resetc@ntr@l\et@tpsaltitude\fi\fi}}
\ctr@ld@f\def\Ps@rcerc#1;#2(#3,#4){\ellBB@x#1;#2,#2(#3,#4,0)%
    \f@gnewpath{\delt@=#2\unit@\delt@=\ptT@ptps\delt@%
    \BdingB@xfalse%
    \PSwrit@cmd{#1}{\repdecn@mb{\delt@}\space #3\space #4\space arc}{\fwf@g}}}
\ctr@ln@m\figdrawarccirc
\ctr@ld@f\def\Q@arccircDD#1;#2(#3,#4){\ifCUR@PS\ifGR@cri%
    \PSc@mment{arccircDD Center=#1 ; Radius=#2 (Ang1=#3, Ang2=#4)}%
    \iffillm@de\Ps@rcerc#1;#2(#3,#4)%
    \f@gfill%
    \else\Ps@rcerc#1;#2(#3,#4)\f@gstroke\fi%
    \PSc@mment{End arccircDD}\fi\fi}
\ctr@ld@f\def\Q@arccircTD#1,#2,#3;#4(#5,#6){{\ifCUR@PS\ifGR@cri\s@uvc@ntr@l\et@tpsarccircTD%
    \PSc@mment{arccircTD Center=#1,P1=#2,P2=#3 ; Radius=#4 (Ang1=#5, Ang2=#6)}%
    \setc@ntr@l{2}\c@lExtAxes#1,#2,#3(#4)\Q@arcellPATD#1,-4,-5(#5,#6)%
    \PSc@mment{End arccircTD}\resetc@ntr@l\et@tpsarccircTD\fi\fi}}
\ctr@ld@f\def\c@lExtAxes#1,#2,#3(#4){%
    \figvectPTD-5[#1,#2]\vecunit@{-5}{-5}\figvectNTD-4[#1,#2,#3]\vecunit@{-4}{-4}%
    \figvectNVTD-3[-4,-5]\delt@=#4\unit@\edef\r@yon{\repdecn@mb{\delt@}}%
    \figpttra-4:=#1/\r@yon,-5/\figpttra-5:=#1/\r@yon,-3/}
\ctr@ln@m\figdrawarccircP
\ctr@ld@f\def\Q@arccircPDD#1;#2[#3,#4]{{\ifCUR@PS\ifGR@cri\s@uvc@ntr@l\et@tpsarccircPDD%
    \PSc@mment{arccircPDD Center=#1; Radius=#2, [P1=#3, P2=#4]}%
    \Ps@ngleparam#1;#2[#3,#4]\ifdim\v@lmin>\v@lmax\advance\v@lmax\DePI@deg\fi%
    \edef\@ngdeb{\repdecn@mb{\v@lmin}}\edef\@ngfin{\repdecn@mb{\v@lmax}}%
    \figdrawarccirc#1;\r@dius(\@ngdeb,\@ngfin)%
    \PSc@mment{End arccircPDD}\resetc@ntr@l\et@tpsarccircPDD\fi\fi}}
\ctr@ld@f\def\Q@arccircPTD#1;#2[#3,#4,#5]{{\ifCUR@PS\ifGR@cri\s@uvc@ntr@l\et@tpsarccircPTD%
    \PSc@mment{arccircPTD Center=#1; Radius=#2, [P1=#3, P2=#4, P3=#5]}%
    \setc@ntr@l{2}\c@lExtAxes#1,#3,#5(#2)\figdrawarcellPP#1,-4,-5[#3,#4]%
    \PSc@mment{End arccircPTD}\resetc@ntr@l\et@tpsarccircPTD\fi\fi}}
\ctr@ld@f\def\Ps@ngleparam#1;#2[#3,#4]{\setc@ntr@l{2}%
    \figvectPDD-1[#1,#3]\vecunit@{-1}{-1}\Figg@tXY{-1}\arct@n\v@lmin(\v@lX,\v@lY)%
    \figvectPDD-2[#1,#4]\vecunit@{-2}{-2}\Figg@tXY{-2}\arct@n\v@lmax(\v@lX,\v@lY)%
    \v@lmin=\rdT@deg\v@lmin\v@lmax=\rdT@deg\v@lmax%
    \v@leur=#2pt\maxim@m{\mili@u}{-\v@leur}{\v@leur}%
    \edef\r@dius{\repdecn@mb{\mili@u}}}
\ctr@ld@f\def\Ps@rcercBz#1;#2(#3,#4){\Ps@rellBz#1;#2,#2(#3,#4,0)}
\ctr@ld@f\def\Ps@rellBz#1;#2,#3(#4,#5,#6){%
    \ellBB@x#1;#2,#3(#4,#5,#6)\BdingB@xfalse%
    \c@lNbarcs{#4}{#5}\v@leur=#4pt\setc@ntr@l{2}\figptell-13::#1;#2,#3(#4,#6)%
    \f@gnewpath\PSwrit@cmd{-13}{\c@mmoveto}{\fwf@g}%
    \s@mme=\z@\bcl@rellBz#1;#2,#3(#6)\BdingB@xtrue}
\ctr@ld@f\def\bcl@rellBz#1;#2,#3(#4){\relax%
    \ifnum\s@mme<\p@rtent\advance\s@mme\@ne%
    \advance\v@leur\delt@\edef\@ngle{\repdecn@mb\v@leur}\figptell-14::#1;#2,#3(\@ngle,#4)%
    \advance\v@leur\delt@\edef\@ngle{\repdecn@mb\v@leur}\figptell-15::#1;#2,#3(\@ngle,#4)%
    \advance\v@leur\delt@\edef\@ngle{\repdecn@mb\v@leur}\figptell-16::#1;#2,#3(\@ngle,#4)%
    \figptscontrolDD-18[-13,-14,-15,-16]%
    \PSwrit@cmd{-18}{}{\fwf@g}\PSwrit@cmd{-17}{}{\fwf@g}%
    \PSwrit@cmd{-16}{\c@mcurveto}{\fwf@g}%
    \figptcopyDD-13:/-16/\bcl@rellBz#1;#2,#3(#4)\fi}
\ctr@ld@f\def\Ps@rell#1;#2,#3(#4,#5,#6){\ellBB@x#1;#2,#3(#4,#5,#6)%
    \f@gnewpath{\v@lmin=#2\unit@\v@lmin=\ptT@ptps\v@lmin%
    \v@lmax=#3\unit@\v@lmax=\ptT@ptps\v@lmax\BdingB@xfalse%
    \PSwrit@cmd{#1}%
    {#6\space\repdecn@mb{\v@lmin}\space\repdecn@mb{\v@lmax}\space #4\space #5\space ellipse}{\fwf@g}}%
    \global\Use@llipsetrue}
\ctr@ln@m\figdrawarcell
\ctr@ld@f\def\Q@arcellDD#1;#2,#3(#4,#5,#6){{\ifCUR@PS\ifGR@cri%
    \PSc@mment{arcellDD Center=#1 ; XRad=#2, YRad=#3 (Ang1=#4, Ang2=#5, Inclination=#6)}%
    \iffillm@de\Ps@rell#1;#2,#3(#4,#5,#6)%
    \f@gfill%
    \else\Ps@rell#1;#2,#3(#4,#5,#6)\f@gstroke\fi%
    \PSc@mment{End arcellDD}\fi\fi}}
\ctr@ld@f\def\Q@arcellTD#1;#2,#3(#4,#5,#6){{\ifCUR@PS\ifGR@cri\s@uvc@ntr@l\et@tpsarcellTD%
    \PSc@mment{arcellTD Center=#1 ; XRad=#2, YRad=#3 (Ang1=#4, Ang2=#5, Inclination=#6)}%
    \setc@ntr@l{2}\figpttraC -8:=#1/#2,0,0/\figpttraC -7:=#1/0,#3,0/%
    \figvectC -4(0,0,1)\figptsrot -8=-8,-7/#1,#6,-4/\Q@arcellPATD#1,-8,-7(#4,#5)%
    \PSc@mment{End arcellTD}\resetc@ntr@l\et@tpsarcellTD\fi\fi}}
\ctr@ln@m\figdrawarcellPA
\ctr@ld@f\def\Q@arcellPADD#1,#2,#3(#4,#5){{\ifCUR@PS\ifGR@cri\s@uvc@ntr@l\et@tpsarcellPADD%
    \PSc@mment{arcellPADD Center=#1,PtAxis1=#2,PtAxis2=#3 (Ang1=#4, Ang2=#5)}%
    \setc@ntr@l{2}\figvectPDD-1[#1,#2]\vecunit@DD{-1}{-1}\v@lX=\ptT@unit@\result@t%
    \edef\XR@d{\repdecn@mb{\v@lX}}\Figg@tXY{-1}\arct@n\v@lmin(\v@lX,\v@lY)%
    \v@lmin=\rdT@deg\v@lmin\edef\Inclin@{\repdecn@mb{\v@lmin}}%
    \figgetdist\YR@d[#1,#3]\Q@arcellDD#1;\XR@d,\YR@d(#4,#5,\Inclin@)%
    \PSc@mment{End arcellPADD}\resetc@ntr@l\et@tpsarcellPADD\fi\fi}}
\ctr@ld@f\def\Q@arcellPATD#1,#2,#3(#4,#5){{\ifCUR@PS\ifGR@cri\s@uvc@ntr@l\et@tpsarcellPATD%
    \PSc@mment{arcellPATD Center=#1,PtAxis1=#2,PtAxis2=#3 (Ang1=#4, Ang2=#5)}%
    \iffillm@de\Ps@rellPATD#1,#2,#3(#4,#5)%
    \f@gfill%
    \else\Ps@rellPATD#1,#2,#3(#4,#5)\f@gstroke\fi%
    \PSc@mment{End arcellPATD}\resetc@ntr@l\et@tpsarcellPATD\fi\fi}}
\ctr@ld@f\def\Ps@rellPATD#1,#2,#3(#4,#5){\let\c@lprojSP=\relax%
    \setc@ntr@l{2}\figvectPTD-1[#1,#2]\figvectPTD-2[#1,#3]\c@lNbarcs{#4}{#5}%
    \v@leur=#4pt\c@lptellP{#1}{-1}{-2}\Figptpr@j-5:/-3/%
    \f@gnewpath\PSwrit@cmdS{-5}{\c@mmoveto}{\fwf@g}{\X@un}{\Y@un}%
    \edef\C@nt@r{#1}\s@mme=\z@\bcl@rellPATD}
\ctr@ld@f\def\bcl@rellPATD{\relax%
    \ifnum\s@mme<\p@rtent\advance\s@mme\@ne%
    \advance\v@leur\delt@\c@lptellP{\C@nt@r}{-1}{-2}\Figptpr@j-4:/-3/%
    \advance\v@leur\delt@\c@lptellP{\C@nt@r}{-1}{-2}\Figptpr@j-6:/-3/%
    \advance\v@leur\delt@\c@lptellP{\C@nt@r}{-1}{-2}\Figptpr@j-3:/-3/%
    \v@lX=\z@\v@lY=\z@\Figtr@nptDD{-5}{-5}\Figtr@nptDD{2}{-3}%
    \divide\v@lX\@vi\divide\v@lY\@vi%
    \Figtr@nptDD{3}{-4}\Figtr@nptDD{-1.5}{-6}\v@lmin=\v@lX\v@lmax=\v@lY%
    \v@lX=\z@\v@lY=\z@\Figtr@nptDD{2}{-5}\Figtr@nptDD{-5}{-3}%
    \divide\v@lX\@vi\divide\v@lY\@vi\Figtr@nptDD{-1.5}{-4}\Figtr@nptDD{3}{-6}%
    \BdingB@xfalse%
    \Figp@intregDD-4:(\v@lmin,\v@lmax)\PSwrit@cmdS{-4}{}{\fwf@g}{\X@de}{\Y@de}%
    \Figp@intregDD-4:(\v@lX,\v@lY)\PSwrit@cmdS{-4}{}{\fwf@g}{\X@tr}{\Y@tr}%
    \BdingB@xtrue\PSwrit@cmdS{-3}{\c@mcurveto}{\fwf@g}{\X@qu}{\Y@qu}%
    \B@zierBB@x{1}{\Y@un}(\X@un,\X@de,\X@tr,\X@qu)%
    \B@zierBB@x{2}{\X@un}(\Y@un,\Y@de,\Y@tr,\Y@qu)%
    \edef\X@un{\X@qu}\edef\Y@un{\Y@qu}\figptcopyDD-5:/-3/\bcl@rellPATD\fi}
\ctr@ld@f\def\c@lNbarcs#1#2{%
    \delt@=#2pt\advance\delt@-#1pt\maxim@m{\v@lmax}{\delt@}{-\delt@}%
    \v@leur=\v@lmax\divide\v@leur45 \p@rtentiere{\p@rtent}{\v@leur}\advance\p@rtent\@ne%
    \s@mme=\p@rtent\multiply\s@mme\thr@@\divide\delt@\s@mme}
\ctr@ld@f\def\figdrawarcellPP#1,#2,#3[#4,#5]{{\ifCUR@PS\ifGR@cri\s@uvc@ntr@l\et@tpsarcellPP%
    \PSc@mment{arcellPP Center=#1,PtAxis1=#2,PtAxis2=#3 [Point1=#4, Point2=#5]}%
    \setc@ntr@l{2}\figvectP-2[#1,#3]\vecunit@{-2}{-2}\v@lmin=\result@t%
    \invers@{\v@lmax}{\v@lmin}%
    \figvectP-1[#1,#2]\vecunit@{-1}{-1}\v@leur=\result@t%
    \v@leur=\repdecn@mb{\v@lmax}\v@leur\edef\AsB@{\repdecn@mb{\v@leur}}
    \c@lAngle{#1}{#4}{\v@lmin}\edef\@ngdeb{\repdecn@mb{\v@lmin}}%
    \c@lAngle{#1}{#5}{\v@lmax}\ifdim\v@lmin>\v@lmax\advance\v@lmax\DePI@deg\fi%
    \edef\@ngfin{\repdecn@mb{\v@lmax}}\figdrawarcellPA#1,#2,#3(\@ngdeb,\@ngfin)%
    \PSc@mment{End arcellPP}\resetc@ntr@l\et@tpsarcellPP\fi\fi}}
\ctr@ld@f\def\c@lAngle#1#2#3{\figvectP-3[#1,#2]%
    \c@lproscal\delt@[-3,-1]\c@lproscal\v@leur[-3,-2]%
    \v@leur=\AsB@\v@leur\arct@n#3(\delt@,\v@leur)#3=\rdT@deg#3}
\ctr@ln@w{newif}\if@rrowratio\@rrowratiotrue
\ctr@ln@w{newif}\if@rrowhfill
\ctr@ln@w{newif}\if@rrowhout
\ctr@ld@f\def\Psset@rrowhe@d#1=#2|{\keln@mun#1|%
    \def\n@mref{a}\ifx\l@debut\n@mref\update@ttr\D@FTarrowheadangle\Q@s@tarrowheadangle{#2}\else
    \def\n@mref{f}\ifx\l@debut\n@mref\update@ttr\D@FTarrowheadfill\Q@s@tarrowheadfill{#2}\else
    \def\n@mref{l}\ifx\l@debut\n@mref\update@ttr\D@FTarrowheadlength\Q@s@tarrowheadlength{#2}\else
    \def\n@mref{o}\ifx\l@debut\n@mref\update@ttr\D@FTarrowheadout\Q@s@tarrowheadout{#2}\else
    \def\n@mref{r}\ifx\l@debut\n@mref\update@ttr\D@FTarrowheadratio\Q@s@tarrowheadratio{#2}\else
    \W@rnmesAttr{figset arrowhead}{#1}\fi\fi\fi\fi\fi}
\ctr@ln@m\@rrowheadangle
\ctr@ln@m\C@AHANG \ctr@ln@m\S@AHANG \ctr@ln@m\UNSS@N
\ctr@ld@f\def\Q@s@tarrowheadangle#1{\edef\@rrowheadangle{#1}{\c@ssin{\C@}{\S@}{#1}%
    \xdef\C@AHANG{\C@}\xdef\S@AHANG{\S@}\v@lmax=\S@ pt%
    \invers@{\v@leur}{\v@lmax}\maxim@m{\v@leur}{\v@leur}{-\v@leur}%
    \xdef\UNSS@N{\the\v@leur}}}
\ctr@ld@f\def\Q@s@tarrowheadfill#1{\expandafter\set@rrowhfill#1:}
\ctr@ld@f\def\set@rrowhfill#1#2:{\if#1n\@rrowhfillfalse\else\@rrowhfilltrue\fi}
\ctr@ld@f\def\Q@s@tarrowheadout#1{\expandafter\set@rrowhout#1:}
\ctr@ld@f\def\set@rrowhout#1#2:{\if#1n\@rrowhoutfalse\else\@rrowhouttrue\fi}
\ctr@ln@m\@rrowheadlength
\ctr@ld@f\def\Q@s@tarrowheadlength#1{\edef\@rrowheadlength{#1}\@rrowratiofalse}
\ctr@ln@m\@rrowheadratio
\ctr@ld@f\def\Q@s@tarrowheadratio#1{\edef\@rrowheadratio{#1}\@rrowratiotrue}
\ctr@ln@m\D@FTarrowheadlength
\ctr@ld@f\def\figresetarrowhead{%
    \Q@s@tarrowheadangle{\D@FTarrowheadangle}%
    \Q@s@tarrowheadfill{\D@FTarrowheadfill}%
    \Q@s@tarrowheadout{\D@FTarrowheadout}%
    \Q@s@tarrowheadratio{\D@FTarrowheadratio}%
    \d@fm@cdim\D@FTarrowheadlength{\D@FTh@rdahlength}
    \Q@s@tarrowheadlength{\D@FTarrowheadlength}}
\ctr@ld@f\def\D@FTarrowheadratio{0.1}
\ctr@ld@f\def\D@FTarrowheadangle{20}
\ctr@ld@f\def\D@FTarrowheadfill{no}
\ctr@ld@f\def\D@FTarrowheadout{no}
\ctr@ld@f\def\D@FTh@rdahlength{8pt}
\ctr@ln@m\figdrawarrow
\ctr@ld@f\def\Q@arrowDD[#1,#2]{{\ifCUR@PS\ifGR@cri\s@uvc@ntr@l\et@tpsarrow%
    \PSc@mment{arrowDD [Pt1,Pt2]=[#1,#2]}\Q@s@tfillmode{no}%
    \Q@arrowheadDD[#1,#2]\setc@ntr@l{2}\figdrawline[#1,-3]%
    \PSc@mment{End arrowDD}\resetc@ntr@l\et@tpsarrow\fi\fi}}
\ctr@ld@f\def\Q@arrowTD[#1,#2]{{\ifCUR@PS\ifGR@cri\s@uvc@ntr@l\et@tpsarrowTD%
    \PSc@mment{arrowTD [Pt1,Pt2]=[#1,#2]}\resetc@ntr@l{2}%
    \Figptpr@j-5:/#1/\Figptpr@j-6:/#2/\let\c@lprojSP=\relax\Q@arrowDD[-5,-6]%
    \PSc@mment{End arrowTD}\resetc@ntr@l\et@tpsarrowTD\fi\fi}}
\ctr@ln@m\figdrawarrowhead
\ctr@ld@f\def\Q@arrowheadDD[#1,#2]{{\ifCUR@PS\ifGR@cri\s@uvc@ntr@l\et@tpsarrowheadDD%
    \if@rrowhfill\def\@hangle{-\@rrowheadangle}\else\def\@hangle{\@rrowheadangle}\fi%
    \if@rrowratio%
    \if@rrowhout\def\@hratio{-\@rrowheadratio}\else\def\@hratio{\@rrowheadratio}\fi%
    \PSc@mment{arrowheadDD Ratio=\@hratio, Angle=\@hangle, [Pt1,Pt2]=[#1,#2]}%
    \Ps@rrowhead\@hratio,\@hangle[#1,#2]%
    \else%
    \if@rrowhout\def\@hlength{-\@rrowheadlength}\else\def\@hlength{\@rrowheadlength}\fi%
    \PSc@mment{arrowheadDD Length=\@hlength, Angle=\@hangle, [Pt1,Pt2]=[#1,#2]}%
    \Ps@rrowheadfd\@hlength,\@hangle[#1,#2]%
    \fi%
    \PSc@mment{End arrowheadDD}\resetc@ntr@l\et@tpsarrowheadDD\fi\fi}}
\ctr@ld@f\def\Q@arrowheadTD[#1,#2]{{\ifCUR@PS\ifGR@cri\s@uvc@ntr@l\et@tpsarrowheadTD%
    \PSc@mment{arrowheadTD [Pt1,Pt2]=[#1,#2]}\resetc@ntr@l{2}%
    \Figptpr@j-5:/#1/\Figptpr@j-6:/#2/\let\c@lprojSP=\relax\Q@arrowheadDD[-5,-6]%
    \PSc@mment{End arrowheadTD}\resetc@ntr@l\et@tpsarrowheadTD\fi\fi}}
\ctr@ld@f\def\Ps@rrowhead#1,#2[#3,#4]{\v@leur=#1\p@\maxim@m{\v@leur}{\v@leur}{-\v@leur}%
    \ifdim\v@leur>\Cepsil@n{
    \PSc@mment{@rrowhead Ratio=#1, Angle=#2, [Pt1,Pt2]=[#3,#4]}\v@leur=\UNSS@N%
    \v@leur=\CUR@width\v@leur\v@leur=\ptpsT@pt\v@leur\delt@=.5\v@leur
    \setc@ntr@l{2}\figvectPDD-3[#4,#3]%
    \Figg@tXY{-3}\v@lX=#1\v@lX\v@lY=#1\v@lY\Figv@ctCreg-3(\v@lX,\v@lY)%
    \vecunit@{-4}{-3}\mili@u=\result@t%
    \ifdim#2pt>\z@\v@lXa=-\C@AHANG\delt@%
     \edef\c@ef{\repdecn@mb{\v@lXa}}\figpttraDD-3:=-3/\c@ef,-4/\fi%
    \edef\c@ef{\repdecn@mb{\delt@}}%
    \v@lXa=\mili@u\v@lXa=\C@AHANG\v@lXa%
    \v@lYa=\ptpsT@pt\p@\v@lYa=\CUR@width\v@lYa\v@lYa=\sDcc@ngle\v@lYa%
    \advance\v@lXa-\v@lYa\gdef\sDcc@ngle{0}%
    \ifdim\v@lXa>\v@leur\edef\c@efendpt{\repdecn@mb{\v@leur}}%
    \else\edef\c@efendpt{\repdecn@mb{\v@lXa}}\fi%
    \Figg@tXY{-3}\v@lmin=\v@lX\v@lmax=\v@lY%
    \v@lXa=\C@AHANG\v@lmin\v@lYa=\S@AHANG\v@lmax\advance\v@lXa\v@lYa%
    \v@lYa=-\S@AHANG\v@lmin\v@lX=\C@AHANG\v@lmax\advance\v@lYa\v@lX%
    \setc@ntr@l{1}\Figg@tXY{#4}\advance\v@lX\v@lXa\advance\v@lY\v@lYa%
    \setc@ntr@l{2}\Figp@intregDD-2:(\v@lX,\v@lY)%
    \v@lXa=\C@AHANG\v@lmin\v@lYa=-\S@AHANG\v@lmax\advance\v@lXa\v@lYa%
    \v@lYa=\S@AHANG\v@lmin\v@lX=\C@AHANG\v@lmax\advance\v@lYa\v@lX%
    \setc@ntr@l{1}\Figg@tXY{#4}\advance\v@lX\v@lXa\advance\v@lY\v@lYa%
    \setc@ntr@l{2}\Figp@intregDD-1:(\v@lX,\v@lY)%
    \ifdim#2pt<\z@\fillm@detrue\figdrawline[-2,#4,-1]
    \else\figptstraDD-3=#4,-2,-1/\c@ef,-4/\s@uvdash{\typ@dash}\Q@s@tdash{\D@FTdash}%
    \figdrawline[-2,-3,-1]\Q@s@tdash{\typ@dash}\fi
    \ifdim#1pt>\z@\figpttraDD-3:=#4/\c@efendpt,-4/\else\figptcopyDD-3:/#4/\fi%
    \PSc@mment{End @rrowhead}}\fi}
\ctr@ld@f\def\sDcc@ngle{0}
\ctr@ld@f\def\Ps@rrowheadfd#1,#2[#3,#4]{{%
    \PSc@mment{@rrowheadfd Length=#1, Angle=#2, [Pt1,Pt2]=[#3,#4]}%
    \setc@ntr@l{2}\figvectPDD-1[#3,#4]\n@rmeucDD{\v@leur}{-1}\v@leur=\ptT@unit@\v@leur%
    \invers@{\v@leur}{\v@leur}\v@leur=#1\v@leur\edef\R@tio{\repdecn@mb{\v@leur}}%
    \Ps@rrowhead\R@tio,#2[#3,#4]\PSc@mment{End @rrowheadfd}}}
\ctr@ln@m\figdrawarrowBezier
\ctr@ld@f\def\Q@arrowBezierDD[#1,#2,#3,#4]{{\ifCUR@PS\ifGR@cri\s@uvc@ntr@l\et@tpsarrowBezierDD%
    \PSc@mment{arrowBezierDD Control points=#1,#2,#3,#4}\setc@ntr@l{2}%
    \if@rrowratio\c@larclengthDD\v@leur,10[#1,#2,#3,#4]\else\v@leur=\z@\fi%
    \Ps@rrowB@zDD\v@leur[#1,#2,#3,#4]%
    \PSc@mment{End arrowBezierDD}\resetc@ntr@l\et@tpsarrowBezierDD\fi\fi}}
\ctr@ld@f\def\Q@arrowBezierTD[#1,#2,#3,#4]{{\ifCUR@PS\ifGR@cri\s@uvc@ntr@l\et@tpsarrowBezierTD%
    \PSc@mment{arrowBezierTD Control points=#1,#2,#3,#4}\resetc@ntr@l{2}%
    \Figptpr@j-7:/#1/\Figptpr@j-8:/#2/\Figptpr@j-9:/#3/\Figptpr@j-10:/#4/%
    \let\c@lprojSP=\relax\ifnum\CUR@proj<\tw@\Q@arrowBezierDD[-7,-8,-9,-10]%
    \else\f@gnewpath\PSwrit@cmd{-7}{\c@mmoveto}{\fwf@g}%
    \if@rrowratio\c@larclengthDD\mili@u,10[-7,-8,-9,-10]\else\mili@u=\z@\fi%
    \p@rtent=\NBz@rcs\advance\p@rtent\m@ne\subB@zierTD\p@rtent[#1,#2,#3,#4]%
    \f@gstroke%
    \advance\v@lmin\p@rtent\delt@
    \v@leur=\v@lmin\advance\v@leur0.33333 \delt@\edef\unti@rs{\repdecn@mb{\v@leur}}%
    \v@leur=\v@lmin\advance\v@leur0.66666 \delt@\edef\deti@rs{\repdecn@mb{\v@leur}}%
    \figptcopyDD-8:/-10/\c@lsubBzarc\unti@rs,\deti@rs[#1,#2,#3,#4]%
    \figptcopyDD-8:/-4/\figptcopyDD-9:/-3/\Ps@rrowB@zDD\mili@u[-7,-8,-9,-10]\fi%
    \PSc@mment{End arrowBezierTD}\resetc@ntr@l\et@tpsarrowBezierTD\fi\fi}}
\ctr@ld@f\def\c@larclengthDD#1,#2[#3,#4,#5,#6]{{\p@rtent=#2\figptcopyDD-5:/#3/%
    \delt@=\p@\divide\delt@\p@rtent\c@rre=\z@\v@leur=\z@\s@mme=\z@%
    \loop\ifnum\s@mme<\p@rtent\advance\s@mme\@ne\advance\v@leur\delt@%
    \edef\T@{\repdecn@mb{\v@leur}}\figptBezierDD-6::\T@[#3,#4,#5,#6]%
    \figvectPDD-1[-5,-6]\n@rmeucDD{\mili@u}{-1}\advance\c@rre\mili@u%
    \figptcopyDD-5:/-6/\repeat\global\result@t=\ptT@unit@\c@rre}#1=\result@t}
\ctr@ld@f\def\Ps@rrowB@zDD#1[#2,#3,#4,#5]{{\Q@s@tfillmode{no}%
    \if@rrowratio\delt@=\@rrowheadratio#1\else\delt@=\@rrowheadlength pt\fi%
    \v@leur=\C@AHANG\delt@\edef\R@dius{\repdecn@mb{\v@leur}}%
    \FigptintercircB@zDD-5::0,\R@dius[#5,#4,#3,#2]%
    \Q@s@tarrowheadlength{\repdecn@mb{\delt@}}\Q@arrowheadDD[-5,#5]%
    \let\n@rmeuc=\n@rmeucDD\figgetdist\R@dius[#5,-3]%
    \FigptintercircB@zDD-6::0,\R@dius[#5,#4,#3,#2]%
    \figptBezierDD-5::0.33333[#5,#4,#3,#2]\figptBezierDD-3::0.66666[#5,#4,#3,#2]%
    \figptscontrolDD-5[-6,-5,-3,#2]\Q@BezierDD1[-6,-5,-4,#2]}}
\ctr@ln@m\figdrawarrowcirc
\ctr@ld@f\def\Q@arrowcircDD#1;#2(#3,#4){{\ifCUR@PS\ifGR@cri\s@uvc@ntr@l\et@tpsarrowcircDD%
    \PSc@mment{arrowcircDD Center=#1 ; Radius=#2 (Ang1=#3,Ang2=#4)}%
    \Q@s@tfillmode{no}\Pscirc@rrowhead#1;#2(#3,#4)%
    \setc@ntr@l{2}\figvectPDD -4[#1,-3]\vecunit@{-4}{-4}%
    \Figg@tXY{-4}\arct@n\v@lmin(\v@lX,\v@lY)%
    \v@lmin=\rdT@deg\v@lmin\v@leur=#4pt\advance\v@leur-\v@lmin%
    \maxim@m{\v@leur}{\v@leur}{-\v@leur}%
    \ifdim\v@leur>\DemiPI@deg\relax\ifdim\v@lmin<#4pt\advance\v@lmin\DePI@deg%
    \else\advance\v@lmin-\DePI@deg\fi\fi\edef\ar@ngle{\repdecn@mb{\v@lmin}}%
    \ifdim#3pt<#4pt\figdrawarccirc#1;#2(#3,\ar@ngle)\else\figdrawarccirc#1;#2(\ar@ngle,#3)\fi%
    \PSc@mment{End arrowcircDD}\resetc@ntr@l\et@tpsarrowcircDD\fi\fi}}
\ctr@ld@f\def\Q@arrowcircTD#1,#2,#3;#4(#5,#6){{\ifCUR@PS\ifGR@cri\s@uvc@ntr@l\et@tpsarrowcircTD%
    \PSc@mment{arrowcircTD Center=#1,P1=#2,P2=#3 ; Radius=#4 (Ang1=#5, Ang2=#6)}%
    \resetc@ntr@l{2}\c@lExtAxes#1,#2,#3(#4)\let\c@lprojSP=\relax%
    \figvectPTD-11[#1,-4]\figvectPTD-12[#1,-5]\c@lNbarcs{#5}{#6}%
    \if@rrowratio\v@lmax=\degT@rd\v@lmax\edef\D@lpha{\repdecn@mb{\v@lmax}}\fi%
    \advance\p@rtent\m@ne\mili@u=\z@%
    \v@leur=#5pt\c@lptellP{#1}{-11}{-12}\Figptpr@j-9:/-3/%
    \f@gnewpath\PSwrit@cmdS{-9}{\c@mmoveto}{\fwf@g}{\X@un}{\Y@un}%
    \edef\C@nt@r{#1}\s@mme=\z@\bcl@rcircTD\f@gstroke%
    \advance\v@leur\delt@\c@lptellP{#1}{-11}{-12}\Figptpr@j-5:/-3/%
    \advance\v@leur\delt@\c@lptellP{#1}{-11}{-12}\Figptpr@j-6:/-3/%
    \advance\v@leur\delt@\c@lptellP{#1}{-11}{-12}\Figptpr@j-10:/-3/%
    \figptscontrolDD-8[-9,-5,-6,-10]%
    \if@rrowratio\c@lcurvradDD0.5[-9,-8,-7,-10]\advance\mili@u\result@t%
    \maxim@m{\mili@u}{\mili@u}{-\mili@u}\mili@u=\ptT@unit@\mili@u%
    \mili@u=\D@lpha\mili@u\advance\p@rtent\@ne\divide\mili@u\p@rtent\fi%
    \Ps@rrowB@zDD\mili@u[-9,-8,-7,-10]%
    \PSc@mment{End arrowcircTD}\resetc@ntr@l\et@tpsarrowcircTD\fi\fi}}
\ctr@ld@f\def\bcl@rcircTD{\relax%
    \ifnum\s@mme<\p@rtent\advance\s@mme\@ne%
    \advance\v@leur\delt@\c@lptellP{\C@nt@r}{-11}{-12}\Figptpr@j-5:/-3/%
    \advance\v@leur\delt@\c@lptellP{\C@nt@r}{-11}{-12}\Figptpr@j-6:/-3/%
    \advance\v@leur\delt@\c@lptellP{\C@nt@r}{-11}{-12}\Figptpr@j-10:/-3/%
    \figptscontrolDD-8[-9,-5,-6,-10]\BdingB@xfalse%
    \PSwrit@cmdS{-8}{}{\fwf@g}{\X@de}{\Y@de}\PSwrit@cmdS{-7}{}{\fwf@g}{\X@tr}{\Y@tr}%
    \BdingB@xtrue\PSwrit@cmdS{-10}{\c@mcurveto}{\fwf@g}{\X@qu}{\Y@qu}%
    \if@rrowratio\c@lcurvradDD0.5[-9,-8,-7,-10]\advance\mili@u\result@t\fi%
    \B@zierBB@x{1}{\Y@un}(\X@un,\X@de,\X@tr,\X@qu)%
    \B@zierBB@x{2}{\X@un}(\Y@un,\Y@de,\Y@tr,\Y@qu)%
    \edef\X@un{\X@qu}\edef\Y@un{\Y@qu}\figptcopyDD-9:/-10/\bcl@rcircTD\fi}
\ctr@ld@f\def\Pscirc@rrowhead#1;#2(#3,#4){{%
    \PSc@mment{circ@rrowhead Center=#1 ; Radius=#2 (Ang1=#3,Ang2=#4)}%
    \v@leur=#2\unit@\edef\s@glen{\repdecn@mb{\v@leur}}\v@lY=\z@\v@lX=\v@leur%
    \resetc@ntr@l{2}\Figv@ctCreg-3(\v@lX,\v@lY)\figpttraDD-5:=#1/1,-3/%
    \figptrotDD-5:=-5/#1,#4/%
    \figvectPDD-3[#1,-5]\Figg@tXY{-3}\v@leur=\v@lX%
    \ifdim#3pt<#4pt\v@lX=\v@lY\v@lY=-\v@leur\else\v@lX=-\v@lY\v@lY=\v@leur\fi%
    \Figv@ctCreg-3(\v@lX,\v@lY)\vecunit@{-3}{-3}%
    \if@rrowratio\v@leur=#4pt\advance\v@leur-#3pt\maxim@m{\mili@u}{-\v@leur}{\v@leur}%
    \mili@u=\degT@rd\mili@u\v@leur=\s@glen\mili@u\edef\s@glen{\repdecn@mb{\v@leur}}%
    \mili@u=#2\mili@u\mili@u=\@rrowheadratio\mili@u\else\mili@u=\@rrowheadlength pt\fi%
    \figpttraDD-6:=-5/\s@glen,-3/\v@leur=#2pt\v@leur=2\v@leur%
    \invers@{\v@leur}{\v@leur}\c@rre=\repdecn@mb{\v@leur}\mili@u
    \mili@u=\c@rre\mili@u=\repdecn@mb{\c@rre}\mili@u%
    \v@leur=\p@\advance\v@leur-\mili@u
    \invers@{\mili@u}{2\v@leur}\delt@=\c@rre\delt@=\repdecn@mb{\mili@u}\delt@%
    \xdef\sDcc@ngle{\repdecn@mb{\delt@}}
    \sqrt@{\mili@u}{\v@leur}\arct@n\v@leur(\mili@u,\c@rre)%
    \v@leur=\rdT@deg\v@leur
    \ifdim#3pt<#4pt\v@leur=-\v@leur\fi%
    \if@rrowhout\v@leur=-\v@leur\fi\edef\cor@ngle{\repdecn@mb{\v@leur}}%
    \figptrotDD-6:=-6/-5,\cor@ngle/\Q@arrowheadDD[-6,-5]%
    \PSc@mment{End circ@rrowhead}}}
\ctr@ln@m\figdrawarrowcircP
\ctr@ld@f\def\Q@arrowcircPDD#1;#2[#3,#4]{{\ifCUR@PS\ifGR@cri%
    \PSc@mment{arrowcircPDD Center=#1; Radius=#2, [P1=#3,P2=#4]}%
    \s@uvc@ntr@l\et@tpsarrowcircPDD\Ps@ngleparam#1;#2[#3,#4]%
    \ifdim\v@leur>\z@\ifdim\v@lmin>\v@lmax\advance\v@lmax\DePI@deg\fi%
    \else\ifdim\v@lmin<\v@lmax\advance\v@lmin\DePI@deg\fi\fi%
    \edef\@ngdeb{\repdecn@mb{\v@lmin}}\edef\@ngfin{\repdecn@mb{\v@lmax}}%
    \figdrawarrowcirc#1;\r@dius(\@ngdeb,\@ngfin)%
    \PSc@mment{End arrowcircPDD}\resetc@ntr@l\et@tpsarrowcircPDD\fi\fi}}
\ctr@ld@f\def\Q@arrowcircPTD#1;#2[#3,#4,#5]{{\ifCUR@PS\ifGR@cri\s@uvc@ntr@l\et@tpsarrowcircPTD%
    \PSc@mment{arrowcircPTD Center=#1; Radius=#2, [P1=#3,P2=#4,P3=#5]}%
    \figgetangleTD\@ngfin[#1,#3,#4,#5]\v@leur=#2pt%
    \maxim@m{\mili@u}{-\v@leur}{\v@leur}\edef\r@dius{\repdecn@mb{\mili@u}}%
    \ifdim\v@leur<\z@\v@lmax=\@ngfin pt\advance\v@lmax-\DePI@deg%
    \edef\@ngfin{\repdecn@mb{\v@lmax}}\fi\Q@arrowcircTD#1,#3,#5;\r@dius(0,\@ngfin)%
    \PSc@mment{End arrowcircPTD}\resetc@ntr@l\et@tpsarrowcircPTD\fi\fi}}
\ctr@ld@f\def\figdrawaxes#1(#2){{\ifCUR@PS\ifGR@cri\s@uvc@ntr@l\et@tpsaxes%
    \PSc@mment{axes Origin=#1 Range=(#2)}\an@lys@xes#2,:\resetc@ntr@l{2}%
    \ifx\t@xt@\empty\ifTr@isDim\Q@@xes#1(0,#2,0,#2,0,#2)\else\Q@@xes#1(0,#2,0,#2)\fi%
    \else\Q@@xes#1(#2)\fi\PSc@mment{End axes}\resetc@ntr@l\et@tpsaxes\fi\fi}}
\ctr@ld@f\def\an@lys@xes#1,#2:{\def\t@xt@{#2}}
\ctr@ln@m\Q@@xes
\ctr@ld@f\def\Q@@xesDD#1(#2,#3,#4,#5){%
    \figpttraC-5:=#1/#2,0/\figpttraC-6:=#1/#3,0/\Q@arrowDD[-5,-6]%
    \figpttraC-5:=#1/0,#4/\figpttraC-6:=#1/0,#5/\Q@arrowDD[-5,-6]}
\ctr@ld@f\def\Q@@xesTD#1(#2,#3,#4,#5,#6,#7){%
    \figpttraC-7:=#1/#2,0,0/\figpttraC-8:=#1/#3,0,0/\Q@arrowTD[-7,-8]%
    \figpttraC-7:=#1/0,#4,0/\figpttraC-8:=#1/0,#5,0/\Q@arrowTD[-7,-8]%
    \figpttraC-7:=#1/0,0,#6/\figpttraC-8:=#1/0,0,#7/\Q@arrowTD[-7,-8]}
\ctr@ln@m\newGr@FN
\ctr@ld@f\def\newGr@FNPDF#1{\s@mme=\Gr@FNb\advance\s@mme\@ne\xdef\Gr@FNb{\number\s@mme}}
\ctr@ld@f\def\newGr@FNDVI#1{\newGr@FNPDF{}\xdef#1{\jobname GI\Gr@FNb.anx}}
\ctr@ld@f\def\figdrawbegin#1{\newGr@FN\DefGIfilen@me\gdef\@utoFN{0}%
    \def\t@xt@{#1}\relax\ifx\t@xt@\empty\GRupdatem@detrue%
    \gdef\@utoFN{1}\Psb@ginfig\DefGIfilen@me\else\expandafter\Psb@ginfigNu@#1 :\fi}
\ctr@ld@f\def\Psb@ginfigNu@#1 #2:{\def\t@xt@{#1}\relax\ifx\t@xt@\empty\def\t@xt@{#2}%
    \ifx\t@xt@\empty\GRupdatem@detrue\gdef\@utoFN{1}\Psb@ginfig\DefGIfilen@me%
    \else\Psb@ginfigNu@#2:\fi\else\Psb@ginfig{#1}\fi}
\ctr@ln@m\PSfilen@me \ctr@ln@m\auxfilen@me
\ctr@ld@f\def\Psb@ginfig#1{\ifCUR@PS\else%
    \edef\PSfilen@me{#1}\edef\auxfilen@me{\jobname.anx}%
    \ifGRupdatem@de\GR@critrue\else\openin\frf@g=\PSfilen@me\relax%
    \ifeof\frf@g\GR@critrue\else\GR@crifalse\fi\closein\frf@g\fi%
    \CUR@PStrue\c@ldefproj\expandafter\setupd@te\D@FTupdate:%
    \ifGR@cri\initb@undb@x%
    \immediate\openout\fwf@g=\auxfilen@me\initpss@ttings\fi%
    \fi}
\ctr@ld@f\def\Gr@FNb{0}
\ctr@ld@f\def\figforTeXFileno{\Gr@FNb}
\ctr@ld@f\def\figforTeXFigno{0 }
\ctr@ld@f\def\figforTeXnextFigno{1 }
\ctr@ld@f\edef\DefGIfilen@me{\jobname GI.anx}
\ctr@ld@f\def\initpss@ttings{\figreset{altitude,arrowhead,curve,general,flowchart,mesh,trimesh}%
    \Use@llipsefalse}
\ctr@ld@f\def\B@zierBB@x#1#2(#3,#4,#5,#6){{\c@rre=\t@n\epsil@n
    \v@lmax=#4\advance\v@lmax-#5\v@lmax=\thr@@\v@lmax\advance\v@lmax#6\advance\v@lmax-#3%
    \mili@u=#4\mili@u=-\tw@\mili@u\advance\mili@u#3\advance\mili@u#5%
    \v@lmin=#4\advance\v@lmin-#3\maxim@m{\v@leur}{-\v@lmax}{\v@lmax}%
    \maxim@m{\delt@}{-\mili@u}{\mili@u}\maxim@m{\v@leur}{\v@leur}{\delt@}%
    \maxim@m{\delt@}{-\v@lmin}{\v@lmin}\maxim@m{\v@leur}{\v@leur}{\delt@}%
    \ifdim\v@leur>\c@rre\invers@{\v@leur}{\v@leur}\edef\Uns@rM@x{\repdecn@mb{\v@leur}}%
    \v@lmax=\Uns@rM@x\v@lmax\mili@u=\Uns@rM@x\mili@u\v@lmin=\Uns@rM@x\v@lmin%
    \maxim@m{\v@leur}{-\v@lmax}{\v@lmax}\ifdim\v@leur<\c@rre%
    \maxim@m{\v@leur}{-\mili@u}{\mili@u}\ifdim\v@leur<\c@rre\else%
    \invers@{\mili@u}{\mili@u}\v@leur=-0.5\v@lmin%
    \v@leur=\repdecn@mb{\mili@u}\v@leur\m@jBBB@x{\v@leur}{#1}{#2}(#3,#4,#5,#6)\fi%
    \else\delt@=\repdecn@mb{\mili@u}\mili@u\v@leur=\repdecn@mb{\v@lmax}\v@lmin%
    \advance\delt@-\v@leur\ifdim\delt@<\z@\else\invers@{\v@lmax}{\v@lmax}%
    \edef\Uns@rAp{\repdecn@mb{\v@lmax}}\sqrt@{\delt@}{\delt@}%
    \v@leur=-\mili@u\advance\v@leur\delt@\v@leur=\Uns@rAp\v@leur%
    \m@jBBB@x{\v@leur}{#1}{#2}(#3,#4,#5,#6)%
    \v@leur=-\mili@u\advance\v@leur-\delt@\v@leur=\Uns@rAp\v@leur%
    \m@jBBB@x{\v@leur}{#1}{#2}(#3,#4,#5,#6)\fi\fi\fi}}
\ctr@ld@f\def\m@jBBB@x#1#2#3(#4,#5,#6,#7){{\relax\ifdim#1>\z@\ifdim#1<\p@%
    \edef\T@{\repdecn@mb{#1}}\v@lX=\p@\advance\v@lX-#1\edef\UNmT@{\repdecn@mb{\v@lX}}%
    \v@lX=#4\v@lY=#5\v@lZ=#6\v@lXa=#7\v@lX=\UNmT@\v@lX\advance\v@lX\T@\v@lY%
    \v@lY=\UNmT@\v@lY\advance\v@lY\T@\v@lZ\v@lZ=\UNmT@\v@lZ\advance\v@lZ\T@\v@lXa%
    \v@lX=\UNmT@\v@lX\advance\v@lX\T@\v@lY\v@lY=\UNmT@\v@lY\advance\v@lY\T@\v@lZ%
    \v@lX=\UNmT@\v@lX\advance\v@lX\T@\v@lY%
    \ifcase#2\or\v@lY=#3\or\v@lY=\v@lX\v@lX=#3\fi\b@undb@x{\v@lX}{\v@lY}\fi\fi}}
\ctr@ld@f\def\PsB@zier#1[#2]{{\f@gnewpath%
    \s@mme=\z@\def\list@num{#2,0}\extrairelepremi@r\p@int\de\list@num%
    \PSwrit@cmdS{\p@int}{\c@mmoveto}{\fwf@g}{\X@un}{\Y@un}\p@rtent=#1\bclB@zier}}
\ctr@ld@f\def\bclB@zier{\relax%
    \ifnum\s@mme<\p@rtent\advance\s@mme\@ne\BdingB@xfalse%
    \extrairelepremi@r\p@int\de\list@num\PSwrit@cmdS{\p@int}{}{\fwf@g}{\X@de}{\Y@de}%
    \extrairelepremi@r\p@int\de\list@num\PSwrit@cmdS{\p@int}{}{\fwf@g}{\X@tr}{\Y@tr}%
    \BdingB@xtrue%
    \extrairelepremi@r\p@int\de\list@num\PSwrit@cmdS{\p@int}{\c@mcurveto}{\fwf@g}{\X@qu}{\Y@qu}%
    \B@zierBB@x{1}{\Y@un}(\X@un,\X@de,\X@tr,\X@qu)%
    \B@zierBB@x{2}{\X@un}(\Y@un,\Y@de,\Y@tr,\Y@qu)%
    \edef\X@un{\X@qu}\edef\Y@un{\Y@qu}\bclB@zier\fi}
\ctr@ln@m\figdrawBezier
\ctr@ld@f\def\Q@BezierDD#1[#2]{\ifCUR@PS\ifGR@cri%
    \PSc@mment{BezierDD N arcs=#1, Control points=#2}%
    \iffillm@de\PsB@zier#1[#2]%
    \f@gfill%
    \else\PsB@zier#1[#2]\f@gstroke\fi%
    \PSc@mment{End BezierDD}\fi\fi}
\ctr@ln@m\et@tpsBezierTD
\ctr@ld@f\def\Q@BezierTD#1[#2]{\ifCUR@PS\ifGR@cri\s@uvc@ntr@l\et@tpsBezierTD%
    \PSc@mment{BezierTD N arcs=#1, Control points=#2}%
    \iffillm@de\PsB@zierTD#1[#2]%
    \f@gfill%
    \else\PsB@zierTD#1[#2]\f@gstroke\fi%
    \PSc@mment{End BezierTD}\resetc@ntr@l\et@tpsBezierTD\fi\fi}
\ctr@ld@f\def\PsB@zierTD#1[#2]{\ifnum\CUR@proj<\tw@\PsB@zier#1[#2]\else\PsB@zier@TD#1[#2]\fi}
\ctr@ld@f\def\PsB@zier@TD#1[#2]{{\f@gnewpath%
    \s@mme=\z@\def\list@num{#2,0}\extrairelepremi@r\p@int\de\list@num%
    \let\c@lprojSP=\relax\setc@ntr@l{2}\Figptpr@j-7:/\p@int/%
    \PSwrit@cmd{-7}{\c@mmoveto}{\fwf@g}%
    \loop\ifnum\s@mme<#1\advance\s@mme\@ne\extrairelepremi@r\p@intun\de\list@num%
    \extrairelepremi@r\p@intde\de\list@num\extrairelepremi@r\p@inttr\de\list@num%
    \subB@zierTD\NBz@rcs[\p@int,\p@intun,\p@intde,\p@inttr]\edef\p@int{\p@inttr}\repeat}}
\ctr@ld@f\def\subB@zierTD#1[#2,#3,#4,#5]{\delt@=\p@\divide\delt@\NBz@rcs\v@lmin=\z@%
    {\Figg@tXY{-7}\edef\X@un{\the\v@lX}\edef\Y@un{\the\v@lY}%
    \s@mme=\z@\loop\ifnum\s@mme<#1\advance\s@mme\@ne%
    \v@leur=\v@lmin\advance\v@leur0.33333 \delt@\edef\unti@rs{\repdecn@mb{\v@leur}}%
    \v@leur=\v@lmin\advance\v@leur0.66666 \delt@\edef\deti@rs{\repdecn@mb{\v@leur}}%
    \advance\v@lmin\delt@\edef\trti@rs{\repdecn@mb{\v@lmin}}%
    \figptBezierTD-8::\trti@rs[#2,#3,#4,#5]\Figptpr@j-8:/-8/%
    \c@lsubBzarc\unti@rs,\deti@rs[#2,#3,#4,#5]\BdingB@xfalse%
    \PSwrit@cmdS{-4}{}{\fwf@g}{\X@de}{\Y@de}\PSwrit@cmdS{-3}{}{\fwf@g}{\X@tr}{\Y@tr}%
    \BdingB@xtrue\PSwrit@cmdS{-8}{\c@mcurveto}{\fwf@g}{\X@qu}{\Y@qu}%
    \B@zierBB@x{1}{\Y@un}(\X@un,\X@de,\X@tr,\X@qu)%
    \B@zierBB@x{2}{\X@un}(\Y@un,\Y@de,\Y@tr,\Y@qu)%
    \edef\X@un{\X@qu}\edef\Y@un{\Y@qu}\figptcopyDD-7:/-8/\repeat}}
\ctr@ld@f\def\NBz@rcs{2}
\ctr@ld@f\def\c@lsubBzarc#1,#2[#3,#4,#5,#6]{\figptBezierTD-5::#1[#3,#4,#5,#6]%
    \figptBezierTD-6::#2[#3,#4,#5,#6]\Figptpr@j-4:/-5/\Figptpr@j-5:/-6/%
    \figptscontrolDD-4[-7,-4,-5,-8]}
\ctr@ln@m\figdrawcirc
\ctr@ld@f\def\Q@circDD#1(#2){\ifCUR@PS\ifGR@cri\PSc@mment{circDD Center=#1 (Radius=#2)}%
    \Q@arccircDD#1;#2(0,360)\PSc@mment{End circDD}\fi\fi}
\ctr@ld@f\def\Q@circTD#1,#2,#3(#4){\ifCUR@PS\ifGR@cri%
    \PSc@mment{circTD Center=#1,P1=#2,P2=#3 (Radius=#4)}%
    \Q@arccircTD#1,#2,#3;#4(0,360)\PSc@mment{End circTD}\fi\fi}
\ctr@ln@m\p@urcent
{\catcode`\%=12\gdef\p@urcent{
\ctr@ld@f\def\PSc@mment#1{\ifGRdebugm@de\immediate\write\fwf@g{\p@urcent\space#1}\fi}
\ctr@ln@m\acc@louv \ctr@ln@m\acc@lfer
{\catcode`\[=1\catcode`\{=12\gdef\acc@louv[{}}
{\catcode`\]=2\catcode`\}=12\gdef\acc@lfer{}]]
\ctr@ld@f\def\PSdict@{\ifUse@llipse%
    \immediate\write\fwf@g{/ellipsedict 9 dict def ellipsedict /mtrx matrix put}%
    \immediate\write\fwf@g{/ellipse \acc@louv ellipsedict begin}%
    \immediate\write\fwf@g{ /endangle exch def /startangle exch def}%
    \immediate\write\fwf@g{ /yrad exch def /xrad exch def}%
    \immediate\write\fwf@g{ /rotangle exch def /y exch def /x exch def}%
    \immediate\write\fwf@g{ /savematrix mtrx currentmatrix def}%
    \immediate\write\fwf@g{ x y translate rotangle rotate xrad yrad scale}%
    \immediate\write\fwf@g{ 0 0 1 startangle endangle arc}%
    \immediate\write\fwf@g{ savematrix setmatrix end\acc@lfer def}%
    \fi\PShe@der{EndProlog}}
\ctr@ld@f\def\Pssetc@rve#1=#2|{\keln@mun#1|%
    \def\n@mref{r}\ifx\l@debut\n@mref\update@ttr\D@FTroundness\Q@s@troundness{#2}\else
    \W@rnmesAttr{figset curve}{#1}\fi}
\ctr@ln@m\curv@roundness
\ctr@ld@f\def\Q@s@troundness#1{\edef\curv@roundness{#1}}
\ctr@ld@f\def\D@FTroundness{0.2} 
\ctr@ln@m\figdrawcurve
\ctr@ld@f\def\Q@curveDD[#1]{{\ifCUR@PS\ifGR@cri\PSc@mment{curveDD Points=#1}%
    \s@uvc@ntr@l\et@tpscurveDD%
    \iffillm@de\Psc@rveDD\curv@roundness[#1]%
    \f@gfill%
    \else\Psc@rveDD\curv@roundness[#1]\f@gstroke\fi%
    \PSc@mment{End curveDD}\resetc@ntr@l\et@tpscurveDD\fi\fi}}
\ctr@ld@f\def\Q@curveTD[#1]{{\ifCUR@PS\ifGR@cri%
    \PSc@mment{curveTD Points=#1}\s@uvc@ntr@l\et@tpscurveTD\let\c@lprojSP=\relax%
    \iffillm@de\Psc@rveTD\curv@roundness[#1]%
    \f@gfill%
    \else\Psc@rveTD\curv@roundness[#1]\f@gstroke\fi%
    \PSc@mment{End curveTD}\resetc@ntr@l\et@tpscurveTD\fi\fi}}
\ctr@ld@f\def\Psc@rveDD#1[#2]{%
    \def\list@num{#2}\extrairelepremi@r\Ak@\de\list@num%
    \extrairelepremi@r\Ai@\de\list@num\extrairelepremi@r\Aj@\de\list@num%
    \f@gnewpath\PSwrit@cmdS{\Ai@}{\c@mmoveto}{\fwf@g}{\X@un}{\Y@un}%
    \setc@ntr@l{2}\figvectPDD -1[\Ak@,\Aj@]%
    \@ecfor\Ak@:=\list@num\do{\figpttraDD-2:=\Ai@/#1,-1/\BdingB@xfalse%
       \PSwrit@cmdS{-2}{}{\fwf@g}{\X@de}{\Y@de}%
       \figvectPDD -1[\Ai@,\Ak@]\figpttraDD-2:=\Aj@/-#1,-1/%
       \PSwrit@cmdS{-2}{}{\fwf@g}{\X@tr}{\Y@tr}\BdingB@xtrue%
       \PSwrit@cmdS{\Aj@}{\c@mcurveto}{\fwf@g}{\X@qu}{\Y@qu}%
       \B@zierBB@x{1}{\Y@un}(\X@un,\X@de,\X@tr,\X@qu)%
       \B@zierBB@x{2}{\X@un}(\Y@un,\Y@de,\Y@tr,\Y@qu)%
       \edef\X@un{\X@qu}\edef\Y@un{\Y@qu}\edef\Ai@{\Aj@}\edef\Aj@{\Ak@}}}
\ctr@ld@f\def\Psc@rveTD#1[#2]{\ifnum\CUR@proj<\tw@\Psc@rvePPTD#1[#2]\else\Psc@rveCPTD#1[#2]\fi}
\ctr@ld@f\def\Psc@rvePPTD#1[#2]{\setc@ntr@l{2}%
    \def\list@num{#2}\extrairelepremi@r\Ak@\de\list@num\Figptpr@j-5:/\Ak@/%
    \extrairelepremi@r\Ai@\de\list@num\Figptpr@j-3:/\Ai@/%
    \extrairelepremi@r\Aj@\de\list@num\Figptpr@j-4:/\Aj@/%
    \f@gnewpath\PSwrit@cmdS{-3}{\c@mmoveto}{\fwf@g}{\X@un}{\Y@un}%
    \figvectPDD -1[-5,-4]%
    \@ecfor\Ak@:=\list@num\do{\Figptpr@j-5:/\Ak@/\figpttraDD-2:=-3/#1,-1/%
       \BdingB@xfalse\PSwrit@cmdS{-2}{}{\fwf@g}{\X@de}{\Y@de}%
       \figvectPDD -1[-3,-5]\figpttraDD-2:=-4/-#1,-1/%
       \PSwrit@cmdS{-2}{}{\fwf@g}{\X@tr}{\Y@tr}\BdingB@xtrue%
       \PSwrit@cmdS{-4}{\c@mcurveto}{\fwf@g}{\X@qu}{\Y@qu}%
       \B@zierBB@x{1}{\Y@un}(\X@un,\X@de,\X@tr,\X@qu)%
       \B@zierBB@x{2}{\X@un}(\Y@un,\Y@de,\Y@tr,\Y@qu)%
       \edef\X@un{\X@qu}\edef\Y@un{\Y@qu}\figptcopyDD-3:/-4/\figptcopyDD-4:/-5/}}
\ctr@ld@f\def\Psc@rveCPTD#1[#2]{\setc@ntr@l{2}%
    \def\list@num{#2}\extrairelepremi@r\Ak@\de\list@num%
    \extrairelepremi@r\Ai@\de\list@num\extrairelepremi@r\Aj@\de\list@num%
    \Figptpr@j-7:/\Ai@/%
    \f@gnewpath\PSwrit@cmd{-7}{\c@mmoveto}{\fwf@g}%
    \figvectPTD -9[\Ak@,\Aj@]%
    \@ecfor\Ak@:=\list@num\do{\figpttraTD-10:=\Ai@/#1,-9/%
       \figvectPTD -9[\Ai@,\Ak@]\figpttraTD-11:=\Aj@/-#1,-9/%
       \subB@zierTD\NBz@rcs[\Ai@,-10,-11,\Aj@]\edef\Ai@{\Aj@}\edef\Aj@{\Ak@}}}
\ctr@ld@f\def\figdrawend{\ifCUR@PS\ifGR@cri\immediate\closeout\fwf@g%
    \immediate\openout\fwf@g=\PSfilen@me\relax%
    \ifPDFm@ke\PSBdingB@x\else%
    \immediate\write\fwf@g{\p@urcent\string!PS-Adobe-2.0 EPSF-2.0}%
    \PShe@der{Creator\string: TeX (fig4tex.tex)}%
    \PShe@der{Title\string: \PSfilen@me}%
    \PShe@der{CreationDate\string: \the\day/\the\month/\the\year}%
    \PSBdingB@x%
    \PShe@der{EndComments}\PSdict@\fi%
    \immediate\write\fwf@g{\c@mgsave}%
    \openin\frf@g=\auxfilen@me\c@pypsfile\fwf@g\frf@g\closein\frf@g%
    \immediate\write\fwf@g{\c@mgrestore}%
    \PSc@mment{End of file.}\immediate\closeout\fwf@g%
    \immediate\openout\fwf@g=\auxfilen@me\immediate\closeout\fwf@g%
    \immediate\write16{File \PSfilen@me\space created.}\fi\fi\CUR@PSfalse\GR@critrue}
\ctr@ld@f\def\PShe@der#1{\immediate\write\fwf@g{\p@urcent\p@urcent#1}}
\ctr@ld@f\def\PSBdingB@x{{\v@lX=\ptT@ptps\c@@rdXmin\v@lY=\ptT@ptps\c@@rdYmin%
     \v@lXa=\ptT@ptps\c@@rdXmax\v@lYa=\ptT@ptps\c@@rdYmax%
     \PShe@der{BoundingBox\string: \repdecn@mb{\v@lX}\space\repdecn@mb{\v@lY}%
     \space\repdecn@mb{\v@lXa}\space\repdecn@mb{\v@lYa}}}}
\ctr@ld@f\def\figdrawfcconnect[#1]{{\ifCUR@PS\ifGR@cri\PSc@mment{fcconnect Points=#1}%
    \Q@s@tfillmode{no}\s@uvc@ntr@l\et@tpsfcconnect\resetc@ntr@l{2}%
    \fcc@nnect@[#1]\resetc@ntr@l\et@tpsfcconnect\PSc@mment{End fcconnect}\fi\fi}}
\ctr@ld@f\def\fcc@nnect@[#1]{\let\N@rm=\n@rmeucDD\def\list@num{#1}%
    \extrairelepremi@r\Ai@\de\list@num\edef\pr@m{\Ai@}\v@leur=\z@\p@rtent=\@ne\c@llgtot%
    \ifcase\fclin@typ@\edef\list@num{[\pr@m,#1,\Ai@}\expandafter\figdrawcurve\list@num]%
    \else\ifdim\fclin@r@d\p@>\z@\Pslin@conge[#1]\else\figdrawline[#1]\fi\fi%
    \v@leur=\@rrowp@s\v@leur\edef\list@num{#1,\Ai@,0}%
    \extrairelepremi@r\Ai@\de\list@num\mili@u=\epsil@n\c@llgpart%
    \advance\mili@u-\epsil@n\advance\mili@u-\delt@\advance\v@leur-\mili@u%
    \ifcase\fclin@typ@\invers@\mili@u\delt@%
    \ifnum\@rrowr@fpt>\z@\advance\delt@-\v@leur\v@leur=\delt@\fi%
    \v@leur=\repdecn@mb\v@leur\mili@u\edef\v@lt{\repdecn@mb\v@leur}%
    \extrairelepremi@r\Ak@\de\list@num%
    \figvectPDD-1[\pr@m,\Aj@]\figpttraDD-6:=\Ai@/\curv@roundness,-1/%
    \figvectPDD-1[\Ak@,\Ai@]\figpttraDD-7:=\Aj@/\curv@roundness,-1/%
    \delt@=\@rrowheadlength\p@\delt@=\C@AHANG\delt@\edef\R@dius{\repdecn@mb{\delt@}}%
    \ifcase\@rrowr@fpt%
    \FigptintercircB@zDD-8::\v@lt,\R@dius[\Ai@,-6,-7,\Aj@]\Q@arrowheadDD[-5,-8]\else%
    \FigptintercircB@zDD-8::\v@lt,\R@dius[\Aj@,-7,-6,\Ai@]\Q@arrowheadDD[-8,-5]\fi%
    \else\advance\delt@-\v@leur%
    \p@rtentiere{\p@rtent}{\delt@}\edef\C@efun{\the\p@rtent}%
    \p@rtentiere{\p@rtent}{\v@leur}\edef\C@efde{\the\p@rtent}%
    \figptbaryDD-5:[\Ai@,\Aj@;\C@efun,\C@efde]\ifcase\@rrowr@fpt%
    \delt@=\@rrowheadlength\unit@\delt@=\C@AHANG\delt@\edef\t@ille{\repdecn@mb{\delt@}}%
    \figvectPDD-2[\Ai@,\Aj@]\vecunit@{-2}{-2}\figpttraDD-5:=-5/\t@ille,-2/\fi%
    \Q@arrowheadDD[\Ai@,-5]\fi}
\ctr@ld@f\def\c@llgtot{\@ecfor\Aj@:=\list@num\do{\figvectP-1[\Ai@,\Aj@]\N@rm\delt@{-1}%
    \advance\v@leur\delt@\advance\p@rtent\@ne\edef\Ai@{\Aj@}}}
\ctr@ld@f\def\c@llgpart{\extrairelepremi@r\Aj@\de\list@num\figvectP-1[\Ai@,\Aj@]\N@rm\delt@{-1}%
    \advance\mili@u\delt@\ifdim\mili@u<\v@leur\edef\pr@m{\Ai@}\edef\Ai@{\Aj@}\c@llgpart\fi}
\ctr@ld@f\def\Pslin@conge[#1]{\ifnum\p@rtent>\tw@{\def\list@num{#1}%
    \extrairelepremi@r\Ai@\de\list@num\extrairelepremi@r\Aj@\de\list@num%
    \figptcopy-6:/\Ai@/\figvectP-3[\Ai@,\Aj@]\vecunit@{-3}{-3}\v@lmax=\result@t%
    \@ecfor\Ak@:=\list@num\do{\figvectP-4[\Aj@,\Ak@]\vecunit@{-4}{-4}%
    \minim@m\v@lmin\v@lmax\result@t\v@lmax=\result@t%
    \det@rm\delt@[-3,-4]\maxim@m\mili@u{\delt@}{-\delt@}\ifdim\mili@u>\Cepsil@n%
    \ifdim\delt@>\z@\figgetangleDD\Angl@[\Aj@,\Ak@,\Ai@]\else%
    \figgetangleDD\Angl@[\Aj@,\Ai@,\Ak@]\fi%
    \v@leur=\PI@deg\advance\v@leur-\Angl@\p@\divide\v@leur\tw@%
    \edef\Angl@{\repdecn@mb\v@leur}\c@ssin{\C@}{\S@}{\Angl@}\v@leur=\fclin@r@d\unit@%
    \v@leur=\S@\v@leur\mili@u=\C@\p@\invers@\mili@u\mili@u%
    \v@leur=\repdecn@mb{\mili@u}\v@leur%
    \minim@m\v@leur\v@leur\v@lmin\edef\t@ille{\repdecn@mb{\v@leur}}%
    \figpttra-5:=\Aj@/-\t@ille,-3/\figdrawline[-6,-5]\figpttra-6:=\Aj@/\t@ille,-4/%
    \figvectNVDD-3[-3]\figvectNVDD-8[-4]\inters@cDD-7:[-5,-3;-6,-8]%
    \ifdim\delt@>\z@\figdrawarccircP-7;\fclin@r@d[-5,-6]\else\figdrawarccircP-7;\fclin@r@d[-6,-5]\fi%
    \else\figdrawline[-6,\Aj@]\figptcopy-6:/\Aj@/\fi
    \edef\Ai@{\Aj@}\edef\Aj@{\Ak@}\figptcopy-3:/-4/}\figdrawline[-6,\Aj@]}\else\figdrawline[#1]\fi}
\ctr@ld@f\def\figdrawfcnode[#1]#2{{\ifCUR@PS\ifGR@cri\PSc@mment{fcnode Points=#1}%
    \s@uvc@ntr@l\et@tpsfcnode\resetc@ntr@l{2}%
    \def\t@xt@{#2}\ifx\t@xt@\empty\def\g@tt@xt{\setbox\Gb@x=\hbox{\Figg@tT{\p@int}}}%
    \else\def\g@tt@xt{\setbox\Gb@x=\hbox{#2}}\fi%
    \v@lmin=\h@rdfcXp@dd\advance\v@lmin\Xp@dd\unit@\multiply\v@lmin\tw@%
    \v@lmax=\h@rdfcYp@dd\advance\v@lmax\Yp@dd\unit@\multiply\v@lmax\tw@%
    \Figv@ctCreg-8(\unit@,-\unit@)\def\list@num{#1}%
    \delt@=\CUR@width bp\divide\delt@\tw@%
    \fcn@de\PSc@mment{End fcnode}\resetc@ntr@l\et@tpsfcnode\fi\fi}}
\ctr@ld@f\def\d@butn@de{\g@tt@xt\v@lX=\wd\Gb@x%
    \v@lY=\ht\Gb@x\advance\v@lY\dp\Gb@x\advance\v@lX\v@lmin\advance\v@lY\v@lmax}
\ctr@ld@f\def\fcn@deE{%
    \@ecfor\p@int:=\list@num\do{\d@butn@de\v@lX=\unssqrttw@\v@lX\v@lY=\unssqrttw@\v@lY%
    \ifdim\thickn@ss\p@>\z@
    \v@lXa=\v@lX\advance\v@lXa\delt@\v@lXa=\ptT@unit@\v@lXa\edef\XR@d{\repdecn@mb\v@lXa}%
    \v@lYa=\v@lY\advance\v@lYa\delt@\v@lYa=\ptT@unit@\v@lYa\edef\YR@d{\repdecn@mb\v@lYa}%
    \arct@n\v@leur(\v@lXa,\v@lYa)\v@leur=\rdT@deg\v@leur\edef\@nglde{\repdecn@mb\v@leur}%
    {\c@lptellDD-2::\p@int;\XR@d,\YR@d(\@nglde)}
    \advance\v@leur-\PI@deg\edef\@nglun{\repdecn@mb\v@leur}%
    {\c@lptellDD-3::\p@int;\XR@d,\YR@d(\@nglun)}%
    \figptstra-6=-3,-2,\p@int/\thickn@ss,-8/\Q@s@tfillmode{yes}%
    \Pss@tspecifSt{color=\DDV@thickcolor}%
    \figdrawline[-2,-3,-6,-5]\figdrawarcell-4;\XR@d,\YR@d(\@nglun,\@nglde,0)%
    \Psrest@reSt{color=\DDV@thickcolor}\fi
    \v@lX=\ptT@unit@\v@lX\v@lY=\ptT@unit@\v@lY%
    \edef\XR@d{\repdecn@mb\v@lX}\edef\YR@d{\repdecn@mb\v@lY}%
    \Q@s@tfillmode{yes}\Pss@tspecifSt{color=\fcbgc@lor}%
    \figdrawarcell\p@int;\XR@d,\YR@d(0,360,0)%
    \Q@s@tfillmode{no}\Psrest@reSt{color=\fcbgc@lor}\figdrawarcell\p@int;\XR@d,\YR@d(0,360,0)}}
\ctr@ld@f\def\fcn@deL{\delt@=\ptT@unit@\delt@\edef\t@ille{\repdecn@mb\delt@}%
    \@ecfor\p@int:=\list@num\do{\Figg@tXYa{\p@int}\d@butn@de%
    \ifdim\v@lX>\v@lY\itis@Ktrue\else\itis@Kfalse\fi%
    \advance\v@lXa-\v@lX\Figp@intreg-1:(\v@lXa,\v@lYa)%
    \advance\v@lXa\v@lX\advance\v@lYa-\v@lY\Figp@intreg-2:(\v@lXa,\v@lYa)%
    \advance\v@lXa\v@lX\advance\v@lYa\v@lY\Figp@intreg-3:(\v@lXa,\v@lYa)%
    \advance\v@lXa-\v@lX\advance\v@lYa\v@lY\Figp@intreg-4:(\v@lXa,\v@lYa)%
    \ifdim\thickn@ss\p@>\z@
    \Figg@tXYa{\p@int}\Q@s@tfillmode{yes}\Pss@tspecifSt{color=\DDV@thickcolor}%
    \c@lpt@xt{-1}{-4}\c@lpt@xt@\v@lXa\v@lYa\v@lX\v@lY\c@rre\delt@%
    \Figp@intregDD-9:(\v@lZ,\v@lYa)\Figp@intregDD-11:(\v@lZa,\v@lYa)%
    \c@lpt@xt{-4}{-3}\c@lpt@xt@\v@lYa\v@lXa\v@lY\v@lX\delt@\c@rre%
    \Figp@intregDD-12:(\v@lXa,\v@lZ)\Figp@intregDD-10:(\v@lXa,\v@lZa)%
    \ifitis@K\figptstra-7=-9,-10,-11/\thickn@ss,-8/\figdrawline[-9,-11,-5,-6,-7]\else%
    \figptstra-7=-10,-11,-12/\thickn@ss,-8/\figdrawline[-10,-12,-5,-6,-7]\fi%
    \Psrest@reSt{color=\DDV@thickcolor}\fi
    \Q@s@tfillmode{yes}\Pss@tspecifSt{color=\fcbgc@lor}\figdrawline[-1,-2,-3,-4]%
    \Q@s@tfillmode{no}\Psrest@reSt{color=\fcbgc@lor}\figdrawline[-1,-2,-3,-4,-1]}}
\ctr@ld@f\def\c@lpt@xt#1#2{\figvectN-7[#1,#2]\vecunit@{-7}{-7}\figpttra-5:=#1/\t@ille,-7/%
    \figvectP-7[#1,#2]\Figg@tXY{-7}\c@rre=\v@lX\delt@=\v@lY\Figg@tXY{-5}}
\ctr@ld@f\def\c@lpt@xt@#1#2#3#4#5#6{\v@lZ=#6\invers@{\v@lZ}{\v@lZ}\v@leur=\repdecn@mb{#5}\v@lZ%
    \v@lZ=#2\advance\v@lZ-#4\mili@u=\repdecn@mb{\v@leur}\v@lZ%
    \v@lZ=#3\advance\v@lZ\mili@u\v@lZa=-\v@lZ\advance\v@lZa\tw@#1}
\ctr@ld@f\def\fcn@deR{\@ecfor\p@int:=\list@num\do{\Figg@tXYa{\p@int}\d@butn@de%
    \advance\v@lXa-0.5\v@lX\advance\v@lYa-0.5\v@lY\Figp@intreg-1:(\v@lXa,\v@lYa)%
    \advance\v@lXa\v@lX\Figp@intreg-2:(\v@lXa,\v@lYa)%
    \advance\v@lYa\v@lY\Figp@intreg-3:(\v@lXa,\v@lYa)%
    \advance\v@lXa-\v@lX\Figp@intreg-4:(\v@lXa,\v@lYa)%
    \ifdim\thickn@ss\p@>\z@
    \Q@s@tfillmode{yes}\Pss@tspecifSt{color=\DDV@thickcolor}%
    \Figv@ctCreg-5(-\delt@,-\delt@)\figpttra-9:=-1/1,-5/%
    \Figv@ctCreg-5(\delt@,-\delt@)\figpttra-10:=-2/1,-5/%
    \Figv@ctCreg-5(\delt@,\delt@)\figpttra-11:=-3/1,-5/%
    \figptstra-7=-9,-10,-11/\thickn@ss,-8/\figdrawline[-9,-11,-5,-6,-7]%
    \Psrest@reSt{color=\DDV@thickcolor}\fi
    \Q@s@tfillmode{yes}\Pss@tspecifSt{color=\fcbgc@lor}\figdrawline[-1,-2,-3,-4]%
    \Q@s@tfillmode{no}\Psrest@reSt{color=\fcbgc@lor}\figdrawline[-1,-2,-3,-4,-1]}}
\ctr@ld@f\def\Pssetfl@wchart#1=#2|{\keln@mtr#1|%
    \def\n@mref{arr}\ifx\l@debut\n@mref\expandafter\keln@mtr\l@suite|%
     \def\n@mref{owp}\ifx\l@debut\n@mref\update@ttr\D@FTfcarrowposition\P@setfcarrowposition{#2}\else
     \def\n@mref{owr}\ifx\l@debut\n@mref\update@ttr\D@FTfcarrowrefpt\P@setfcarrowrefpt{#2}\else
     \W@rnmesAttr{figset flowchart}{#1}\fi\fi\else%
    \def\n@mref{bgc}\ifx\l@debut\n@mref\update@ttr\D@FTfcbgcolor\P@setfcbgcolor{#2}\else
    \def\n@mref{lin}\ifx\l@debut\n@mref\update@ttr\D@FTfcline\P@setfcline{#2}\else
    \def\n@mref{pad}\ifx\l@debut\n@mref\update@ttr\D@FTfcxpadding\P@setfcxpadding{#2}%
                                       \update@ttr\D@FTfcypadding\P@setfcypadding{#2}\else
    \def\n@mref{rad}\ifx\l@debut\n@mref\update@ttr\D@FTfcradius\P@setfcradius{#2}\else
    \def\n@mref{sha}\ifx\l@debut\n@mref\update@ttr\D@FTfcshape\P@setfcshape{#2}\else
    \def\n@mref{thi}\ifx\l@debut\n@mref\expandafter\keln@mtr\l@suite|%
     \def\n@mref{ckc}\ifx\l@debut\n@mref\update@ttr\D@FTref\P@setfcthickcolor{#2}\else
     \def\n@mref{ckn}\ifx\l@debut\n@mref\update@ttr\D@FTfcthickness\P@setfcthickness{#2}\else
     \W@rnmesAttr{figset flowchart}{#1}\fi\fi\else%
    \def\n@mref{xpa}\ifx\l@debut\n@mref\update@ttr\D@FTfcxpadding\P@setfcxpadding{#2}\else
    \def\n@mref{ypa}\ifx\l@debut\n@mref\update@ttr\D@FTfcypadding\P@setfcypadding{#2}\else
    \W@rnmesAttr{figset flowchart}{#1}\fi\fi\fi\fi\fi\fi\fi\fi\fi}
\ctr@ln@m\@rrowp@s
\ctr@ld@f\def\P@setfcarrowposition#1{\edef\@rrowp@s{#1}}
\ctr@ln@m\@rrowr@fpt
\ctr@ld@f\def\P@setfcarrowrefpt#1{\setfcr@fpt#1|}
\ctr@ld@f\def\setfcr@fpt#1#2|{\if#1e\def\@rrowr@fpt{1}\else\def\@rrowr@fpt{0}\fi}
\ctr@ln@m\fcbgc@lor
\ctr@ld@f\def\P@setfcbgcolor#1{\edef\fcbgc@lor{#1}}
\ctr@ln@m\fclin@typ@
\ctr@ld@f\def\P@setfcline#1{\setfccurv@#1|}
\ctr@ld@f\def\setfccurv@#1#2|{\if#1c\def\fclin@typ@{0}\else\def\fclin@typ@{1}\fi}
\ctr@ln@m\fclin@r@d
\ctr@ld@f\def\P@setfcradius#1{\edef\fclin@r@d{#1}}
\ctr@ln@m\fcn@de \ctr@ln@m\fcsh@pe
\ctr@ln@m\h@rdfcXp@dd \ctr@ln@m\h@rdfcYp@dd
\ctr@ld@f\def\P@setfcshape#1{\setfcshap@#1|}
\ctr@ld@f\def\setfcshap@#1#2|{%
    \if#1e\let\fcn@de=\fcn@deE\def\h@rdfcXp@dd{4pt}\def\h@rdfcYp@dd{4pt}%
     \edef\fcsh@pe{ellipse}\else%
    \if#1l\let\fcn@de=\fcn@deL\def\h@rdfcXp@dd{4pt}\def\h@rdfcYp@dd{4pt}%
     \edef\fcsh@pe{lozenge}\else%
          \let\fcn@de=\fcn@deR\def\h@rdfcXp@dd{6pt}\def\h@rdfcYp@dd{6pt}%
     \edef\fcsh@pe{rectangle}\fi\fi}
\ctr@ln@m\DDV@thickcolor
\ctr@ld@f\def\P@setfcthickcolor#1{\edef\DDV@thickcolor{#1}}
\ctr@ln@m\thickn@ss
\ctr@ld@f\def\P@setfcthickness#1{\edef\thickn@ss{#1}}
\ctr@ln@m\Xp@dd
\ctr@ld@f\def\P@setfcxpadding#1{\edef\Xp@dd{#1}}
\ctr@ln@m\Yp@dd
\ctr@ld@f\def\P@setfcypadding#1{\edef\Yp@dd{#1}}
\ctr@ld@f\def\figdrawline[#1]{{\ifCUR@PS\ifGR@cri\PSc@mment{line Points=#1}%
    \let\figdrawlign@=\Pslign@P\Pslin@{#1}\PSc@mment{End line}\fi\fi}}
\ctr@ld@f\def\figdrawlineF#1{{\ifCUR@PS\ifGR@cri\PSc@mment{lineF Filename=#1}%
    \let\figdrawlign@=\Pslign@F\Pslin@{#1}\PSc@mment{End lineF}\fi\fi}}
\ctr@ld@f\def\figdrawlineC(#1){{\ifCUR@PS\ifGR@cri\PSc@mment{lineC}%
    \let\figdrawlign@=\Pslign@C\Pslin@{#1}\PSc@mment{End lineC}\fi\fi}}
\ctr@ld@f\def\Pslin@#1{\iffillm@de\figdrawlign@{#1}%
    \f@gfill%
    \else\figdrawlign@{#1}\ifx\derp@int\premp@int%
    \f@gclosestroke%
    \else\f@gstroke\fi\fi}
\ctr@ld@f\def\Pslign@P#1{\def\list@num{#1}\extrairelepremi@r\p@int\de\list@num%
    \edef\premp@int{\p@int}\f@gnewpath%
    \PSwrit@cmd{\p@int}{\c@mmoveto}{\fwf@g}%
    \@ecfor\p@int:=\list@num\do{\PSwrit@cmd{\p@int}{\c@mlineto}{\fwf@g}%
    \edef\derp@int{\p@int}}}
\ctr@ld@f\def\Pslign@F#1{\s@uvc@ntr@l\et@tPslign@F\setc@ntr@l{2}\openin\frf@g=#1\relax%
    \ifeof\frf@g\message{*** File #1 not found !}\end\else%
    \read\frf@g to\tr@c\edef\premp@int{\tr@c}\expandafter\extr@ctCF\tr@c:%
    \f@gnewpath\PSwrit@cmd{-1}{\c@mmoveto}{\fwf@g}%
    \loop\read\frf@g to\tr@c\ifeof\frf@g\mored@tafalse\else\mored@tatrue\fi%
    \ifmored@ta\expandafter\extr@ctCF\tr@c:\PSwrit@cmd{-1}{\c@mlineto}{\fwf@g}%
    \edef\derp@int{\tr@c}\repeat\fi\closein\frf@g\resetc@ntr@l\et@tPslign@F}
\ctr@ln@m\extr@ctCF
\ctr@ld@f\def\extr@ctCFDD#1 #2:{\v@lX=#1\unit@\v@lY=#2\unit@\Figp@intregDD-1:(\v@lX,\v@lY)}
\ctr@ld@f\def\extr@ctCFTD#1 #2 #3:{\v@lX=#1\unit@\v@lY=#2\unit@\v@lZ=#3\unit@%
    \Figp@intregTD-1:(\v@lX,\v@lY,\v@lZ)}
\ctr@ld@f\def\Pslign@C#1{\s@uvc@ntr@l\et@tPslign@C\setc@ntr@l{2}%
    \def\list@num{#1}\extrairelepremi@r\p@int\de\list@num%
    \edef\premp@int{\p@int}\f@gnewpath%
    \expandafter\Pslign@C@\p@int:\PSwrit@cmd{-1}{\c@mmoveto}{\fwf@g}%
    \@ecfor\p@int:=\list@num\do{\expandafter\Pslign@C@\p@int:%
    \PSwrit@cmd{-1}{\c@mlineto}{\fwf@g}\edef\derp@int{\p@int}}%
    \resetc@ntr@l\et@tPslign@C}
\ctr@ld@f\def\Pslign@C@#1 #2:{{\def\t@xt@{#1}\ifx\t@xt@\empty\Pslign@C@#2:
    \else\extr@ctCF#1 #2:\fi}}
\ctr@ld@f\def\Pssetm@sh#1=#2|{\keln@mde#1|%
    \def\n@mref{co}\ifx\l@debut\n@mref\update@ttr\D@FTref\P@setmeshcolor{#2}\else
    \def\n@mref{da}\ifx\l@debut\n@mref\update@ttr\D@FTref\P@setmeshdash{#2}\else
    \def\n@mref{di}\ifx\l@debut\n@mref\update@ttr\D@FTmeshdiag\Q@s@tmeshdiag{#2}\else
    \def\n@mref{wi}\ifx\l@debut\n@mref\update@ttr\D@FTref\P@setmeshwidth{#2}\else
    \W@rnmesAttr{figset mesh}{#1}\fi\fi\fi\fi}
\ctr@ln@m\c@ntrolmesh
\ctr@ld@f\def\Q@s@tmeshdiag#1{\edef\c@ntrolmesh{#1}}
\ctr@ld@f\def\D@FTmeshdiag{0}    
\ctr@ln@m\DDV@meshcolor
\ctr@ld@f\def\P@setmeshcolor#1{\edef\DDV@meshcolor{#1}}
\ctr@ln@m\DDV@meshdash
\ctr@ld@f\def\P@setmeshdash#1{\edef\DDV@meshdash{#1}}
\ctr@ln@m\DDV@meshwidth
\ctr@ld@f\def\P@setmeshwidth#1{\edef\DDV@meshwidth{#1}}
\ctr@ld@f\def\figdrawmesh#1,#2[#3,#4,#5,#6]{{\ifCUR@PS\ifGR@cri%
    \PSc@mment{mesh N1=#1, N2=#2, Quadrangle=[#3,#4,#5,#6]}\s@uvc@ntr@l\et@tpsmesh%
    \Pss@tspecifSt{color=\DDV@meshcolor,dash=\DDV@meshdash,width=\DDV@meshwidth}%
    \setc@ntr@l{2}%
    \ifnum#1>\@ne\Psmeshp@rt#1[#3,#4,#5,#6]\fi%
    \ifnum#2>\@ne\Psmeshp@rt#2[#4,#5,#6,#3]\fi%
    \ifnum\c@ntrolmesh>\z@\Psmeshdi@g#1,#2[#3,#4,#5,#6]\fi%
    \ifnum\c@ntrolmesh<\z@\Psmeshdi@g#2,#1[#4,#5,#6,#3]\fi%
    \Psrest@reSt{color=\DDV@meshcolor,dash=\DDV@meshdash,width=\DDV@meshwidth}%
    \figdrawline[#3,#4,#5,#6,#3]\PSc@mment{End mesh}\resetc@ntr@l\et@tpsmesh\fi\fi}}
\ctr@ld@f\def\Psmeshp@rt#1[#2,#3,#4,#5]{{\l@mbd@un=\@ne\l@mbd@de=#1\loop%
    \ifnum\l@mbd@un<#1\advance\l@mbd@de\m@ne\figptbary-1:[#2,#3;\l@mbd@de,\l@mbd@un]%
    \figptbary-2:[#5,#4;\l@mbd@de,\l@mbd@un]\figdrawline[-1,-2]\advance\l@mbd@un\@ne\repeat}}
\ctr@ld@f\def\Psmeshdi@g#1,#2[#3,#4,#5,#6]{\figptcopy-2:/#3/\figptcopy-3:/#6/%
    \l@mbd@un=\z@\l@mbd@de=#1\loop\ifnum\l@mbd@un<#1%
    \advance\l@mbd@un\@ne\advance\l@mbd@de\m@ne\figptcopy-1:/-2/\figptcopy-4:/-3/%
    \figptbary-2:[#3,#4;\l@mbd@de,\l@mbd@un]%
    \figptbary-3:[#6,#5;\l@mbd@de,\l@mbd@un]\Psmeshdi@gp@rt#2[-1,-2,-3,-4]\repeat}
\ctr@ld@f\def\Psmeshdi@gp@rt#1[#2,#3,#4,#5]{{\l@mbd@un=\z@\l@mbd@de=#1\loop%
    \ifnum\l@mbd@un<#1\figptbary-5:[#2,#5;\l@mbd@de,\l@mbd@un]%
    \advance\l@mbd@de\m@ne\advance\l@mbd@un\@ne%
    \figptbary-6:[#3,#4;\l@mbd@de,\l@mbd@un]\figdrawline[-5,-6]\repeat}}
\ctr@ln@m\figdrawnormal
\ctr@ld@f\def\Q@normalDD#1,#2[#3,#4]{{\ifCUR@PS\ifGR@cri%
    \PSc@mment{normal Length=#1, Lambda=#2 [Pt1,Pt2]=[#3,#4]}%
    \s@uvc@ntr@l\et@tpsnormal\resetc@ntr@l{2}\figptendnormal-6::#1,#2[#3,#4]%
    \figptcopyDD-5:/-1/\figdrawarrow[-5,-6]%
    \PSc@mment{End normal}\resetc@ntr@l\et@tpsnormal\fi\fi}}
\ctr@ld@f\def\figreset#1{\trtlis@rg{#1}{\Psreset@}}
\ctr@ld@f\def\Psreset@#1|{\def\t@xt@{#1}\ifx\t@xt@\empty\P@resetg@n
    \else\keln@mde#1|%
    \def\n@mref{al}\ifx\l@debut\n@mref%
        \figset altitude(blcolor=default,bldash=default,blwidth=default,%
        sqcolor=default,sqdash=default,sqwidth=default)\else
    \def\n@mref{ar}\ifx\l@debut\n@mref\figresetarrowhead\else
    \def\n@mref{cu}\ifx\l@debut\n@mref\figset curve(roundness=\D@FTroundness)\else
    \def\n@mref{ge}\ifx\l@debut\n@mref\P@resetg@n\else
    \def\n@mref{fl}\ifx\l@debut\n@mref%
        \figset flowchart(arrowp=\D@FTfcarrowposition,arrowr=\D@FTfcarrowrefpt,%
	bgcolor=\D@FTfcbgcolor,line=\D@FTfcline,radius=\D@FTfcradius,%
	shape=\D@FTfcshape,thickcolor=default,thickness=\D@FTfcthickness,%
	xpadd=\D@FTfcxpadding,ypadd=\D@FTfcypadding)\else
    \def\n@mref{me}\ifx\l@debut\n@mref\figset mesh(diag=\D@FTmeshdiag,%
        color=default,dash=default,width=default)\else
    \def\n@mref{tr}\ifx\l@debut\n@mref%
        \figset trimesh(color=default,dash=default,width=default)\else
    \W@rnmeskwd{figreset}{#1}\fi\fi\fi\fi\fi\fi\fi\fi}
\ctr@ld@f\def\P@resetg@n{\figset (color=\D@FTcolor,dash=\D@FTdash,fill=\D@FTfill,%
    join=\D@FTjoin,width=\D@FTwidth)}
\ctr@ld@f\def\figset#1(#2){\def\t@xt@{#1}\ifx\t@xt@\empty\trtlis@rg{#2}{\Pssetg@n}
    \else\keln@mde#1|%
    \def\n@mref{al}\ifx\l@debut\n@mref\trtlis@rg{#2}{\Psset@lti}\else
    \def\n@mref{ar}\ifx\l@debut\n@mref\trtlis@rg{#2}{\Psset@rrowhe@d}\else
    \def\n@mref{cu}\ifx\l@debut\n@mref\trtlis@rg{#2}{\Pssetc@rve}\else
    \def\n@mref{fl}\ifx\l@debut\n@mref\trtlis@rg{#2}{\Pssetfl@wchart}\else
    \def\n@mref{ge}\ifx\l@debut\n@mref\trtlis@rg{#2}{\Pssetg@n}\else
    \def\n@mref{me}\ifx\l@debut\n@mref\trtlis@rg{#2}{\Pssetm@sh}\else
    \def\n@mref{pr}\ifx\l@debut\n@mref\ifCUR@PS\W@rnmesIgn{figset proj(...)}%
     \else\trtlis@rg{#2}{\Figsetpr@j}\fi\else
    \def\n@mref{tr}\ifx\l@debut\n@mref\trtlis@rg{#2}{\Pssettrim@sh}\else
    \def\n@mref{wr}\ifx\l@debut\n@mref\let\M@cro=\Figsetwr@te\trtlis@rgtok{#2,|}\else
    \W@rnmeskwd{figset}{#1}\fi\fi\fi\fi\fi\fi\fi\fi\fi\fi\ignorespaces}
\ctr@ld@f\def\figsetdefault#1(#2){\ifCUR@PS\W@rnmesIgn{figsetdefault}\else%
    \def\t@xt@{#1}\ifx\t@xt@\empty\trtlis@rg{#2}{\Pssd@g@n}\else\keln@mun#1|
    \def\n@mref{a}\ifx\l@debut\n@mref\trtlis@rg{#2}{\Pssd@@rrowhe@d}\else
    \def\n@mref{c}\ifx\l@debut\n@mref\trtlis@rg{#2}{\Pssd@c@rve}\else
    \def\n@mref{g}\ifx\l@debut\n@mref\trtlis@rg{#2}{\Pssd@g@n}\else
    \def\n@mref{f}\ifx\l@debut\n@mref\trtlis@rg{#2}{\Pssd@fl@wchart}\else
    \def\n@mref{m}\ifx\l@debut\n@mref\trtlis@rg{#2}{\Pssd@m@sh}\else
    \W@rnmeskwd{figsetdefault}{#1}\fi\fi\fi\fi\fi\fi\initpss@ttings\fi}
\ctr@ld@f\def\Pssd@g@n#1=#2|{\keln@mun#1|%
    \def\n@mref{c}\ifx\l@debut\n@mref\edef\D@FTcolor{#2}\else
    \def\n@mref{d}\ifx\l@debut\n@mref\edef\D@FTdash{#2}\else
    \def\n@mref{f}\ifx\l@debut\n@mref\edef\D@FTfill{#2}\else
    \def\n@mref{j}\ifx\l@debut\n@mref\edef\D@FTjoin{#2}\else
    \def\n@mref{u}\ifx\l@debut\n@mref\edef\D@FTupdate{#2}\Q@s@tupdate{#2}\else
    \def\n@mref{w}\ifx\l@debut\n@mref\edef\D@FTwidth{#2}\else
    \W@rnmesAttr{figsetdefault}{#1}\fi\fi\fi\fi\fi\fi}
\ctr@ld@f\def\Pssd@@rrowhe@d#1=#2|{\keln@mun#1|%
    \def\n@mref{a}\ifx\l@debut\n@mref\edef\D@FTarrowheadangle{#2}\else
    \def\n@mref{f}\ifx\l@debut\n@mref\edef\D@FTarrowheadfill{#2}\else
    \def\n@mref{l}\ifx\l@debut\n@mref\y@tiunit{#2}\ifunitpr@sent%
     \edef\D@FTh@rdahlength{#2}\else\edef\D@FTh@rdahlength{#2pt}%
     \message{*** \BS@ figsetdefault (..., #1=#2, ...) : unit is missing, pt is assumed.}%
     \fi\else
    \def\n@mref{o}\ifx\l@debut\n@mref\edef\D@FTarrowheadout{#2}\else
    \def\n@mref{r}\ifx\l@debut\n@mref\edef\D@FTarrowheadratio{#2}\else
    \W@rnmesAttr{figsetdefault arrowhead}{#1}\fi\fi\fi\fi\fi}
\ctr@ld@f\def\Pssd@c@rve#1=#2|{\keln@mun#1|%
    \def\n@mref{r}\ifx\l@debut\n@mref\edef\D@FTroundness{#2}\else%
    \W@rnmesAttr{figsetdefault curve}{#1}\fi}
\ctr@ld@f\def\Pssd@fl@wchart#1=#2|{\keln@mtr#1|%
    \def\n@mref{arr}\ifx\l@debut\n@mref\expandafter\keln@mtr\l@suite|%
     \def\n@mref{owp}\ifx\l@debut\n@mref\edef\D@FTfcarrowposition{#2}\else
     \def\n@mref{owr}\ifx\l@debut\n@mref\edef\D@FTfcarrowrefpt{#2}\else
                     \W@rnmesAttr{figsetdefault flowchart}{#1}\fi\fi\else%
    \def\n@mref{bgc}\ifx\l@debut\n@mref\edef\D@FTfcbgcolor{#2}\else
    \def\n@mref{lin}\ifx\l@debut\n@mref\edef\D@FTfcline{#2}\else
    \def\n@mref{pad}\ifx\l@debut\n@mref\edef\D@FTfcxpadding{#2}%
                    \edef\D@FTfcypadding{#2}\else
    \def\n@mref{rad}\ifx\l@debut\n@mref\edef\D@FTfcradius{#2}\else
    \def\n@mref{sha}\ifx\l@debut\n@mref\edef\D@FTfcshape{#2}\else
    \def\n@mref{thi}\ifx\l@debut\n@mref\expandafter\keln@mtr\l@suite|%
     \def\n@mref{ckn}\ifx\l@debut\n@mref\edef\D@FTfcthickness{#2}\else
                     \W@rnmesAttr{figsetdefault flowchart}{#1}\fi\else%
    \def\n@mref{xpa}\ifx\l@debut\n@mref\edef\D@FTfcxpadding{#2}\else
    \def\n@mref{ypa}\ifx\l@debut\n@mref\edef\D@FTfcypadding{#2}\else
    \W@rnmesAttr{figsetdefault flowchart}{#1}\fi\fi\fi\fi\fi\fi\fi\fi\fi}
\ctr@ld@f\def\D@FTfcarrowposition{0.5}
\ctr@ld@f\def\D@FTfcarrowrefpt{start}
\ctr@ld@f\def\D@FTfcbgcolor{1}
\ctr@ld@f\def\D@FTfcline{polygon}
\ctr@ld@f\def\D@FTfcradius{0}
\ctr@ld@f\def\D@FTfcshape{rectangle}
\ctr@ld@f\def\D@FTfcthickness{0}
\ctr@ld@f\def\D@FTfcxpadding{0}
\ctr@ld@f\def\D@FTfcypadding{0}
\ctr@ld@f\def\Pssd@m@sh#1=#2|{\keln@mun#1|%
    \def\n@mref{d}\ifx\l@debut\n@mref\edef\D@FTmeshdiag{#2}\else%
    \W@rnmesAttr{figsetdefault mesh}{#1}\fi}
\ctr@ln@w{newif}\iffillm@de
\ctr@ld@f\def\Q@s@tfillmode#1{\expandafter\setfillm@de#1:}
\ctr@ld@f\def\setfillm@de#1#2:{\if#1n\fillm@defalse\else\fillm@detrue\fi}
\ctr@ld@f\def\D@FTfill{no}     
\ctr@ln@w{newif}\ifGRupdatem@de
\ctr@ld@f\def\Q@s@tupdate#1{\ifCUR@PS\W@rnmesIgn{figset (update=...)}%
    \else\expandafter\setupd@te#1:\fi}
\ctr@ld@f\def\setupd@te#1#2:{\if#1n\GRupdatem@defalse\else\GRupdatem@detrue\fi}
\ctr@ld@f\def\D@FTupdate{no}     
\ctr@ln@m\CUR@color \ctr@ln@m\CUR@colorc@md
\ctr@ld@f\def\s@uvcolor#1{\edef#1{\CUR@color}}
\ctr@ld@f\def\D@FTcolor{0}       
\ctr@ld@f\def\Pssetc@lor#1{\ifGR@cri\result@tent=\@ne\expandafter\c@lnbV@l#1 :%
    \def\CUR@color{}\def\CUR@colorc@md{}%
    \ifcase\result@tent\or\Q@s@tgray{#1}\or\or\Q@s@trgb{#1}\or\Q@s@tcmyk{#1}\fi\fi}
\ctr@ln@m\CUR@colorc@mdStroke
\ctr@ld@f\def\Q@s@tcmyk#1{\ifGR@cri\def\CUR@color{#1}\def\CUR@colorc@md{\c@msetcmykcolor}%
    \def\CUR@colorc@mdStroke{\c@msetcmykcolorStroke}%
    \ifCUR@PS\PSc@mment{setcmyk Color=#1}\us@primarC@lor\fi\fi}
\ctr@ld@f\def\Q@s@trgb#1{\ifGR@cri\def\CUR@color{#1}\def\CUR@colorc@md{\c@msetrgbcolor}%
    \def\CUR@colorc@mdStroke{\c@msetrgbcolorStroke}%
    \ifCUR@PS\PSc@mment{setrgb Color=#1}\us@primarC@lor\fi\fi}
\ctr@ld@f\def\Q@s@tgray#1{\ifGR@cri\def\CUR@color{#1}\def\CUR@colorc@md{\c@msetgray}%
    \def\CUR@colorc@mdStroke{\c@msetgrayStroke}%
    \ifCUR@PS\PSc@mment{setgray Gray level=#1}\us@primarC@lor\fi\fi}
\ctr@ln@m\fillc@md
\ctr@ld@f\def\us@primarC@lor{\immediate\write\fwf@g{\d@fprimarC@lor}%
    \let\fillc@md=\prfillc@md}
\ctr@ld@f\def\prfillc@md{\d@fprimarC@lor\space\c@mfill}
\ctr@ld@f\def\c@lnbV@l#1 #2:{\def\t@xt@{#1}\relax\ifx\t@xt@\empty\c@lnbV@l#2:
    \else\c@lnbV@l@#1 #2:\fi}
\ctr@ld@f\def\c@lnbV@l@#1 #2:{\def\t@xt@{#2}\ifx\t@xt@\empty%
    \def\t@xt@{#1}\ifx\t@xt@\empty\advance\result@tent\m@ne\fi
    \else\advance\result@tent\@ne\c@lnbV@l@#2:\fi}
\ctr@ld@f\def\Blackcmyk{0 0 0 1}
\ctr@ld@f\def\Whitecmyk{0 0 0 0}
\ctr@ld@f\def\Cyancmyk{1 0 0 0}
\ctr@ld@f\def\Magentacmyk{0 1 0 0}
\ctr@ld@f\def\Yellowcmyk{0 0 1 0}
\ctr@ld@f\def\Redcmyk{0 1 1 0}
\ctr@ld@f\def\Greencmyk{1 0 1 0}
\ctr@ld@f\def\Bluecmyk{1 1 0 0}
\ctr@ld@f\def\Graycmyk{0 0 0 0.50}
\ctr@ld@f\def\BrickRedcmyk{0 0.89 0.94 0.28} 
\ctr@ld@f\def\Browncmyk{0 0.81 1 0.60} 
\ctr@ld@f\def\ForestGreencmyk{0.91 0 0.88 0.12} 
\ctr@ld@f\def\Goldenrodcmyk{ 0 0.10 0.84 0} 
\ctr@ld@f\def\Marooncmyk{0 0.87 0.68 0.32} 
\ctr@ld@f\def\Orangecmyk{0 0.61 0.87 0} 
\ctr@ld@f\def\Purplecmyk{0.45 0.86 0 0} 
\ctr@ld@f\def\RoyalBluecmyk{1. 0.50 0 0} 
\ctr@ld@f\def\Violetcmyk{0.79 0.88 0 0} 
\ctr@ld@f\def\Blackrgb{0 0 0}
\ctr@ld@f\def\Whitergb{1 1 1}
\ctr@ld@f\def\Redrgb{1 0 0}
\ctr@ld@f\def\Greenrgb{0 1 0}
\ctr@ld@f\def\Bluergb{0 0 1}
\ctr@ld@f\def\Cyanrgb{0 1 1}
\ctr@ld@f\def\Magentargb{1 0 1}
\ctr@ld@f\def\Yellowrgb{1 1 0}
\ctr@ld@f\def\Grayrgb{0.5 0.5 0.5}
\ctr@ld@f\def\Chocolatergb{0.824 0.412 0.118}
\ctr@ld@f\def\DarkGoldenrodrgb{0.722 0.525 0.043}
\ctr@ld@f\def\DarkOrangergb{1 0.549 0}
\ctr@ld@f\def\Firebrickrgb{0.698 0.133 0.133}
\ctr@ld@f\def\ForestGreenrgb{0.133 0.545 0.133}
\ctr@ld@f\def\Goldrgb{1 0.843 0}
\ctr@ld@f\def\HotPinkrgb{1 0.412 0.706}
\ctr@ld@f\def\Maroonrgb{0.690 0.188 0.376}
\ctr@ld@f\def\Pinkrgb{1 0.753 0.796}
\ctr@ld@f\def\RoyalBluergb{0.255 0.412 0.882}
\ctr@ld@f\def\Pssetg@n#1=#2|{\keln@mun#1|%
    \def\n@mref{c}\ifx\l@debut\n@mref\update@ttr\D@FTcolor\Pssetc@lor{#2}\else
    \def\n@mref{d}\ifx\l@debut\n@mref\update@ttr\D@FTdash\Q@s@tdash{#2}\else
    \def\n@mref{f}\ifx\l@debut\n@mref\update@ttr\D@FTfill\Q@s@tfillmode{#2}\else
    \def\n@mref{j}\ifx\l@debut\n@mref\update@ttr\D@FTjoin\Q@s@tjoin{#2}\else
    \def\n@mref{u}\ifx\l@debut\n@mref\update@ttr\D@FTupdate\Q@s@tupdate{#2}\else
    \def\n@mref{w}\ifx\l@debut\n@mref\update@ttr\D@FTwidth\Q@s@twidth{#2}\else
    \W@rnmesAttr{figset}{#1}\fi\fi\fi\fi\fi\fi}
\ctr@ln@m\CUR@dash
\ctr@ld@f\def\s@uvdash#1{\edef#1{\CUR@dash}}
\ctr@ld@f\def\D@FTdash{1}        
\ctr@ld@f\def\Q@s@tdash#1{\ifGR@cri\edef\CUR@dash{#1}\ifCUR@PS\expandafter\Pssetd@sh#1 :\fi\fi}
\ctr@ld@f\def\Pssetd@shI#1{\PSc@mment{setdash Index=#1}\ifcase#1%
    \or\immediate\write\fwf@g{[] 0 \c@msetdash}
    \or\immediate\write\fwf@g{[6 2] 0 \c@msetdash}
    \or\immediate\write\fwf@g{[4 2] 0 \c@msetdash}
    \or\immediate\write\fwf@g{[2 2] 0 \c@msetdash}
    \or\immediate\write\fwf@g{[1 2] 0 \c@msetdash}
    \or\immediate\write\fwf@g{[2 4] 0 \c@msetdash}
    \or\immediate\write\fwf@g{[3 5] 0 \c@msetdash}
    \or\immediate\write\fwf@g{[3 3] 0 \c@msetdash}
    \or\immediate\write\fwf@g{[3 5 1 5] 0 \c@msetdash}
    \or\immediate\write\fwf@g{[6 4 2 4] 0 \c@msetdash}
    \fi}
\ctr@ld@f\def\Pssetd@sh#1 #2:{{\def\t@xt@{#1}\ifx\t@xt@\empty\Pssetd@sh#2:
    \else\def\t@xt@{#2}\ifx\t@xt@\empty\Pssetd@shI{#1}\else\s@mme=\@ne\def\debutp@t{#1}%
    \an@lysd@sh#2:\ifodd\s@mme\edef\debutp@t{\debutp@t\space\finp@t}\def\finp@t{0}\fi%
    \PSc@mment{setdash Pattern=#1 #2}%
    \immediate\write\fwf@g{[\debutp@t] \finp@t\space\c@msetdash}\fi\fi}}
\ctr@ld@f\def\an@lysd@sh#1 #2:{\def\t@xt@{#2}\ifx\t@xt@\empty\def\finp@t{#1}\else%
    \edef\debutp@t{\debutp@t\space#1}\advance\s@mme\@ne\an@lysd@sh#2:\fi}
\ctr@ln@m\CUR@width
\ctr@ld@f\def\s@uvwidth#1{\edef#1{\CUR@width}}
\ctr@ld@f\def\D@FTwidth{0.4}     
\ctr@ld@f\def\Q@s@twidth#1{\ifGR@cri\edef\CUR@width{#1}\ifCUR@PS%
    \PSc@mment{setwidth Width=#1}\immediate\write\fwf@g{#1 \c@msetlinewidth}\fi\fi}
\ctr@ln@m\CUR@join
\ctr@ld@f\def\s@uvjoin#1{\edef#1{\CUR@join}}
\ctr@ld@f\def\D@FTjoin{miter}   
\ctr@ld@f\def\Q@s@tjoin#1{\ifGR@cri\edef\CUR@join{#1}\ifCUR@PS\expandafter\Pssetj@in#1:\fi\fi}
\ctr@ld@f\def\Pssetj@in#1#2:{\PSc@mment{setjoin join=#1}%
    \if#1r\def\t@xt@{1}\else\if#1b\def\t@xt@{2}\else\def\t@xt@{0}\fi\fi%
    \immediate\write\fwf@g{\t@xt@\space\c@msetlinejoin}}
\ctr@ld@f\def\Pss@tspecifSt#1{\trtlis@rg{#1}{\Pss@tspecifSt@}}
\ctr@ld@f\def\Pss@tspecifSt@#1=#2|{\keln@mun#1|%
    \def\n@mref{c}\ifx\l@debut\n@mref\def\n@mref{#2}\ifx\n@mref\D@FTref\else%
     \s@uvcolor{\typ@color}\Pssetc@lor{#2}\fi\else
    \def\n@mref{d}\ifx\l@debut\n@mref\def\n@mref{#2}\ifx\n@mref\D@FTref\else%
     \s@uvdash{\typ@dash}\Q@s@tdash{#2}\fi\else
    \def\n@mref{j}\ifx\l@debut\n@mref\def\n@mref{#2}\ifx\n@mref\D@FTref\else%
     \s@uvjoin{\typ@join}\Q@s@tjoin{#2}\fi\else
    \def\n@mref{w}\ifx\l@debut\n@mref\def\n@mref{#2}\ifx\n@mref\D@FTref\else%
     \s@uvwidth{\typ@width}\Q@s@twidth{#2}\fi\else
    \W@rnmeskwd{Pss@tspecifSt}{#1}\fi\fi\fi\fi}
\ctr@ld@f\def\Psrest@reSt#1{\trtlis@rg{#1}{\Psrest@reSt@}}
\ctr@ld@f\def\Psrest@reSt@#1=#2|{\keln@mun#1|%
    \def\n@mref{c}\ifx\l@debut\n@mref\def\n@mref{#2}\ifx\n@mref\D@FTref\else%
     \Pssetc@lor{\typ@color}\fi\else
    \def\n@mref{d}\ifx\l@debut\n@mref\def\n@mref{#2}\ifx\n@mref\D@FTref\else%
     \Q@s@tdash{\typ@dash}\fi\else
    \def\n@mref{j}\ifx\l@debut\n@mref\def\n@mref{#2}\ifx\n@mref\D@FTref\else%
     \Q@s@tjoin{\typ@join}\fi\else
    \def\n@mref{w}\ifx\l@debut\n@mref\def\n@mref{#2}\ifx\n@mref\D@FTref\else%
     \Q@s@twidth{\typ@width}\fi\else
    \W@rnmeskwd{Psrest@reSt}{#1}\fi\fi\fi\fi}
\ctr@ld@f\def\Pssettrim@sh#1=#2|{\keln@mde#1|%
    \def\n@mref{co}\ifx\l@debut\n@mref\update@ttr\D@FTref\P@settmeshcolor{#2}\else
    \def\n@mref{da}\ifx\l@debut\n@mref\update@ttr\D@FTref\P@settmeshdash{#2}\else
    \def\n@mref{wi}\ifx\l@debut\n@mref\update@ttr\D@FTref\P@settmeshwidth{#2}\else
    \W@rnmesAttr{figset trimesh}{#1}\fi\fi\fi}
\ctr@ln@m\DDV@tmeshcolor
\ctr@ld@f\def\P@settmeshcolor#1{\edef\DDV@tmeshcolor{#1}}
\ctr@ln@m\DDV@tmeshdash
\ctr@ld@f\def\P@settmeshdash#1{\edef\DDV@tmeshdash{#1}}
\ctr@ln@m\DDV@tmeshwidth
\ctr@ld@f\def\P@settmeshwidth#1{\edef\DDV@tmeshwidth{#1}}
\ctr@ld@f\def\figdrawtrimesh#1[#2,#3,#4]{{\ifCUR@PS\ifGR@cri%
    \PSc@mment{trimesh Type=#1, Triangle=[#2,#3,#4]}%
    \s@uvc@ntr@l\et@tpstrimesh\ifnum#1>\@ne%
    \Pss@tspecifSt{color=\DDV@tmeshcolor,dash=\DDV@tmeshdash,width=\DDV@tmeshwidth}%
    \setc@ntr@l{2}%
    \Pstrimeshp@rt#1[#2,#3,#4]\Pstrimeshp@rt#1[#3,#4,#2]\Pstrimeshp@rt#1[#4,#2,#3]%
    \Psrest@reSt{color=\DDV@tmeshcolor,dash=\DDV@tmeshdash,width=\DDV@tmeshwidth}%
    \fi\figdrawline[#2,#3,#4,#2]%
    \PSc@mment{End trimesh}\resetc@ntr@l\et@tpstrimesh\fi\fi}}
\ctr@ld@f\def\Pstrimeshp@rt#1[#2,#3,#4]{{\l@mbd@un=\@ne\l@mbd@de=#1\loop\ifnum\l@mbd@de>\@ne%
    \advance\l@mbd@de\m@ne\figptbary-1:[#2,#3;\l@mbd@de,\l@mbd@un]%
    \figptbary-2:[#2,#4;\l@mbd@de,\l@mbd@un]\figdrawline[-1,-2]%
    \advance\l@mbd@un\@ne\repeat}}
\initpr@lim\initpss@ttings\initPDF@rDVI
\ctr@ln@w{newbox}\figBoxA
\ctr@ln@w{newbox}\figBoxB
\ctr@ln@w{newbox}\figBoxC
\catcode`\@=12

\newtheorem{theorem}{Theorem}[section]

\newtheorem{remark}[theorem]{Remark}
\numberwithin{equation}{section}

\numberwithin{equation}{section}
\numberwithin{theorem}{section}

\title{Propagation of one and two-dimensional discrete waves under finite difference approximation}

\author{Umberto Biccari\textsuperscript{1,2}}  
\address{\textsuperscript{1}\,DeustoTech, University of Deusto, 48007 Bilbao, Basque Country, Spain.}
\address{\textsuperscript{2}\,Facultad de Ingenier\'ia, Universidad de Deusto, Avenida de las Universidades 24, 48007 Bilbao, Basque Country, Spain, +34 944139003 - 3282.}
\thanks{This project has received funding from the European Research Council (ERC) under the European Union's Horizon 2020 research and innovation programme (grant agreement No. 694126-DyCon). The work of the first and of the third author was partially supported by the Grants MTM2014-52347 and MTM2017-92996 of MINECO (Spain) and by the Grant FA9550-18-1-0242 of AFOSR. The work of the second author was partially supported by CNCS-UEFISCDI Grant No. PN-III-P4-ID-PCE-2016-0035. The work of the third author was partially supported by the Grant ICON of the French ANR.} 

\author{Aurora Marica\textsuperscript{3}\,}
\address{\textsuperscript{3}\,Politehnica University of Bucharest, Faculty of Applied Sciences, Department of Mathematics and Computer Sciences, Splaiul Independentei 313, Bucharest, 060042, Romania.}

\author{Enrique Zuazua\textsuperscript{1,2,4,5}}
\address{\textsuperscript{4}\, Departamento de Matem\'aticas, Universidad Aut\'onoma de Madrid, 28049 Madrid, Spain.} 
\address{\textsuperscript{5}\, Sorbonne Universit\'es, UPMC Univ Paris 06, CNRS UMR 7598, Laboratoire Jacques-Louis Lions, F-75005, Paris, France.}
\email{umberto.biccari@deusto.es, mauramihaela@yahoo.com, enrique.zuazua@deusto.es}


\begin{document}
	
\bibliographystyle{acm}
 
\begin{abstract}
	We analyze the propagation properties of the numerical versions of one and two-dimensional wave equations, semi-discretized in space by finite difference schemes. We focus on high-frequency solutions whose propagation can be described, both at the continuous and semi-discrete level, by micro-local tools. We do it both for uniform and non-uniform numerical grids and also for constant coefficients and variable ones. The energy of continuous and semi-discrete high-frequency solutions propagates along bi-characteristic rays, but their dynamics differ from the continuous to the semi-discrete setting, because of the different nature of the corresponding Hamiltonians. One of the main objectives of this paper is to illustrate through accurate numerical simulations that, in agreement with the micro-local theory, numerical high-frequency solutions can bend in an unexpected manner, as a result of the accumulation of the local effects introduced by the heterogeneity of the numerical grid. These effects are enhanced in the multi-dimensional case where the interaction and combination of such behaviors in the various space directions may produce, for instance, the \textit{rodeo effect},  i. e. waves that are trapped by the numerical grid in closed loops, without ever getting to the exterior boundary. Our analysis allows explaining all such pathological behaviors. Moreover, the discussion in this paper also contributes to the existing theory about the necessity of filtering high-frequency numerical components when dealing with control and inversion problems for waves, which is based very much in the theory of rays and, in particular, on the fact that they  can be observed when reaching the exterior boundary of the domain, a key property that can be lost through numerical discretization.
\end{abstract}
    		
\maketitle 

\section{Introduction}

The analysis of propagation properties of numerical waves obtained through a finite difference discretization on uniform or non-uniform meshes is a topic which has been extensively investigated in the literature. Among other contributions, we mention the works  \cite{trefethen1982wave,trefethen1982group,vichnevetsky1980propagation,vichnevetsky1981energy,vichnevetsky1981propagation,vichnevetsky1987wave,vichnevetsky1982fourier}. 

In this paper, we are interested in discussing several aspects of wave propagation in a computational framework, and in the comparison with the usual behavior of the continuous models. In more detail, we consider here wave-like equations and their numerical approximation by a finite difference scheme, with the intent to illustrate theoretically and computationally the dynamics that this discretization introduces, and to comment the main differences with respect to the expected comportment of the original continuous PDE. 

Our analysis will address both constant and variable coefficients models, in one and two space dimensions. Moreover, our approach is based on the study of the propagation of high-frequency Gaussian beam solutions (that is, solutions originated from highly concentrated and oscillating initial data), both in continuous and discrete media. 

In a continuous setting, these kind of techniques date back to the works of H\"ormander (\cite{hormander1985analysis}), and they have been later extended by several authors, with the developments of tools like micro-local defect measures (introduced independently by G\'erard in \cite{gerard1991microlocal} and by Tartar in \cite{tartar1990h}, in the context of nonlinear partial differential equations and of homogenization, respectively) or Wigner measures (\cite{lions1993mesures,markowich1994wigner,wigner1932quantum})

More recently, this same question has been addressed also at the numerical level, and we can nowadays mention contributions on the extension of micro-local techniques to the study of the propagation properties for discrete waves (\cite{macia2002propagacion,marica2015propagation}).

Roughly speaking, the idea at the basis of this techniques is that the energy of Gaussian beam solutions propagates along bi-characteristic rays, which are obtained from the Hamiltonian system associated to the symbol of the operator under consideration. At the continuous level, if the coefficients of the equation are constants, these mentioned rays are straight lines and travel with a uniform velocity. In the case of variable coefficients, instead, the heterogeneity of the medium where waves propagate produces the bending of the rays and, consequently, the increasing or the decreasing of their velocity. 

On the other hand, the finite difference space semi-discretization of the equation may introduce different dynamics, with a series of unexpected propagation properties at high frequencies, that substantially differ from the expected behavior of the continuous equation. For instance, one can generate spurious solutions traveling at arbitrarily small velocities (\cite{trefethen1982group}) which, therefore, show lack of propagation in space. As we shall see, this phenomenon is related with the particular nature of the discrete group velocity which, differently from the continuous equation, may vanish at certain frequencies. In addition, the introduction of a non-uniform mesh for the discretization of the equation makes the situation even more intricate. For instance, as indicated in \cite{marica2015propagation,vichnevetsky1987wave}, for some numerical grids the rays of geometric optics may present internal reflections, meaning that the waves change direction without hitting the boundary. 

All these pathologies are purely numerical, and they are related to changes in the Hamiltonian system giving the equations of the rays. Actually, as we will discuss with more details later, they can be rigorously explained by the behavior of the corresponding phase portrait.

With this considerations in mind, our main goal in the present article is to exhibit the kind of spurious effects one may encounter when performing numerical approximations by means of finite difference schemes. Our principal contribution is the detection of these singular phenomena in the numerical propagation of the waves. 

More than on the presentation of the mathematical details, in this article we focus on a self contained review of the propagation properties of the numerical solutions of discrete wave equations, with the support of several numerical simulations showing the aforementioned pathologies. For a more complete and exhaustive discussion on this topic, the interested reader may refer also to \cite{macia2002propagacion,macia2002lack,marica2015propagation}. We start by considering the following one-dimensional system 
\begin{align}\label{wave_1d}
	\begin{cases}
		\partial^2_tu-\partial^2_xu = 0, & (x,t)\in (-1,1)\times(0,T)
		\\
		u(-1,t) = u(1,t) = 0, & t\in(0,T)
		\\
		u(x,0) = u^0(x),\;\;\;\partial_tu(x,0) = u^1(x), & x\in(-1,1),
	\end{cases}
\end{align}

A classical property of the above equation is that its total energy, is conserved in time. In other words, if we define the energy associated to \eqref{wave_1d} as  the quantity
\begin{align}\label{energy_1d_const}
	\mathcal E(u,\partial_t u) :=\frac 12 \int_{-1}^1 \left(|\partial_t u|^2 + |\partial_x u|^2\right)\,dx,
\end{align}
we then have 
\begin{align*}
	\mathcal E(u,\partial_t u) = \mathcal E(u^0,u^1),\;\;\;\forall\,t\in(0,T).
\end{align*}
Moreover, when $x\in\RR$, according to the d'Alembert formula 
\begin{align}\label{waves_dal}
	u(x,t) = \frac 12 \Big(u^0(x+t) + u^0(x-t)\Big) + \frac 12 \int_{x-t}^{x+t} u^1(z)\,dz
\end{align}
the solution can be uniquely decomposed into two waves, each one of them propagating along one of the straight characteristics $x\pm t$. In addition, it is well-known (see \cite{bardos1992sharp,burq1997condition,burq1998controle,macia2002propagacion,macia2002lack,ralston1982gaussian,rauch2005polynomial}) that the energy of initial data presenting high-frequency oscillation and/or concentration effects propagates along these characteristic rays. 

As we mentioned, these concentration features can be observed also at the numerical level but, in this framework, the high-frequency solutions exhibit several pathological behaviors. In order to detect and describe these propagation properties, we are going to discretize \eqref{wave_1d} through a finite difference scheme for the space variable and a classical leapfrog method for the time integration. 

Moreover, the analysis that we are going to present can be easily extended also to the case of variable coefficients and to  multi-dimensional problems. In fact, in the second part of this paper we will briefly discuss the following one-dimensional variable coefficients wave equation
\begin{align}\label{wave_var_1d_intro}
	\begin{cases}
		\rho(x)\partial^2_tu-\partial_x\big(\sigma(x)\partial_xu\big) = 0, & (x,t)\in (-1,1)\times(0,T)
		\\
		u(-1,t) = u(1,t) = 0, & t\in(0,T)
		\\
		u(x,0) = u^0(x),\;\;\;\partial_tu(x,0) = u^1(x), & x\in(-1,1),
	\end{cases}
\end{align}
where $\rho$ and $\sigma$ are chosen to be $L^\infty(\RR)$-functions with the strict hyperbolicity assumption $\rho(x)\geq\rho^*>0$ and $\sigma(x)\geq\sigma^*>0$. 

Even if for \eqref{wave_var_1d_intro} there is no explicit formula for the solutions analogous to \eqref{waves_dal}, it still holds that  highly concentrated and oscillatory initial data lead to Gaussian wave packet-type solutions concentrated along one of the bi-characteristics (which are not straight lines anymore), and their energy localized outside any neighborhood of the ray vanishes as the wavelength parameter tends to zero.

Once again, these concentration features can be observed also at the numerical level. Moreover, also in this case the space semi-discretization of \eqref{wave_var_1d_intro} may introduce dynamics different than the ones that we would expect from the analysis of the continuous equation. We will present some of these phenomena later in the manuscript. 

Lastly, the final part of the present work will be devoted to a brief discussion of the following two-dimensional system 
\begin{align}\label{wave_var_2d_intro}
	\begin{cases}
		\rho(\Bf{z})\partial^2_tu-\textrm{div}_{\Bf{z}}\big(\sigma(\Bf{z})\nabla_{\Bf{z}} u\big) = 0, & (\Bf{z},t)\in \Omega\times(0,+\infty)
		\\
		\left. u\right|_{\partial\Omega} = 0, & t\in(0,+\infty)
		\\
		u(\Bf{z},0) = u^0(\Bf{z}),\;\;\;\partial_tu(\Bf{z},0) = u^1(\Bf{z}), & \Bf{z}\in\Omega,
	\end{cases}
\end{align}
where we indicate $\Omega:=(-1,1)^2$ and in which the same assumptions on the coefficients $\rho$ and $\sigma$ will be considered. 

For the sake of completeness, let us conclude this introduction by mentioning that our study is motivated by control and inverse problems. Indeed, it is well-known that boundary controllability and identifiability properties of solutions of wave equations hold because of the fact that the energy is driven by characteristics that reach a subregion of the domain or of its boundary where the controllers or observers are placed. This, in particular, allows the \textit{observability} of the solutions (namely, the possibility to obtain estimates of the total energy in terms of the energy concentrated on the support of the control along time), which is classically known to be equivalent to control properties. 

In the framework of wave-like processes, observability is guaranteed by the so-called geometric control condition (GCC), requiring all rays of geometric optics to enter the control region during the control time. This condition has been proved to be sufficient and almost necessary in \cite{bardos1992sharp} (see also \cite{burq1997condition,burq1998controle}). In particular, the necessity of GCC is related to the fact that, around each ray that does not meet the  observation/control region, one can always build concentrated Gaussian beams making the observability inequality impossible.

When the wave equation is approximated by finite difference methods on uniform meshes, the control for the discretized model does not necessarily yields a good approximation of the control for the original continuous problem. Instead, observability/controllability may be lost under numerical discretization as the mesh size tends to zero, due to the existence of high-frequency spurious solutions for which the group velocity vanishes. In particular, it is by now well-known that the discretization of the equation generates wave packets traveling with a group velocity which, at high frequencies, is of the order of $1/h$, $h$ being the space mesh size. These high-frequency solutions are such that the energy concentrated in the control region is asymptotically smaller than the total energy, and they produce the exponential blow-up of the observability constant as $h\to 0$. A deeper discussion on this topic is out of the scope of the present paper and, therefore, it will not be considered here. The interested reader may refer to the well extended bibliography on the control properties for discrete waves \cite{castro2006boundary,ervedoza2015numerical,ervedoza2012wave,ervedoza2013numerical,glowinski1992ensuring,glowinski1990numerical,glowinski1990approach,ignat2009convergence,infante1999boundary,loreti2008ingham,macia2002propagacion,macia2002lack,micu2002uniform,negreanu2004convergence,zuazua1999boundary,zuazua2005propagation}).

This paper is organized as follows. In Section \ref{1d_sec}, we discuss the one-dimensional case. In more detail, in Section \ref{discretization_sec} we are going to introduce the finite difference discretization of the one-dimensional model \eqref{wave_1d}. In Section \ref{ham_sec}, we will present and discuss the Hamiltonian system, giving the equations for the bi-characteristic rays, while Section \ref{numerical_sec} will be devoted to show and comment the results of our numerical simulations. We will complete this first part of the paper with Section \ref{var_coeff_sec}, which is devoted to the case of a variable-coefficients wave equation. Finally, in Section \ref{2d_sec}, we will briefly consider the extension of our analysis to two-dimensional problems.

\section{One-dimensional wave equation}\label{1d_sec}

We discuss here the propagation properties of the finite difference approximation of the solution to the one-dimensional wave equation \eqref{wave_1d}. We start by introducing the numerical scheme that we will employ.
 
\subsection{Semi-discrete approximation of \eqref{wave_1d}}\label{discretization_sec}

Let us introduce the grid and the space semi-discrete approximation of \eqref{wave_1d} that we will employ in our analysis. 

Let $g:[-1,1]\to [-1,1]$ be a diffeomorphism. For $N\in\NN^*$, we set $h=2/(N+1)$ as the step size of the uniform mesh 
\begin{align*}
	\mathcal G^h:= \Big\{x_j := -1+jh,\,j=0,\ldots, N+1\,\Big\},
\end{align*}
and we consider the non-uniform grid 
\begin{align*}
	\Ggh:=\Big\{g_j:= g(x_j),\,x_j\in\mathcal G^h\,\Big\}
\end{align*} 
obtained by transforming $\mathcal G^h$ through the map $g$. 
\begin{figure}[h]
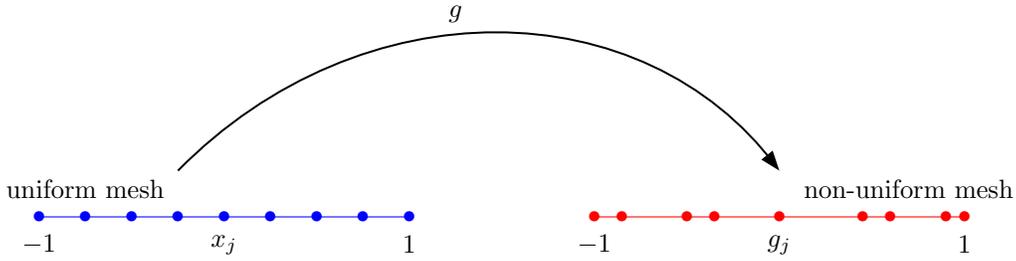

	\figinit{0.7pt}
	
	\figpt 1:(-250,0) \figpt 2:(-50,0)
	\figpt 3:(50,0)  \figpt 4:(250,0)
	
	\figpt 5:(-225,0) \figpt 6:(-200,0)
	\figpt 7:(-175,0) \figpt 8:(-150,0)
	\figpt 9:(-125,0) \figpt 10:(-100,0)
	\figpt 11:(-75,0) 	
	
	\figpt 12:(65,0) \figpt 13:(100,0)
	\figpt 14:(115,0) \figpt 15:(150,0)
	\figpt 16:(195,0) \figpt 17:(210,0)
	\figpt 18:(240,0)
	
	\figpt 19:(-175,25) \figpt 20:(-75,125)
	\figpt 21:(75,125)  \figpt 22:(150,25)
	
	\figpt 23:(-150,-15) \figpt 24:(150,-15)
	\figpt 25:(-25,110) \figpt 26:(-225,15)
	\figpt 27:(220,15) \figpt 28:(-250,-15)
	\figpt 29:(-50,-15) \figpt 30:(50,-15)
	\figpt 31:(250,-15)
	 
	
	\figdrawbegin{}
	\figset (color=\Bluergb)
	\figdrawline[1,2]
	\figset (color=\Redrgb)
	\figdrawline[3,4]
	\figset (color=default)
	\figset arrowhead(fillmode=switch)
	\figset (width=0.75)
	\figdrawarrowBezier[19,20,21,22]
	
	\figdrawend
	
	\figvisu{\figBoxA}{}{
		\figwritec[1,2,5,6,7,8,9,10,11]{$\color{blue}\bullet$}	
		\figwritec[3,4,12,13,14,15,16,17,18]{$\color{red}\bullet$}	
		\figwritec [23]{$x_j$}
		\figwritec [24]{$g_j$}
		\figwritec [25]{$g$}
		\figwritec [26]{uniform mesh}
		\figwritec [27]{non-uniform mesh}
		\figwritec [28,30]{$-1$}
		\figwritec [29,31]{$1$}
	}
	\centerline{\box\figBoxA}
	\caption{The diffeomorphism $g$ transform the nodes $x_j$ of the uniform mesh $\mathcal G^h$ into the nodes $g_j$ of the non-uniform one $\mathcal G^h_g$.}\label{mesh_fig}
\end{figure}

Along all this paper, we focus on the case in which $g$ is regular. More precisely, we assume that $g\in C^2(\RR)$ with $0<g_d^-\leq|g'(x)|\leq g_d^+<+\infty$ and $|g''(x)|\leq g_{dd}<+\infty$ for some given constant $g_d^-,g_d^+,g_{dd}>0$ and for all $x\in\RR$. We stress that this also includes the case of a uniform grid, that is $g(x) = x$. 

Notice that the above assumptions yield that $1/g'$ belongs to the space $C^{0,1}(\RR)$ of the Lipschitz continuous functions. As we will see in the next section, this Lipschitz regularity assumption will be important when introducing the Hamiltonian system for the bi-characteristic rays. 

We set ${g_{j+1/2}:= g(x_{j+1/2})}$ to be the image through the map $g$ of the midpoints ${x_{j+1/2}:=-1+(j+1/2)h}$, and we denote by 
\begin{align*}
	&h_{j+1/2} :=g_{j+1}-g_j,\;\;\; j=0,\ldots,N
	\\
	&h_{j-1/2} :=g_j-g_{j-1},\;\;\;j=1,\ldots,N+1
	\\
	&h_j :=\frac{h_{j+1/2}-h_{j-1/2}}{2},\;\;\;j=1,\ldots,N
\end{align*}
the heterogeneous mesh sizes. The semi-discretization of the wave equation \eqref{wave_var_1d} on the non-uniform grid $\Ggh$ is then given as follows:
\begin{align}\label{wave_discr_1d}
	\begin{cases}
		\displaystyle h_ju_j''(t) - \left(\frac{u_{j+1}(t)-u_j(t)}{h_{j+1/2}}-\frac{u_j(t)-u_{j-1}(t)}{h_{j-1/2}}\right) = 0, & j=1,\ldots,N,\;\;t\in(0,T)
		\\
		u_0(t) = u_{N+1}(t) = 0, & t\in(0,T)
		\\
		u_j(0) = u_j^0, \;\;\; u_j'(0) = u_j^1, & j=1,\ldots,N.
	\end{cases}
\end{align}

Here we indicate $u_j(t):=u(g_j,t)$. Moreover, to avoid possible confusions, when needed, we will denote the solutions of \eqref{wave_var_discr_1d} by $\textbf{u}^{h,g}(t):=(u^g_j(t))_{j=0}^{N+1}$ for emphasizing the dependence from the diffeomorphic transformation $g$. But, in general, we will simply write down $\mathbf{u}^h(t):=(u_j(t))_{j=0}^{N+1}$, $\mathbf{u}^{0,h}:=(u^0_j)_{j=0}^{N+1}$ and $\mathbf{u}^{1,h}:=(u^1_j)_{j=0}^{N+1}$ for the solutions and for the initial data.

It has been shown in \cite[Proposition 3.4]{marica2015propagation} that, under our regularity assumptions on the mesh function $g$, the numerical approximation \eqref{wave_discr_1d} with initial data $(\mathbf{u}^{0,h},\mathbf{u}^{1,h})$ is a convergent scheme of order $O(h)$ for the wave equation \eqref{wave_1d} in the appropriate $\ell^2$ setting. Moreover, it can be readily checked that, as it happens for the continuous model, system \eqref{wave_discr_1d} enjoys the property of energy conservation, that is,
\begin{align*}
	\mathcal E^{h,g}(\mathbf{u}^{h,g}(t),\partial_t \mathbf{u}^{h,g}(t)):= \frac 12 \sum_{j=1}^N h_j|\partial_t u_j^g(t)|^2 + \frac 12 \sum_{j=1}^N h_{j+1/2}\left|\frac{u_{j+1}^g(t)-u_j^g(t)}{h_{j+1/2}}\right|^2  = \mathcal E^{h,g}(\mathbf{u}^{0,h},\mathbf{u}^{1,h}).
\end{align*}

\subsection{The Hamiltonian system}\label{ham_sec}

As we mentioned in the introduction, it is classically known (see \cite{bardos1992sharp}) that at the basis of the analysis of the observability properties of the continuous wave equation in any space dimension there are the propagation properties of bi-characteristic rays, leading to the so-called Geometric Control Condition (GCC). This same analysis can be repeated also at the numerical level, and  the observability of the discrete wave equation \eqref{wave_discr_1d} can be discussed through the study of the propagation of discrete bi-characteristic rays. 

Rays of geometric optics are defined as the projections on the physical space $(x,t)$ of the bi-characteristic rays given by the Hamiltonian corresponding to the principal part of the operator. In the case of the one-dimensional wave equation \eqref{wave_1d}, this Hamiltonian is given by 
\begin{align}\label{ham_cont_1d}
	\mathcal H_c(x,t,\xi,\tau) = -\tau^2 + \xi^2,
\end{align}
and the bi-characteristic rays are the curves $s\mapsto (x(s),t(s),\xi(s),\tau(s))$ solving the first order ODE system
\begin{align}\label{ham_syst_cont_1d}
	\begin{cases}
		\dot{x}(s) = \partial_\xi \mathcal H_c(x(s),t(s),\xi(s),\tau(s)) = 2\xi(s)
		\\
		\dot{t}(s) = \partial_\tau \mathcal H_c(x(s),t(s),\xi(s),\tau(s)) = -2\tau(s)
		\\
		\dot{\xi}(s) = -\partial_x \mathcal H_c(x(s),t(s),\xi(s),\tau(s)) = 0
		\\
		\dot{\tau}(s) = -\partial_t \mathcal H_c(x(s),t(s),\xi(s),\tau(s)) = 0.
	\end{cases}
\end{align}

In \eqref{ham_syst_cont_1d}, we indicate with $\dot{x}$ the time derivative of the variable $x$, for differentiating it from the notation $'$ that we are using for general derivatives. Moreover, to system \eqref{ham_syst_cont_1d} we associate an initial datum $(x(0),t(0),\xi(0),\tau(0)) = (x_0,0,\xi_0,\tau_0)$ such that $\mathcal H_c(x_0,0,\xi_0,\tau_0)=0$. 

It is immediately seen that the above system can be explicitly solved, and we thus obtain that for any $\xi_0$ there are two characteristics starting from the point $x_0$ which are given by the straight lines $x^\pm(t)=x_0\mp t$. Moreover, each one of these characteristics reaches the boundary of the interval $(-1,1)$ in a uniform time not depending on the frequency $\xi_0$. Then, when hitting the boundary, the ray reflects according to the Descartes-Snell's law. 

At the discrete level, the situation changes and it is more delicate. As a first thing, one needs to find a suitable counterpart of the Hamiltonian \eqref{ham_cont_1d}. This was done in \cite{macia2002propagacion} for the discrete wave equation in an infinite lattice
(i.e. with no boundary) and, more recently, the analysis was extended in \cite{marica2015propagation} to non-uniform meshes obtained by means of a $C^2$ diffeomorphism.

To be more precise, in \cite{marica2015propagation} the authors constructed solutions of \eqref{wave_discr_1d} with frequencies of the order of $1/h$, localized around the rays of geometric optics given by the projection on the physical space of the bi-characteristic curves provided by the Hamiltonian 
\begin{align}\label{ham_discr_1d}
	\mathcal H(y,t,\xi,\tau) = -g'(y)\tau^2 + 4\sin^2\left(\frac \xi2\right)\frac{1}{g'(y)},\;\;\; y=g^{-1}(x),
\end{align}
or equivalently
\begin{align}\label{ham_discr_1d_2}
	\mathcal H(y,t,\xi,\tau) = -\tau^2 + \frac{1}{g'(y)^2}\omega(\xi)^2, 
\end{align}
with
\begin{align*}
	\omega(\xi):=2\sin\left(\frac \xi2\right).
\end{align*}

We observe some changes in this discrete symbol with respect to the continuous one. Firstly, it depends on the space variable  $y=g^{-1}(x)$ corresponding to the uniform grid. Note also the appearance of the factor $1/g'(x)$ accompanying each space derivative, which is also due to the grid transformation.  Moreover, the Fourier symbol $\xi^2$ of the second-order space derivative has been replaced by the corresponding symbol $4\sin^2(\xi/2)$ of the three-point finite difference approximation of the Laplacian. Lastly, observe that in the symbol the parameter $h$ disappears. This is due to the fact that we are analyzing the propagation of waves of wavelength of the order of $h$.

Set $c_g:=1/g'$. It is easy to see that from the discrete Hamiltonian \eqref{ham_discr_1d} the ODE system for the rays that one obtains is the following
\begin{align}\label{ham_syst_discr_1d}
	\begin{cases}
		\dot{y}(s) = 2c_g(y(s))^2\omega(\xi(s))\partial_\xi \omega(\xi(s))
		\\
		\dot{t}(s) = -2\tau(s)
		\\
		\dot{\xi}(s) = -2c_g(y(s))\partial_yc_g(y(s))\omega(\xi(s))^2
		\\
		\dot{\tau}(s) = 0,
	\end{cases}
\end{align}
with an initial condition at $s=0$ given by $(y_0,0,\xi_0,\tau_0)$ satisfying $\mathcal H(y_0,0,\xi_0,\tau_0) = 0$. 

In \eqref{ham_syst_discr_1d}, $\partial_\xi \omega(\xi)$ is the \textit{group velocity}, i.e. the speed at which the energy associated with wave number $\xi$ moves. As it is natural to expect, this quantity has a fundamental role in the analysis of the propagation of numerical solutions.  

To better understand the dynamics of the bi-characteristic rays remark that, for all $s$, $\tau(s) = \tau_0$. Besides, $\mathcal H(y(s),t(s),\xi(s),\tau(s)) = 0$, so that for all $s$ we encounter two different possible values for $\tau_0$, namely
\begin{align*}
	\tau_0^\pm = \pm c_g(y(s))|\omega(\xi(s))|.
\end{align*}

Now, as $dt/ds$ does not vanish, the Inverse Function Theorem allows to parametrize the curve $s\mapsto (y(s),t(s),\xi(s),\tau_0^\pm)$ by $t\mapsto (y(t),t,\xi(t),\tau_0^\pm)$, hence obtaining
\begin{align}\label{ham_syst_discr_1d_tau}
	\begin{cases}
		\displaystyle\dot{y}^\pm(t) = -\frac{1}{\tau_0^\pm} c_g(y^\pm(t))^2\omega(\xi^\pm(t))\partial_\xi \omega(\xi^\pm(t))
		\\[10pt]
		\displaystyle\dot{\xi}^\pm(t) = \frac{1}{\tau_0^\pm}c_g(y^\pm(t))\partial_yc_g(y^\pm(t))\omega(\xi^\pm(t))^2
		\\[10pt]
		y^\pm(0) = y_0,\;\;\;\xi^\pm(0)=\xi_0,
	\end{cases}
\end{align}
or, equivalently,
\begin{align}\label{ham_syst_discr_1d_t}
	\begin{cases}
		\displaystyle\dot{y}^\pm(t) = \mp c_g(y^\pm(t))\partial_\xi \omega(\xi^\pm(t))
		\\
		\displaystyle\dot{\xi}^\pm(t) = \pm\partial_yc_g(y^\pm(t))\omega(\xi^\pm(t))
		\\
		y^\pm(0) = y_0,\;\;\;\xi^\pm(0)=\xi_0
	\end{cases}
\end{align}
Since by assumption $c_g(\cdot)>0$, thanks to \eqref{ham_syst_discr_1d_t} we immediately see that
\begin{align*}
	|\dot{y}^\pm(t)| = c_g(y^\pm(t))\left|\partial_\xi\omega(\xi^\pm(t))\right|.
\end{align*}

In view of that, the velocity of the rays vanishes if, and only if, $\partial_\xi(\omega)=0$ for some $\xi$. When $\omega(\xi)=\xi$, corresponding to the continuous case, this cannot happen. Actually, in this case the rays travel with velocity one until they hit the boundary, where they are reflected. On the other hand, when $\omega(\xi) = 2\sin(\xi/2)$, corresponding to the finite difference discretization that we are considering, we immediately see that $\partial_\xi \omega(\xi) = \cos(\xi/2)$ vanishes for $\xi = (2k+1)\pi$, $k\in\ZZ$. As our simulations will highlight, in this case we have the phenomenon of non-propagating waves. Notice also that the possibility that the velocity of the discrete rays is zero is independent of the choice of the mesh and of the variable coefficients. In other words, no matter whether the coefficients are constant or variables, and no matter what mesh we select for solving our problem, there will always be certain frequencies for which the group velocity of the numerical waves vanishes.  

\subsection{Numerical results}\label{numerical_sec}

We present several numerical simulations showing the propagation of the solutions to \eqref{wave_1d} with highly concentrated and oscillating initial data.

\begin{figure}[!h]
	\centering 
	\begin{minipage}{0.4\textwidth}
		\includegraphics[scale=1]{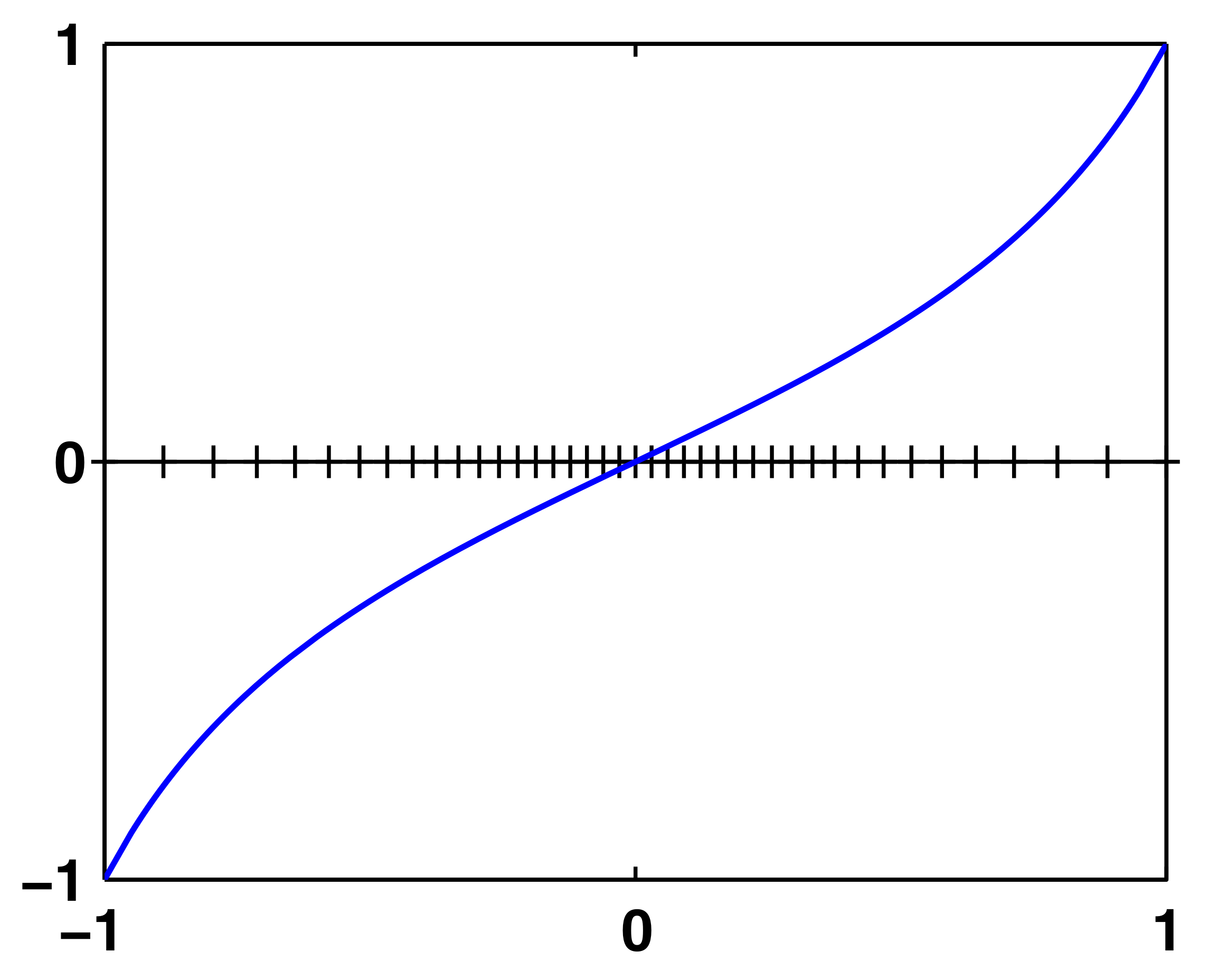}
		\caption*{Grid $\textbf{g}^{\,h,1}$}
	\end{minipage}
	\hspace{0.5cm}
	\begin{minipage}{0.4\textwidth}
		\vspace{-0.1cm}
		\includegraphics[scale=1]{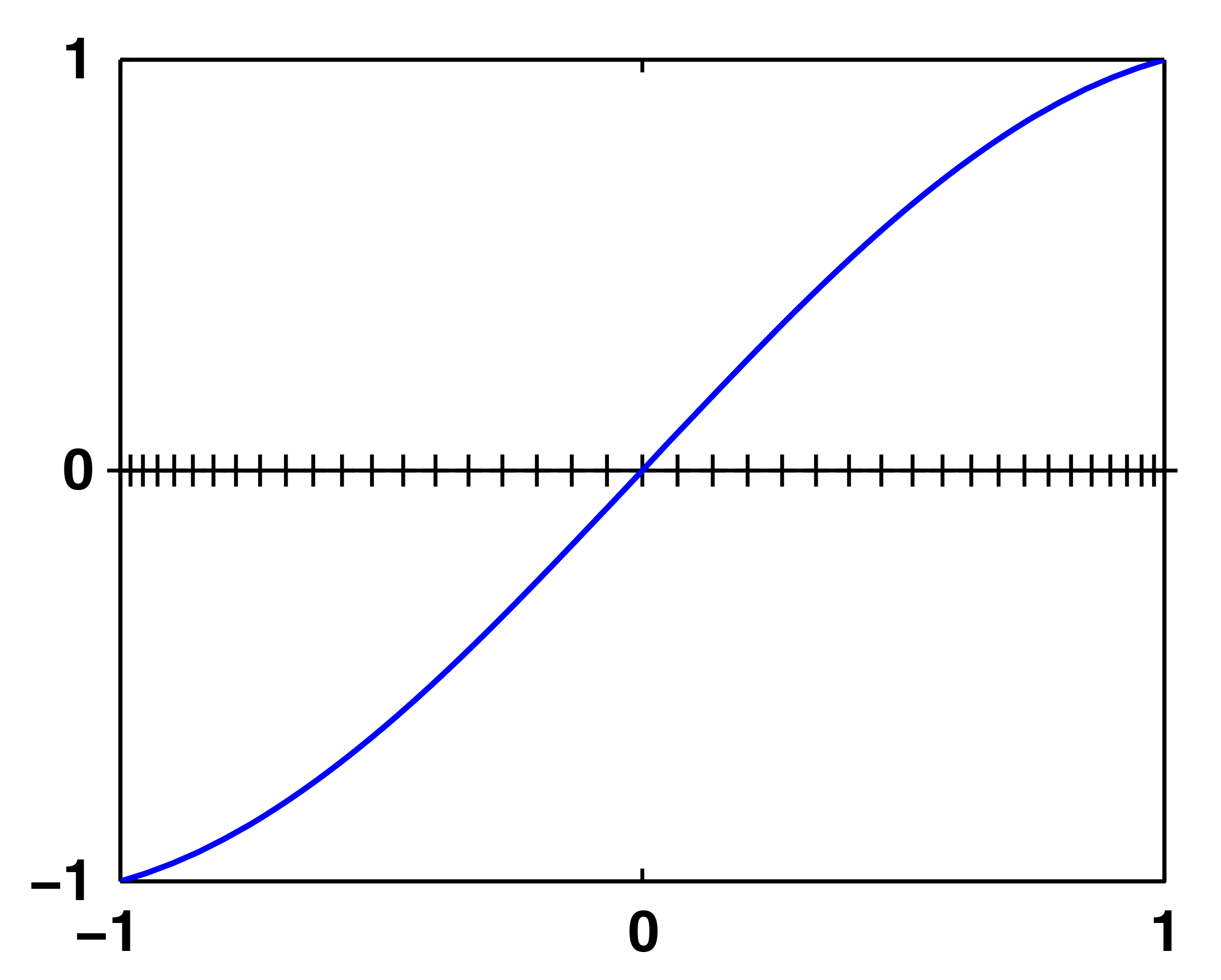}
		\caption*{Grid $\textbf{g}^{\,h,2}$}
	\end{minipage}
	\caption{Representation of the two non-uniform grids $\textbf{g}^{\,h,1}$ and $\textbf{g}^{\,h,2}$. By cross markers we indicate the grid points and by solid line the application $g_i$ generating them.}\label{meshes_fig}
\end{figure}

Set $\textbf{x}^h$ to be the uniform mesh of size $h=2/(N+1)$. We consider two non-uniform grids produced by the transformations 
\begin{align}\label{mesh_fun_1d}
	g_1:= \tan\left(\frac \pi4 x\right)\;\;\;\textrm{ and} \;\;\; g_2 (x):= 2\sin\left(\frac \pi6 x\right),
\end{align}
yielding a gradual refinement at the center $x = 0$ of the space interval and at the two endpoints $x=\pm 1$, respectively (see Figure \ref{meshes_fig}). In what follows, we will indicate
\begin{align*}
	\textbf{g}^{\,h,1}:= g_1(\textbf{x}^h)=\tan\left(\frac \pi4 \textbf{x}^h\right)\;\;\;\textrm{ and} \;\;\; \textbf{g}^{\,h,2}:=g_2(\textbf{x}^h)= 2\sin\left(\frac \pi6 \textbf{x}^h\right).
\end{align*}
	
The time discretization is done by means of a leap-frog scheme $(u^{n+1}-2u^n+u^{n-1})/(\delta t)^2$ with Courant-Fiedrichs-Lewy (CFL) condition $\delta t = 0.1\cdot h$ (which is necessary since the method is explicit). Moreover, the initial data are constructed starting from the following Gaussian profile
\begin{align}\label{in_data_wave}
	G_\gamma(x) = e^{-\frac \gamma2 \big(g^{-1}(x)-g^{-1}(x_0)\big)^2}e^{\,i\frac{\xi_0}{h} g^{-1}(x)},\;\;\;\gamma:=h^{-0.9}.
\end{align}

In more detail, we will take $\textbf{u}^{0,h}=G_\gamma(\textbf{g}^{\,h,i})$, $i=1,2$, and $\textbf{u}^{1,h}=G_\gamma'(\textbf{g}^{\,h,i})$. Moreover, in what follows we will always assume $\xi\in[0,2\pi]$.

In our further discussion, it will be fundamental a deeper analysis of the Hamiltonian system associated to \eqref{wave_discr_1d} which, we recall, is given by (see also \eqref{ham_syst_discr_1d_t})
\begin{align*}
	\begin{cases}
		\displaystyle\dot{y}^\pm(t) = \mp c_g(y^\pm(t))\cos\left(\frac{\xi^\pm(t)}{2}\right)
		\\[8pt]
		\displaystyle\dot{\xi}^\pm(t) = \pm 2\partial_yc_g(y^\pm(t))\sin\left(\frac{\xi^\pm(t)}{2}\right)
		\\[8pt]
		y^\pm(0) = y_0,\;\;\;\xi^\pm(0)=\xi_0.
	\end{cases}
\end{align*}

Recall also that in the above system the variable $y=g^{-1}(x)$ is the one corresponding to the refined mesh. In the original variable $x=g(y)$ associated to the uniform partition of the space interval $[-1,1]$, \eqref{ham_syst_discr_1d_t} reads instead as 
\begin{align*}
	\begin{cases}
		\displaystyle\dot{x}^\pm(t) = \mp a_g(x^\pm(t))\cos\left(\frac{\xi^\pm(t)}{2}\right)
		\\[8pt]
		\displaystyle\dot{\xi}^\pm(t) = \pm 2b_g(x^\pm(t))\sin\left(\frac{\xi^\pm(t)}{2}\right)
		\\[8pt]
		x^\pm(0) = x_0,\;\;\;\xi^\pm(0)=\xi_0.
	\end{cases}
\end{align*}
with $a_g(\cdot):=(g'c_g)(g^{-1}(\cdot))$, $b_g(\cdot):=c_g'(g^{-1}(\cdot))$ and, clearly, $x_0=g(y_0)$.

Notice that, independently of the choice of the function $g$, we always have $a_g\equiv 1$. Moreover, for each one of the functions that we are considering for our mesh refinement also $b_g$ can be computed explicitly. In more detail, we have
\begin{displaymath}
	\begin{array}{lll}
		\displaystyle g(y) = \tan\left(\frac \pi4 y\right) &\Rightarrow &\displaystyle b_g(x) = -\frac{2x}{x^2+1}
		\\[10pt]
		\displaystyle g(y) = 2\sin\left(\frac \pi6 y\right) &\Rightarrow &\displaystyle b_g(x) = \frac{x}{4-x^2},
	\end{array}
\end{displaymath}
and the Hamiltonian systems become
\begin{align}\label{ham_g1}
	g(y) = \tan\left(\frac \pi4 y\right) \;\;\;\Rightarrow \;\;\;
		\begin{cases}
			\displaystyle\dot{x}^\pm(t) = \mp \cos\left(\frac{\xi^\pm(t)}{2}\right)
			\\[9pt]
			\displaystyle\dot{\xi}^\pm(t) = \mp \frac{4x^\pm(t)}{x^\pm(t)^2+1}\sin\left(\frac{\xi^\pm(t)}{2}\right)
			\\[9pt]
			x^\pm(0) = x_0,\;\;\;\xi^\pm(0)=\xi_0
		\end{cases}
\end{align}
and
\begin{align}\label{ham_g2}
	g(y) = 2\sin\left(\frac \pi6 y\right) \;\;\;\Rightarrow \;\;\;
		\begin{cases}
			\displaystyle\dot{x}^\pm(t) = \mp \cos\left(\frac{\xi^\pm(t)}{2}\right)
			\\[9pt]
			\displaystyle\dot{\xi}^\pm(t) = \pm \frac{2x^\pm(t)}{4-x^\pm(t)^2}\sin\left(\frac{\xi^\pm(t)}{2}\right)
			\\[9pt]
			x^\pm(0) = x_0,\;\;\;\xi^\pm(0)=\xi_0.
		\end{cases}
\end{align}

\subsection{Discrete phase portraits and their interpretation}\label{pp_sec}

\begin{figure}[!h]
	\centering
	\includegraphics[scale=0.5]{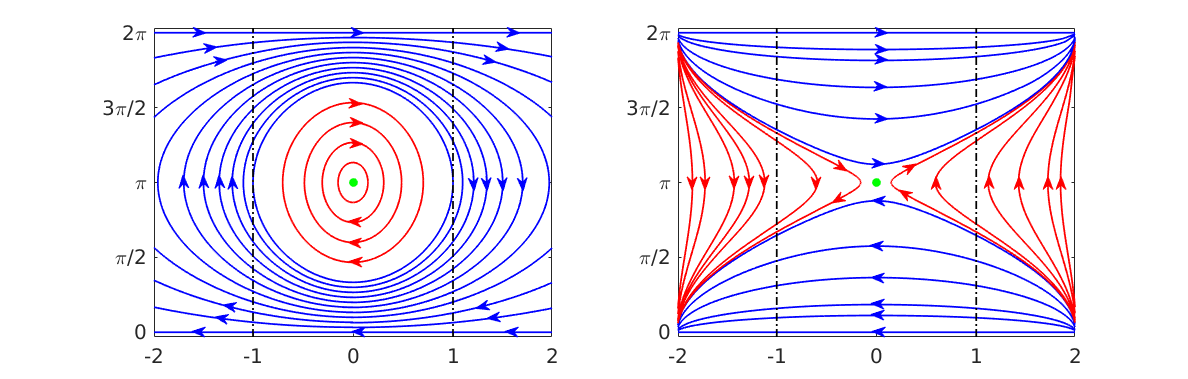}
	\caption{Phase portrait of the Hamiltonian system for the numerical wave equation and the grid transformations $\textbf{g}^{\,h,1}$ (left) and $\textbf{g}^{\,h,2}$ (right). We put $x^\pm(t)$, $\xi^\pm(t)$ on the horizontal/vertical direction.}\label{plot_phaseportrait}
\end{figure}

To help us in our further discussion, we include in Figure \ref{plot_phaseportrait} the phase diagrams corresponding to \eqref{ham_g1} and \eqref{ham_g2}, i.e. to the Hamiltonian system associated to \eqref{wave_discr_1d} on the two non-uniform meshes that we selected. The case of a uniform mesh is not displayed there since, the related discrete symbol being independent of the variable $x$, the orbits would simply be straight lines parallel to the horizontal axis. 

A couple of additional comments are needed here. First of all, in both phase diagrams we have a unique equilibrium in the point $P_{\textrm{e}}:=(x_{\textrm{e}},\xi_{\textrm{e}})= (0,\pi)$ (the green one). Nevertheless, the nature of this equilibrium changes when changing the mesh function. This is easily seen by considering the linearization of \eqref{ham_g1} and \eqref{ham_g2} around $P_{\textrm{e}}$, which is given by the linear systems
\begin{align*}
	\left(\begin{array}{c}
		\dot{x}^\pm \\ \dot{\xi}^\pm
		\end{array}\right) = A_{g_i}\left(\begin{array}{c}
		x^\pm \\ \xi^\pm
		\end{array}\right),\;\;\;i=1,2,
\end{align*}
with
\begin{align*}
	A_{g_1}:=\left(\begin{array}{cc}
		0 & 1/2 \\ -4 & 0
		\end{array}\right) \;\;\;\textrm{ and }\;\;\;
	A_{g_2}:=\left(\begin{array}{cc}
		0 & 1/2 \\ 1/2 & 0
		\end{array}\right).
\end{align*}

In the first case, it can be readily checked that the eigenvalues of $A_{g_1}$ are $\lambda_{1,2} = \pm i\sqrt{2}$, i.e. they are purely imaginary. In view of that, we can conclude that $P_{\textrm{e}}$ is a center. On the other hand, for the other mesh refinement that we are considering, the eigenvalues of $A_{g_2}$ are $\lambda_{1,2} = \pm 1/2$, that is, they are purely real with opposite signs. This implies that this times the fixed point is a saddle.

The second observation is related to the range of frequencies that we chose to include in our phase diagrams. In principle, the most suitable choice for the domain of the phase variable in this finite difference setting would be $\xi\in[-\pi,\pi]$ this being related essentially to the $2\pi$-periodicity of the discrete Wigner transform and to the fact that the solutions of our semi-discrete wave equation may be written as linear combinations of monochromatic waves given by the complex exponentials 
\begin{align*}
	e^{\pm ij\left(\frac{\sqrt{\lambda_j}}{j\pi}t-x\right)}, \;\;\; j\in 1\ldots N,
\end{align*}
$\lambda_j:= (4/h^2)\sin\left(j\pi h/2\right)$ being the eigenvalues of the one-dimensional finite difference Dirichlet Laplacian (see \cite{macia2002propagacion,marica2015propagation,zuazua2005propagation} and the references therein). In view of that, the relevant range of frequencies for the semi-discrete waves is $\xi\in[0,\pi]$, and considering $\xi\in[-\pi,\pi]$ in the phase portraits takes into account the two branches of the associated bi-characteristic rays.

Despite of these considerations, in Figure \ref{plot_phaseportrait} we chose to take in to account a range of frequencies $\xi\in[0,2\pi]$, since we believe that in this way it is more visible the stable/unstable nature of the equilibrium points. Consequently, our phase diagrams have to be interpreted as showing in the upper part $\xi\in[\pi,2\pi]$ what would actually correspond to $\xi\in[-\pi,0]$.
To better understand this fact, let us follow, for instance, one of the blue orbits in the lower part of the portraits, corresponding to some given initial frequency $\xi_0$ close to zero. The wave corresponding to this trajectory starts propagating to the left, until it hits the boundary of the physical domain $(-1,1)$ (indicated with the two black dotted lines) with a frequency $\xi_1>\xi_0$. Then, it reflects according to the laws of geometric optics, and starts propagating to the right, along a trajectory in the phase portrait whose initial frequency is $\xi_2=2\pi-\xi_1$. Notice that, if we were considering a phase portrait on $[-\pi,\pi]$ instead of on $[0,2\pi]$, the initial frequency of the reflected trajectory would instead be $\xi_2=-\xi_1$. 

We now present and analyze the results of our simulations. In all this part, if not specified, we always consider a time horizon $T=5s$. Moreover, since the solution to \eqref{wave_discr_1d} starting from an initial datum as in \eqref{in_data_wave} is complex, in our plots we decided to show its modulus.

We start by observing that, at low frequencies (that is, $\xi\in[0,\pi]$ but not too close to $\xi=\pi$), the numerical solutions of \eqref{wave_discr_1d} behave basically like the solution of the continuous model \eqref{wave_1d}: it starts traveling to the left along the straight characteristic line $x+t$ and, after having hit the boundary, it reflects following the Descartes-Snell's law and continues propagating, this time to the left along the other branch of the characteristic ($x-t$).
\begin{SCfigure}[1][h]
	\centering 
	\includegraphics[scale=0.4]{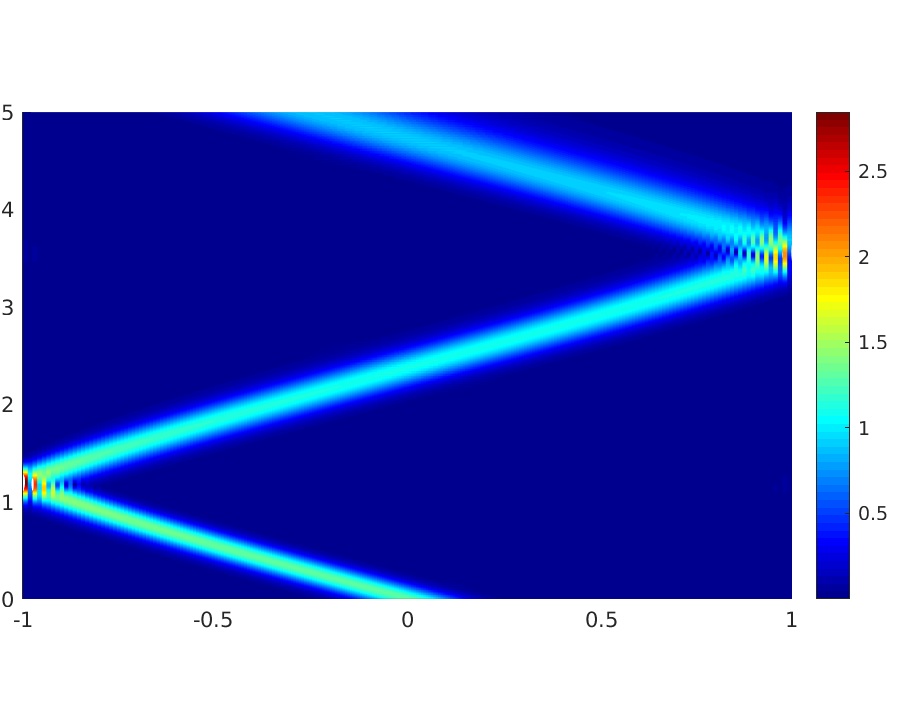}
	\caption{Propagation of a Gaussian wave packet for the finite difference scheme of equation \eqref{wave_discr_1d} on the non-uniform mesh $\textbf{g}^{\,h,1}$ at frequency $\xi_0 = \pi/4$. We put on the horizontal axis the space domain $(-1,1)$ and on the vertical on the time domain $(0,T)$.}\label{plot_01_wave}
\end{SCfigure} 

\begin{figure}[!h]
	\centering
	\begin{minipage}{0.3\textwidth}
		\includegraphics[scale=0.32]{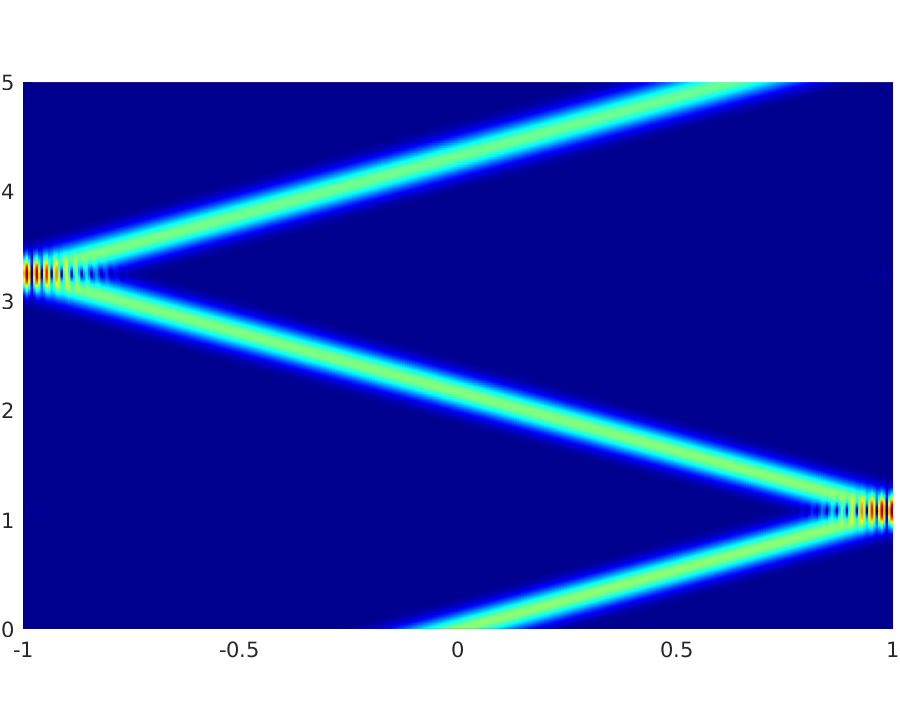}
	\end{minipage}
	\hspace{0.4cm}
	\begin{minipage}{0.3\textwidth}
		\includegraphics[scale=0.32]{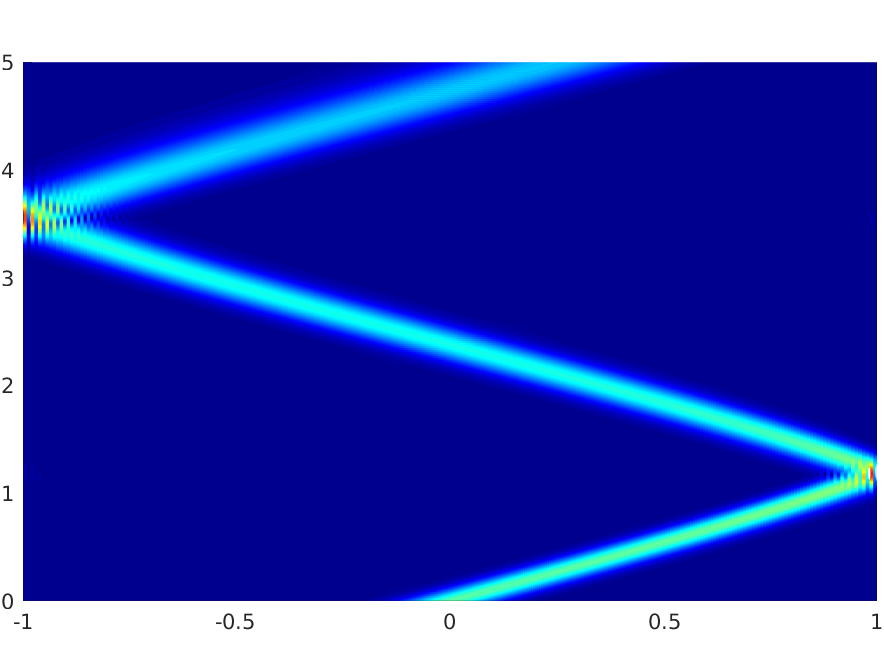}
	\end{minipage}
	\hspace{0.4cm}
	\begin{minipage}{0.3\textwidth}
		\includegraphics[scale=0.32]{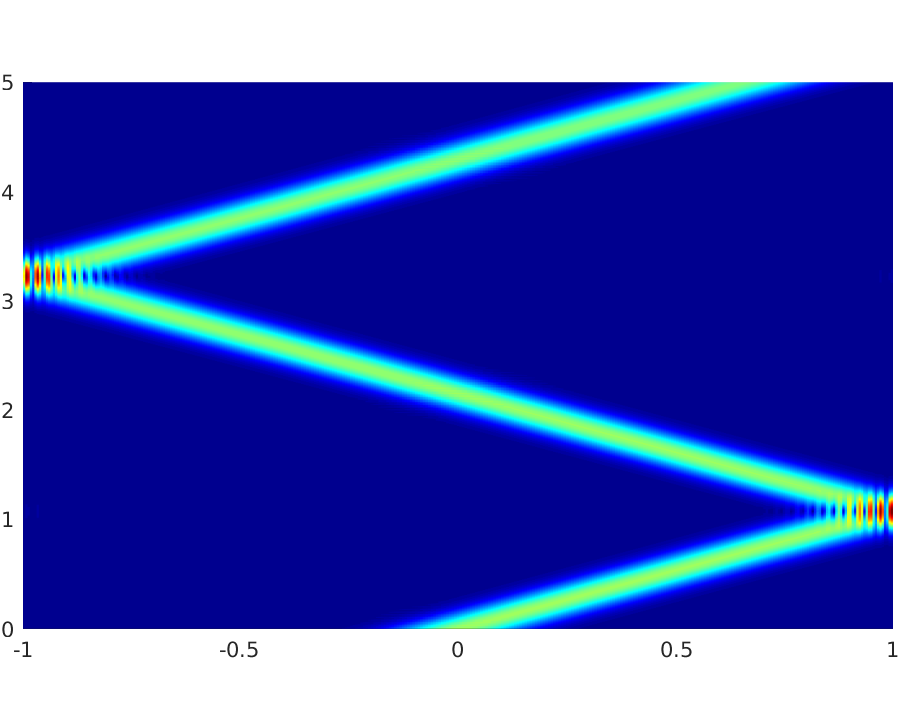}
	\end{minipage}	
	\caption{Numerical solutions with $x_0=0$, $\xi_0=7\pi/4 = 2\pi-\pi/4$ and uniform mesh (left), non-uniform mesh $\textbf{g}^{\,h,1}$ (middle) and non-uniform mesh $\textbf{g}^{\,h,2}$ (right).}\label{plot_05_wave}
\end{figure}

Moreover, it can be seen in Figure \ref{plot_05_wave} that an analogous but specular behavior is encountered also for frequencies $\xi\in[\pi,2\pi]$ sufficiently close to $\xi=2\pi$. This is in accordance with the discussion in the last part of Section \ref{pp_sec}.

When increasing the frequency, instead, the situation changes and we encounter several interesting phenomena and pathologies: 
\begin{itemize}
	\item[$\bullet$] Non-propagating waves, corresponding to equilibrium (fixed) points on the phase diagram. 
		
	\item[$\bullet$] The so-called umklapp or U-process, also known as internal reflection (see \cite{vichnevetsky1981propagation}), consisting in the reflection of waves without touching only one or both the endpoints of the space interval.
\end{itemize}

All these phenomena can be justified through the discussion that we introduced before and by looking at the phase portrait of the rays.

\subsubsection{Non-propagating waves}
\begin{figure}[!h]
	\centering 
	\begin{minipage}{0.3\textwidth}
		\includegraphics[scale=0.32]{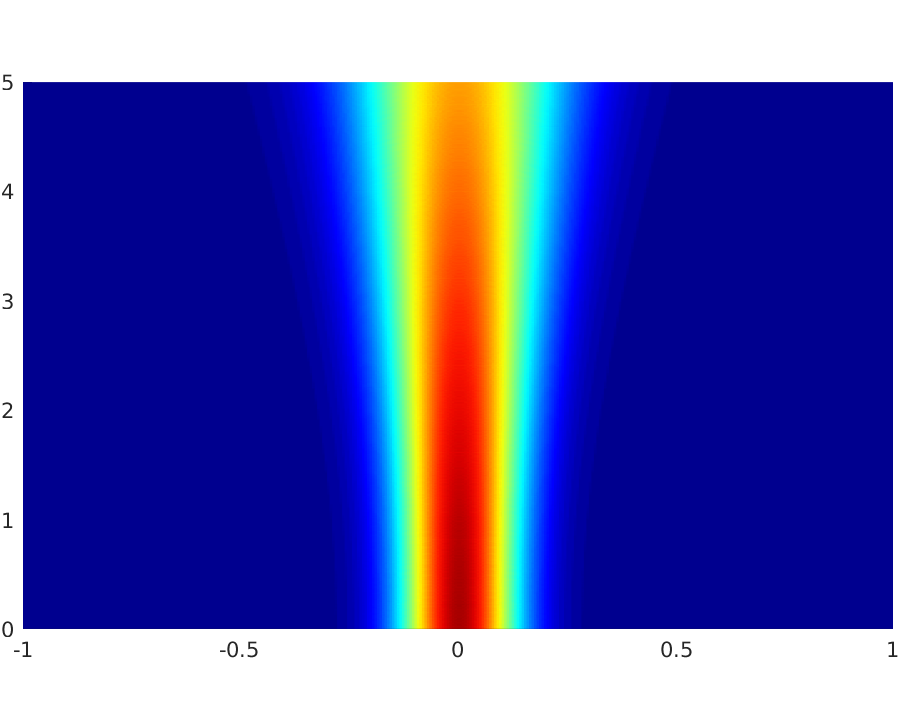}
	\end{minipage}
	\hspace{0.4cm}
	\begin{minipage}{0.3\textwidth}
		\includegraphics[scale=0.325]{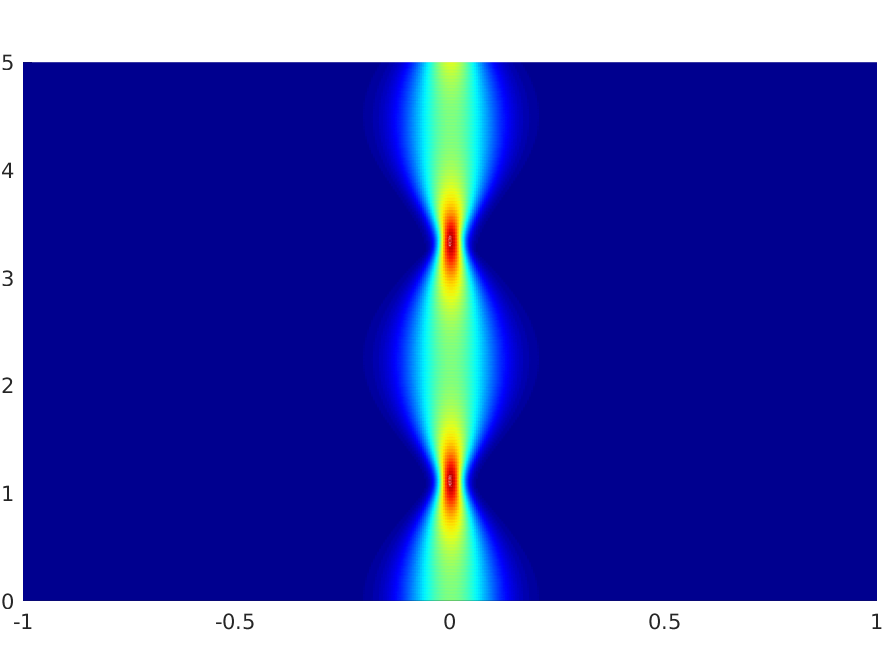}
	\end{minipage}
	\hspace{0.4cm}
	\begin{minipage}{0.3\textwidth}
		\includegraphics[scale=0.32]{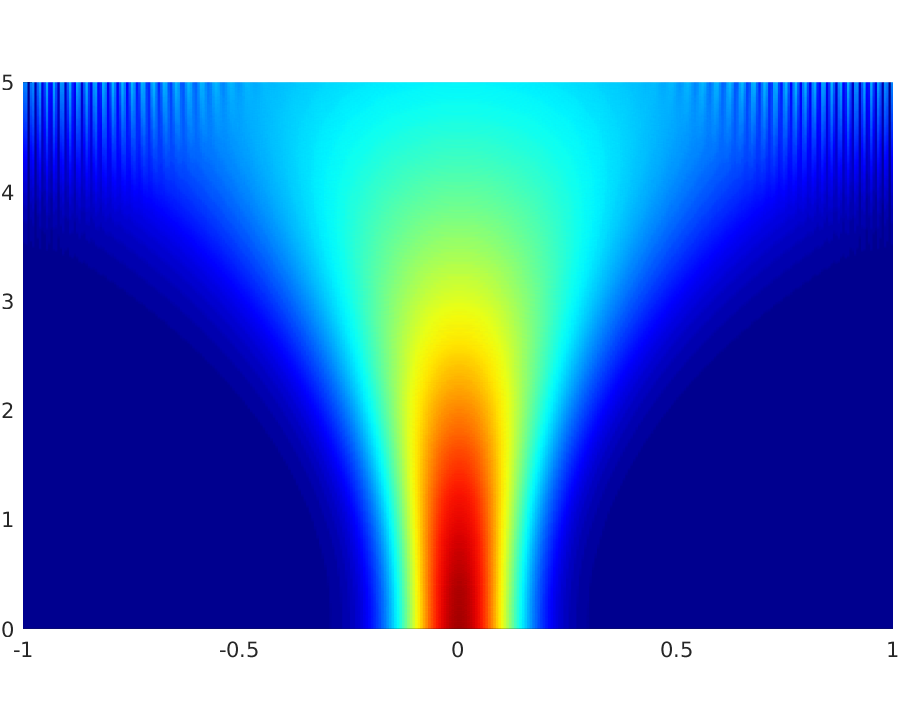}
	\end{minipage}	
	\caption{Numerical solutions with $x_0=0$, $\xi_0=\pi$ and uniform mesh (left), non-uniform mesh $\textbf{g}^{\,h,1}$ (middle) and non-uniform mesh $\textbf{g}^{\,h,2}$ (right).}\label{plot_02_wave}
\end{figure}

In Figure \ref{plot_02_wave} we observe waves that do not propagate. As we can see from the plots, and as we were mentioning before, this happens both with uniform and non-uniform meshes. 
$\newline$

The justification to this fact is that, as we saw in the previous section, for $\xi=\pi$ we have $\partial_\xi\omega(\xi)=0$ and, therefore, the velocity of the rays vanishes. 

However, we notice a big difference between the two plots corresponding to a non-uniform grid, concerning the dispersion along the ray. When solving using the mesh $\textbf{g}^{\,h,1}$, the wave remains concentrated along the ray, while when using the mesh $\textbf{g}^{\,h,2}$ the wave is very dispersive. This fact finds explanation in the analysis of the phase portrait associated to \eqref{ham_g1} and \eqref{ham_g2}. In both cases the non propagating wave corresponds to the only equilibrium point $P_{\textrm{e}}:=(x_{\textrm{e}},\xi_{\textrm{e}})= (0,\pi)$ (the green one) on the corresponding phase diagrams in Figure \ref{plot_phaseportrait}. Nevertheless, as we observed in Section \ref{pp_sec}, the nature of this equilibrium changes when changing the mesh function: for the grid $\textbf{g}^{\,h,1}$, the point $P_{\textrm{e}}$ is a center and, therefore, it produces a non-dispersive wave. On the other hand, for the other mesh refinement that we are considering, the fixed point is a saddle and, being this an unstable equilibrium, it generates the dispersion displayed in Figure \ref{plot_02_wave}.

\subsubsection{Internal reflection}

\begin{figure}[!h]
	\centering
	\begin{minipage}{0.3\textwidth}
		\includegraphics[scale=0.32]{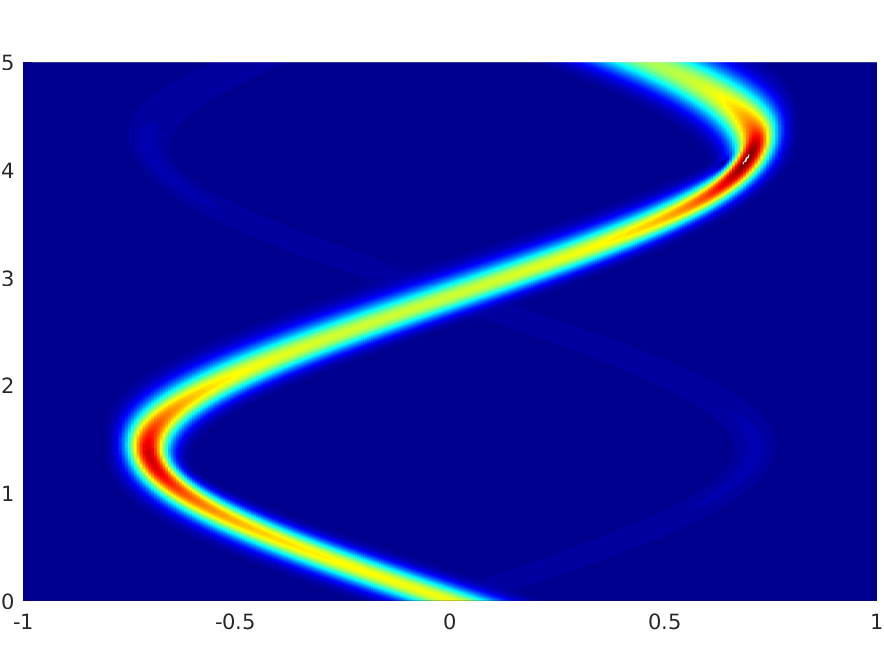}
	\end{minipage}
	\hspace{0.5cm}
	\begin{minipage}{0.3\textwidth}
		\includegraphics[scale=0.32]{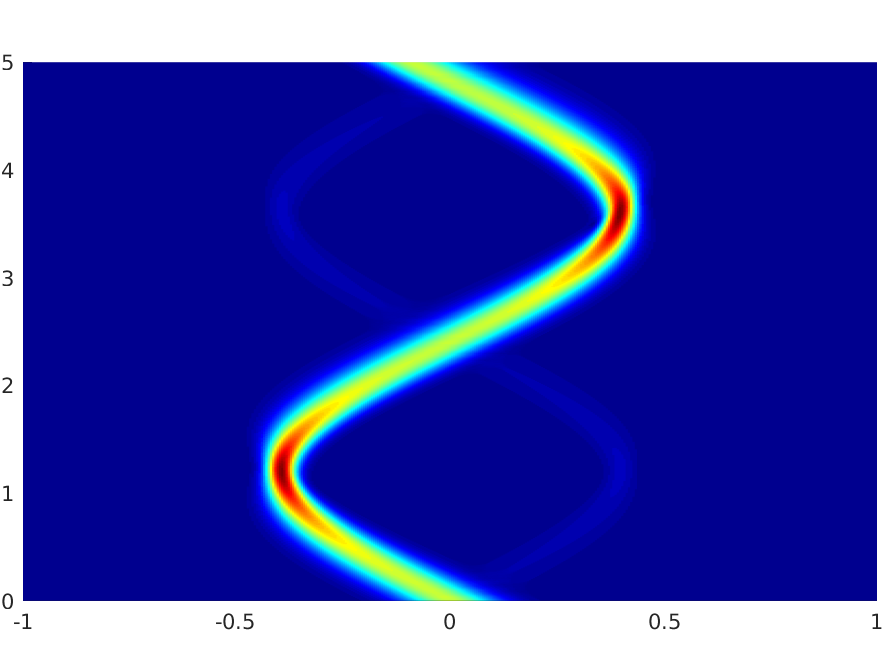}
	\end{minipage}	
	\hspace{0.5cm}
	\begin{minipage}{0.3\textwidth}
		\includegraphics[scale=0.32]{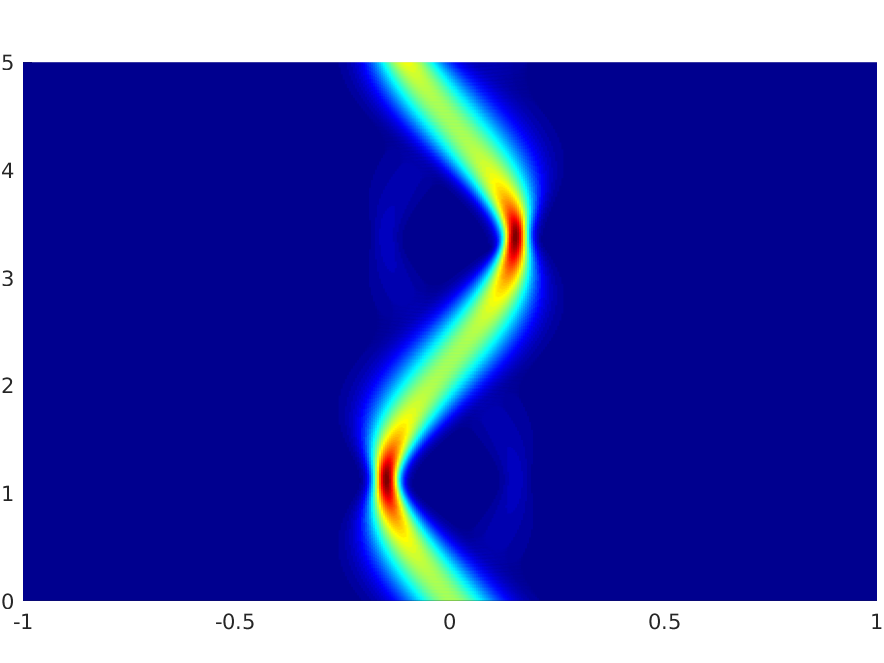}
	\end{minipage}	
	\caption{Numerical solutions corresponding to the mesh $\textbf{g}^{\,h,1}$ with $x_0=0$ and $\xi_0=7\pi/15$ (left), $\xi_0=10\pi/15$ (middle) and $\xi_0=13\pi/15$ (right).}\label{plot_03_wave}
\end{figure}

\begin{figure}[h]
	\centering
	\subfloat[$x_0=1/2$]{ 
			\includegraphics[scale=0.32]{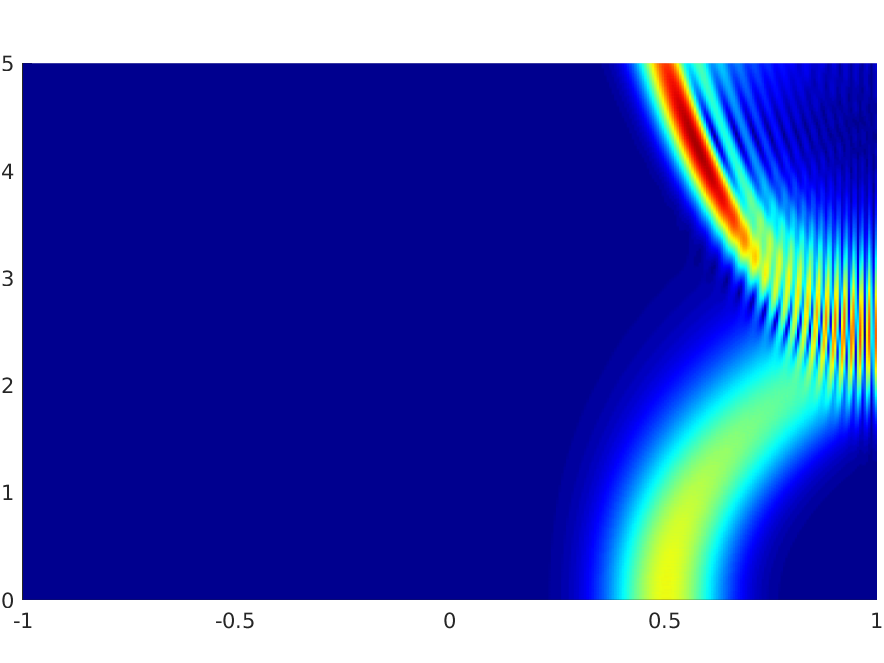}
	}
	\hspace{0.5cm}
	\subfloat[$x_0=-1/2$]{
			\includegraphics[scale=0.32]{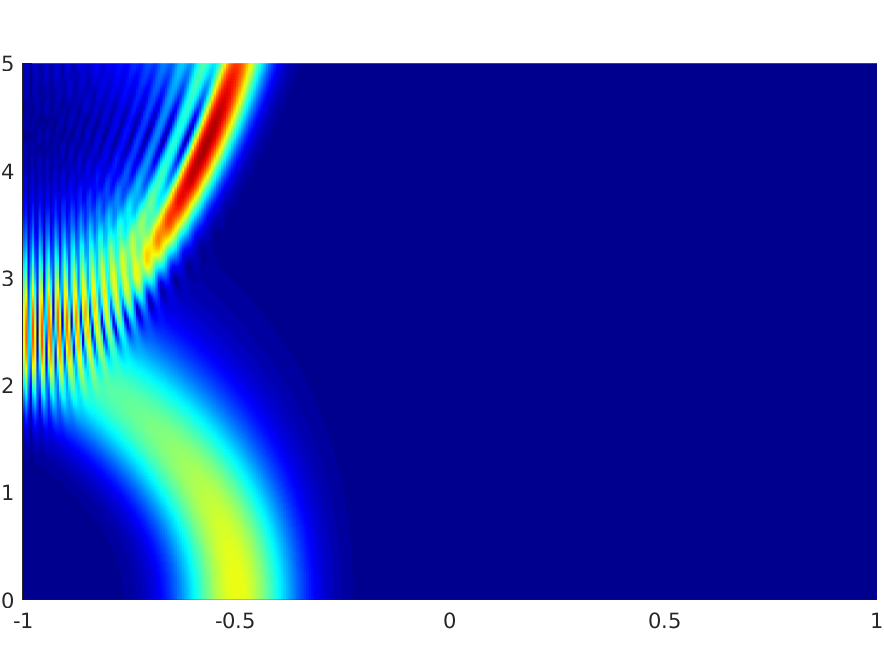}
	}	
	\caption{Numerical solutions corresponding to the mesh $\textbf{g}^{\,h,2}$ with $x_0=\pm 1/2$ and $\xi_0=\pi$.}\label{plot_04_wave}
\end{figure}

It is well known that, for the continuous case or for numerical waves on uniform meshes concentrated on frequencies where the group velocity is not trivial, all the generalized rays are straight lines reflecting at both endpoints. Instead, when the mesh is non-uniform, we observe that certain solutions to our wave equation do not preserve this behavior. In particular, our simulations show the following two pathologies:
\begin{enumerate}
	\item[(i)] waves oscillating in the interior of the computational domain and reflecting without touching the boundary (see
	Figure \ref{plot_03_wave});
	\item[(ii)] waves that oscillate in the interior of the domain and reflect touching the boundary only at one of the endpoints (see Figure \ref{plot_04_wave}).
\end{enumerate}

These trapped rays correspond to trajectories which remain always in the red area of the phase portraits in Figure \ref{plot_phaseportrait}. More precisely, the situation (i) appears when employing the mesh $\textbf{g}^{\,h,1}$, and it corresponds to periodic orbits in the phase diagram which are completely included in the region between the two dotted black vertical asymptotes indicating the computational domain $[-1,1]$. 

Notice that, as the frequency increases towards $\xi=\pi$, the amplitude of the oscillation decreases, the limit case being the non-propagating wave of Figure \ref{plot_02_wave}. 

Moreover, this phenomena can be explained if we consider that the mesh $\textbf{g}^{\,h,1}$ is an expanding one, meaning that the step size increases approaching the endpoints of the domain. Consequently, the group velocity $1/h$ of the high-frequency waves decreases while moving away from $x=0$. If this group velocity vanishes before the wave has reached the boundary, then this results in a process of internal reflection. 

Furthermore, the non-uniformity of the mesh size $h$ is also responsible of the increasing and decreasing of the magnitude of the solution during its propagation (roughly speaking, the change of color from yellow to red in our plots). Indeed, the amplitude of the wave is the one of the Gaussian profile of the initial datum, which is given by the constant $\gamma$ in \eqref{in_data_wave} and depends on $h$. In particular, on the mesh $\textbf{g}^{\,h,1}$ that we are using in Figure \ref{plot_03_wave}, while approaching the boundary $h$ increases. Therefore, the support of the ray shrinks and, due to energy conservation, the high of the corresponding wave has to increase since the same amount of energy is now concentrated in a smaller region. Then, when moving again towards the center of the physical domain, where $h$ is smaller, the support of the ray increases again and the magnitude of the wave decreases. 

Let us now conclude by briefly discussing the situation (ii), which appears when employing the mesh $\textbf{g}^{\,h,2}$. Recall that in this case $P_{\textrm{e}}$ is a saddle point, which is characterized by the fact that the space around it is divided into four sectors by two curves (the separatrices) passing through the equilibrium. In view of that, the red curves always remain trapped either in the region $x\in[0,1]$ (Figure \ref{plot_04_wave}a) or $x\in[-1,0]$ (Figure \ref{plot_04_wave}b). Moreover, notice that this time the mesh is coarser around $x=0$. On the one hand, this this implies that the velocity of high-frequency waves may vanish approaching the center of the domain, thus generating again some internal reflection phenomena. On the other hand, this reflects again on the increasing of the magnitude of the wave, according to the discussion that we presented above. 
%

\subsection{The case of variable coefficients}\label{var_coeff_sec}

We discuss briefly here the case of the variable coefficients wave equation. In more detail, the model that we are going to consider is the following
\begin{align}\label{wave_var_1d}
	\begin{cases}
		\rho(x)\partial^2_tu-\partial_x\big(\sigma(x)\partial_xu\big) = 0, & (x,t)\in (-1,1)\times(0,T)
		\\
		u(-1,t) = u(1,t) = 0, & t\in(0,T)
		\\
		u(x,0) = u^0(x),\;\;\;\partial_tu(x,0) = u^1(x), & x\in(-1,1),
	\end{cases}
\end{align}
where, we recall, $\rho$ and $\sigma$ are chosen to be $L^\infty(\RR)$-functions with the strict hyperbolicity assumption $\rho(x)\geq\rho^*>0$ and $\sigma(x)\geq\sigma^*>0$. 

As for the constant-coefficients wave equation \eqref{wave_1d}, also for \eqref{wave_var_1d} we have that the total energy given by the quantity
\begin{align*}
	\mathcal E_{\rho,\sigma} (u,\partial_t u) :=\frac 12 \int_{-1}^1\left(\rho(x)|u_t|^2+\sigma(x)|u_x|^2\right)\,dxdt
\end{align*}
is conserved in time, that is,
\begin{align*}
	\mathcal E_{\rho,\sigma} (u,\partial_t u) = \mathcal E_{\rho,\sigma} (u^0,u^1), \forall\, t\in(0,T).
\end{align*}

Moreover, even if solutions to \eqref{wave_var_1d} cannot be represented explicitly through a formula like \eqref{waves_dal}, we know that wave packets originated by highly concentrated and oscillatory initial data remain concentrated along one of the bi-characteristic lines, and their energy localized outside any neighborhood of the ray vanishes as the wavelength parameter tends to zero. 

We recall that the principal symbol associated to the one-dimensional variable coefficients wave equation \eqref{wave_var_1d} is given by 
\begin{align}\label{ham_cont_var_1d}
	\mathcal H_c(x,t,\xi,\tau) = -\rho(x)\tau^2 + \sigma(x)\xi^2,
\end{align}
and that the bi-characteristic rays are the curves $s\mapsto (x(s),t(s),\xi(s),\tau(s))$ solving the first order ODE system
\begin{align}\label{ham_syst_cont_var_1d}
	\begin{cases}
		\dot{x}(s) = 2\sigma(x(s))\xi(s)
		\\
		\dot{t}(s) = -2\rho(x(s))\tau(s)
		\\
		\dot{\xi}(s) = \rho'(x(s))\tau^2(s) - \sigma'(x(s))\xi^2(s)
		\\
		\dot{\tau}(s) = 0.
	\end{cases}
\end{align}

Moreover, the initial datum $(x(0),t(0),\xi(0),\tau(0)) = (x_0,0,\xi_0,\tau_0)$ associated to system \eqref{ham_syst_cont_var_1d} is chosen so that $\mathcal H_c(x_0,0,\xi_0,\tau_0)=0$. Then, we immediately see that $\tau(s)=\tau_0^\pm=\pm c(x(s))\xi(s)$ for all $s$, where with $c(\cdot)$ we denote the function 
\begin{align}\label{function_c}
	c(\cdot):=\sqrt{\frac{\sigma(\cdot)}{\rho(\cdot)}}.
\end{align}
Notice that now the bi-characteristics are not straight lines anymore, since $\dot{\xi}(s)$ does not vanish. 

Once again, the concentration features of the solutions to \eqref{wave_var_1d} can be observed also at the numerical level, and this will be the concern of the remaining of the present section. The space discretization that we are going to employ is totally analogous to the finite differences scheme on non-uniform mesh that we used in Section \ref{discretization_sec}. In particular, we still focus on the case in which $g$ is regular ($g\in C^2(\RR)$ with $0<g_d^-\leq|g'(x)|\leq g_d^+<+\infty$ and $|g''(x)|\leq g_{dd}<+\infty$ for some given constant $g_d^-,g_d^+,g_{dd}>0$ and for all $x\in\RR$). Notice that this, joint with the hypothesis on the coefficients $\rho$ and $\sigma$, implies that $c_g:=c(g)/g'$ belongs to the space $C^{0,1}(\RR)$ of the Lipschitz continuous function. As we will see, this Lipschitz regularity assumption will be very important when introducing the Hamiltonian system for the bi-characteristic rays. Moreover, notice that when $\rho=\sigma\equiv 1$ the function $c_g$ becomes the one that we introduced in the case of constant coefficients.

Also the semi-discretization of the wave equation \eqref{wave_var_1d} on the non-uniform grid $\Ggh$ is analogous to the constant coefficients one, and it is given as
\begin{align}\label{wave_var_discr_1d}
	\begin{cases}
		\displaystyle h_j\rho(g_j)u_j''(t) - \left(\sigma(g_{j+1/2})\frac{u_{j+1}(t)-u_j(t)}{h_{j+1/2}}-\sigma(g_{j-1/2})\frac{u_j(t)-u_{j-1}(t)}{h_{j-1/2}}\right) = 0, & j=1,\ldots,N,\;\;t\in(0,T)
		\\
		u_0(t) = u_{N+1}(t) = 0, & t\in(0,T)
		\\
		u_j(0) = u_j^0, \;\;\; u_j'(0) = u_j^1, & j=1,\ldots,N.
	\end{cases}
\end{align}

Here, as before, we indicate $u_j(t):=u(g_j,t)$. Moreover, to avoid possible confusions, when needed, we will denote the solutions of \eqref{wave_var_discr_1d} by $\textbf{u}^{h,g}(t):=(u^g_j(t))_{j=0}^{N+1}$ for emphasizing the dependence from the diffeomorphic transformation $g$. But, in general, we will simply write down $\mathbf{u}^h(t):=(u_j(t))_{j=0}^{N+1}$, $\mathbf{u}^{0,h}:=(u^0_j)_{j=0}^{N+1}$ and $\mathbf{u}^{1,h}:=(u^1_j)_{j=0}^{N+1}$ for the solutions and for the initial data.

It has been shown in \cite[Proposition 3.4]{marica2015propagation} that, under our regularity assumptions on the coefficients $\rho$ and $\sigma$ and on the mesh function $g$, the numerical approximation \eqref{wave_var_discr_1d} with initial data $(\mathbf{u}^{0,h},\mathbf{u}^{1,h})$ is a convergent scheme of order $O(h)$ for the wave equation \eqref{wave_var_1d} in the appropriate $\ell^2$ setting. Moreover, it can be readily checked that, as it happens for the continuous model, system \eqref{wave_var_discr_1d} enjoys the property of energy conservation, that is,
\begin{align*}
	\mathcal E^{h,g}(\mathbf{u}^{h,g}(t),\partial_t \mathbf{u}^{h,g}(t)) &:= \frac 12 \sum_{j=1}^N h_j\rho(g_j)|\partial_t u_j^g(t)|^2 + \frac 12 \sum_{j=1}^N h_{j+1/2}\sigma(g_{j+1/2})\left|\frac{u_{j+1}^g(t)-u_j^g(t)}{h_{j+1/2}}\right|^2 
	\\
	&= \mathcal E^{h,g}(\mathbf{u}^{0,h},\mathbf{u}^{1,h}).
\end{align*}

Concerning now the propagation of the discrete bi-characteristic rays, it has been shown again in \cite{marica2015propagation} that the principal symbol associated to \eqref{wave_var_discr_1d} is given by 
\begin{align}\label{ham_discr_var_1d}
	\mathcal H(y,t,\xi,\tau) = -g'(y)\rho(g(y))\tau^2 + 4\sin^2\left(\frac \xi2\right)\frac{\sigma(g(y))}{g'(y)},\;\;\; y=g^{-1}(x),
\end{align}
or equivalently
\begin{align}\label{ham_discr_var_1d_2}
	\mathcal H(y,t,\xi,\tau) = -\tau^2 + c_g(y)^2\omega(\xi)^2, 
\end{align}
with
\begin{align*}
	\omega(\xi):=2\sin\left(\frac \xi2\right).
\end{align*}
Then, proceeding as in Section \ref{ham_sec} it is easy to obtain the following ODE system for the rays 
\begin{align}\label{ham_syst_discr_var_1d_t}
	\begin{cases}
		\displaystyle\dot{y}^\pm(t) = \mp c_g(y^\pm(t))\partial_\xi \omega(\xi^\pm(t))
		\\
		\displaystyle\dot{\xi}^\pm(t) = \pm\partial_yc_g(y^\pm(t))\omega(\xi^\pm(t))
		\\
		y^\pm(0) = y_0,\;\;\;\xi^\pm(0)=\xi_0
	\end{cases}
\end{align}

At this level, the Lipschitz regularity of the function $c_g$ becomes fundamental in order to guarantee existence and uniqueness of solutions to the above Cauchy problem. Moreover, since by assumption $c_g(y)>0$, thanks to \eqref{ham_syst_discr_var_1d_t} we immediately see that
\begin{align*}
	|\dot{y}^\pm(t)| = c_g(y^\pm(t))\left|\partial_\xi\omega(\xi^\pm(t))\right|.
\end{align*}

In view of that, as it was for the constant coefficients case, the velocity of the rays vanishes for $\xi = (2k+1)\pi$, $k\in\ZZ$, since these are the zeros of $\partial_\xi \omega(\xi) = \cos(\xi/2)$. Consequently, for these frequencies we will have again the phenomenon of non-propagating waves. 

\subsubsection*{Numerical results}

We conclude this section presenting and discussing several simulations for equation \eqref{wave_var_discr_1d}. As before, we chose an initial datum as in \eqref{in_data_wave}, and we used a leapfrog scheme for the time integration. Moreover, we employed both a uniform mesh and the non-uniform ones that we already used in Section \ref{numerical_sec}.

In what follows, for keeping our discussion lighter an more clear, we decided to consider only a specific case of variable coefficients, namely
\begin{align}\label{var_coeff}
	\rho(x)\equiv 1 \;\;\;\textrm{ and } \;\;\; \sigma(x)=1+A\cos^2(\kappa\pi x), \;\;\; A>0,\;\kappa\in\NN^*.
\end{align}

A complete dissertation with general coefficients would be, in our opinion, very lengthy and, possibly, difficult to follow. By focusing only on a particular election of $\rho$ and $\sigma$, instead, we believe that we can describe more clearly several phenomena in the propagation of the discrete rays. Actually, the pathologies that we discussed in the previous section are related to the mesh employed and to the trigonometric nature of the discrete group velocity, more than to the choice of the coefficients.   

Let us start with the propagation at low frequencies. Here we will only consider simulations on a uniform mesh since, when the frequency of the initial datum is small enough, the effects associated to the change of the grid are very mild and cannot be really appreciated in the plots. 

\begin{figure}[!h]
	\centering 
	\begin{minipage}{0.33\textwidth}
		\includegraphics[scale=0.32]{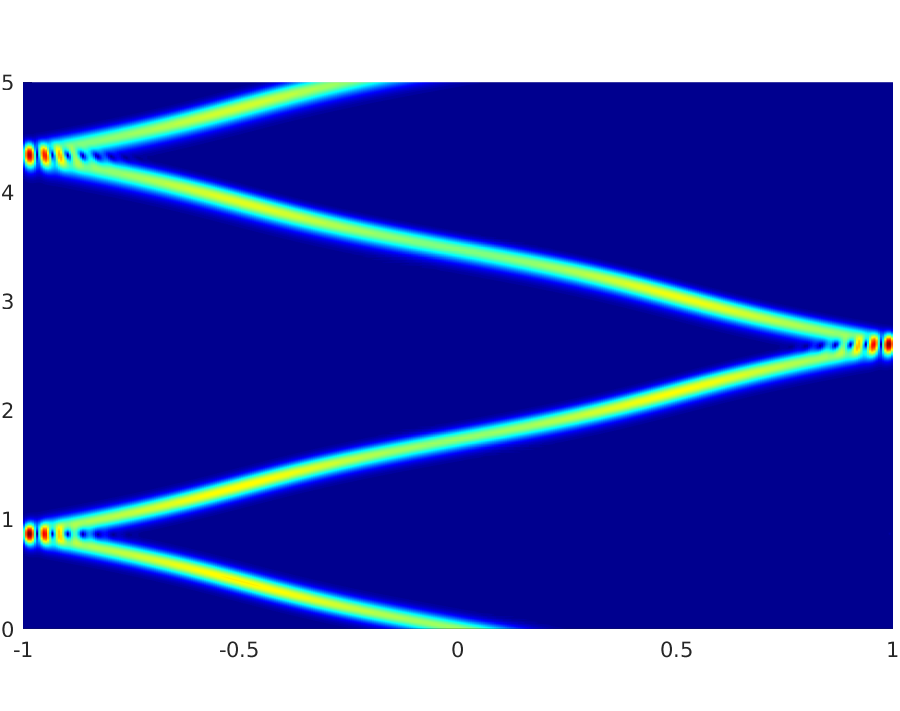}
	\end{minipage}
	\hspace{0.5cm}
	\begin{minipage}{0.33\textwidth}
		\includegraphics[scale=0.32]{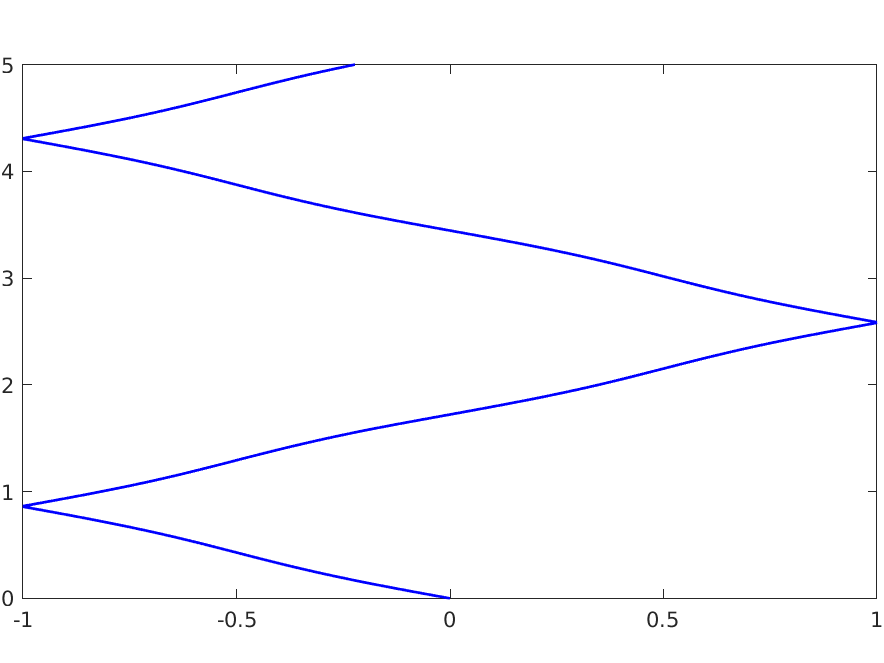}
	\end{minipage}	
	
	\begin{minipage}{0.33\textwidth}
		\includegraphics[scale=0.32]{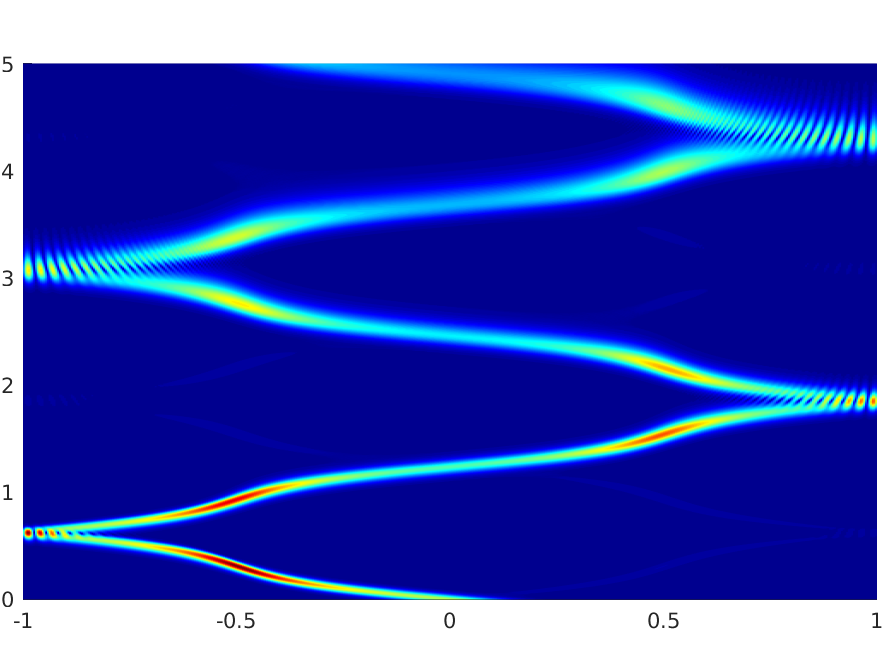}
	\end{minipage}
	\hspace{0.5cm}
	\begin{minipage}{0.33\textwidth}
		\includegraphics[scale=0.32]{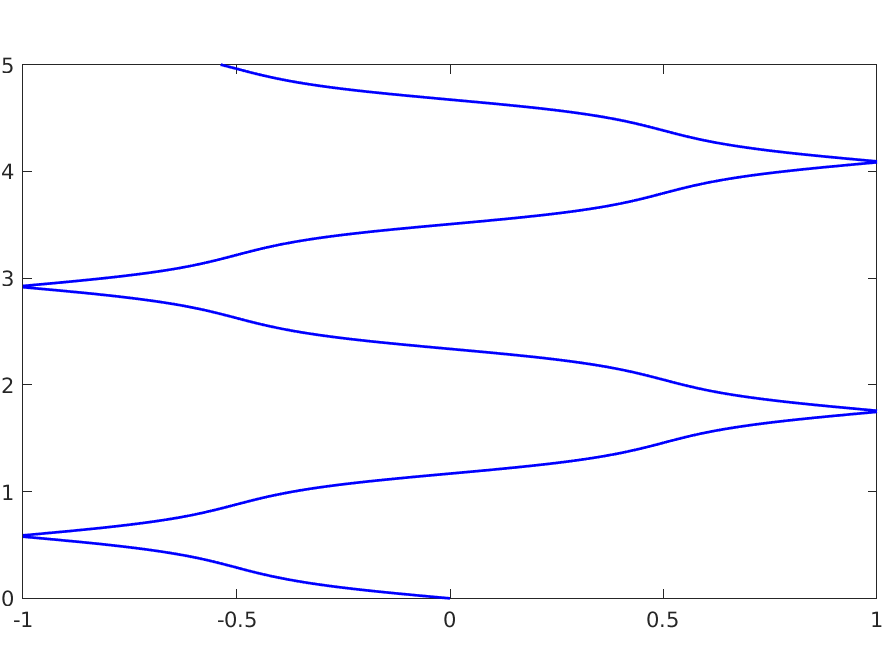}
	\end{minipage}
	\caption{Numerical solutions to \eqref{wave_var_1d} on a uniform mesh, corresponding to an initial datum with $x_0=0$ and $\xi_0=\pi/7$, and to a coefficient $\sigma$ as in \eqref{var_coeff} with $\kappa=1$ and $A=1$ (top) or $A=7$ (bottom). On the left we plot the wave, while on the right we plot the corresponding ray.}\label{plot_01_wave_var}
\end{figure}

We can see in Figure \ref{plot_01_wave_var} how the waves travel along characteristics and reach the boundary, where they are reflected according to the Descartes-Snell's law. This behavior is consistent with what we already observed for the constant coefficients equation, and it is totally expected since we know that at low frequency the discretization does no introduce relevant pathologies.

Recall that this time the rays are not straight lines anymore, $\dot{\xi}^\pm(t)$ in \eqref{ham_syst_discr_var_1d_t} being not zero. In fact, in our plots it is possible to appreciate a slight bending, generated by the variable coefficient $\sigma$. For having a better perception of that, in what follows we are plotting on the left the wave and on the right the ray, solution of the Hamiltonian system.

Moreover, we can observe that the parameters $A$ and $\kappa$ in the coefficient $\sigma$ affect the shape of the rays. First of all, $A$ measures the amplitude of the oscillation of the coefficient $\sigma$ and, therefore, it is natural to expect that increasing this parameter also the bending of the ray should increase. This fact can be clearly seen in Figure \ref{plot_01_wave_var}.

Moreover, a further justification of this phenomenon comes from an easy analysis of the Hamiltonian system giving the equations of the bi-characteristics $(x^\pm(t),\xi^\pm(t))$. 

For simplifying the presentation, here we will only consider one branch of the ray, that we will denote simply $(x(t),\xi(t))$. It can be readily checked that these curves can be obtained by solving the following ODE system:
\begin{align}\label{ham_syst_var_low}
	\begin{cases}
		\displaystyle\dot{x}(t) = -\sqrt{1+A\cos^2(\kappa\pi x(t))}\cos\left(\frac{\xi(t)}{2}\right)
		\\[10pt]
		\displaystyle\dot{\xi}(t) = -A\kappa\pi\frac{\sin(2\kappa\pi x(t))}{\sqrt{1+A\cos^2(\kappa\pi x(t))}}\sin\left(\frac{\xi(t)}{2}\right)
		\\[10pt]
		x(0) = x_0,\;\xi(0)=\xi_0.
	\end{cases}
\end{align}

Now, from \eqref{ham_syst_var_low} we can obtain informations on the shape of the characteristic ray in terms of the parameters $A$ and $\kappa$. In what follows, we will use the notation $x_{A,\kappa}(t)$ for the trajectory of the ray, for highlighting its dependence on the two parameters. 

Let us start by fixing a value for $\kappa$, say $\kappa=1$, and by studying the behavior of $x_{A,1}(t)$ with respect to $A$. First of all, we notice that, according to \eqref{ham_syst_var_low}, we have
\begin{align*}
	A_1\geq A_2 \;\;\Rightarrow\;\; |\dot{x}_{A_1,1}(t)|\geq |\dot{x}_{A_2,1}(t)|.
\end{align*}

In other words, the velocity of $x_{A,1}(t)$ is an increasing function of $A$. In view of that, the ray corresponding to $A_2$ has a lower inclination with respect to the one corresponding to $A_1$ and it needs less time for reaching the boundary. This can be observed in Figure \ref{plot_01_wave_var}.

As a second thing, from \eqref{ham_syst_var_low} we can also compute the curvature of the rays, which is measured from the second derivative $\ddot{x}_{A,1}(t)$. In particular, we have
\begin{align*}
	\ddot{x}_{A,1}(t) = -\frac{A\pi}{2}\sin(2\pi x(t)).
\end{align*}
Then, we immediately see that also the modulus of the curvature is an increasing function of $A$, that is, 
\begin{align*}
	A_1\geq A_2 \;\;\Rightarrow\;\; |\ddot{x}_{A_1,1}(t)|\geq |\ddot{x}_{A_2,1}(t)|.
\end{align*}

As a consequence, when $A$ grows the ray increases its bending while, for low values of $A$, it is almost a straight line. Also this feature appears in our simulations.

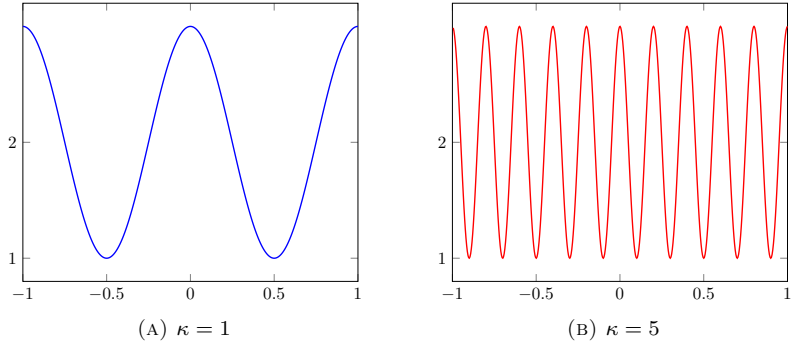
\begin{figure}[h]
	\begin{center}
		\subfloat[$\kappa=1$]{
			\begin{tikzpicture}[scale=0.65]
			\pgfplotsset{settings/.style={legend cell align = {left},
					legend pos =  south east,
					legend style = {draw = none}}}
			\begin{axis}[settings,xmin=-1,xmax=1, ytick={1,2}]
			
			\addplot[domain=-1:1, samples=700,thick,color=blue]{1+2*cos(deg(pi*x))^2};
			\end{axis}
			\end{tikzpicture}
		}\hspace{0.5cm}
		\subfloat[$\kappa=5$]{
			\begin{tikzpicture}[scale=0.65]
			\pgfplotsset{settings/.style={legend cell align = {left},
					legend pos =  south east,
					legend style = {draw = none}}}
			\begin{axis}[settings,xmin=-1,xmax=1, ytick={1,2}]
			
			\addplot[domain=-0.995:1, samples=500,thick,color=red]{1+2*cos(deg(5*pi*x))^2};
			\end{axis}
			\end{tikzpicture}
		}
		\caption{Function $\sigma(x)=1+2\cos^2(\kappa\pi x)$ for different values of $\kappa\in\NN^*$.}\label{sigma_fig}
	\end{center}
\end{figure}

Let us now fix a value for $A$ (say $A=2$), and let us briefly analyze what happens when varying the parameter $\kappa\in\NN^*$. Observe that the role of $\kappa$ is only to modulate the frequency of the oscillations in the coefficient, since it is evident that $\sigma$ is a periodic function of period $T=2\kappa$ (see Figure \ref{sigma_fig}).

\begin{figure}[!h]
	\centering 
	\begin{minipage}{0.3\textwidth}
		\includegraphics[scale=0.32]{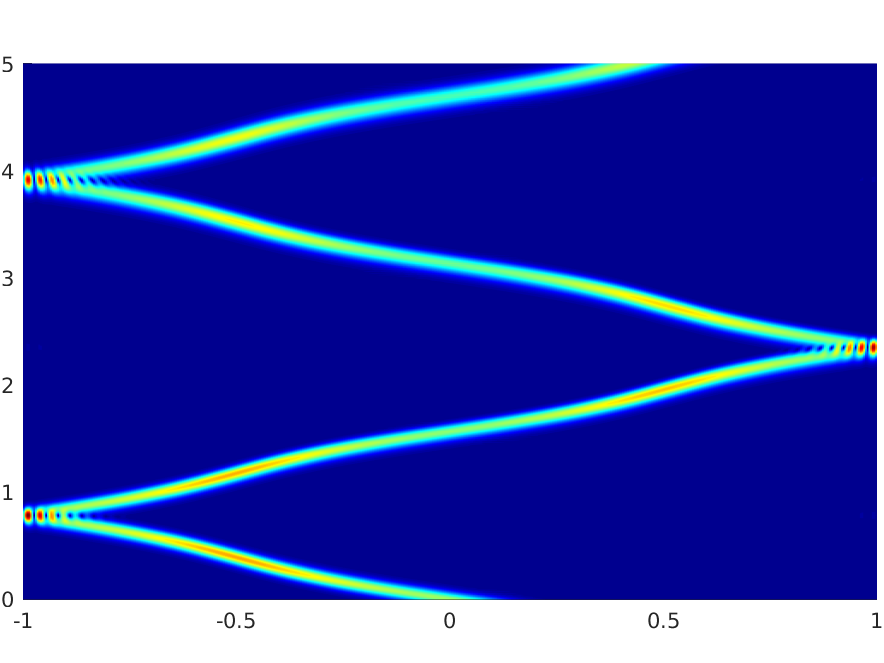}
	\end{minipage}
	\hspace{0.5cm}
	\begin{minipage}{0.3\textwidth}
		\includegraphics[scale=0.32]{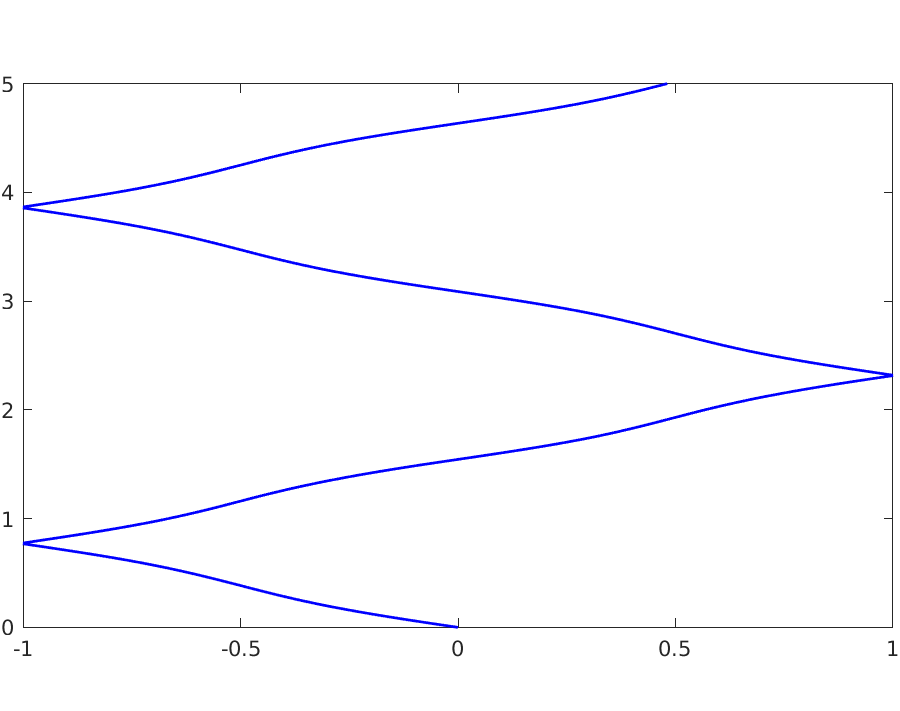}
	\end{minipage}
	
	\begin{minipage}{0.3\textwidth}
		\includegraphics[scale=0.32]{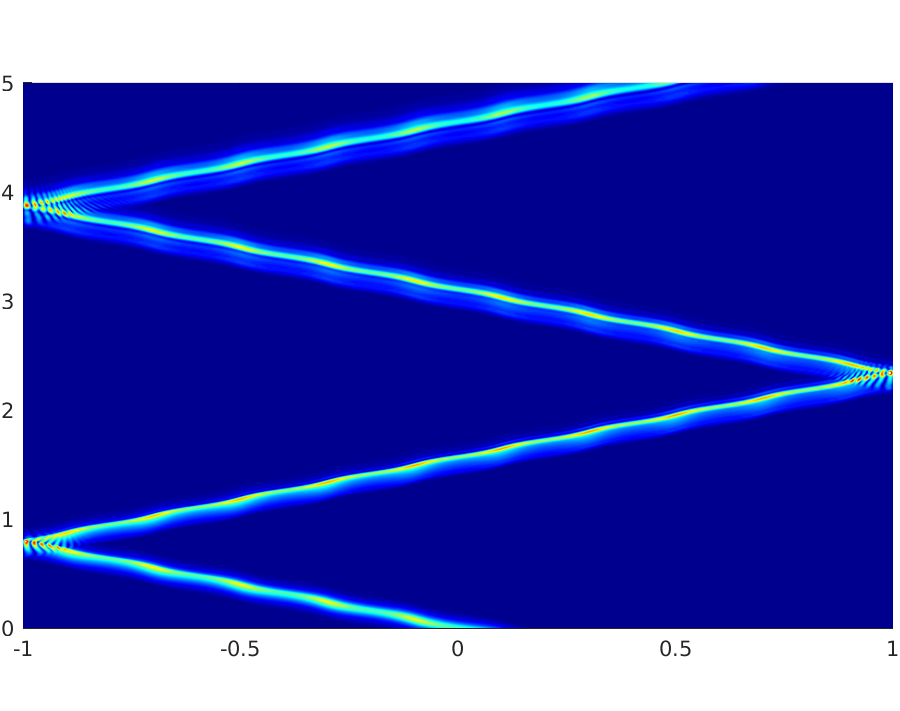}
	\end{minipage}
	\hspace{0.5cm}
	\begin{minipage}{0.3\textwidth}
		\includegraphics[scale=0.32]{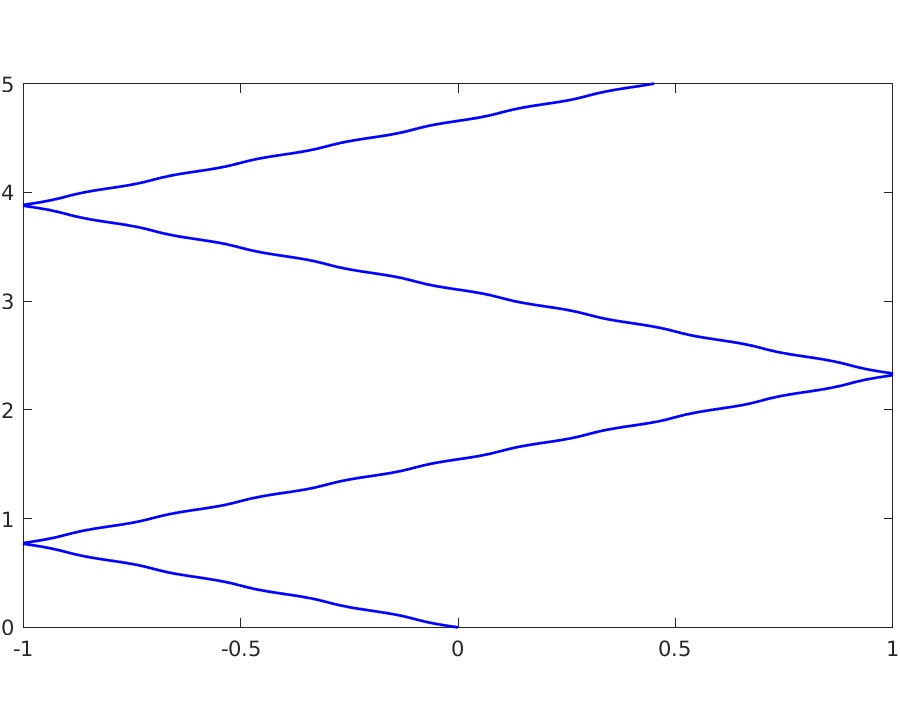}
	\end{minipage}	
	\caption{Numerical solutions to \eqref{wave_var_1d} on a uniform mesh, corresponding to an initial datum with $x_0=0$ and $\xi_0=\pi/7$, and to a coefficient $\sigma$ as in \eqref{var_coeff} with $A=2$ and $\kappa=1$ (top) or $\kappa=5$ (bottom). On the left we plot the wave, while on the right we plot the corresponding ray.}\label{plot_wave_var_freq}
\end{figure}

This reflects also on the Hamiltonian system \eqref{ham_syst_var_low}, in particular on the trajectory in the physical space given by the curve $x(t)$. Actually, it can be readily checked that, as $\kappa$ grows, also the number of equilibrium points in the corresponding phase portrait raises, and this produces an increasing in the oscillations of the orbits. Consequently, we can observe in Figure \ref{plot_wave_var_freq} how, augmenting $\kappa$, also the number of oscillations of the ray augments. Moreover, the high oscillation of the coefficient add also some numerical dispersion phenomena in our simulations. 

\begin{remark}
In this phenomena of changing of oscillation in the rays when increasing or decreasing $\kappa$, we find analogies with the works of Pjatnicki{\u\i} and, later, Allaire on the homogenization of rapidly oscillating coefficients wave equations (see, e.g., \cite{allaire2003dispersive,allaire2005homogenization,pjatnickiui1982limit} and the references therein). In particular, in \cite{pjatnickiui1982limit} the limit behavior of the domain of dependence of a hyperbolic equation with rapidly oscillating coefficients in the form
\begin{align}\label{wave_hom}
	\frac{\partial^2}{\partial t^2}u^\varepsilon - \frac{\partial}{\partial x_i}a_{i,j}\left(\frac x\varepsilon\right)\frac{\partial}{\partial x_j}u^\varepsilon = 0,\;\;\;\varepsilon>0
\end{align}
is analyzed through the study of the bi-characteristic flow. As a matter of fact, we know that this domain of dependence is given by a conical region delimited by the rays. 

It is shown there that the solutions to the Hamiltonian systems associated to \eqref{wave_hom} and to its homogenized version are at a distance of the order of $\varepsilon$ and that, consequently, the domains of dependence for the two systems are at a distance of the order of $\sqrt{\varepsilon}$. Moreover, the cone for the non-homogenized equation is always slightly wider than the cone for the homogenized one. Clearly, the shape of the domain of dependence of \eqref{wave_hom} is affected by the changing of the parameter $\varepsilon$ and by the the effects that this produces in terms of oscillations in the ray. Finally, we mention that these kind of phenomena are related to the behavior of the spectrum of the elliptic operator associated to \eqref{wave_hom}, which has been studied by Castro and Zuazua in \cite{castro1996remark,castro2000low,castro2002concentration} with control purposes.
\end{remark}

Let us now conclude this section by analyzing the behavior of the solutions to \eqref{wave_var_discr_1d} for high frequencies, where the same kind of pathological behaviors that we detected in the simulations with constant coefficients appear again. In what follows, we will always assume $A=1$ and $\kappa=1$ in the coefficient $\sigma$.

\begin{figure}[!h]
	\centering
	\subfloat{ 
		\includegraphics[scale=0.35]{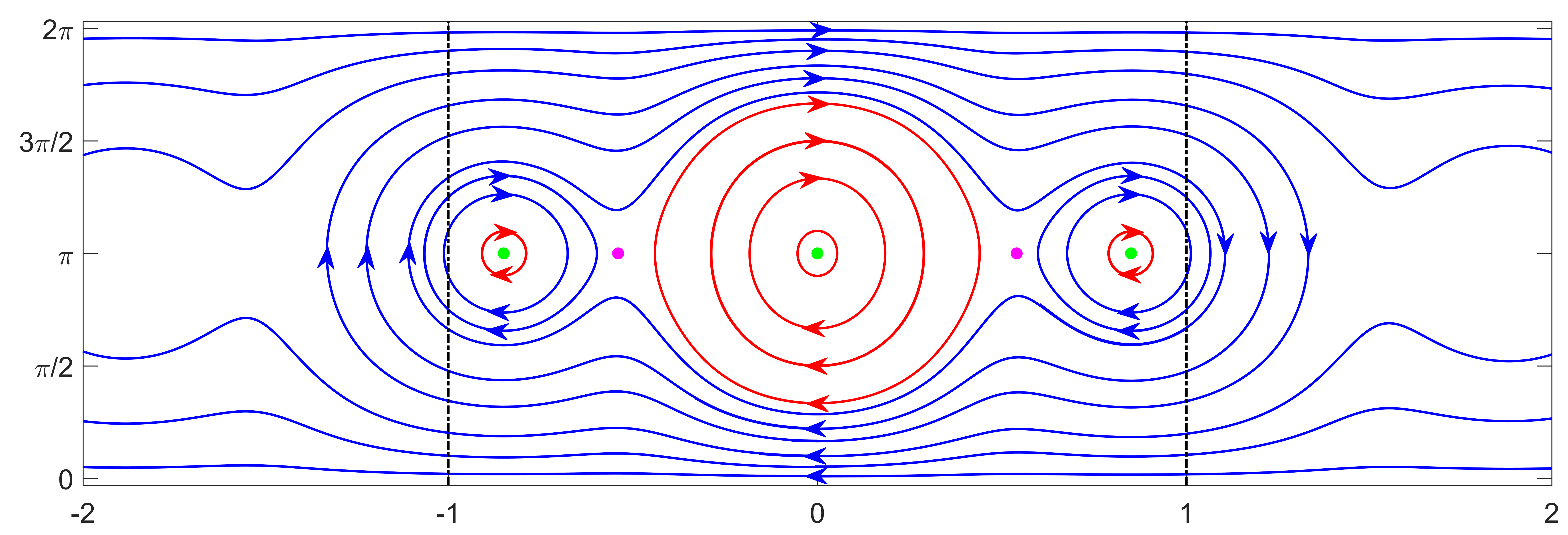}
	}
	
	\subfloat{
		\hspace{0.05cm}\includegraphics[scale=0.35]{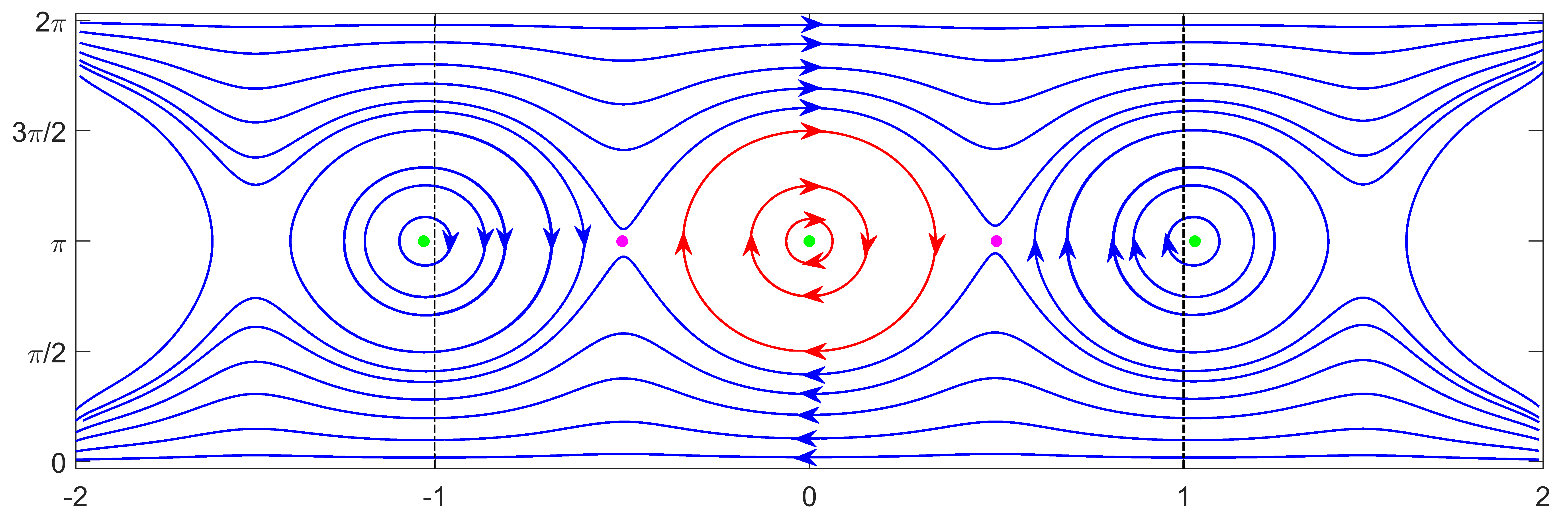}
	}	
	\caption{Phase portrait of the Hamiltonian system for the numerical wave equation and the grid transformations $\textbf{g}^{\,h,1}$ (top) and $\textbf{g}^{\,h,2}$ (bottom). We put $x^\pm(t)$, $\xi^\pm(t)$ on the horizontal/vertical direction.}\label{plot_ham_syst_var}
\end{figure}

To help us in our presentation, we include in Figure \ref{plot_ham_syst_var} the phase portraits corresponding to the solution of \eqref{wave_var_discr_1d} obtained through the employment of the non-uniform meshes $\textbf{g}^{\,h,1}$ and $\textbf{g}^{\,h,2}$ that we introduced before.

We can notice how the introduction of an oscillating variable coefficient changes the dynamics of these phase portraits by generating several equilibria (while, for the constant coefficient case, we had only one). Moreover, it can be easily seen through a linearization around those equilibria that the points marked in green are centers, while the other ones marked in pink are saddles. 

The phenomena that we already detected in Section \ref{numerical_sec} show up also in this case of a variable coefficient wave equation. Actually, our simulations turn out to be quite similar in several aspects to the ones that we displayed for the constant coefficients case. Therefore, in order not to extend our discussion more than necessary, in what follows we decided to include only the plots that we believe add something really new and interesting. Nonetheless, we will still describe all the high-frequency pathologies taking as a reference the phase portraits in Figure \ref{plot_ham_syst_var}.

\begin{figure}[!h]
	\centering 
	\begin{minipage}{0.3\textwidth}
		\includegraphics[scale=0.32]{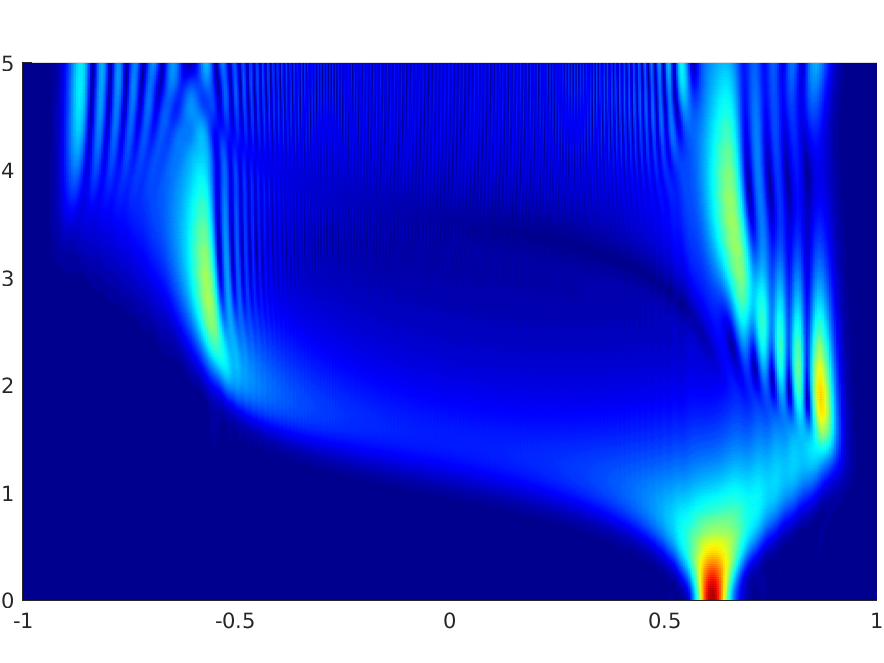}
	\end{minipage}
	\hspace{0.5cm}
	\begin{minipage}{0.3\textwidth}
		\includegraphics[scale=0.32]{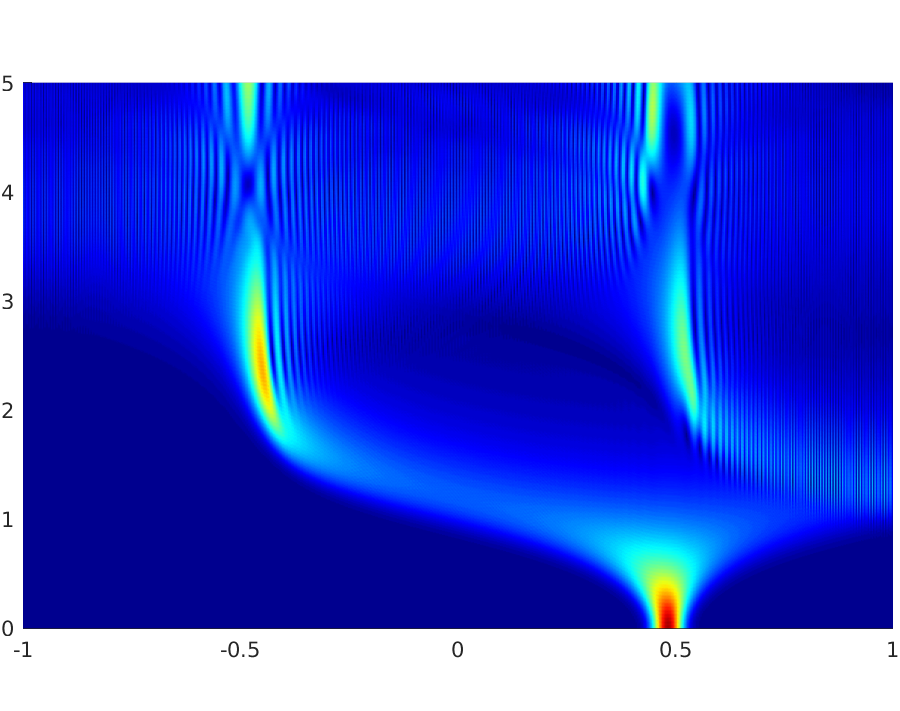}
	\end{minipage}
	\caption{Numerical solutions to \eqref{wave_var_1d} on a the non-uniform mesh $\textbf{g}^{\,h,1}$ (left) and $\textbf{g}^{\,h,2}$ (right), with $\sigma$ as in \eqref{var_coeff} with $A=1$ and $\kappa=1$ and with an initial datum corresponding to one unstable equilibrium (in pink) in the phase portrait.}\label{plot_wave_var_dispersive}
\end{figure}

First of all, we have that waves originated by an initial datum corresponding to one of the stable fixed points show a lack of propagation, similarly to what we already observed in Figure \ref{plot_02_wave}. The difference here is that, this time, for the mesh $\textbf{g}^{\,h,1}$ we have three stable fixed points inside the physical domain $(-1,1)$. In view of that, we have three different initial positions which, at frequency $\xi_0=\pi$, generate non propagating waves. On the other hands, initial data corresponding to one of the unstable fixed point produce solutions that, apart from showing absence of propagation, present also a huge dispersion (see Figure \ref{plot_wave_var_dispersive}). Moreover, we can notice differences with the nice symmetric shape of the dispersive wave observed in Figure \ref{plot_02_wave}. This because these solutions, as soon as they move away from the unstable equilibrium point, are quite immediately affected by the orbits around the stable ones, thus generating the comeback effects that can be appreciated in the plots.

As a second and last phenomenon, we can observe that around each one of the stable equilibria there are several trajectories (the red ones) which are completely contained in the central regions of the phase portraits, delimited by $x=\pm 1$. In view of that, solutions corresponding to these frequencies will show again the pathology of internal reflection, related to the fact that their group velocity vanishes before reaching the boundary of the physical domain.

\subsection{Conclusions}

Summarizing, our analysis shows several interesting phenomena in the propagation of numerical solutions of one-dimensional wave equations obtained under finite difference discretizations. 

At low frequencies, the discrete solutions behave similarly as the continuous one, propagating along the bi-characteristics with a positive velocity $\omega'(\xi)$. Nevertheless, the fact that this discrete velocity is given by a trigonometric function of the frequency yields pathological behaviors as $\xi$ increases. Already when employing a uniform mesh for our approximation, we see non-propagating waves, obtained when $\omega'(\xi)=0$. In addition, the introduction of non-uniform meshes generates also phenomena of internal reflection, in which the discrete waves remain trapped in some region of the domain, without the possibility to reach either one or both the boundary points.

Finally, the aforementioned phenomena can be observed both for a constant coefficient wave equation and when considering the case of variable coefficients (in the simulations that we presented, $\rho(x)\equiv 1$ and $\sigma(x) = 1+A\cos^2(\kappa\pi x)$, $A>0$, $\kappa\in\NN^*$).

\section{Two-dimensional wave equation}\label{2d_sec}

We briefly discuss here the two dimensional case. In more detail, the main model analyzed in this section is the following one
\begin{align}\label{wave_var_2d}
	\begin{cases}
		\rho(\Bf{z})\partial^2_tu-\textrm{div}_\Bf{z}\big(\sigma(\Bf{z})\nabla_\Bf{z} u\big) = 0, & (\Bf{z},t)\in \Omega\times(0,T)
		\\
		\left. u\right|_{\partial\Omega} = 0, & t\in(0,T)
		\\
		u(\Bf{z},0) = u^0(\Bf{z}),\;\;\;\partial_tu(\Bf{z},0) = u^1(\Bf{z}), & \Bf{z}\in\Omega,
	\end{cases}
\end{align}
where we indicate $\Bf{z}:=(x,y)$ and $\Omega:=(-1,1)^2$. Moreover, as for the one-dimensional equation \eqref{wave_var_1d}, we will consider $\rho,\sigma\in L^\infty(\Omega)$ with the hyperbolicity conditions $\rho(\Bf{z})\geq\rho^*>0$ and $\sigma(\Bf{z})\geq\sigma^*>0$. The total energy of the solutions written below is conserved in time, i.e. 
\begin{align}\label{energy_2d}
	\mathcal E_{\rho,\sigma}(u,\partial_t u) :=\frac 12 \int_{-1}^1\int_{-1}^1 \left(\rho(\Bf{z})|\partial_t u|^2 + \sigma(\Bf{z})|\nabla u|^2\right)\,d\Bf{z} = \mathcal E_{\rho,\sigma}(u^0,u^1).
\end{align}

Furthermore, we still have that the energy of solutions originated from highly oscillating and/or concentrated initial data propagates along the characteristic rays, obtained by solving the corresponding Hamiltonian system (\cite{ralston1982gaussian,rauch2005polynomial}). As we will see from our simulation, this concentration effect is preserved also at the numerical level. Nevertheless, we will observe that the discrete solution to \eqref{wave_var_2d} not always propagates as one would expect from the continuous model. In certain cases, instead, similar phenomena to the ones that we already observed for the one-dimensional equation appear. 

In what follows, we describe the numerical finite difference scheme that we shall employ for approximating \eqref{wave_var_2d}, and we present the results of our simulations. Since there will be many similarities with the one-dimensional analysis of Section \ref{1d_sec}, we will keep our discussion as short as possible, giving only the essential details.

\subsection{Finite difference approximation of \eqref{wave_var_2d} and Hamiltonian system}\label{2d_fd_sec} 

Let us introduce the grid and the finite difference discretization of the 2-d wave equation \eqref{wave_var_2d}. In the simulations that we will present in the next section, we will focus only on the case of constant coefficients $\rho=\sigma\equiv 1$. Nevertheless, 
for the sake of completeness in what follows we are going to consider the general case  of variable coefficients. 

Let $g_1,g_2:\Omega\to\Omega$ be two diffeomorphisms of the domain $\Omega$, and set $\Bf{g}:=(g_1,g_2)$. For $M,N\in\NN^*$ given, let $h_x:= 2/(M+1)$ and $h_y:=2/(N+1)$, and set $\Bf{h} = (h_x,h_y)$. Consider the uniform grid of $\Omega$
\begin{align*}
	\Bf{G}^h :=\left\{\Bf{z}_{j,k}:=(x_j,y_k) = (-1+jh_x,-1+kh_y), j=0,\ldots,M+1, k=0,\ldots,N+1\right\},
\end{align*}
and the non-uniform one 
\begin{align*}
	\Bf{G}^h_\Bf{g} :=\left\{\Bs{\omega}_{j,k}:=(\upsilon_j,\zeta_k)=(g_1(x_j),g_2(y_k))\right\},
\end{align*}
obtained by transforming $\Bf{G}^h$ through the map $\Bf{g}$. 

We denote $\Bs{\omega}:=(\upsilon,\zeta)$, and we set $\upsilon_{j+1/2} = g_1(x_{j+1/2})$ and $\zeta_{k+1/2} = g_2(y_{k+1/2})$, with $x_{j+1/2} = -1+(j + 1/2)h_x$ and $y_{k+1/2} = -1+(k + 1/2)h_y.$ The finite difference approximation of \eqref{wave_var_2d} is then given in the following form:
\begin{align}\label{wave_var_2d_discr}
	\begin{cases}
		\rho_{j,k}\partial^2_tu_{j,k}-\textrm{div}_{\Bs{\omega}}^{\Bf{g},\Bf{h}}\big(\sigma\nabla_{\Bs{\omega}}^{\Bf{g},\Bf{h}} u\big)_{j,k} = 0, & j=1,\ldots,M,\; k=1,\ldots,N,\;t\in(0,T)
		\\
		u_{0,0} = u_{M+1,N+1} = 0, & t\in(0,T)
		\\
		u_{j,k}(0) = u_{j,k}^0,\;\;\;\partial_tu_{j,k}(0) = u^1_{j,k}, & j=1,\ldots,M,\; k=1,\ldots,N.
	\end{cases}
\end{align}

Here, $\rho_{j,k}:=\rho(\upsilon_j,\zeta_k)$, and $\textrm{div}_{\Bs{\omega}}^{\Bf{g},\Bf{h}}\big(\sigma\nabla_{\Bs{\omega}}^{\Bf{g},\Bf{h}}\cdot \big)$ is the five-points finite difference approximation of the $\textrm{div}_{\Bf{z}}(\sigma\nabla_{\Bf{z}}\cdot)$-operator on the non-uniform grid $\Bf{G}^h_\Bf{g}$, which is defined as
\begin{align*}
	\textrm{div}_{\Bs{\omega}}^{\Bf{g},\Bf{h}}\big(\sigma\nabla_{\Bs{\omega}}^{\Bf{g},\Bf{h}}u\big)_{j,k}:= & \,\frac{\sigma_{j+1/2,\,k}\frac{u_{j+1,\,k}-u_{j,\,k}}{\upsilon_{j+1}-\upsilon_j}-\sigma_{j-1/2,\,k}\frac{u_{j,\,k}-u_{j-1,\,k}}{\upsilon_j-\upsilon_{j-1}}}{\frac{\upsilon_{j+1}-\upsilon_{j-1}}{2}} 
	\\
	&+ \frac{\sigma_{j,\,k+1/2}\frac{u_{j,\,k+1}-u_{j,\,k}}{\zeta_{k+1}-\zeta_k}-\sigma_{j,\,k-1/2}\frac{u_{j,\,k}-u_{j,\,k-1}}{\zeta_k-\zeta_{k-1}}}{\frac{\zeta_{k+1}-\zeta_{k-1}}{2}}.
\end{align*}

Also, $\sigma_{j\pm 1/2,\,k}:=\sigma(\upsilon_{j\pm 1/2},\zeta_k)$ and $\sigma_{j,\,k\pm 1/2}:=\sigma(\upsilon_j,\zeta_{k\pm 1/2})$. The total energy of the solutions of \eqref{wave_var_2d_discr} is conserved in time:
\begin{align*}
	\mathcal E_{\rho,\sigma}^{\Bf{g},\Bf{h}}(\Bf{u}^{\Bf{h}}(t), &\partial_t \Bf{u}^{\Bf{h}}(t))
	\\
	:= &\, \frac 18 \sum_{j=1}^N \sum_{k=1}^M \rho_{j,k}(\upsilon_{j+1}-\upsilon_j)(\zeta_{k+1}-\zeta_k)|\partial_t u_{j,k}|^2 
	\\
	&+ \frac 14 \sum_{j=1}^N \sum_{k=1}^M \left(\frac{\zeta_{k+1}-\zeta_k}{\upsilon_{j+1}-\upsilon_j}\sigma_{j+1/2,\,k}\,|u_{j+1,k}-u_{j,k}|^2 + \frac{\zeta_{j+1}-\zeta_j}{\upsilon_{k+1}-\upsilon_k}\sigma_{j,\,k+1/2}\,|u_{j,k+1}-u_{j,k}|^2\right).
\end{align*}

Moreover, as we mentioned, solutions originated from highly oscillating and/or concentrated initial data propagates along the characteristic rays, and these localized solution can be employed for studying the motion of the discrete waves. Hence, before presenting our numerical results, let us quickly discuss the Hamiltonian system.

Following the Wigner measure approach, the principal symbol associated to the finite difference scheme \eqref{wave_var_2d_discr} is (see, e.g., \cite{macia2002propagacion})
\begin{align}\label{symbol_2d_discr}
	\mathcal P(x,y,t,\xi,\eta,\tau) := \tau^2 - \Lambda (x,y,\xi,\eta)
\end{align}
with 
\begin{align}\label{symbol_Lambda}
	\Lambda (x,y,\xi,\eta):= \frac{\sigma(x,y)}{\rho(x,y)}\left(4\sin^2\left(\frac \xi2\right)\frac{1}{g_1'(x)^2} + 4\sin^2\left(\frac \eta2\right)\frac{1}{g_2'(y)^2}\right).
\end{align} 

We note that, once again, this symbol depends on the space variable  $\Bf{z} = \Bf{g}^{-1}(\Bs{\omega})$ corresponding to the uniform grid. As a consequence, the variable coefficients $\rho$ and $\sigma$ have to be composed with $\Bf{g}$. Note also the appearance of the factors $1/g_1'(x)$ and $1/g_2'(x)$ accompanying each space derivative, which is also due to the grid transformation. Moreover, note that the Fourier symbols $\xi^2$ and $\eta^2$ of second-order space derivative have been replaced by the corresponding symbols $4\sin^2(\xi/2)$ and $4\sin^2(\eta/2)$ of the five-points finite difference approximation of the Laplacian.

Set $\Bs{\omega}_e:= (\upsilon,\zeta,t)$, $\boldsymbol{\theta}:= (\xi,\eta)$ and $\boldsymbol{\theta}_e:= (\xi,\eta,\tau)$ (the subscript $e$ stands for \textit{extended}). The null bi-characteristic rays $(\Bf{z}_e(s),\Bs{\theta}_e(s))$ are solutions of the following Hamiltonian system of non-linear ODEs:
\begin{align}\label{ham_syst_2d_discr}
	\begin{cases}
		\dot{\Bf{z}}_e(s) = \nabla_{\Bs{\theta}_e}\mathcal P(\Bf{z}_e(s),\Bs{\theta}_e(s))
		\\
		\dot{\Bs{\theta}}_e(s) = -\nabla_{\Bf{z}_e}\mathcal P(\Bf{z}_e(s),\Bs{\theta}_e(s)),
	\end{cases}
\end{align}
with initial conditions $\Bf{z}_e(0)=\Bf{z}_e^0:=(x_0,y_0,t_0)$ and $\Bs{\theta}_e(0)=\Bs{\theta}_e^0:=(\xi_0,\eta_0,\tau_0)$. 

Note that the principal symbol is conserved along these null bi-characteristic rays (i.e., $\mathcal P(\Bf{z}_e(s),\Bs{\theta}_e(s)) = \mathcal P(\Bf{z}_e^0,\Bs{\theta}_e^0)=0$ and that the $\tau$ component of the Hamiltonian system \eqref{ham_syst_2d_discr} does not depend on $s$ (since the principal symbol $\mathcal P$ does not depend on $t$). Then, there are two solutions $\tau_0$ of the equation $\mathcal P(\Bf{z}(s),t(s),\boldsymbol{\theta}(s),\tau_0)=0$, $\tau_0^\pm:=\pm\sqrt{\Lambda(\Bf{z}(s),\Bs{\theta}(s))}$, and, correspondingly, two families of solutions of \eqref{ham_syst_2d_discr}. By considering $\Bf{z}$ and $\Bs{\theta}$ as functions of the time variable $t$ instead of $s$ (this can be done since $dt/ds$ does not vanish), the two families of characteristic rays $(\Bf{z}^\pm(t),\Bs{\theta}^\pm(t))$ solve the Hamiltonian system
\begin{align}\label{ham_syst_2d_discr_t}
	\begin{cases}
		\dot{\Bf{z}}^\pm(t) = \pm\nabla_{\Bs{\theta}}\sqrt{\Lambda(\Bf{z}^\pm(t),\Bs{\theta}^\pm(t))}
		\\
		\dot{\Bs{\theta}}^\pm(t) = \mp\nabla_{\Bf{z}}\sqrt{\Lambda(\Bf{z}^\pm(t),\Bs{\theta}^\pm(t))}.
	\end{cases}
\end{align}

When $\sigma/\rho$ is constant (assume $\sigma = \rho\equiv 1$ for simplicity), then \eqref{ham_syst_2d_discr_t} can be decoupled into two Hamiltonian systems corresponding to the variables $(x,\xi)$ and $(y,\eta)$. Indeed, note firstly that the following quantities are conserved in time along the characteristics:
\begin{align*}
	r_0:=\sqrt{\Lambda(\Bf{z}^\pm(t),\Bs{\theta}^\pm(t))},\;\;\; r_1:=\lambda_1(x^\pm(t),\xi^\pm(t)),\;\;\; r_2:=\lambda_2(y^\pm(t),\eta^\pm(t)),
\end{align*}
where
\begin{align}
	\lambda_1(x,\xi):=2\sin\left(\frac \xi2\right)\frac{1}{g_1'(x)}\;\;\;\textrm{ and } \;\;\; \lambda_2(y,\eta):=2\sin\left(\frac \eta2\right)\frac{1}{g_2'(y)}.
\end{align}
The two Hamiltonian systems corresponding to each direction are as follows:
\begin{align*}
	\begin{cases}
		\displaystyle\dot{x}^\pm(t) = \pm \frac{r_1}{r_0}\partial_\xi\lambda_1(x^\pm(t),\xi^\pm(t))
		\\[8pt]
		\displaystyle\dot{\xi}^\pm(t) = \pm \frac{r_1}{r_0}\partial_x\lambda_1(x^\pm(t),\xi^\pm(t))
	\end{cases}
	\;\;\;\textrm{ and }\;\;\;
	\begin{cases}
		\displaystyle\dot{y}^\pm(t) = \pm \frac{r_2}{r_0}\partial_\eta\lambda_2(y^\pm(t),\eta^\pm(t))
		\\[8pt]
		\displaystyle\dot{\eta}^\pm(t) = \pm \frac{r_2}{r_0}\partial_{y}\lambda_2(y^\pm(t),\eta^\pm(t)).
	\end{cases}
\end{align*}
Then, the original variables $(x:= g_1^{-1}(\upsilon),\xi)$ and $(y:= g_2^{-1}(\zeta),\xi)$ satisfy the following ODE systems:
\begin{align}\label{ham_syst_2d_x}
	\begin{cases}
		\displaystyle\dot{x}^\pm(t) = \mp \frac{r_1}{r_0}g_1'(g_1^{-1}(x^\pm(t)))\partial_\xi\lambda_1(g_1^{-1}(x^\pm(t)),\xi^\pm(t))
		\\[8pt]
		\displaystyle\dot{\xi}^\pm(t) = \mp \frac{r_1}{r_0}\partial_x\lambda_1(g_1^{-1}(x^\pm(t)),\xi^\pm(t))
	\end{cases}
\end{align}
and
\begin{align}\label{ham_syst_2d_y}
	\begin{cases}
		\displaystyle\dot{y}^\pm(t) = \mp \frac{r_2}{r_0}g_2'(g_2^{-1}(y^\pm(t)))\partial_\eta\lambda_2(g_2^{-1}(y^\pm(t)),\eta^\pm(t))
		\\[8pt]
		\displaystyle\dot{\eta}^\pm(t) = \mp \frac{r_2}{r_0}\partial_y\lambda_2(g_2^{-1}(y^\pm(t)),\eta^\pm(t)).
	\end{cases}
\end{align}

In \eqref{ham_syst_2d_x} and \eqref{ham_syst_2d_y}, $\partial_\xi \lambda_1$ and $\partial_\eta \lambda_2$ are the two components of the group velocity. They describe the speed at which solutions associated with wave number $(\xi,\eta)$ move in the corresponding direction. Notice that
\begin{align*}
	\partial_\xi \lambda_1(x,\xi) = \cos\left(\frac \xi2\right)\frac{1}{g_1'(x)} = 0\;\;\; \textrm{ iff } \;\;\;\xi = (2k+1)\pi,\;\;\;k\in\ZZ
	\\
	\partial_\eta \lambda_2(y,\eta) = \cos\left(\frac \eta2\right)\frac{1}{g_2'(y)} = 0\;\;\; \textrm{ iff } \;\;\;\eta = (2k+1)\pi,\;\;\;k\in\ZZ,
\end{align*} 
and that this is independent on the choice of $g_1$ and $g_2$. Therefore, no matter what mesh we select for our discretization, there are certain frequencies $(\xi,\eta)$ at which the group velocity vanishes in at least one component, thus producing a lack of propagation of the wave in the corresponding direction. This fact will be pointed out by our simulations.

\subsection{Numerical results}

\begin{figure}[!h]
	\centering
	\begin{minipage}{0.42\textwidth}
		\includegraphics[scale=0.55]{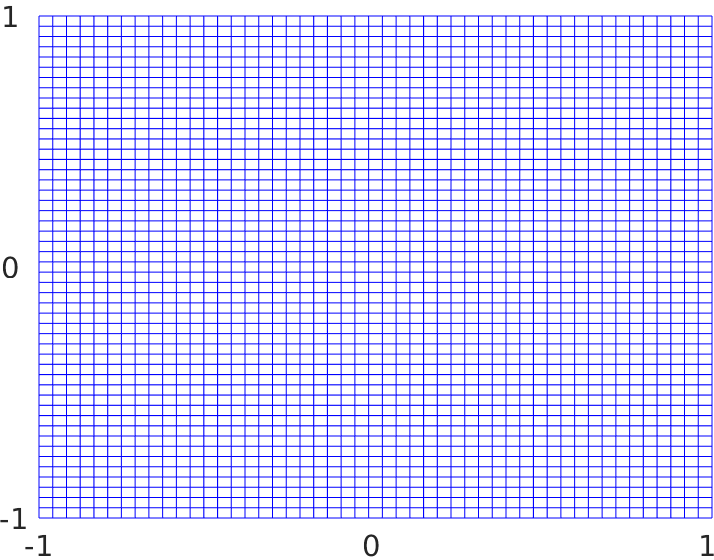}
	\end{minipage}
	\hspace{0.7cm}
	\begin{minipage}{0.42\textwidth}
		\includegraphics[scale=0.515]{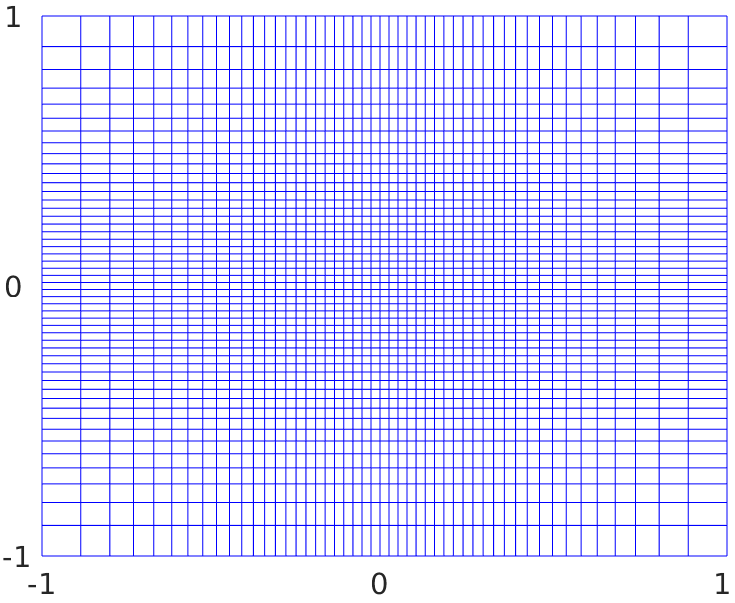}
	\end{minipage}	
	\caption{Uniform mesh on $\Omega$ and its refinement through the function $\Bf{g}$.}\label{mesh_2d_fig}
\end{figure}

We present and discuss here several simulations for the two-dimensional wave equation \eqref{wave_var_2d_discr}. In what follows, analogously to the one-dimensional case that we discussed before, we are going to consider as mesh functions
\begin{align}\label{mesh_fun_2d}
	g_1(x)=g_2(x)=\tan\left(\frac \pi4 x\right)=:g(x),
\end{align}
yielding a gradual refinement of the grid around the point $(0,0)$ (see Figure \ref{mesh_2d_fig}). 
Moreover, as we were mentioning before, we focus here only on the case $\rho=\sigma\equiv 1$. 

For the numerical resolution of our equation, instead of the leapfrog scheme that we used in the one-dimensional case, we are going to follow a different approach. Taking advantage of the fact that the equation has constant coefficients, and that both the domain and the mesh are symmetric in the $x$ and in the $y$ variables, we are going to compute the solution in Fourier series, that is, 
\begin{align}\label{sol_fourier}
	\Bf{u}^\Bf{h} = \sum_{j=1}^M\sum_{k=1}^N \beta_{j,k}\Bs{\Phi}_{j,k}(\Bs{\omega})\,e^{\,it\sqrt{\lambda_{j,k}}}.
\end{align}

In \eqref{sol_fourier}, $\{\Bs{\Phi}_{j,k},\lambda_{j,k}\}$ are the eigenvector and the eigenvalues of the discrete Laplacian $-\Delta_{\Bs{\omega}}$ on the refined mesh $\Bf{G}_{\Bf{g}}^h$, i.e., 
\begin{align*}
	-\Delta_{\Bs{\omega}} \Bs{\Phi}_{j,k} = \lambda_{j,k} \Bs{\Phi}_{j,k}, \;\;\; j=1,\ldots,M,\;\; k=1,\ldots,N,
\end{align*}
while $\beta_{j,k}$ are the corresponding Fourier coefficients of the initial datum $\Bf{u}^{0,\Bf{h}}$. Moreover, due to symmetry reasons, we have that the eigenvectors $\Bs{\Phi}_{j,k}(\Bs{\omega})$ are, actually, in separated variables and that the eigenvalues $\lambda_{j,k}$ are given by the sum of the eigenvalues of the corresponding 1-D problems in the $\upsilon$ and in the $\zeta$ directions. In other words, we have $\Bs{\Phi}_{j,k}(\Bs{\omega})=\Bs{\Psi}_j(\upsilon)\Bs{\Upsilon}_k(\zeta)$ and $\lambda_{j,k} = \mu_j+\nu_k$, with 
\begin{align*}	
	-\Delta_{\upsilon} \Bs{\Psi}_j = \mu_j \Bs{\Psi}_j \;\;\; \textrm{ and }\;\;\; -\Delta_{\zeta} \Bs{\Upsilon}_k = \nu_k \Bs{\Upsilon}_k,\;\;\; j=1,\ldots,M,\;\; k=1,\ldots,N.
\end{align*}

We stress that this approach that we just described can be adopted since we are limiting our analysis only to a very particular case (constant coefficients and square domain). If one would treat, instead, the variable coefficients wave equation on $(-1,1)^2$, the discretization shall be done, for instance, joining the scheme that we presented in Section \ref{2d_fd_sec} and a leapfrog method for the time integration.

We present below several simulations obtained with the methodology just described. We considered a non-uniform mesh as in Figure \ref{mesh_2d_fig}, with an equal number of points in both directions (i.e. $M=N$) and, for simplicity we denote $h_x=h=h_y$. Moreover, this time we show the plots in the space domain $(-1,1)^2$ (in the one dimensional case, we were showing the space-time domain $(-1,1)\times(0,T)$). We will then indicate in each case the time horizon of our simulations.

As for the one-dimensional case before, we choose an initial datum given by a Gaussian wave packet concentrated at $(x_0,y_0)$ and oscillating at the wave number $(\xi_0/h,\eta_0/h)$, namely
\begin{align}\label{in_data_2d}
	u^0(x,y) = \exp\bigg[-\gamma\left((x-x_0)^2+(y-y_0)^2\right)\bigg]\exp\left[\,i\left(\frac{x\xi_0}{h} + \frac{y\eta_0}{h}\right)\right],
\end{align}
where $\gamma:=h^{-0.9}$. Moreover, we will observe several analogies with what we showed in Section \ref{numerical_sec}. 

\begin{figure}[!h]
	\centering
	\begin{minipage}{0.3\textwidth}
		\includegraphics[scale=0.16]{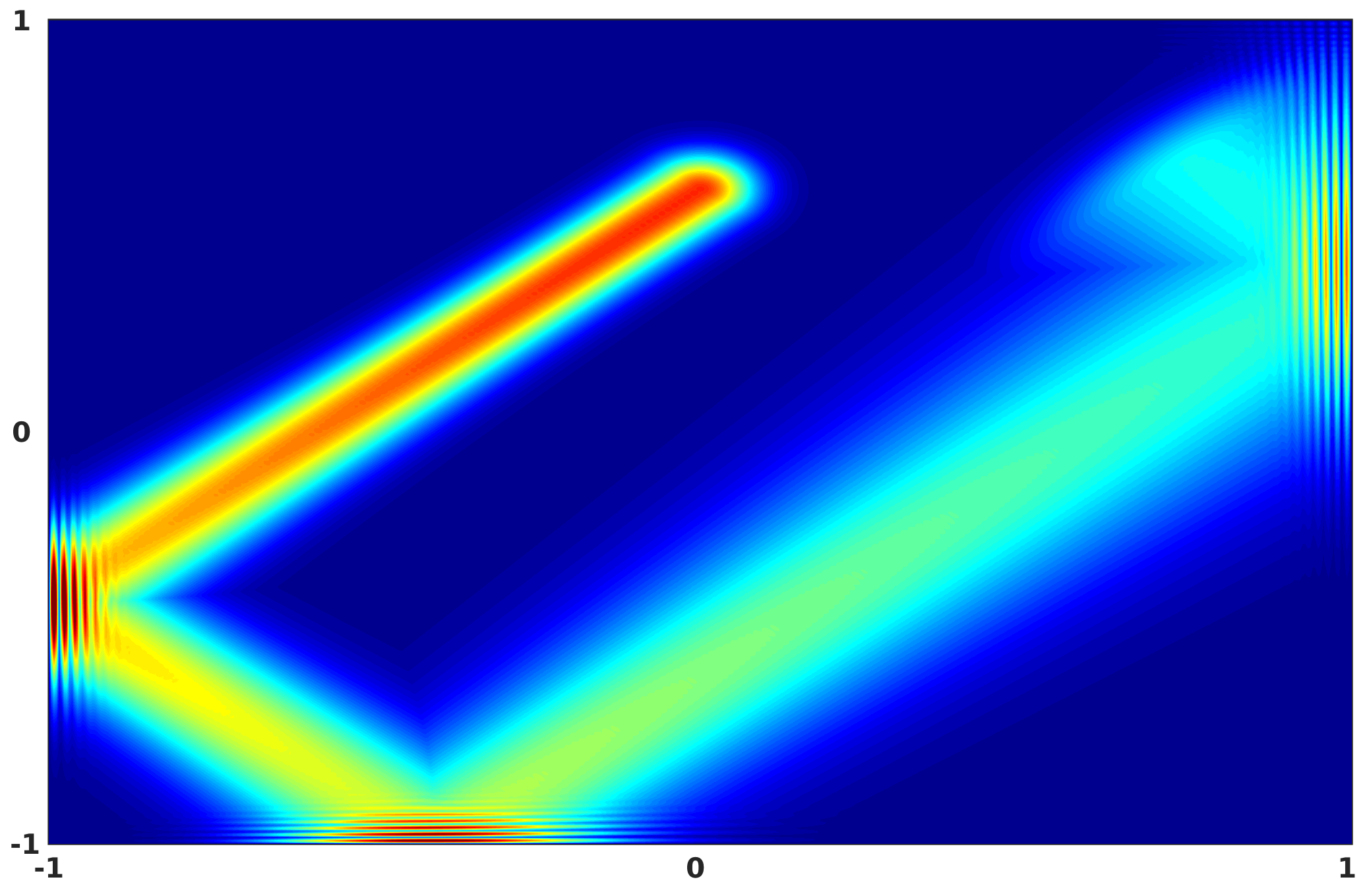}
	\end{minipage}	
	\hspace{1.5cm}
	\begin{minipage}{0.3\textwidth}
		\includegraphics[scale=0.16]{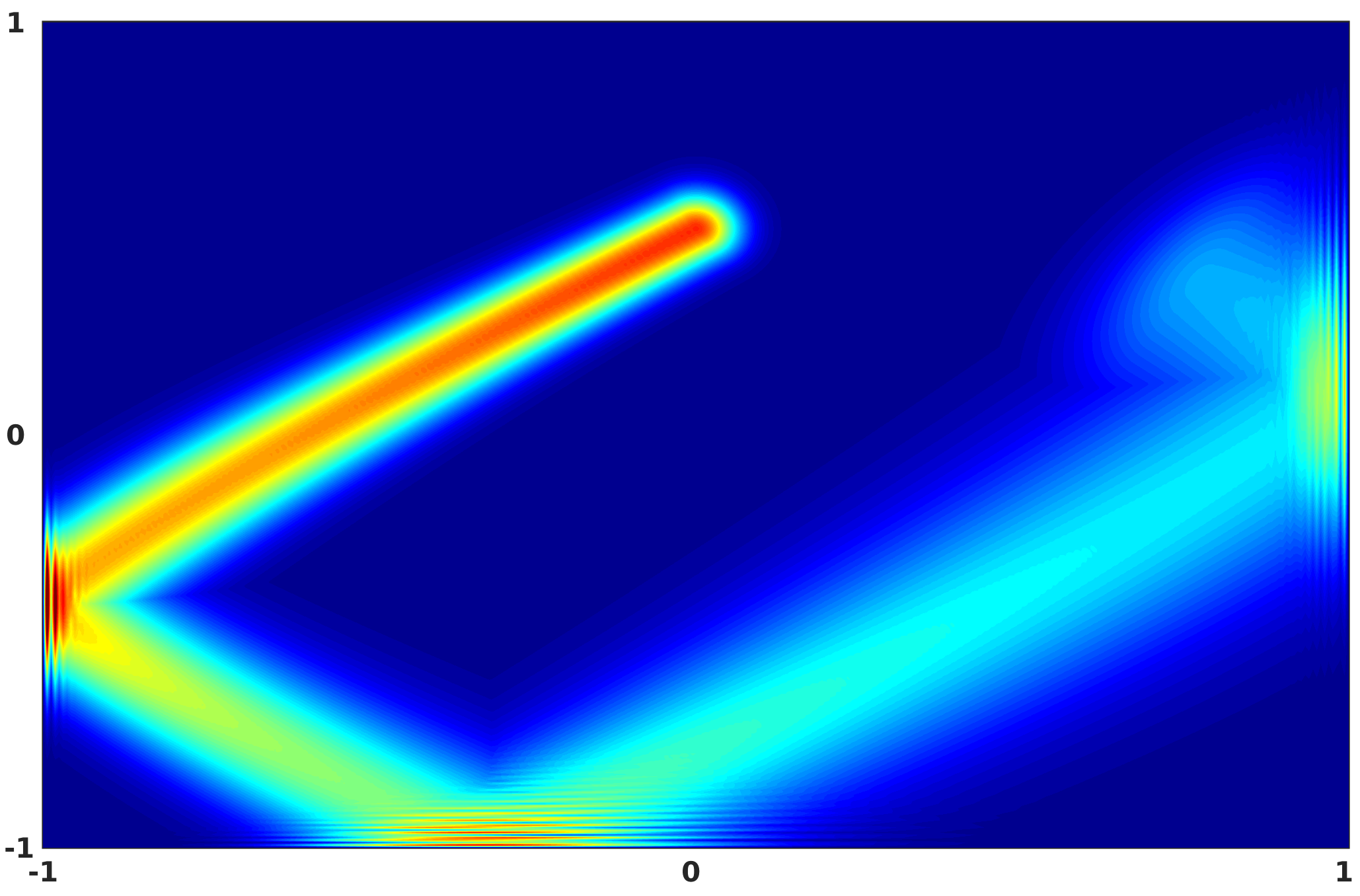}
	\end{minipage}
	\caption{Numerical solutions with initial datum \eqref{in_data_2d} and parameters $(x_0,y_0,\xi_0,\eta_0) = (0,1/2,\pi/4,\pi/4)$. The discretization is done on a uniform mesh (left) and on a non-uniform one obtained through the mesh function $\Bf{g}$ (right). The time horizon is $T=5s$ in both cases.}\label{plot_wave2d_low1}
\end{figure}

It can be observed in Figure \ref{plot_wave2d_low1} that, at low frequencies, the solution remains concentrated and propagates along straight characteristics which reach the boundary, where there is reflection according to the Descartes-Snell's law. This independently on whether we use a uniform or a non-uniform mesh. 

\begin{figure}[!h]
	\centering
	\begin{minipage}{0.3\textwidth}
		\includegraphics[scale=0.16]{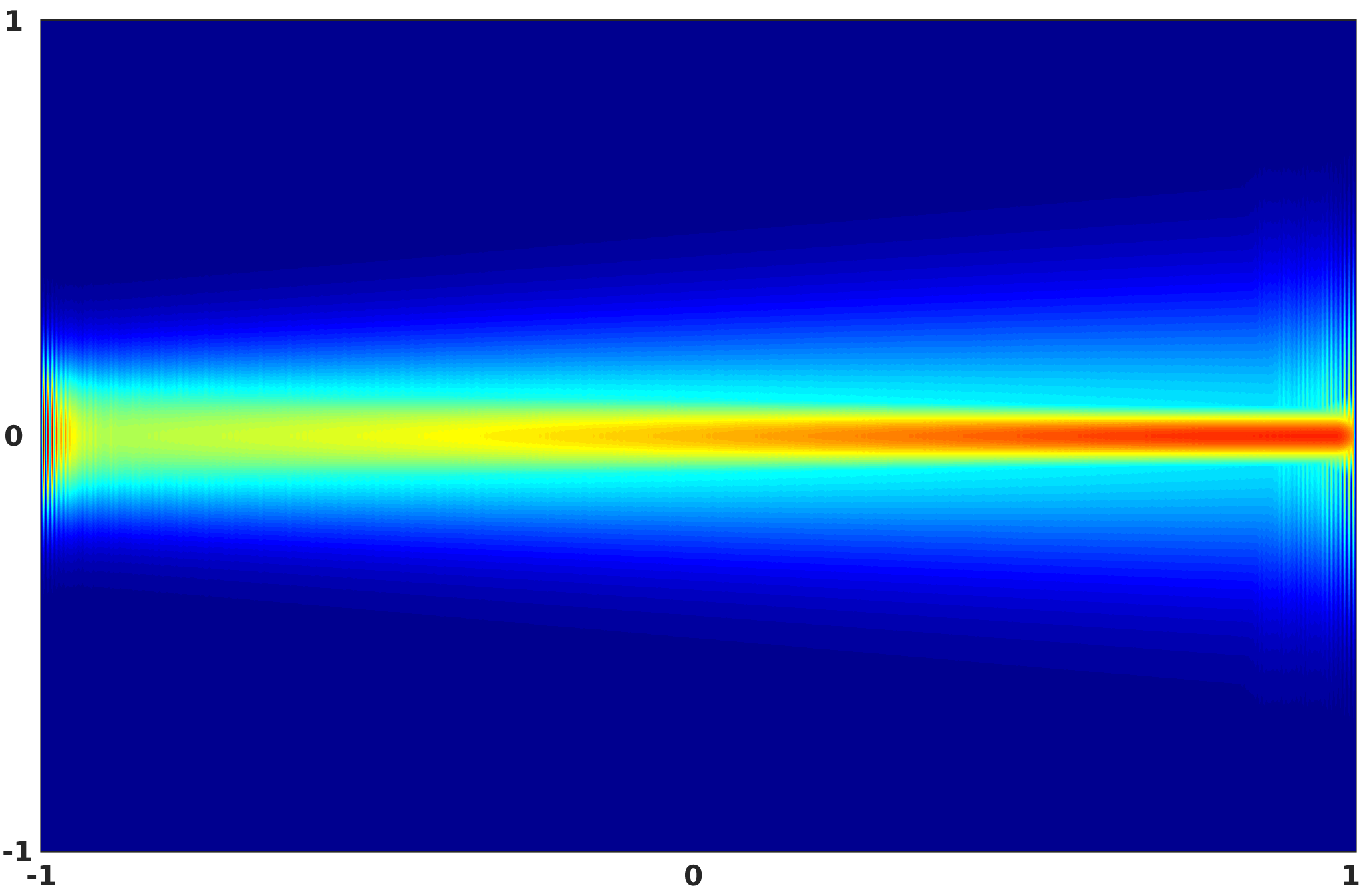}
	\end{minipage}	
	\hspace{1.5cm}
	\begin{minipage}{0.3\textwidth}
		\includegraphics[scale=0.16]{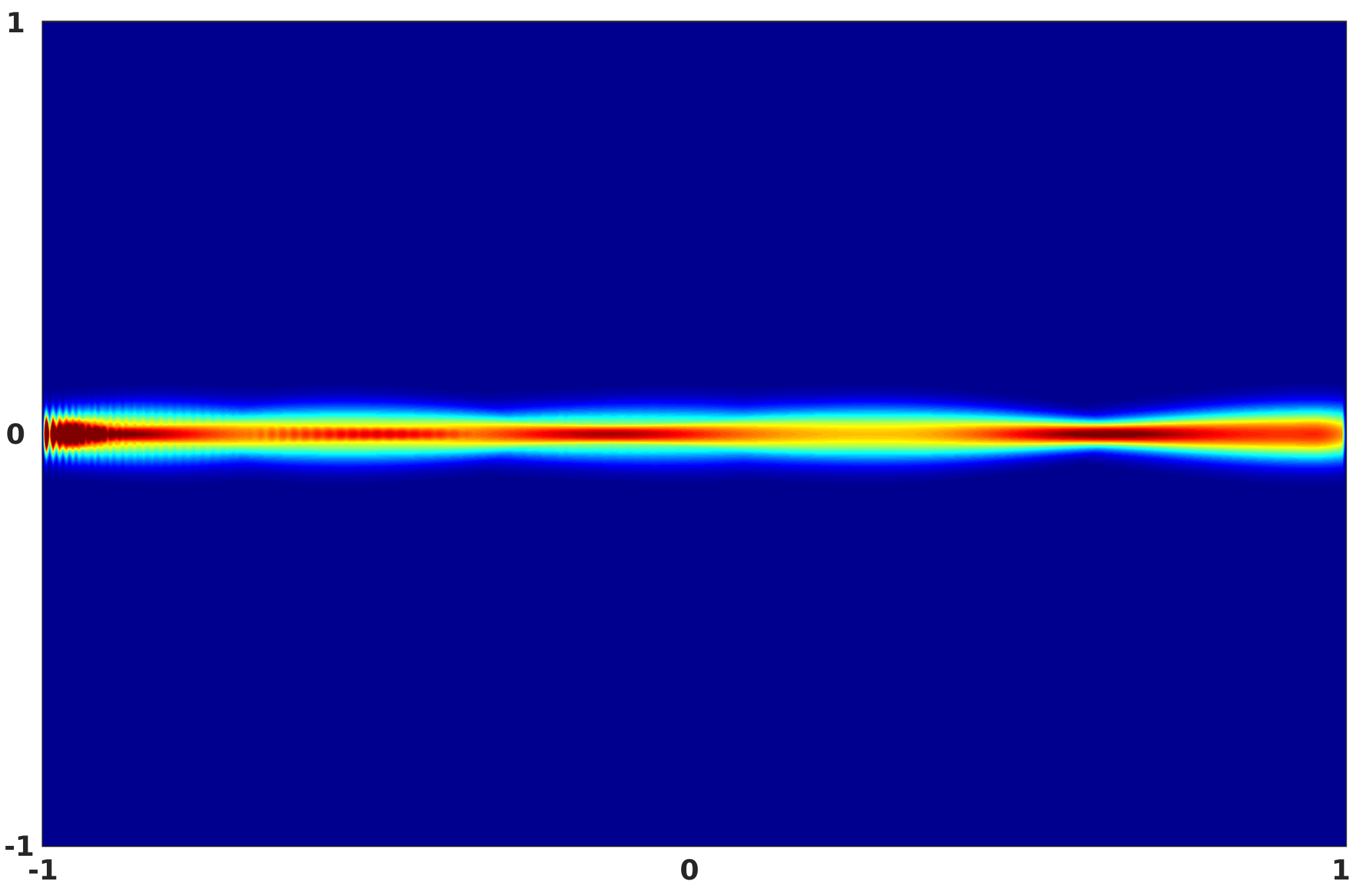}
	\end{minipage}
	\caption{Numerical solutions with initial datum \eqref{in_data_2d} and parameters $(x_0,y_0,\xi_0,\eta_0) = (1,0,\pi/2,\pi)$. The discretization is done on a uniform mesh (left) and on a non-uniform one obtained through the mesh function $\Bf{g}$ (right). The time horizon is $T=10s$ in both cases.}\label{plot_wave2d_np1}
\end{figure}

\begin{figure}[!h]
	\centering
	\begin{minipage}{0.3\textwidth}
		\includegraphics[scale=0.16]{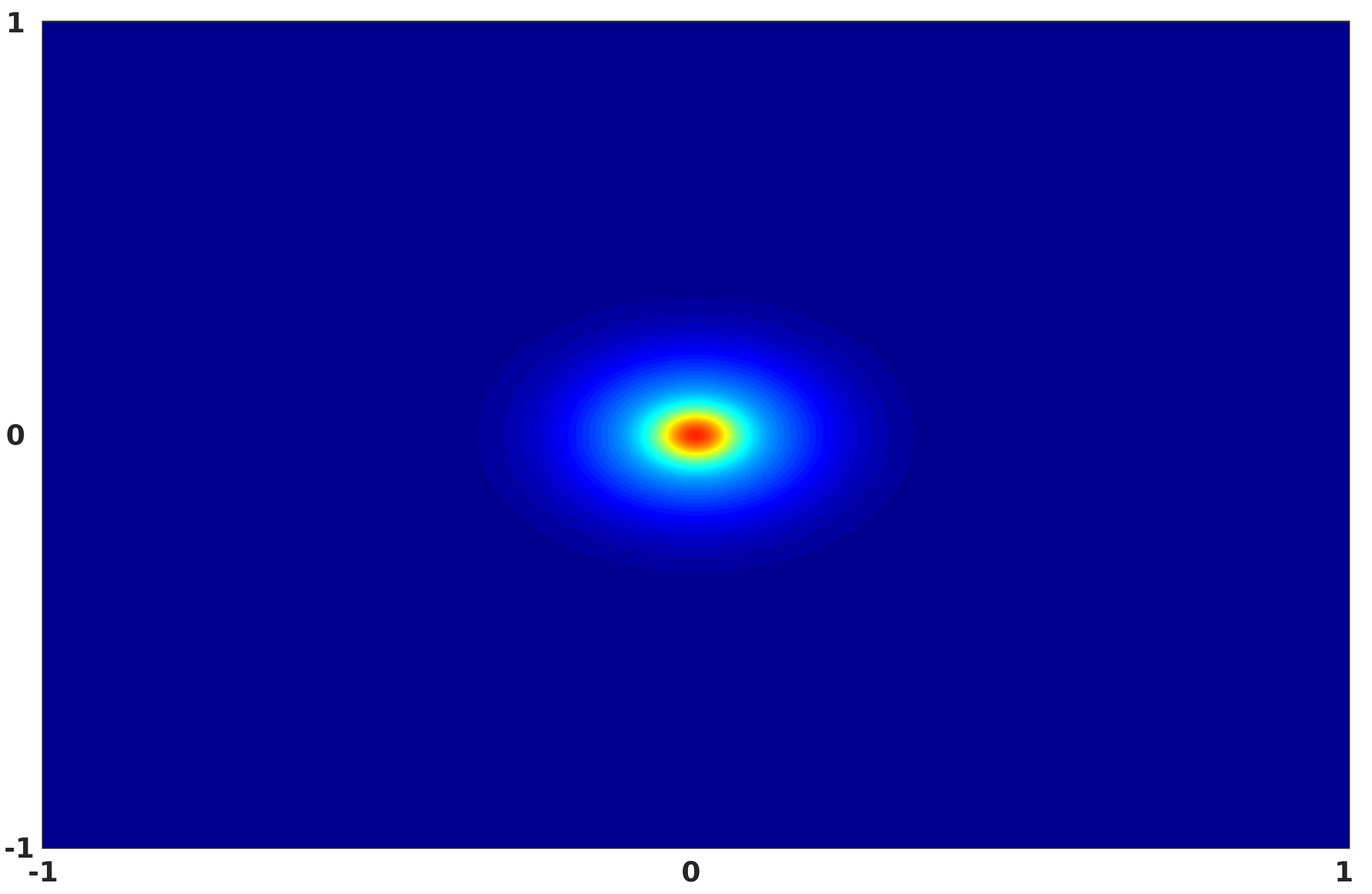}
	\end{minipage}	
	\hspace{1.5cm}
	\begin{minipage}{0.3\textwidth}
		\includegraphics[scale=0.16]{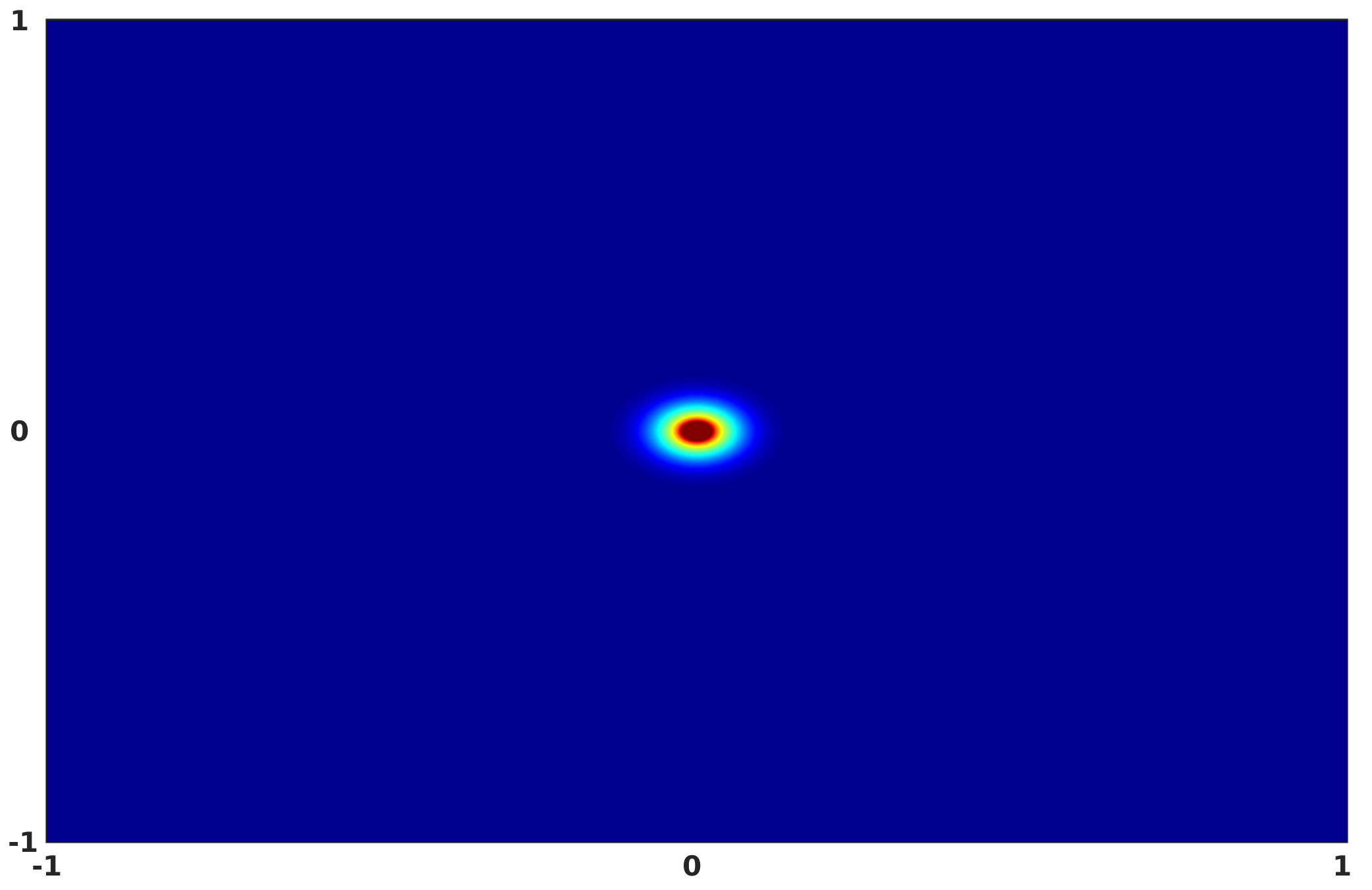}
	\end{minipage}
	\caption{Numerical solutions with initial datum \eqref{in_data_2d} and parameters $(x_0,y_0,\xi_0,\eta_0) = (0,0,\pi,\pi)$. The discretization is done on a uniform mesh (left) and on a non-uniform one obtained through the mesh function $\Bf{g}$ (right). The time horizon is $T=10s$ in both cases.}\label{plot_wave2d_np3}
\end{figure}

Nevertheless, increasing the frequencies similar phenomena as in the one-dimensional case show up. First of all, in Figures \ref{plot_wave2d_np1} and \ref{plot_wave2d_np3} we observe again waves that do not propagate. 
The justification to this fact is that, for the frequency considered there, we have that either $\partial_\xi\lambda_1$ or $\partial_\eta\lambda_2$ (or even both) vanishes, i.e. the velocity of propagation of the rays is zero in one or both spatial direction. Once again, this fact is not related to the particular mesh that we are employing. On the other hand, it is due to the changes in the Hamiltonian when passing from the continuous to the discrete setting, and to the trigonometric nature of the discrete velocity.

In Figure \ref{plot_wave2d_np1}, we considered an initial datum with parameters such that $\partial_\xi\lambda_2$ is zero, and this produces a loss o propagation in the vertical direction. Therefore, the wave remains trapped bouncing between the two sides $x=-1$ and $x=1$. 

Specularly, if one chooses parameters that annul $\partial_\xi\lambda_1$, the resulting wave shows no propagation in the horizontal direction. Finally, in Figure \ref{plot_wave2d_np3}, we considered an initial datum with parameters such that both $\partial_\xi\lambda_1$ and $\partial_\xi\lambda_2$ vanish. In view of that, a wave starting from such initial datum cannot move, and remains trapped around the point $(x,y)=(0,0)$ for any time.

This phenomena has been already observed and discussed in \cite[Chapter 4]{macia2002lack}, using the approach of Wigner measures. Furthermore, this lack of propagation finds explanation also in the analysis of the phase portrait associated to \eqref{ham_syst_2d_discr}. Indeed, consider for instance a discretization on non-uniform mesh obtained through the function in \eqref{mesh_fun_2d}. It is easily seen that, in this particular case, the Hamiltonian systems for the rays $(x^\pm(t),\xi^\pm(t))$ and $(y^\pm(t),\eta^\pm(t))$ become, respectively,
\begin{align}\label{ham_syst_2d_eq1}
	\begin{cases}
		\displaystyle\dot{x}^\pm(t) = \mp \frac{4}{r_0\pi}\sin(\xi^\pm(t))\frac{1}{x^\pm(t)^2+1}
		\\[8pt]
		\displaystyle\dot{\xi}^\pm(t) = \mp \frac{32}{r_0\pi}\sin^2\left(\frac{\xi^\pm(t)}{2}\right)\frac{x^\pm(t)}{(x^\pm(t)^2+1)^2}
	\end{cases}
\end{align}
and
\begin{align}\label{ham_syst_2d_eq2}
	\begin{cases}
		\displaystyle\dot{y}^\pm(t) = \mp \frac{4}{r_0\pi}\sin(\eta^\pm(t))\frac{1}{y^\pm(t)^2+1}
		\\[8pt]
		\displaystyle\dot{\eta}^\pm(t) = \mp \frac{32}{r_0\pi}\sin^2\left(\frac{\eta^\pm(t)}{2}\right)\frac{y^\pm(t)}{(y^\pm(t)^2+1)^2}.
	\end{cases}
\end{align}

Both \eqref{ham_syst_2d_eq1} and \eqref{ham_syst_2d_eq2} admit a unique equilibrium in the point $P_e:=(0,\pi)$. Moreover, a linearization around that point shows that it is a center. In view of that, as it is shown by our simulations: 
\begin{itemize}
	\item[$\bullet$] when considering the initial datum $(x_0,y_0,\xi_0,\eta_0)=(0,y_0,\pi,\eta_0)$, the corresponding solution does not propagates in the vertical direction.
	\item[$\bullet$] when considering the initial datum $(x_0,y_0,\xi_0,\eta_0)=(x_0,0,\xi_0,\pi)$, the corresponding solution does not propagates in the horizontal direction.
	\item[$\bullet$] when considering the initial datum $(x_0,y_0,\xi_0,\eta_0)=(0,0,\pi,\pi)$, the corresponding solution does not propagates neither in the vertical nor in the horizontal direction.
\end{itemize}

Finally, also in the two-dimensional case they appear phenomena of internal reflection, which are purely numerical and are related to the employment of a non-uniform mesh. 

\begin{figure}[h]
	\centering
	\begin{minipage}{0.3\textwidth}
		\includegraphics[scale=0.16]{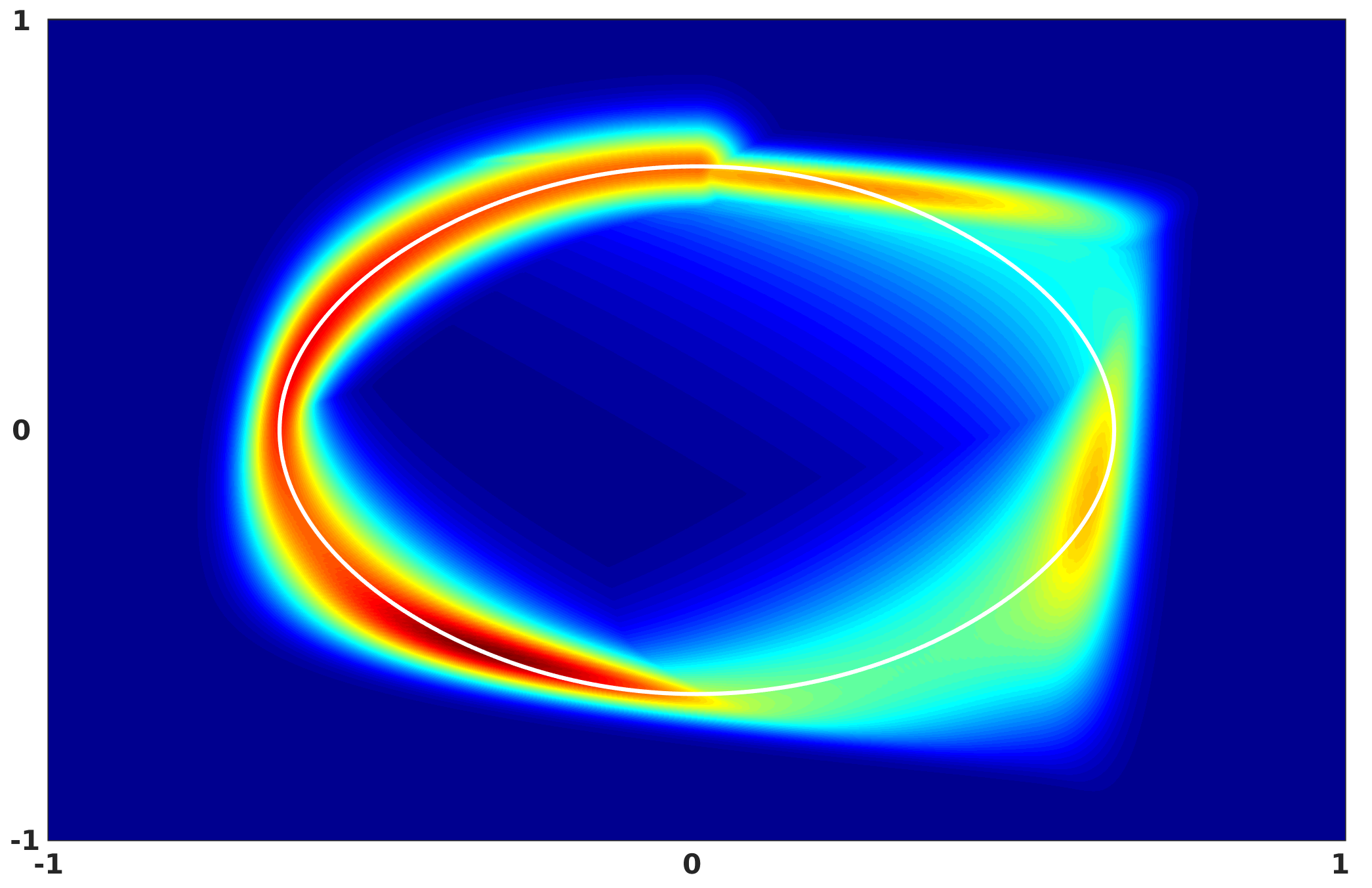}\caption*{(a)}
	\end{minipage}
	\hspace{1.5cm}
	\begin{minipage}{0.3\textwidth}
		\includegraphics[scale=0.16]{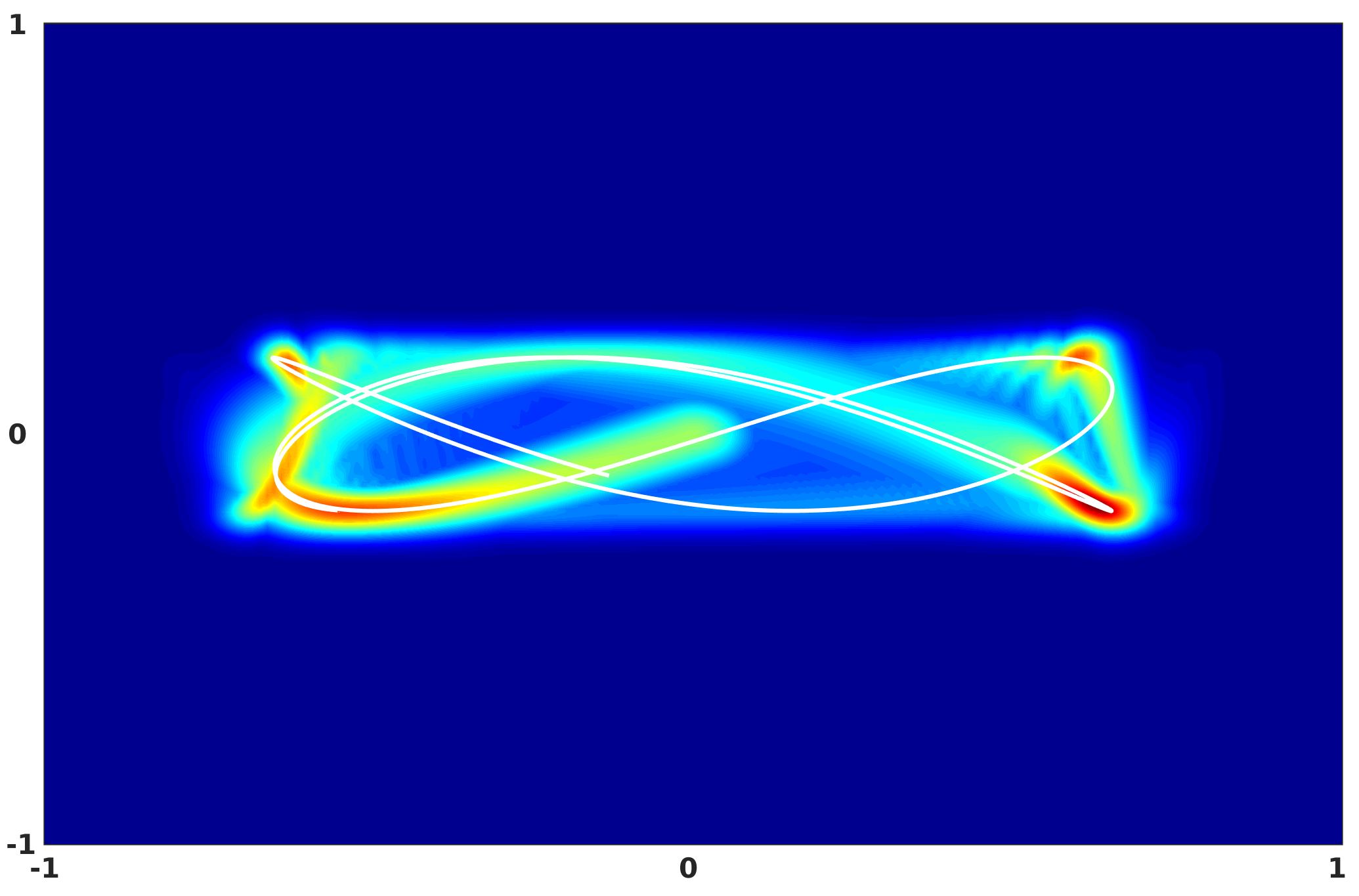}\caption*{(b)}
	\end{minipage}
	\\
	\begin{minipage}{0.3\textwidth}
		\includegraphics[scale=0.16]{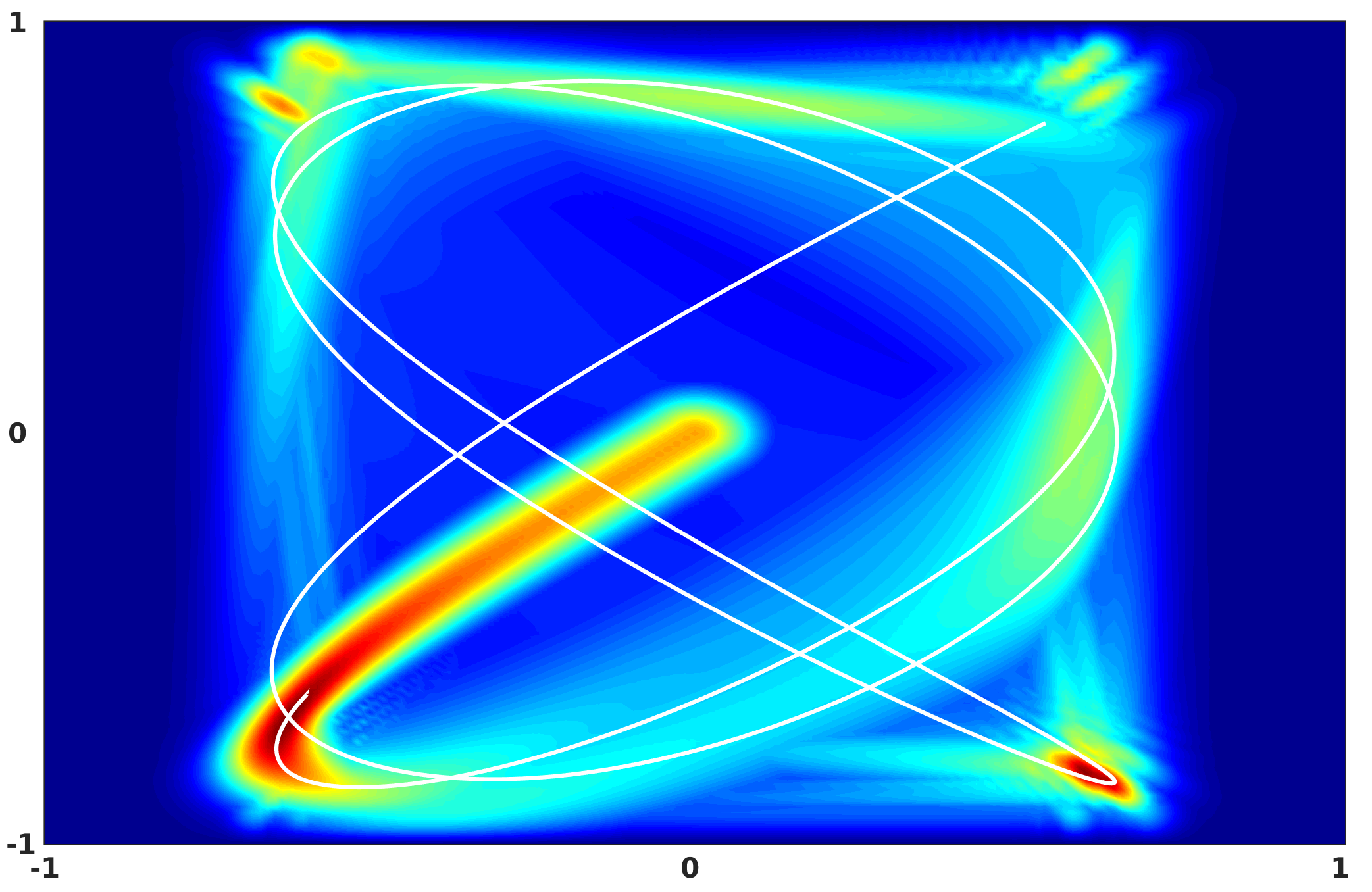}\caption*{(c)}
	\end{minipage}
	\hspace{1.5cm}
	\begin{minipage}{0.3\textwidth}
		\includegraphics[scale=0.16]{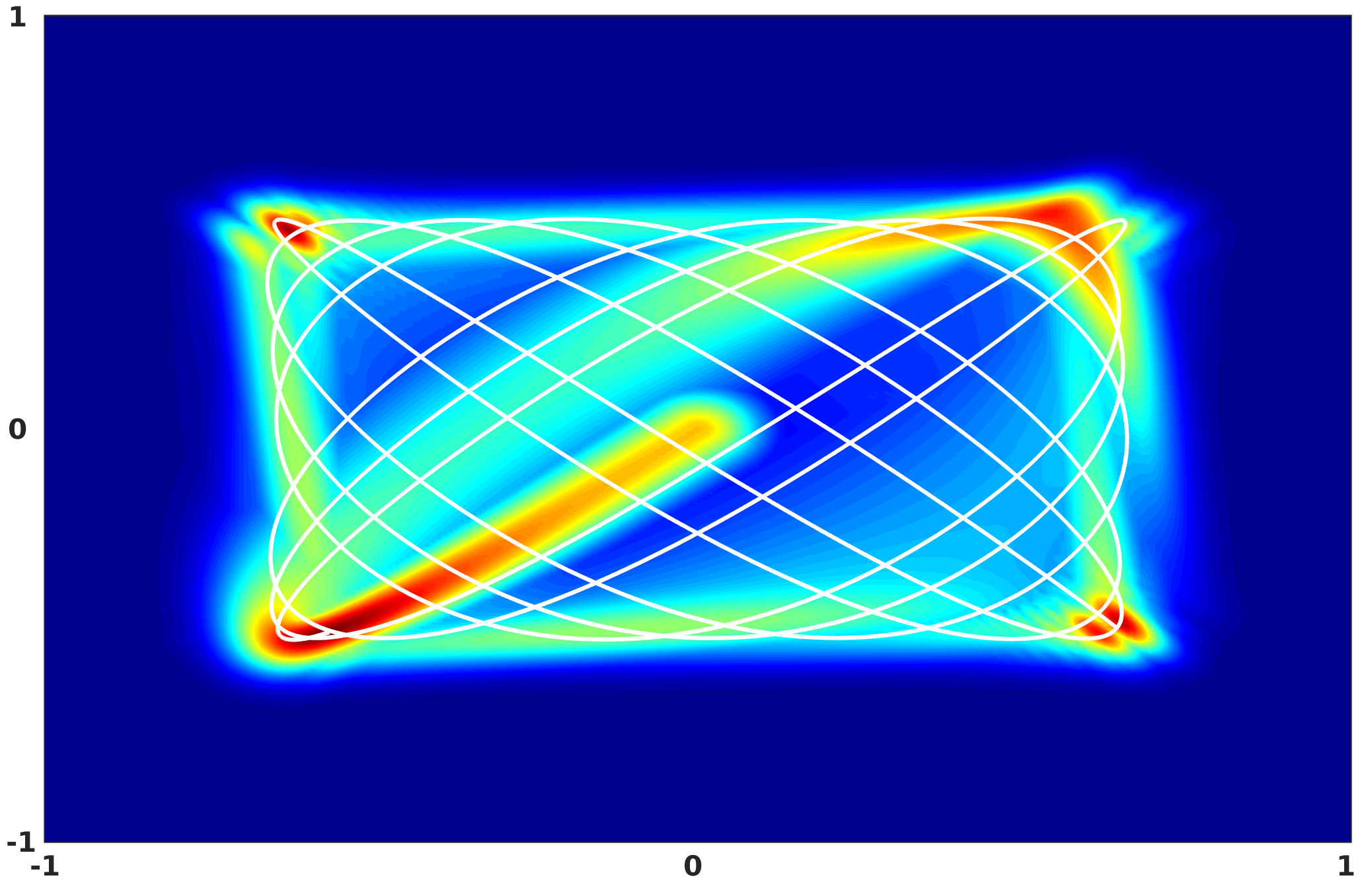}\caption*{(d)}
	\end{minipage}
	\caption{Numerical solutions corresponding to the initial datum \eqref{in_data_2d} with parameters given in Table \ref{in_data_table}.}\label{plot_wave2d_u3}
\end{figure}

In Figure \ref{plot_wave2d_u3}, we consider an initial data with parameters $(x_0,y_0,\xi_0,\eta_0)$ which are listed in Table \ref{in_data_table}.
\begin{table}[!h]
\centering 
\begin{tabular}{|c|c|c|c|c|c|}
	\hline
	& $x_0$ & $y_0$ & $\xi_0$ & $\eta_0$ & $T$
	\\
	\hline
	Figure \ref{plot_wave2d_u3}a & $0$ & $\tan(\arccos(\sqrt[4]{1/2}))$ & $\pi/2$ & $\pi$ & $8s$
	\\
	\hline
	Figure \ref{plot_wave2d_u3}b & $0$ & $0$ & $\pi/2$ & $5\pi/6$ & $21s$
	\\
	\hline
	Figure \ref{plot_wave2d_u3}c & $0$ & $0$ & $\pi/2$ & $7\pi/18$ & $37s$
	\\
	\hline
	Figure \ref{plot_wave2d_u3}d & $0$ & $0$ & $\pi/2$ & $7\pi/12$ & $118s$
	\\
	\hline
\end{tabular}\caption{Initial data $(x_0,y_0,\xi_0,\eta_0)$ used in the simulations.}\label{in_data_table}
\end{table}

The RGB surface represents the evolution of a Gaussian wave packet. Moreover, in order to better appreciate the concentration along the rays (which is in part lost due to some numerical dispersion effect), in white we plot the component $x(t)$ against $y(t)$ for $t\in [0,T]$.   

We note that the projection of the characteristic rays $(x^\pm(t),y^\pm(t),\xi^\pm(t),\eta^\pm(t))$ corresponding to these
initial data on the physical space are Lissajous-type curves (see \cite{von1994practical}), which means that in each direction the solutions $(x^\pm(t), \xi^\pm(t))$ and $(y^\pm(t),\eta^\pm(t))$ of \eqref{ham_syst_2d_eq1} and \eqref{ham_syst_2d_eq2} are orbits around the fixed point $(0,\pi)$. They give rise in the physical space to curves $x^\pm(t)$ and $y^\pm(t)$ presenting each the umklapp phenomenon and a oscillatory shapes with periodicity respectively (see \cite[Section 3.2]{ervedoza2015numerical})
\begin{align*}
	T_1(x_0) = \frac{2r_0}{r_1}\int_{-x^*}^{x^*} \frac{g'(z)}{\sqrt{1-\left(\frac{g'(z)}{g'(x^*)}\right)^2}}\,dz\;\;\;\textrm{ and }\;\;\; T_2(y_0) = \frac{2r_0}{r_1}\int_{-y^*}^{y^*} \frac{g'(z)}{\sqrt{1-\left(\frac{g'(z)}{g'(y^*)}\right)^2}}\,dz,
\end{align*}
where $x^*$ and $y^*$ are such that $g'(x_0) = g'(x^*)\sin(\xi_0/2)$ and $g'(y_0) = g'(y^*)\sin(\eta_0/2)$. Lastly, this internal reflection phenomena are once again related to the coarseness of the mesh around the boundary, which makes the velocity of the high-frequency waves decrease (until, eventually, change sign) in one ore both the spatial directions while moving away from the center of the domain. 

\subsection{Conclusions}

Summarizing, our analysis shows that, also in the two-dimensional case, the discretization of \eqref{wave_var_2d} introduces several interesting phenomena in the propagation of numerical solutions. In more detail, the solutions of the discrete two-dimensional wave equation \eqref{wave_var_2d_discr} can be classified essentially in three groups. On the one hand, we have low-frequency solutions that are able to travel freely until they reach the boundary of the domain, where they are reflected (Figure \ref{plot_wave2d_low1}). On the other hand, we also have high-frequency waves which either remain trapped bouncing indefinitely between two parallel sides of the domain (Figure \ref{plot_wave2d_np1}) or are kept confined in the interior of the domain, without possibility of reaching the boundary (Figures \ref{plot_wave2d_np3} and \ref{plot_wave2d_u3}). Some of these pathologies appear both with uniform and non-uniform meshes, and are a consequence of the trigonometric nature of the discrete group velocity. In particular, when one or both the components of this group velocity become zero, the waves show lack of propagation either in the horizontal direction, in the vertical one or in both. Lastly, the phenomena of internal reflection are once again due to the introduction of non-uniform meshes, which generate fictitious numerical boundaries when the grid passes from fine to coarse.

\section{Final comments and open problems}

In this article, we analyzed the propagation properties of the finite difference approximation of one and two-dimensional wave equations, both on uniform and non-uniform numerical grids.

Starting from the observation that the energy of continuous and semi-discrete high-frequency solutions propagates along bi-characteristic rays, we showed that the dynamics may change from the continuous to the semi-discrete setting, because of the different nature of the corresponding Hamiltonians. In particular, as a result of the accumulation of the local effects introduced by the heterogeneity of the employed grid, numerical high-frequency solutions can bend in a singular and unexpected manner. Moreover, this phenomenon has to be added to the well known numerical dispersion effect, producing the high-frequency discrete group velocity to vanish, even in uniform grids. 

Overall, in the one-dimensional case, the result of these pathologies are slowly propagating numerical high-frequency components that never get to the exterior boundary of the domain, a fact which is against  the behavior of the continuous solutions. Besides, these effects are enhanced in the multi-dimensional case where the interaction and combination of such behaviors in the various space directions may produce, for instance, waves that are trapped by the numerical grid in closed loops, without ever getting to the exterior boundary.

Our analysis allows to explain all such unusual behaviors and illustrates that the effect of the non-uniformity of the numerical mesh is similar to the one that the heterogeneity of the coefficients introduces on the propagation of continuous waves, making the bi-characteristic rays bend. One of the main objectives of this paper has been to illustrate all these possible pathologies through accurate numerical simulations. 

Furthermore, our results constitute a warning both for adaptivity and for the treating of control and inverse problems. 

In broad terms, the goal of adaptivity is to refine a mesh on the support of the solution, keeping it coarse where the solution has little oscillations and energy. Our analysis shows that, in this context, adaptivity has to be performed with some attention. Indeed, if one is not careful enough when refining the mesh, they can be produced spurious effects due to the fact that waves feel the fictitious numerical boundaries that are generated when the grid passes from fine to coarse.

Finally, the results of this paper are also a signal that the dangers of uniform meshes in the study of numerical control and inverse problems may be enhanced when the mesh is non-uniform. In more details, we showed that heterogeneity of the grid introduces added trapping effects, which need to be avoided in order to prove convergence in the context of controllability, stabilization or inversion algorithms.

We now conclude our discussion by introducing some suggestions of future research, which are related with the topics addressed in the present paper.

\begin{enumerate}
	\item \textbf{General analysis of non-uniform grids.} In this paper, we limited our analysis only to the selection of some particular example of non-uniform mesh that may be employed for the space semi-discretization of a wave equation. On the other hand, it may be interesting to address a general discussion on the link between the choice of a mesh and the effects that this produces in terms of wave propagation. For our numerical implementation, we focused only on two classes of non-uniform grids: one fine at the center of the computational domain and coarse at the extrema and, specularly, one fine at the extrema but coarse at the center. What would happen if, for instance, we used a mesh generated by a periodic function, with repeated regions of coarseness and refinement? This periodic structure of the grid would be reflected also on the propagation of the bi-characteristics? Let us mention also that, for the approximation of the two-dimensional wave equation, we employed the same refining function in both the space directions, but we could instead have decided to combine two difference kind of refinement for the $x$ and for the $y$ axis. In that case, we can expect that new classes of phenomena arise, as a consequence of the non uniformity of the mesh. 
		
	\item \textbf{Irregular meshes.} The analysis addressed in this paper remains valid while the micro-local description of the phenomena we are interested in may be performed. For this, we need at least to make sure that the Hamiltonian system of bi-characteristic rays admits solutions, which requires a minimal regularity of the coefficients and of the grid. On the other hand, it would be certainly interesting to explore the possibility of employing irregular non-uniform meshes, generated by the transformation of a uniform one through a singular map. It is natural to expect that other new phenomena and pathologies will appear, and a systematic analysis of this extra possible behaviors is worth to be developed. To some extents, this is also related to the behavior of one-dimensional continuous waves in the lack of $C^1$-regularity of the coefficients, which exhibit unexpected concentration phenomena of the high-frequency that contradict all the well-known propagation and dispersion properties of waves in homogeneous media (see \cite{castro2002concentration}). 
	
	\item \textbf{Other numerical schemes.} In this paper, we have considered only the case of finite difference approximations for the wave equation. On the other hand, it would be interesting to consider also other numerical schemes such as mixed finite element, which have been introduced, for instance, in \cite{castro2006boundary,ervedoza2010observability,glowinski1989mixed} with control purposes. According to our analysis, the first step would be to identify the principal symbol associated to the FE discretization of the variable coefficients wave equation in one and two space dimension and both on uniform and non-uniform meshes. From there, the study of the corresponding Hamiltonian system would again allow to describe the propagation properties of the numerical solutions. To the best of our knowledge, this is a topic that is yet to be addressed.
	
	\item \textbf{Filtering mechanisms on non-uniform meshes.} We mentioned that our analysis is motivated by the study of control properties for discrete wave equations. As we discussed in the introduction to this work, the presence of high-frequency spurious numerical solutions traveling at velocity of the order of $h$ yields the exponential blow-up of the observability constant and, consequently, the divergence of numerical controls. For overcoming this problem, several different approaches have been proposed, all brought together by the common idea of filtering the problematic high frequencies (see \cite{ervedoza2012wave} and the references therein). In addition, we mention that the very same methodology applies in the context of inverse problems for the discrete wave equation, as it has been discussed, e.g., in \cite{baudouin2013convergence,baudouin2015stability}. These filtering mechanisms have been fully analyzed in the context of a uniform mesh discretization but, to the best of our knowledge, nothing is known concerning their extension to the non-uniform mesh case. On the other hand, to address this kind of analysis also for non-uniform meshes is a question worth to be considered.	 
	
	\item \textbf{Fully discrete wave equations.} The analysis of the principal symbol and of the corresponding Hamiltonian system that we addressed in this article has been developed starting from semi-discrete models. On the other hand, the same study can be performed also for fully discrete approximations of the wave equation, both in one and two dimensions. The equation being second order in time, it is natural to expect that the changes that one would observe in the Hamiltonian are analogous to the ones that already appeared when passing from the continuous to the semi-discrete system. To analyze how these changes affect the propagation of numerical waves, it is certainly an interesting issue.
\end{enumerate}

\bibliography{biblio}

\begin{thebibliography}{10}

\bibitem{allaire2003dispersive}
{\sc Allaire, G.}
\newblock Dispersive limits in the homogenization of the wave equation.
\newblock In {\em Ann. Fac. Sci. Toulouse Math.\/} (2003), vol.~12,
  Universit{\'e} Paul Sabatier, Institut de Math{\'e}matiques, pp.~415--431.

\bibitem{allaire2005homogenization}
{\sc Allaire, G., and Pjatnicki{\u\i}, A.}
\newblock Homogenization of the {S}chr{\"o}dinger equation and effective mass
  theorems.
\newblock {\em Commun. Math. Phys. 258}, 1 (2005), 1--22.

\bibitem{bardos1992sharp}
{\sc Bardos, C., Lebeau, G., and Rauch, J.}
\newblock Sharp sufficient conditions for the observation, control, and
  stabilization of waves from the boundary.
\newblock {\em SIAM J. Control Optim. 30}, 5 (1992), 1024--1065.

\bibitem{baudouin2013convergence}
{\sc Baudouin, L., and Ervedoza, S.}
\newblock Convergence of an inverse problem for a 1-{D} discrete wave equation.
\newblock {\em SIAM J. Control Optim. 51}, 1 (2013), 556--598.

\bibitem{baudouin2015stability}
{\sc Baudouin, L., Ervedoza, S., and Osses, A.}
\newblock Stability of an inverse problem for the discrete wave equation and
  convergence results.
\newblock {\em J. Math. Pures Appl. 103}, 6 (2015), 1475--1522.

\bibitem{burq1997condition}
{\sc Burq, N., and G{\'e}rard, P.}
\newblock Condition n{\'e}cessaire et suffisante pour la contr{\^o}labilit{\'e}
  exacte des ondes.
\newblock {\em C.R. Acad. Sci. Paris S{\'e}r. I 325}, 7 (1997), 749--752.

\bibitem{burq1998controle}
{\sc Burq, N., and Schlenker, J.-M.}
\newblock Contr{\^o}le de l'{\'e}quation des ondes dans des ouverts comportant
  des coins.
\newblock {\em Bull. Soc. Math. France 126}, 4 (1998), 601.

\bibitem{castro2006boundary}
{\sc Castro, C., and Micu, S.}
\newblock Boundary controllability of a linear semi-discrete 1-{D} wave
  equation derived from a mixed finite element method.
\newblock {\em Numer. Math. 102}, 3 (2006), 413--462.

\bibitem{castro1996remark}
{\sc Castro, C., and Zuazua, E.}
\newblock A remark on the spectral asymptotic analysis in homogenization.
\newblock {\em C. R. Acad. Sci. Ser. I 322}, 11 (1996), 1043--1047.

\bibitem{castro2000low}
{\sc Castro, C., and Zuazua, E.}
\newblock Low frequency asymptotic analysis of a string with rapidly
  oscillating density.
\newblock {\em SIAM J. Appl. Math. 60}, 4 (2000), 1205--1233.

\bibitem{castro2002concentration}
{\sc Castro, C., and Zuazua, E.}
\newblock Concentration and lack of observability of waves in highly
  heterogeneous media.
\newblock {\em Arch. Rat. Mech. Anal. 164}, 1 (2002), 39--72.

\bibitem{ervedoza2010observability}
{\sc Ervedoza, S.}
\newblock Observability properties of a semi-discrete 1d wave equation derived
  from a mixed finite element method on nonuniform meshes.
\newblock {\em ESAIM: Control Optim. Calc. Var. 16}, 2 (2010), 298--326.

\bibitem{ervedoza2015numerical}
{\sc Ervedoza, S., Marica, A., and Zuazua, E.}
\newblock Numerical meshes ensuring uniform observability of one-dimensional
  waves: construction and analysis.
\newblock {\em IMA J. Numer. Anal. 36}, 2 (2015), 503--542.

\bibitem{ervedoza2012wave}
{\sc Ervedoza, S., and Zuazua, E.}
\newblock The wave equation: control and numerics.
\newblock {\em Control and stabilization of PDEs, P.M. Cannarsa and J.M. Coron
  eds., Lecture Notes in Mathematics 2048\/} (2012), 245--340.

\bibitem{ervedoza2013numerical}
{\sc Ervedoza, S., and Zuazua, E.}
\newblock On the numerical approximation of exact controls for waves.
\newblock {\em Springer Briefs in Mathematics XVII\/} (2013).

\bibitem{gerard1991microlocal}
{\sc G{\'e}rard, P.}
\newblock Microlocal defect measures.
\newblock {\em Comm. Partial Differential Equations 16}, 11 (1991), 1761--1794.

\bibitem{glowinski1992ensuring}
{\sc Glowinski, R.}
\newblock Ensuring well-posedness by analogy; {S}tokes problem and boundary
  control for the wave equation.
\newblock {\em J. Comput. Phys. 103}, 2 (1992), 189--221.

\bibitem{glowinski1990numerical}
{\sc Glowinski, R., and Chin-Hsien, L.}
\newblock On the numerical implementation of the {H}ilbert uniqueness method
  for the exact boundary controllability of the wave equation.
\newblock {\em C. R. Acad. Sc. S{\'e}r. 1 Math. 311}, 2 (1990), 135--142.

\bibitem{glowinski1989mixed}
{\sc Glowinski, R., Kinton, W., and Wheeler, M.~F.}
\newblock A mixed finite element formulation for the boundary controllability
  of the wave equation.
\newblock {\em Internat. J.Numer. Methods Engrg 27}, 3 (1989), 623--635.

\bibitem{glowinski1990approach}
{\sc Glowinski, R., Li, C.-H., and Lions, J.-L.}
\newblock A numerical approach to the exact boundary controllability of the
  wave equation (i) {D}irichlet controls: {D}escription of the numerical
  methods.
\newblock {\em Japan J. Appl. Math. 7}, 1 (1990), 1--76.

\bibitem{hormander1985analysis}
{\sc Hormander, L.}
\newblock The analysis of linear partial differential operators, {V}ols {I},
  {III}, 1985.

\bibitem{ignat2009convergence}
{\sc Ignat, L.~I., and Zuazua, E.}
\newblock Convergence of a two-grid algorithm for the control of the wave
  equation.
\newblock {\em J. Europ. Math. Soc. 11}, 2 (2009), 351--391.

\bibitem{infante1999boundary}
{\sc Infante, J.~A., and Zuazua, E.}
\newblock Boundary observability for the space semi-discretizations of the 1-d
  wave equation.
\newblock {\em Math. Model. Numer. Anal. 33}, 2 (1999), 407--438.

\bibitem{lions1993mesures}
{\sc Lions, P.-L., and Paul, T.}
\newblock Sur les mesures de {W}igner.
\newblock {\em Rev. Mat. Iberoamericana 9}, 3 (1993), 553--618.

\bibitem{loreti2008ingham}
{\sc Loreti, P., and Mehrenberger, M.}
\newblock An {I}ngham type proof for a two-grid observability theorem.
\newblock {\em ESAIM: Control Optim. Calc. Var. 14}, 3 (2008), 604--631.

\bibitem{macia2002propagacion}
{\sc Maci{\`a}, F.}
\newblock {\em Propagaci{\'o}n y control de vibraciones en medios discretos y
  continuos}.
\newblock PhD thesis, PhD Thesis, Universidad Complutense de Madrid, 2002.

\bibitem{macia2002lack}
{\sc Mac\'ia, F., and Zuazua, E.}
\newblock On the lack of observability for wave equations: a gaussian beam
  approach.
\newblock {\em Asympt. Anal. 32}, 1 (2002), 1--26.

\bibitem{marica2015propagation}
{\sc Marica, A., and Zuazua, E.}
\newblock Propagation of 1d waves in regular discrete heterogeneous media: a
  {W}igner measure approach.
\newblock {\em Found. Comp. Math. 15}, 6 (2015), 1571--1636.

\bibitem{markowich1994wigner}
{\sc Markowich, P., Mauser, N., and Poupaud, F.}
\newblock A {W}igner-function approach to (semi) classical limits: {E}lectrons
  in a periodic potential.
\newblock {\em J. Math. Phys. 35}, 3 (1994), 1066--1094.

\bibitem{micu2002uniform}
{\sc Micu, S.}
\newblock Uniform boundary controllability of a semi-discrete 1-{D} wave
  equation.
\newblock {\em Numer. Math. 91}, 4 (2002), 723--768.

\bibitem{negreanu2004convergence}
{\sc Negreanu, M., and Zuazua, E.}
\newblock Convergence of a multigrid method for the controllability of a 1-d
  wave equation.
\newblock {\em C.R. Acad. Sci. Paris 338}, 5 (2004), 413--418.

\bibitem{pjatnickiui1982limit}
{\sc Pjatnicki{\u\i}}.
\newblock On the limit behavior of the domain of dependence of a hyperbolic
  equation with rapidly oscillating coefficients.

\bibitem{ralston1982gaussian}
{\sc Ralston, J.}
\newblock Gaussian beams and the propagation of singularities.
\newblock {\em Studies in partial differential equations 23\/} (1982),
  206--248.

\bibitem{rauch2005polynomial}
{\sc Rauch, J., Zhang, X., and Zuazua, E.}
\newblock Polynomial decay for a hyperbolic-parabolic coupled system☆.
\newblock {\em J. Math. Pures Appl. 84}, 4 (2005), 407--470.

\bibitem{tartar1990h}
{\sc Tartar, L.}
\newblock H-measures, a new approach for studying homogenisation, oscillations
  and concentration effects in partial differential equations.
\newblock {\em Proc. Roy. Soc. Edinburgh Sec. A 115}, 3-4 (1990), 193--230.

\bibitem{trefethen1982wave}
{\sc Trefethen, L.~N.}
\newblock {\em Wave propagation and stability for finite difference schemes}.
\newblock PhD thesis, PhD Thesis.

\bibitem{trefethen1982group}
{\sc Trefethen, L.~N.}
\newblock Group velocity in finite difference schemes.
\newblock {\em SIAM rev. 24}, 2 (1982), 113--136.

\bibitem{vichnevetsky1980propagation}
{\sc Vichnevetsky, R.}
\newblock Propagation properties of semi-discretizations of hyperbolic
  equations.
\newblock {\em Math. Comp. Simul. 22}, 2 (1980), 98--102.

\bibitem{vichnevetsky1981energy}
{\sc Vichnevetsky, R.}
\newblock Energy and group velocity in semi discretizations of hyperbolic
  equations.
\newblock {\em Math. Comp. Simul. 23}, 4 (1981), 333--343.

\bibitem{vichnevetsky1981propagation}
{\sc Vichnevetsky, R.}
\newblock Propagation through numerical mesh refinement for hyperbolic
  equations.
\newblock {\em Math. Comp. Simul. 23}, 4 (1981), 344--353.

\bibitem{vichnevetsky1987wave}
{\sc Vichnevetsky, R.}
\newblock Wave propagation and reflection in irregular grids for hyperbolic
  equations.
\newblock {\em Appl. Numer. Math. 3\/} (1987), 133--166.

\bibitem{vichnevetsky1982fourier}
{\sc Vichnevetsky, R., and Bowles, J.~B.}
\newblock {\em Fourier analysis of numerical approximations of hyperbolic
  equations}, vol.~5.
\newblock Siam, 1982.

\bibitem{von1994practical}
{\sc Von~Seggern, D.~H.}
\newblock {\em Practical handbook of curve design and generation}.
\newblock CRC Press, 1994.

\bibitem{wigner1932quantum}
{\sc Wigner, E.}
\newblock On the quantum correction for thermodynamic equilibrium.
\newblock {\em Phys. Rev. 40}, 5 (1932), 749.

\bibitem{zuazua1999boundary}
{\sc Zuazua, E.}
\newblock Boundary observability for the finite-difference space
  semi-discretizations of the 2-{D} wave equation in the square.
\newblock {\em J. Math. Pures Appl. 78}, 5 (1999), 523--563.

\bibitem{zuazua2005propagation}
{\sc Zuazua, E.}
\newblock Propagation, observation, control and numerical approximation of
  waves.
\newblock {\em SIAM Rev. 47}, 2 (2005), 197--243.

\end{thebibliography}
\end{document}